\documentclass[reqno,11pt]{amsart}
\usepackage{amsfonts}
\usepackage{times,amssymb, amsmath,latexsym,amsfonts,amsbsy, amsthm}
\usepackage{color,esint,pgfplots}
\usepackage{cases}
\usepackage{mathtools}
\usepackage{hyperref}
\usepackage{enumitem}
\usepackage{setspace}
\usepackage{mathrsfs}
\usepackage{stmaryrd}
\usepackage{graphicx}
\allowdisplaybreaks[4]
\usepackage[left=2cm, right=2cm, top=3cm, bottom=2cm]{geometry}
\newcommand{\picdis}[1]{}
\newcommand{\p}{\partial}

\newcommand{\mfu}{\mathfrak{u}}
\newcommand{\mfv}{\mathfrak{v}}

\newcommand{\udr}{\underline{\rho}}
\newcommand{\udu}{\underline{u}}
\newcommand{\udv}{\underline{v}}

\renewcommand{\[}{\begin{equation}}
\renewcommand{\]}{\end{equation}}

\renewcommand{\div}{\operatorname{div}}

\newcommand{\eqdef}{\overset{\mbox{\tiny{def}}}{=}}
\newcommand{\eps}{\epsilon}

\newtheorem{theorem}{Theorem}[section]

\newtheorem{lemma}[theorem]{Lemma}
\newtheorem{proposition}[theorem]{Proposition}
\newtheorem{corollary}[theorem]{Corollary}

\theoremstyle{definition}
\newtheorem{remark}{Remark}[section]
\begin{document}
\date{\today}
\title[Structural stability of subsonic boundary layers]{Structural stability of boundary layers in the entire subsonic regime}

\author{Shengxin Li}
\address{Department of Applied Mathematics, The Hong Kong Polytechnic University, Hong Kong, China. }
\email{sh1li@polyu.edu.hk}
\author{Tong Yang}
\address{Institute for Math \& AI, Wuhan University, Wuhan, China. }
\email{tongyang@whu.edu.cn}
\author{Zhu Zhang}
\address{Department of Applied Mathematics, The Hong Kong Polytechnic University, Hong Kong, China. }
\email{zhuama.zhang@polyu.edu.hk}

\begin{abstract}
Despite the physical importance, there are limited mathematical theories for the compressible Navier-Stokes equations with strong boundary layers. This is mainly due to the absence of a stream function structure, unlike the extensively studied incompressible fluid dynamics in two dimensions. This paper aims to establish the structural stability of boundary layer profiles in the form of shear flow for the two-dimensional steady compressible Navier-Stokes equations. Our estimates are uniform across the entire subsonic regime, where the Mach number $m\in (0,1)$. As a byproduct, we provide the first result concerning the low Mach number limit in the presence of Prandtl boundary layers. The proof relies on the quasi-compressible-Stokes iteration introduced in \cite{YZ23}, along with a subtle analysis of the interplay between density and velocity variables in different frequency regimes, and the identification of cancellations in higher-order estimates.
\end{abstract}
\numberwithin{equation}{section}

\setcounter{secnumdepth}{3}

\date{\today}
%\thanks{}
\maketitle

\tableofcontents
\thispagestyle{empty}

%%      ---------------------------------------------------------------------
%%      ------------------- TABLE OF CONTENTS (OPTIONAL) --------------------
%%      ---------------------------------------------------------------------

% \tableofcontents

%%      ---------------------------------------------------------------------
%%      ---------------------------- BODY OF PAPER --------------------------
%%      ---------------------------------------------------------------------

%%      Please input or insert the body of your paper here.

\section{Introduction}
\subsection{Background}
We are interested in the two-dimensional steady compressible Navier-Stokes equations for isentropic flow in half-plane $\Omega=\{(x, y)\mid (x, y)\in\mathbb{T}_L\times \mathbb{R}_+\}:$
\begin{align}\label{1.1}
\begin{cases}
\nabla\cdot(\rho^\nu{\bf{u}}^\nu)=0,\\
\nabla\cdot(\rho^\nu{\bf{u}}^\nu\otimes {\bf{u}}^\nu)+\nabla P(\rho^\nu)-\mu\nu\Delta {\bf{u}}^\nu-\lambda\nu\nabla(\nabla\cdot {\bf{u}}^\nu)=\rho^\nu{\bf{F}}^\nu,\\
{\bf{u}}^\nu|_{y=0}={\bf{0}}.
\end{cases}
\end{align}
In these equations, $\rho^\nu$, ${\bf{u}}^\nu=(u^\nu, v^\nu)$ and $P=P(\rho^\nu)$ denote the density, velocity field, and pressure respectively, and ${\bf{F}}^\nu$ is a given external force. The parameter $\nu$ stands for the recipocal of the Reynolds number. The constants $\mu>0$ and $\lambda\geq 0$ represent the rescaled shear and bulk viscosities. For simplicity, we set $\mu=1$ throughout paper.  The tangential variable $x$ is defined on the torus $\mathbb{T}_L$ with length $L$.

Understanding the asymptotic behavior of fluid variables $(\rho^\nu,u^\nu,v^\nu)$ in the limit $\nu\rightarrow 0^+$ is a fundamental problem in fluid mechanics. The interaction between fluids and rigid wall generates complex and singular structures, such as the creation of large vorticity, and emergence of boundary layers. To understand these structures, Prandtl introduced an asymptotic theory of boundary layers in 1904. According to Prandtl's theory, the fluid domain has the following scenario:  away from the boundary, the viscosity is negligible and the fluid can be effectively described by the Euler system; near the boundary, a thin layer forms with a thickness propotional to $\nu^{\frac12}$. Within this layer, the density $\rho^\nu$ remains approximately as the Euler solution. However, the velocity fields exhibit the  following asymptotic behavior:
$$(u^\nu,v^\nu)\approx (U(x,y/\sqrt{\nu}),V(x,y/\sqrt{\nu})). 
$$
Here the boundary layer profile $(U,V)$ satisfies the Prandtl equation.

 The mathematical theories for Prandtl equation started from the foundational work by Oleinik \cite{O63}, where the local-in-space solution to steady system was constructed under the assumption of monotonicity of the horizontal velocity component at the boundary. However, extending this result to more general settings, such as global-in-space or unsteady flows, presents significant challenges. These challenges arise from various instability mechanisms inherent in the Prandtl model. For instance, E-Engquist \cite{EE} constructed blow up solutions. G\'erard-Varet and Dormy \cite{GD} explored instabilities that occur in non-monotone flows. Further studies by Dalibard-Masmoudi \cite{DM} analyze behavior of boundary layer flow near its separation.  Collot-Ghoul-Ibrahim-Masmoudi \cite{CGIM} have investigated singularity formations in the Prandtl equation. Despite these challenges, the well-posedness theories for Prandtl equation have been established by imposing specific structural conditions on the flow, as discussed in \cite{AWXY15,I,IM22,MW15}, or for high-regularity data, as explored in \cite{DG,SC98}, among many others. 

Although the Prandtl equations can be solved under specific conditions, the validity of boundary layer expansion is still very difficult due to intrinsic instabilities in the Navier-Stokes system at high Reynolds numbers. Two main destabilizing effects have been identified so far: The first one is induced by the inflection point within the boundary layer. Such a profile can lead to strong ill-posedness in the Navier-Stokes system below the analytic regularity. For the nonlinear instabilities associated to these profiles, we refer to the works of Grenier \cite{G00} and Grenier-Nguyen \cite{GN19}. The second one is the Tollmien-Schlichting instability that was rigorously justified by Grenier-Guo-Nguyen \cite{GGN16}.
That is,  even inviscid-stable profiles (monotone and concave) can lose their stabilities with small viscosity. Due to these instability mechanisms, justification of the Prandtl ansatz can only be achieved in the function spaces of high regularities. In this direction, the first result was obtained by Sammartino-Caflisch \cite{SC98}, who justified the boundary layer expansion in the analytic setting. Subsequent research in the analytic framework has been explored in various situations. For example, Maekawa \cite{M14} explored scenarios where the vorticity is supported away from the boundary; Wang-Wang-Zhang \cite{WWZ17} justified the Prandtl expansion based on co-normal estimates;  Nguyen-Nguyen \cite{NN18} investigated  the inviscid limit without relying on boundary layer expansions; Kukavica-Vicol-Wang \cite{KVW20} demonstrated that analyticity is only needed near the boundary. If the regularity is below analyticity, G\'erard-Varet-Maekawa-Masmoudi \cite{GMM18,GMM20} justified the Prandtl ansatz in the critical Gevrey 3/2 space, and  Chen-Wu-Zhang \cite{CWZ22} proved the  $L^\infty$-stability up to initial time. We also mention the work by  Lopes Filho et al. \cite{LMN08} and Mazzucato-Taylor \cite{MT08} for the stability results under some symmetry conditions.

Boundary layers exhibit better behavior in steady flow because of  the absence of certain instability mechanisms so that the justification of boundary layer expansions can be proved in the Sobolev framework. In this context, two research directions have been identified. The first approach is to treat horizontal variable $x$ like the ``time'' variable. By imposing appropriate boundary conditions on horizontal boundaries, one can justify the Prandtl expansion for general $x$-dependent boundary layer profiles. Notable contributions in this area include the local-in-$x$ results by Guo-Nguyen \cite{GN19} addressing moving boundaries, as well as Guo-Iyer \cite{GI} which focuses on no-slip boundary conditions. Furthermore, Iyer-Masmoudi \cite{IM20} has studied asymptotic behavior for large $x$ of Blasius solutions.  The second approach, initiated by G\'erard-Varet-Maekawa \cite{GM19}, involves imposing periodic conditions in $x$ direction to investigate the structural stability properties of shear boundary layer profiles. We also mention the work by  Chen-Wu-Zhang \cite{CWZ23} for a result in large toroidal domains subject to a spectral condition.

In this paper, we study stability properties of Prandtl boundary layers for  compressible fluid. Unlike boundary layers in incompressible fluid, stability mechanisms of compressible boundary layers are significantly influenced by the Mach number. For subsonic flows, the Tollmien-Schlichting (T-S) type wave induces destabilization, as demonstrated by the second and third authors \cite{YZ23} for Mach number $m<1/\sqrt{3}$, and by Masmoudi-Wang-Wu-Zhang \cite{MWWZ24} for the entire subsonic regime. When the Mach number is greater than one,  compressibility effect becomes more significant that  leads to an inviscid-type instability, which is known as the Mack second mode.  This instability mechanism has been justified by Masmoudi-Wang-Wu-Zhang \cite{MWWZ241}.  Apart from the work by Wang-Wang-Zhang \cite{WWZ23}, where the inviscid limit in the analytic framework was studied, 
to our knowledge, there is no stability result without  analytic regularity. Given its complexity, a first attempt in this direction is to study the steady subsonic boundary layer flows. {\it The motivation of this paper is twofold. The first one is to establish the structural stability of boundary layers across the entire subsonic regime $m\in (0,1)$ in the Sobolev setting. And the second one is about the low Mach number limit in the presence of Prandtl boundary layers.}
\subsection{Main results}

Motivated by G\'erard-Varet-Maekawa \cite{GM19}, we set the $x$ variable in the torus $\mathbb{T_L}$, and study the structural stability of the following laminar boundary layer flow:
\begin{align}
(\rho_s, {\bf{u}}_s)\eqdef(1, U_s(Y), 0),\quad Y=\frac{y}{\sqrt{\nu}},~\text{with } U_s(0)=0,~\lim_{Y\rightarrow +\infty}U_s(Y)=1.\label{BL}
\end{align}
This is an exact solution to \eqref{1.1} with external force ${\bf{F}}^\nu=(-\partial_Y^2U_s, 0)$.  According to Prandtl's theory, the pressure does not exhibit a leading order boundary layer, and consequently, the density has no boundary layer. Therefore, without loss of generality, we assume that the background density is equal to one. Moreover, the background shear flow $U_s(Y)\in C^3({\mathbb{R}_+})\cap C^1(\overline{\mathbb{R}_+})$ is assumed to satisfy following structural conditions: 
\begin{align}\label{2.1}
	&U_s(0)=0,\quad \partial_YU_s(0)=1,\quad U_s(Y)>0, \forall Y>0,\quad \lim\limits_{Y\to\infty}U_s(Y)=1;\\
	&\sup\limits_{Y\ge 0}(1+Y)^s\left(|1-U_s(Y)|+\sum_{k=1}^3|\partial_Y^k U_s(Y)|\right)<\infty,~\text{for }s>3.\label{2.2}
\end{align}

To study stability properties of the above boundary layer profile, we denote the perturbation  by $(\rho, {\bf{u})}=(\rho^\nu, {\bf{u}}^\nu)-(\rho_s, {\bf{U}}_s)$, which is  induced by the  perturbation of external force
 \begin{align}
	(F_{{\rm ext},1}, F_{{\rm ext},2})=(F^\nu_1, F_2^\nu)-(-\partial_Y^2U_s,0).\label{2.2-1}
\end{align}
The equations for $(\rho, {\bf{u})}$ read
\begin{align}\label{1.5}
	\begin{cases}
		U_s\partial_x\rho+\partial_xu+\partial_yv=-\partial_x(\rho u)-\partial_y(\rho v),\\
		U_s\partial_xu+v\partial_yU_s+m^{-2}\partial_x\rho-\nu\Delta u-\lambda\nu\partial_x(\partial_xu+\partial_yv)+\nu\rho\partial_y^2U_s=N_u(\rho,u,v),\\
		U_s\partial_xv+m^{-2}\partial_y\rho-\nu\Delta v-\lambda\nu\partial_y(\partial_xu+\partial_yv)=N_v(\rho,u,v),\\
		u|_{y=0}=v|_{y=0}=0,
	\end{cases}
\end{align}
where  the nonlinear operators $N_u$ and $N_v$ in momentum equations are defined by
\begin{align}
	N_u(\rho,u,v)=&-\partial_x\left[(1+\rho)u^2\right]-U_s\partial_x(\rho u)-\partial_y\left[(1+\rho)uv\right]\nonumber\\
	&-\rho v\partial_yU_s-\left[P^\prime(1+\rho)-P^\prime(1)\right]\partial_x\rho+(1+\rho)F_{{\rm ext},1},\label{N1}\\
	N_v(\rho,u,v)=&-\partial_x\left[(1+\rho)uv\right]-U_s\partial_x(\rho v)-\partial_y\left[(1+\rho)v^2\right]\nonumber\\
	&-\left[P^\prime(1+\rho)-P^\prime(1)\right]\partial_y\rho+(1+\rho)F_{{\rm ext},2}.\label{N2}
\end{align}
In \eqref{1.5}, there is an important dimensionless parameter for compressible fluid,  the Mach number, which is defined as $m\eqdef\frac{1}{\sqrt{P'(1)}}$. This number represents the ratio of the tangential velocity field to the sound speed at the far field of the boundary layer profile. In this paper, we aim to study the structural stability of boundary layer profile \eqref{BL} in the entire subsonic regime $m\in (0,1)$.

Now we introduce some notations. For any $g\in L^2(\Omega)$, we denote by $\mathcal{P}_n (n\in\mathbb{Z})$ the orthogonal projection on the $n$-th Fourier mode:
\begin{align}%\label{1.3}
(\mathcal{P}_ng)(x, y)=g_n(y)e^{i\hat{n}x},\quad \hat{n}=\frac{n}{L},\quad g_n(y)=\frac{1}{2\pi L}\int_0^{2\pi L} g(x, y)e^{-i\hat{n}x}\;dx.\nonumber
\end{align}
We also denote by $g_{\neq}= (I-\mathcal{P}_0)g$  the non-zero modes of $g$. The solution space is given by
\begin{align}
\mathfrak{X}=\bigg\{ (\rho, u, v)\bigg| &\rho\in H^3(\mathbb{T}_L\times \mathbb{R}_+), v_0\in L^1(\mathbb{R}_+)\cap H^4(\mathbb{R}_+), u_0\in L^{\infty}(\mathbb{R}_+)\cap \dot{H^4}(\mathbb{R}_+),\nonumber\\
&(u_{\neq},v_{\neq})\in H^4(\mathbb{T}_L\times \mathbb{R}_+)^2,~
u|_{y=0}=v|_{y=0}=0,\|(\rho,u,v)\|_{\mathfrak{X}}<\infty\bigg\}.\label{SS}
\end{align}
Here the solution norm $\|\cdot\|_{\mathfrak{X}}$ is defined by
\begin{align}
	\|(\rho,u,v)\|_{\mathfrak{X}}\eqdef & ~\|v_0\|_{L^1_y}+\|(m^{-2}\rho_0,u_0,v_0)\|_{L^\infty_y}+\|(m^{-2}\rho_0,v_0)\|_{L^2_y}+\|(m^{-2}\rho_{\neq},u_{\neq},v_{\neq})\|_{L^2(\Omega)}\nonumber\\
	& +\nu^{\frac12}\|\nabla_{x,y}(m^{-2}\rho,u,v)\|_{L^2(\Omega)}+\nu^{\frac{13}{8}}\|\nabla_{x,y}^2(m^{-2}\rho,u,v)\|_{L^2(\Omega)}\nonumber\\
	&+\nu^{\frac{21}{8}} \|\nabla_{x,y}^3(m^{-2}\rho,u,v)\|_{L^2(\Omega)}+\nu^{\frac{29}{8}} \|\nabla_{x,y}^4(u,v)\|_{L^2(\Omega)}.
	\label{S1}
\end{align}
We also define the weighted norm for the perturbation of external force \eqref{2.2-1}
\begin{align}
	\|(F_{{\rm ext},1},F_{{\rm ext},2})\|_{w}\eqdef  &\|(1+y)^s(F_{{\rm ext},1},F_{{\rm ext},2})\|_{L^2(\Omega)}+\nu^{\frac{13}{8}}\|\nabla_{x,y}(F_{{\rm ext},1},F_{{\rm ext},2})\|_{L^2(\Omega)}\nonumber\\
	&
	+\nu^{\frac{21}{8}}\|\nabla^2_{x,y}(F_{{\rm ext},1},F_{{\rm ext},2})\|_{L^2(\Omega)},\qquad s>2.\label{F}
\end{align}
And for simplicity, we assume that $(F_{\rm ext,1},F_{\rm ext,2})$ has no zero mode: $\mathcal{P}_0F_{\rm ext,1}=\mathcal{P}_0F_{\rm ext,2}=0.$ The general case can be treated by imposing appropriate smallness and decay conditions on $(\mathcal{P}_0F_{\rm ext,1},\mathcal{P}_0F_{\rm ext,2})$.

%Now we are ready to  state main results in the paper.

\begin{theorem}[Structural Stablity]\label{T1.1}
Let the Mach number $m\in (0,1)$. Suppose that $U_s$ satisfies structural conditions \eqref{2.1}-\eqref{2.2}. Then there exists a positive constant $L_0$, such that for any $L\in (0, L_0)$ and $0<\nu\ll1$, if  $\|(F_{{\rm ext},1}, F_{{\rm ext},2})\|_{w}\leq \nu^{\frac{9}{8}+}$, the system \eqref{1.5} admits a unique solution $(\rho, u, v)\in \mathfrak{X}$ satisfying
\begin{align}\label{1.4}
\|(\rho,u,v)\|_{\mathfrak{X}}\leq C\|(F_{{\rm ext},1}, F_{{\rm ext},2})\|_{w}.
\end{align}
Moreover, the solution satisfies the following zero-mass condition
\begin{align}
	\iint_{\Omega}\rho(x,y)dx dy=0.\label{z}
\end{align}
\end{theorem}
\begin{remark}\label{rmk1}
	\begin{itemize}
		\item[(i)] By \eqref{1.4}, $L^\infty$ inequality \eqref{lw} and  interpolation, the following $W^{1,\infty}$ and $W^{2,\infty}$ estimates hold for density $\rho$ and velocity field $(u,v)$ respectively:
		\begin{align}
			\|(m^{-2}\rho,u,v)\|_{L^\infty(\Omega)}+\nu^{\frac{9}{8}}\|\nabla_{x,y}(m^{-2}\rho,u,v)\|_{L^\infty(\Omega)}+\nu^{\frac{17}{8}}\|\nabla_{x,y}^2(u,v)\|_{L^\infty(\Omega)}\leq \nu^{\frac58+}.\label{1.4-1}
		\end{align} 
		See \eqref{Lw1}-\eqref{Lw3} for details.  According to \eqref{F}, we do not require uniform-in-$\nu$ $L^\infty$ bound on the  derivatives of external forces. Thus, the $L^\infty$ bounds of $(\nabla_{x,y}\rho,\nabla_{x,y}u,\nabla_{x,y}v)$ and $(\nabla_{x,y}^2u,\nabla_{x,y}^2v)$ can be of the orders of   $\nu^{-(\frac12)-}$ and $\nu^{-(\frac{3}{2})-}$ respectively.
		\item[(ii)] The subsonic assumption $m\in (0,1)$ is essential for the stability analysis. When $m>1$, the Mack unstable modes arise that  lead to  ill-posedness at the inviscid level. We refer to Masmoudi-Wang-Wu-Zhang \cite{MWWZ241} for the mathematical justification.
		\item [(iii)]	The solution $(\rho,u,v)$ satisfies an additional boundary condition $\div_{x, y}({u},{v})|_{y=0}=0$. This can be derived from the equation $\eqref{1.5}_1$. In fact,  we have the identity
		$(1+\rho)\div_{x,y}(u,v)|_{y=0}=-U_s\partial_x\rho|_{y=0}-u\partial_x\rho|_{y=0}-v\partial_y\rho|_{y=0}=0.$ Then the result holds from $\|\rho\|_{L^\infty}\ll 1$ due to \eqref{1.4-1}. This boundary condition is crucial for establishing higher order estimates.
		\item [(iv)] 	The mass constraint \eqref{z} is not only physical, but also crucial for the higher-order estimates, see Proposition \ref{L7.1}. This condition is ensured by the analysis of zero mode in Section \ref{S2}.
		\item[(v)] As one can see from the proof, the length of torus $L$ can be large if the amplitude of boundary layer profile $\|\partial_YU_s\|_{L^\infty}$ is sufficiently small.
	\end{itemize}
\end{remark}
\bigbreak
It is noteworthy that the bound \eqref{1.4} is uniform for the Mach number $m\in [0,m_0],$ for any $m_0\in (0,1)$. This uniform bound enables us to derive the low Mach number limit. Formally, as $m\rightarrow 0^+$, we obtain the following steady incompressible Navier-Stokes equations 
\begin{align}\label{1.5-1}
	\begin{cases}
		\partial_xu^{\rm in}+\partial_yv^{\rm in}=0,\\
		U_s\partial_xu^{\rm in}+v^{\rm in}\partial_yU_s+\partial_xP^{\rm in}-\nu\Delta u^{\rm in}=-u^{\rm in}\partial_xu^{\rm in}-v^{\rm in}\partial_yu^{\rm in}+F_{\rm ext,1},\\
		U_s\partial_xv^{\rm in}+\partial_yP^{\rm in}-\nu\Delta v^{\rm in}=-u^{\rm in}\partial_xv^{\rm in}-v^{\rm in}\partial_yv^{\rm in}+F_{\rm ext,2},\\
		u^{\rm in}|_{y=0}=v^{\rm in}|_{y=0}=0,
	\end{cases}
\end{align}
where $(u^{\rm in},v^{\rm in})$ is the perturbation around the boundary layer profile $(U_s(Y),0)$ and $P^{\rm in}$ is the pressure.  

For any Mach number $m\in (0,1)$, let $(\rho^m,u^m,v^m)$ be the solution to steady compressible Navier-Stokes equations \eqref{1.5}  obtained in Theorem \ref{T1.1}. We define a lower order norm
$$	
\begin{aligned}
\|(\rho^m,u^m,v^m)\|_{2}\eqdef & ~\|v_0^m\|_{L^1_y}+\|(m^{-2}\rho_0^m,u_0^m,v_0^m)\|_{L^\infty_y}+\|(m^{-2}\rho_0^m,v_0^m)\|_{L^2_y}+\|(m^{-2}\rho_{\neq}^m,u_{\neq}^m,v_{\neq}^m)\|_{L^2(\Omega)}\nonumber\\
& +\nu^{\frac12}\|\nabla_{x,y}(m^{-2}\rho^m,u^m,v^m)\|_{L^2(\Omega)}+\nu^{\frac{13}{8}}\|\nabla_{x,y}^2(m^{-2}\rho^m,u^m,v^m)\|_{L^2(\Omega)},
\end{aligned}
$$
where $(\rho^m_0,u^m_0,v^m_0)$ and $(\rho^m_{\neq},u^m_{\neq },v^m_{\neq})$ are zero and non-zero modes of $(\rho^m,u^m,v^m)$ respectively.

\begin{theorem}[Low Mach number limit]\label{T1.2}
Let $m_0\in (0,1)$. There exists a solution $(P^{\rm in},u^{\rm in},v^{\rm in})\in H^3(\Omega)\times \left(L^\infty(\Omega)\cap\dot{H^4}(\Omega)\right)\times H^4(\Omega)$ to the incompressible Navier-Stokes system \eqref{1.5-1}, such that for any $m\in (0,m_0)$, the solution $(\rho^m,u^m,v^m)$ obtained in Theorem \ref{T1.1} satisfies the following bound
	\begin{align}\label{T1.2-1}
		\|(\rho^m-m^2P^{\rm in},u^m-u^{\rm in},v^m-v^{\rm in})\|_2\leq Cm^2\nu^{-\frac{5}{8}-}\|(F_{\rm ext,1},F_{\rm ext,2})\|_{w}.
	\end{align}
In particular,
\begin{align}
	\|(m^{-2}\rho^m-P^{\rm in},u^m-u^{\rm in},v^m-v^{\rm in})\|_{L^\infty(\Omega)}\leq Cm^2\nu^{-\frac98-}\|(F_{\rm ext,1},F_{\rm ext,2})\|_{w}.\label{T1.2-2}
\end{align}
\end{theorem}
\begin{remark}
	\begin{itemize}
		\item[(i)]  In the low Mach number limit $m\rightarrow 0^+$, we can recover a steady solution to the incompressible Navier-Stokes equations. This solution features a boundary layer profile $(1,U_s(Y),0)$ as the leading order term. Therefore, Theorem \ref{T1.2} provides the first result concerning the low Mach number limit in the presence of Prandtl boundary layers for the Navier-Stokes system. 
		\item[(ii)] Solvability of \eqref{1.5-1} in the space $H^1(\Omega)$ has been established by G\'erard-Varet and Maekawa in \cite{GM19}. In our analysis, we require higher regularity to perform the low Mach number limit.  
	\end{itemize}
  
\end{remark}
%\begin{remark}
%	\begin{itemize}
%	\item[(i)] The smallness of $m_0$ does not depend on the viscosity $\nu$.
%	\end{itemize}
%\end{remark}

\subsection{ Strategy of Proof} In Sections \ref{S2}--\ref{S7}, we consider the following linearized Navier-Stokes system in $\Omega$:
\begin{align}\label{1.6}
	\begin{cases}
		U_s\partial_x\rho+\div_{x,y}(u,v)=g_{\rho},\\
		U_s\partial_xu+v\partial_yU_s+m^{-2}\partial_x\rho-\nu\Delta u-\lambda\nu\partial_x\div_{x,y}(u,v)+\rho\nu\partial_y^2U_s=g_u,\\
		U_s\partial_xv+m^{-2}\partial_y\rho-\nu\Delta v-\lambda\nu\partial_y\div_{x,y}(u,v)=g_v,\\
		u|_{y=0}=v|_{y=0}=0,
	\end{cases}
\end{align}
where $(g_{\rho},g_u,g_v)$ is a  given inhomogeneous source term. To take advantage of shear structure of $U_s$, we study \eqref{1.6} at each Fourier mode $(\rho_n,u_n,v_n)$ of $(\rho,u,v)$, which satisfies
\begin{align}\label{1.7}
	\begin{cases}
		i\hat{n} U_s\rho_n+i\hat{n} u_n+\partial_yv_n=g_{\rho, n},\\
		i\hat{n}U_su_n+v_n\partial_yU_s+(i\hat{n} m^{-2}+\nu\partial_y^2U_s)\rho_n-\nu(\partial_y^2-\hat{n}^2) u_n-\lambda\nu i\hat{n}(i\hat{n}u_n+\partial_yv_n)=g_{u, n},\\
		i\hat{n} U_sv_n+m^{-2}\partial_y\rho_n-\nu(\partial_y^2-\hat{n}^2) v_n-\lambda\nu\partial_y(i\hat{n} u_n+\partial_yv_n)=g_{v, n},
	\end{cases}
\end{align}
with no-slip boundary conditions
\begin{align}
	u_n|_{y=0}=v_n|_{y=0}=0. \label{1.7-1}
\end{align}
Here $(g_{\rho, n}, g_{u, n}, g_{v, n})$ is the $n$-th Fourier mode of $(g_\rho, g_u, g_v)$.  We divide the analysis into two cases.\\

\underline{ \it Case 1. Zero mode:} when $n=0$, the ODE system \eqref{1.7} can be solved directly. A crucial aspect of this solution is to ensure the zero mass condition \eqref{z}. Inspired by Esposito et al. \cite{EGKM}, we introduce a parameter $\delta>0$, and incorporate a penalty term into linear system  \eqref{3.14}. Assuming the zero mass condition $\iint_{\Omega} g_{\rho}(x,y)dx dy=0$ for the inhomogeneous source term $g_\rho$, which is also satisfied by the nonlinear term, we can deduce that
$$\int_0^\infty\rho_{0,\delta}(y)dy=0
$$
holds for any $\delta$. By establishing uniform-in-$\delta$ estimate and subsequently  taking the limit as $\delta\rightarrow 0^+$, we obtain the solution $(\rho_0,u_0,v_0)$, which satisfies zero mass condition \eqref{3.13},  cf. Theorem \ref{thmo}.\\

\underline{\it Case 2. Non-zero modes:}  for convenience, we introduce the rescaled variables $(X,Y)=\nu^{-\frac{1}{2}}(x,y)$. Then the linear system \eqref{1.7} reads
\begin{align}\label{4.1}
	\begin{cases}
		i\alpha U_s\rho_n+\div_{\alpha}(u_n,v_n)=f_{\rho, n},\\
		i\alpha U_su_n+v_n\partial_YU_s+i\alpha m^{-2}\rho_n-\sqrt{\nu}\Delta_\alpha u_n-\lambda\sqrt{\nu} i\alpha\div_{\alpha}(u_n,v_n)=f_{u, n},\\
		i\alpha U_s v_n+m^{-2}\partial_Y\rho_n-\sqrt{\nu}\Delta_\alpha v_n-\lambda\sqrt{\nu}\partial_Y\div_{\alpha}(u_n,v_n)=f_{v, n},
	\end{cases}
\end{align}
with no-slip boundary conditions
\begin{align}
	u_{n}|_{Y=0}=v_{n}|_{Y=0}=0.\label{4.1-1}
\end{align}
Here $\alpha=\hat{n}\sqrt{\nu}$ is the rescaled frequency, $\Delta_\alpha=\partial_Y^2-\alpha^2$ and $\div_{\alpha}(u,v)=i\alpha u+\partial_Yv$ are Fourier transform of Laplacian and divergence operator respectively, and $(f_{\rho,n},f_{u,n},f_{v,n})$ is the rescaled homogeneous source term satisfying
\begin{align}\label{gs}
	(f_{\rho,n},f_{u,n},f_{v,n})=\sqrt{\nu}(g_{\rho,n},g_{u,n},g_{v,n}).
\end{align}
For simplicity, we denote \eqref{4.1} by $\mathcal{L}(\rho_n,u_n,v_n)=(f_{\rho,n},f_{u,n},f_{v,n})$.
\bigbreak
For the incompressible flow, the linearized Navier-Stokes system can be reduced into a {\it scalar} Orr-Sommerfeld equation through vorticity. This approach turns out to be powerful  for analyzing stability properties of boundary layers in various settings, as shown in \cite{GM19,GMM18,GMM20,GI,GGN16,IM20,M14}. However, this formulation is not directly applicable for compressible flow. The main obstacle is the strong coupling among density, divergence and vorticity fields of the fluid, which prevents the reduction of the linear system into a single equation analogous to the Orr-Sommerfeld equation.

To address this difficulty, in \cite{YZ23} the second and third authors of this paper introduced a quasi-compressible-Stokes iteration approach for the isentropic fluid, which corresponds to Rayleigh-Airy iteration  \cite{GM19,GMM18,GGN16} used in incompressible fluid. However, it is not trivial to directly extend this approach to the problem in this paper. In the stability analysis, it is essential to consider {\it all frequency ranges}, rather than focusing on a specific regime as  in \cite{YZ23}. 

Now we briefly comment difficulties and state key ideas for the analysis of linear system \eqref{4.1} with no-slip boundary conditions \eqref{4.1-1}.
\bigbreak
\underline{ (i) {\it Quasi-compressible approximation.}}  The first step is to introduce an approximate system, the so-called quasi-compressible system, to \eqref{4.1}. This is achieved by replacing the physical viscosity by an artificial one suitably:
	$$
\underbrace{\left(\begin{array}{ccc}
		0&0&0\\
		0&\sqrt{\nu}\partial_Y^2-(1+\lambda)\sqrt{\nu} \alpha^2&\lambda\sqrt{\nu} i\alpha\partial_{Y}\\
		0&\lambda\sqrt{\nu} i\alpha\partial_Y&(1+\lambda)\sqrt{\nu}\partial_Y^2+\sqrt{\nu}\alpha^2
	\end{array}\right)}_{\text{physical viscosity}}\Rightarrow \underbrace{\left(\begin{array}{ccc}
		0&0&0\\
		\sqrt{\nu}\Delta_\alpha (U_s\cdot)&\sqrt{\nu}\Delta_\alpha&0\\
		0&0&\sqrt{\nu}\Delta_\alpha
	\end{array}\right)}_{\text{artificial viscosity}}.
$$
 We refer to \eqref{4.3} for its precise definition. Surprisingly, fluid variables in this approximate system exhibit a weaker coupling compared to the original Navier-Stokes system \eqref{4.1}. In fact, it
can be reduced to a scalar equation \eqref{4.122-1}, which shares a similar structure as the classical Orr-Sommerfeld equation when  Mach number is less than one (subsonic regime). 
%Consequently, this formulation enables us to apply ``Incompressible techniques''  to effectively analyze stability problems in compressible fluids.
 Section \ref{S3} is dedicated to solving this equation. Due to the similarity between \eqref{4.122-1} with its incompressible counterpart, we can employ the Rayleigh-Airy iteration in the steady setting developed by G\'erard-Varet-Maekawa \cite{GM19}. In the compressible setting, some higher order derivatives appear as commutators in the Airy approximation. Foturnately, these terms are typically accompanied  with suitable multiplied factors that vanish on the boundary. Thus, we can control them by the co-normal bounds of Airy solution; see Remark \ref{rmk3.1}.  Finally, it is worth pointing out that the smallness of $L$ is essential for closing the Rayleigh-Airy iteration; see Proposition \ref{P4.4}.

 After obtaining the solution to Orr-Sommerfeld equation, we  recover the fluid quantities of the inhomogeneous quasi-compressible system in Section \ref{S4}. From Corollary \ref{P4.9}, we can see that  compressible components: $\div_{\alpha}(\mathfrak{u},\mathfrak{v})$  and $\varrho$ have better estimates than the incompressible component. This is due to the absence of sublayer related to these two components. We stress that this stronger estimate is crucial in the proof of convergence of quasi-compressible-Stokes iteration.
\bigbreak
(ii) \underline{\it Quasi-compressible-Stokes iteration.} Even though the quasi-compressible system $\mathcal{L}_Q$ is an approximation to linear Navier-Stokes equations \eqref{4.1} at a formal level, it cannot be directly applied to construct the original linear solution for the following two reasons: (I) $\mathcal{L}_Q$ only works when the inhomogeneous source term $f_{\rho,n}$ in the continuity equation vanishes; (II) it induces a small yet singular error:
$$E_{Q}(\varrho,\mathfrak{u},\mathfrak{v})\eqdef \mathcal{L}(\varrho,\mathfrak{u},\mathfrak{v})-\mathcal{L}_Q(\varrho,\mathfrak{u},\mathfrak{v})=\sqrt{\nu}U_s\Delta_{\alpha}\rho+\text{lower order terms}.
$$
The regularity of $\Delta_{\alpha}\varrho$ is only  $L^2(\mathbb{R}_+)$, which is not enough to solve \eqref{4.1} through iteration. This is because the quasi-compressible system need to be solved at the level of vorticity, and $H^1(\mathbb{R}_+)$ regularity is neccessary.

In Section \ref{S6.1}, we introduce the Stokes regularizing system (denoted by $\mathcal{L}_S$) to handle the inhomogeneous source term $f_{\rho,n}$, and to smooth out the error term generated by the quasi-compressible system. This system is derived by removing the stretching term $v_n\partial_YU_s$ from the equation \eqref{4.1}, and  by replacing the no-slip boundary condition by the slip one.  Unlike the time-dependent case studied in \cite{YZ23}, the Stokes system in the steady setting does not yield a closed basic elliptic estimate:
$$\sqrt{\nu}\|(\partial_Yu,\partial_Yv)\|_{L^2}^2+\alpha\|\sqrt{U_s}(u,v)\|_{L^2}^2\leq \alpha m^{-2}\|\sqrt{U_s}\rho\|_{L^2}^2+\cdots;
$$
as shown in equations \eqref{5.10} and \eqref{5.11}. However, by taking advantage of the slip boundary condition, we can derive the following estimate for the density:
$$m^{-4}\|\rho\|_{L^2}^2\leq \|\sqrt{U_s}(u,v)\|_{L^2}^2+\cdots,
$$
cf. \eqref{rho}. This allows us to close the estimate because $0\leq U_s\leq 1$, and the Mach number $m$ is strictly less than one.

In Section \ref{S6.2}, we apply the quasi-compressible-Stokes iteration to solve the linearized Navier-Stokes system \eqref{4.1}. The procedure is outlined as follows: firstly, we solve the Stokes regularization system \eqref{5.1} that carries  the inhomogenuity in the continuity equation:
$$
\begin{aligned}
&\mathcal{L}_S(\udr,\udu,\udv)=(f_{\rho,n},f_{u,n},f_{v,n}),\\
&\partial_Y\udu|_{Y=0}=\udv|_{Y=0}=0.
\end{aligned}
$$
In the second step, we introduce a solution to following quasi-compressible system 
$$
\begin{aligned}
	&\mathcal{L}_Q(\varrho,\mathfrak{u},\mathfrak{v})=(0,\udv\partial_YU_s,0),\\
	&\mathfrak{v}|_{Y=0}=0.
\end{aligned}
$$
to correct the stretching term generated in Step 1. It is important to note that the error term introduced by $\mathcal{L}_{Q}$ is small in $\nu$ that  can be regularized by the Stokes system. Therefore, by alternately iterating these two steps, we can solve the linear Navier-Stokes system \eqref{4.1}, {\it with only one boundary condition $v_n|_{Y=0}=0$}. Again, we emphasize that the convergence of quasi-compressible-Stokes iteration is shown for all frequencies within the regime $\hat{n}\lesssim \nu^{-\frac{3}{4}}$. \\

(iii) \underline{\it Boundary layer corrections.} The iteration approach is not sufficient for solving the linearized Navier-Stokes system \eqref{4.1} with no-slip boundary condition \eqref{4.1-1}. There are two main difficulties: (I) In the low and middle frequency regime, the quasi-compressible iteration can preserve only one boundary condition $v|_{Y=0}=0.$ (II) In the high frequency regime $\hat{n}\gg \nu^{-\frac34}$, even though the stretching term can be controlled by the tangential diffusion, it is still very hard to directly solve the compressible Navier-Stokes system \eqref{4.1} under  no-slip boundary conditions \eqref{4.1-1}. The main reason is that, in constrast to the incompressible case, where the pressure term can be eliminated due to divergence free condition, here we have to handle the density/pressure term, which is subject to a non-local boundary condition that is not known a priori. A natural attempt is to solve the same system using slip boundary conditions instead of no-slip one, as done in the incompressible case \cite{CLWZ20,GMM20}.

Therefore, to overcome the  difficulties (I) and (II), we need to introduce boundary layer corrections to recover no-slip boundary conditions {\it for all frequencies }. In particular, in the high-frequency regime $n\gg \nu^{-\frac34}$,  it is difficult to  find a non-trivial profile from \eqref{4.1} by the dimensional analysis.

The key idea is to use the quasi-compressible formulation again. Our strategy is summarized as follows.

\begin{itemize}
	\item[] Step 1. We construct a homogeneous solution $\phi_{H}$ to the compressible Orr-Sommerfeld equation
	\begin{equation}
		\begin{aligned}
			&{\rm OS_{CNS}}(\phi_H)=i\eps\Lambda(\Delta_\alpha\phi_{H})+U_s\Lambda(\phi_H)-\phi_{H}\partial_Y(A^{-1}\partial_YU_s)=0,\\
			&\phi_{H}|_{Y=0}=0,~\partial_Y\phi_{H}|_{Y=0}\neq 0.\nonumber
		\end{aligned} 
	\end{equation}
	In the low and middle frequency regime $\hat{n}\lesssim \nu^{-\frac34}$, the structure of the compressible Orr-Sommerfeld equation resembles that of the classical one. This similarity allows us to construct the homogeneous solution $\phi_H$ by using techniques developed for the incompressible case, such as those by G\'erard-Varet-Maekawa \cite{GM19} and Chen-Wu-Zhang \cite{CWZ23}. However, these profiles are not accurate in the high frequency regime $\hat{n}\gg \nu^{-\frac34}$. Instead, in this regime, we construct $\phi_H$ around a specific profile $\phi_{f,app,h}=\frac{Ye^{-\alpha Y}}{2\alpha}$, which is a non-trivial solution to the following equation
	$$\Delta_\alpha^2\phi=0,~\phi|_{Y=0}=0.
	$$
	\item[] Step 2. We obtain a homogeneous solution $(\varrho_{H},\mathfrak{u}_{H},\mathfrak{v}_{H})$ to the quasi-compressible system \eqref{4.151} in terms of $\phi_H$. Note that $\mathfrak{u}_H$ and $\mathfrak{v}_H$ satisfy  boundary conditions:
	$$\mathfrak{v}_H|_{Y=0}=0,~ \mathfrak{u}_{H}|_{Y=0}\neq 0.$$
	\item[] Step 3. We construct a boundary layer profile $(\rho_{n,b},u_{n,b},v_{n,b})$ near the homogeneous quasi-compressible solution $(\varrho_{H},\mathfrak{u}_{H},\mathfrak{v}_{H})$ by stability estimate established in Section \ref{S6}.
\end{itemize}

After obtaining $(\rho_{n,sl},u_{n,sl},v_{n,sl})$ with slip boundary conditions in Section \ref{S6}, and boundary layer corrections $(\rho_{n,b},u_{n,b},v_{n,b})$ in Section \ref{S7}, we have the linear stability result Theorem \ref{thmno} in Section \ref{S7.4}.\\

Boostraping the linear stability to the nonlinear stability is still challenging since the non-linear term $\nabla (\rho\bf{u})$ in the continuity equation causes derivative loss. Inspired by Kagei and Nishida \cite{KN19}, we introduce a modified linear system  \eqref{7.1}. To solve it, we need to establish higher order estimates on solutions. Our strategy is as follows.
\begin{itemize}
	\item[] Step 1. We obtain the estimates on $x$-derivatives of $\div_{x,y}(u,v)$.
	\begin{align}
		&\nu\|\nabla_{x,y}\partial_{x}^2(u,v)\|_{L^2(\Omega)}+\nu\|\partial_x^2\div_{x,y}(u,v)\|_{L^2(\Omega)}\nonumber\\
		&\qquad\leq o(1)m^{-2}\|\nabla_{x,y}^2\rho\|_{L^2(\Omega)}+\text{lower order norms}+\text{source}.\nonumber
	\end{align}

	\item[] Step 2. We estimate $\|\Delta_{x,y}\rho\|_{L^2}$ and full derivatives of divengence field $\|\nabla_{x,y}^2\div_{x,y}(u,v)\|_{L^2(\Omega)}$ in terms of
	 $\|\partial_x^2\div_{x,y}(u,v)\|_{L^2(\Omega)}$. 
	 From the momentum equations \eqref{7.1}, we get
	$$m^{-2}\Delta_{x,y}\rho-\nu(1+\lambda)\Delta_{x,y}\div(u,v)=\text{lower order norms}+\text{source}.
	$$
	The key point is to use continuity equation to decouple $\Delta_{x,y}\rho$ and $\Delta_{x,y}\div(u,v)$. In fact, taking $\partial_y^2$ to $\eqref{7.1}_1$ yields
	$$\partial_y^2\div_{x,y}(u,v)=-(U_s+\tilde{u})\partial_{x}^2\partial_y\rho-\tilde{v}\partial_y^3\rho+\cdots.
	$$
	Combining these two equations, we deduce
	$$m^{-2}\Delta_{x,y}\rho=-\nu(1+\lambda)\partial_x^2\div_{x,y}(u,v)-\nu(1+\lambda)(U_s+\tilde{u})\partial_{x}^2\partial_y\rho-\tilde{v}\partial_y^3\rho+\cdots,
	$$
	where the first term has been bounded in Step 1, and high order derivatives in the second and third terms vanish after integrating by parts because of the  no-slip boundary conditions $\tilde{u}|_{y=0}=\tilde{v}|_{y=0}=0$.  Bounds on full derivatives $\|\nabla_{x,y}^2\div_{x,y}(u,v)\|_{L^2(\Omega)}$ are derived by using the additional boundary condition $\div_{x,y}(u,v)|_{y=0}=0$, cf.  Remark \ref{rmk1}(iv)  and $W^{2,2}$-estimates for Laplacian.\\
	
	\item[] Step 3. To bound full derivatives of $(\rho,u,v)$, we rewrite the momentum equations in \eqref{7.1} into a Stokes equation, with source term depending only on highest derivative of divergence field and some lower order terms, cf. \eqref{St}. We use a cut-off technique to divide the fluid domain $\Omega$ into two parts: (I) In a bounded domain, the higher-order estimate for Stokes equations developed in monographs \cite{GA,SH} can be applied; (II) For the remainded part where the vorticity is away from $y=0$,  elliptic estimates in the whole space can be applied, cf.  Lemma \ref{lem9.1} for  details.
\end{itemize}

Finally, we prove main theorems \ref{T1.1} about the nonlinear stability and \ref{T1.2} about the low Mach number limit in Sections \ref{S8.3} and \ref{S8.4} respectively.
\bigbreak
{\bf{Notations:}} For any function $f\in L^1(\mathbb{R}_+)$, we define  $\mathcal{I}[f]:=\int_0^y f(y^\prime)\;dy^\prime$ and  $\partial_y^{-1}f:=-\int_y^\infty f(y^\prime)\;dy^\prime$.  We use the notation $A\lesssim B$ to indicate that $A\le CB$ for some generic constant $C>0$, which is independent of both the viscosity $\nu$ and the Mach number $m$.  The notation $a+$ is used to denote any number greater than $a$, while $a-$ denotes any number less than $a$. The standard $L^p$-norm with respect to variable $y$ is denoted by $\|\cdot\|_{L^p_y}$.

\section{Linear stability of zero mode}\label{S2}
When $n=0$, the linear problem \eqref{1.7} is reduced to the following  ODE system.
\begin{align}\label{3.1}
\begin{cases}
\partial_yv_0=g_{\rho, 0},\\
v_0\partial_yU_s+\nu\rho_0\partial_y^2U_s-\nu\partial_y^2u_0=g_{u, 0},\\
m^{-2}\partial_y\rho_0-(1+\lambda)\nu\partial_y^2v_0=g_{v, 0},\\
u_0|_{y=0}=v_0|_{y=0}=0,
\end{cases}
\end{align}
where $(g_{\rho,0},g_{u,0},g_{v,0})$ is the zero mode of $(g_{\rho},g_u,g_v)$. Solvability of \eqref{3.1} is given in the following theorem.

\begin{theorem}\label{thmo}
Let $g_{\rho, 0}\in H^1(\mathbb{R}_+)\cap L^1(\mathbb{R}_+)\cap L^\infty(\mathbb{R}_+)$, $\mathcal{I}(g_{\rho,0})$ and  $\partial^{-1}_yg_{v,0}\in L^1(\mathbb{R}_+)\cap L^\infty(\mathbb{R}_+)$, $\partial_y^{-1}g_{u, 0}\in L^1(\mathbb{R}_+)\cap L^2(\mathbb{R}_+)$, and $(g_{u, 0}, g_{v, 0})\in L^2(\mathbb{R}_+)^2$. There exists a unique solution $(\rho_0, u_0, v_0)$ to \eqref{3.1} satisfying
\begin{align}
	\|v_0\|_{L^p_y}&\le C\|\mathcal{I}[g_{\rho, 0}]\|_{L^p_y},~~\forall p\in [1,\infty]\label{3.5}\\ 
	\|\partial_yv_0\|_{L^2_y}&\le C\|g_{\rho, 0}\|_{L^2_y},\label{3.6}\\
	\|\partial_y^2v_0\|_{L^2_y}&\le C\|\partial_yg_{\rho, 0}\|_{L^2_y},\label{3.7}\\
m^{-2}\|\rho_0\|_{L^p_y}&\le C\nu\|g_{\rho, 0}\|_{L^p_y}+C\|\partial_y^{-1}g_{v, 0}\|_{L^p_y},~~\forall p\in [1,\infty],\label{3.8}\\ 
m^{-2}\|\partial_g\rho_0\|_{L^2_y}&\le C{\nu}\|\partial_yg_{\rho, 0}\|_{L^2_y}+C\|g_{v, 0}\|_{L^2_y},\label{3.9}\\
\|u_0\|_{L^\infty_y}&\le C\nu^{-\frac12}\|\mathcal{I}[g_{\rho, 0}]\|_{L^\infty_y}+C{\nu}\|g_{\rho, 0}\|_{L^\infty_y}\nonumber\\
&\qquad+C\nu^{-1}\|\partial_y^{-1}g_{u, 0}\|_{L^1_y}+C\|\partial_y^{-1}g_{v, 0}\|_{L^\infty_y},\label{3.10}\\
\|\partial_yu_0\|_{L^2_y}&\le C\nu^{-\frac34}\|\mathcal{I}[g_{\rho, 0}]\|_{L^\infty_y}+C{\nu}^{\frac34}\|g_{\rho, 0}\|_{L^\infty_y}\nonumber\\
&\qquad+C\nu^{-1}\|\partial_y^{-1}g_{u, 0}\|_{L^2_y}+C\nu^{-\frac14}\|\partial_y^{-1}g_{v, 0}\|_{L^\infty_y},\label{3.11}\\
\|\partial_y^2u_0\|_{L^2_y}&\le C\nu^{-\frac54}\|\mathcal{I}[g_{\rho, 0}]\|_{L^\infty_y}+C{\nu}^{\frac14}\|g_{\rho, 0}\|_{L^\infty_y}\nonumber\\
&\qquad+C\nu^{-1}\|g_{u, 0}\|_{L^2_y}+C\nu^{-\frac34}\|\partial_y^{-1}g_{v, 0}\|_{L^\infty_y}.\label{3.12}
\end{align}
If in addition $\int_{0}^\infty g_{\rho,0}(y)d y=0$, we have
\begin{align}
	\int_{0}^\infty \rho_{0}(y)d y=0.\label{3.13}
\end{align}
\end{theorem}
\begin{remark}
	\begin{itemize}
		\item[(i)] Different from the incompressible case studied in \cite{GM19},  $v_0$ does not vanish. The $L^1$-estimate of $v_0$ is crucial for closing the nonlinear estimates.
		\item[(ii)] The zero mass condition \eqref{3.13} is crucial for obtaining higher-order estimates, cf. Proposition \ref{L7.1}.
	\end{itemize}	
\end{remark}

\begin{proof} Inspired by Esposito et al. \cite{EGKM},
	 we introduce a penalty term in the first equation of \eqref{3.1} to ensure the zero mass condition \eqref{3.13},
	\begin{align}\label{3.14}
		\begin{cases}
			\delta\rho_{0,\delta}+\partial_yv_{0,\delta}=g_{\rho, 0},\\
			v_{0,\delta}\partial_yU_s+\nu\rho_{0,\delta}\partial_y^2U_s-\nu\partial_y^2u_{0,\delta}=g_{u, 0},\\
			m^{-2}\partial_y\rho_{0,\delta}-(1+\lambda){\nu}\partial_y^2v_{0,\delta}=g_{v, 0},\\
			u_{0,\delta}|_{y=0}=v_{0,\delta}|_{y=0}=0,
		\end{cases}
	\end{align}
where $\delta\in (0,1)$ is an arbitrary constant. The system \eqref{3.14} can be solved explicitly as follows. Integrating $\eqref{3.14}_3$ yields
\begin{align}\label{3.2}
m^{-2}\rho_{0,\delta}-(1+\lambda)\nu\partial_yv_{0,\delta}=\partial_y^{-1}g_{v,0}.
\end{align}
Then by multiplying $\eqref{3.14}_1$ by $(1+\lambda)\nu$ and adding it to \eqref{3.2}, we solve $\rho_{0,\delta}$ as follows.
\begin{align}
	m^{-2}\rho_{0,\delta}=\frac{1}{1+\delta m^2(1+\lambda)\nu}\left[\partial_y^{-1}g_{v,0}+(1+\lambda)\nu g_{\rho,0}\right].\label{3.4}
\end{align}
The formula \eqref{3.4} yields following uniform-in-$\delta$ estimate
\begin{align}
	m^{-2}\|\rho_{0,\delta}\|_{L^p_y}&\leq \|\partial_{y}^{-1}g_{v,0}\|_{L^p_y}+C\nu\|g_{\rho,0}\|_{L^p_y},\label{3.15}\\
	m^{-2}\|\partial_y\rho_{0,\delta}\|_{L^2_y}&\leq \|g_{v,0}\|_{L^2_y}+C{\nu}\|\partial_yg_{\rho,0}\|_{L^2_y}.\label{3.16}
\end{align}
For $v_{0,\delta}$, multplying $\eqref{3.14}_1$ and \eqref{3.2} by $m^{-2}$ and $-\delta$ respectively, and then integrating  the sum from $0$ to $y$, we obtain %can solve $v_{0,\delta}$ as follows:
\begin{align}
	v_{0,\delta}=\frac{1}{1+\delta m^{2}(1+\lambda)\nu}\left[\mathcal{I}(g_{\rho,0})-\delta m^2\mathcal{I}(\partial_y^{-1}g_{v,0})\right].\label{3.17}
\end{align}
From \eqref{3.17}, we deduce that $v_{0,\delta}\in L^p_{loc}(\mathbb{R}_+)\cap \dot{H}^2(\mathbb{R}_+)$ (Note that the operator $\mathcal{I}(\cdot)$ introduces a growth in $y$).
Hence, % From \eqref{3.17}, 
we derive following bounds on $v_{0,\delta}$:
\begin{align}
	\|\partial_yv_{0,\delta}\|_{L^2_y}&\leq \|g_{\rho,0}\|_{L^2_y}+\delta\|\partial_y^{-1}g_{v,0}\|_{L^2_y},\label{3.18}\\
	\|\partial_y^2v_{0,\delta}\|_{L^2_y}&\leq \|\partial_yg_{\rho,0}\|_{L^2_y}+\delta\|g_{v,0}\|_{L^2_y},\label{3.19}\\
	\|\frac{v_{0,\delta}}{1+y}\|_{L^\infty_y}&\leq C\|\mathcal{I}(g_{\rho,0})\|_{L^\infty_y}+C\delta\|\partial_y^{-1}g_{v,0}\|_{L^\infty_y}\label{3.20},\\
	\|v_{0,\delta}\mathbf{1}_{0\leq y\leq N}\|_{L^p_y}&\leq \|\mathcal{I}(g_{\rho,0})\|_{L^p_y}+\delta N^{\frac1p}\|\partial_y^{-1}g_{v,0}\|_{L^\infty_y}.\label{3.21}
\end{align}

Next we estimate $u_{0,\delta}$. Integrating $\eqref{3.14}_2$ yields
%,  we get the following explicit formula 
\begin{align}
	u_{0,\delta}&=\frac{1}{\nu}\int_0^y\;dy^\prime\int_{y^\prime}^{+\infty}(g_{u, 0}-v_{0,\delta}\partial_yU_s-{\nu}\rho_{0,\delta}\partial_y^2U_s)(y^{\prime\prime})\;dy^{\prime\prime}.\nonumber
\end{align}
By boundary layer structure and decay properties \eqref{2.2} of $U_s$, we deduce the following bounds on $u_{0,\delta}$:
\begin{align}
	\|u_{0,\delta}\|_{L^\infty_y}&\leq C\nu^{-1}\left(\|\partial_{y}^{-1}g_{u,0}\|_{L^1_y}+\nu^{\frac12}\|\frac{v_{0,\delta}}{1+y}\|_{L^\infty_y}+{\nu}\|\rho_{0,\delta}\|_{L^\infty_y}\right)\nonumber\\
	&\leq C\nu^{-1}\left(\|\partial_{y}^{-1}g_{u,0}\|_{L^1_y}+\nu^{\frac12}\|\mathcal{I}(g_{\rho,0})\|_{L^\infty_y}+\nu\|\partial_y^{-1}g_{v,0}\|_{L^\infty_y}+{\nu}^2\|g_{\rho,0}\|_{L^\infty_y}+\delta\nu^{\frac12}\|\partial_y^{-1}g_{v,0}\|_{L^\infty_y}\right),\label{3.22}\\
	\|\partial_yu_{0,\delta}\|_{L^2_y}&\leq C\nu^{-1}\left( \|\partial_{y}^{-1}g_{u,0}\|_{L^2_y}+\nu^{\frac14}\|\frac{v_{0,\delta}}{1+y}\|_{L^\infty_y}+\nu^{\frac34}\|\rho_{0,\delta}\|_{L^\infty_y} \right)\nonumber\\
	&\leq C\nu^{-1}\left(\|\partial_{y}^{-1}g_{u,0}\|_{L^2_y}+\nu^{\frac14}\|\mathcal{I}(g_{\rho,0})\|_{L^\infty_y}+\nu^{\frac34}\|\partial_y^{-1}g_{v,0}\|_{L^\infty_y}+{\nu^{\frac74}}\|g_{\rho,0}\|_{L^\infty_y}+\delta\nu^{\frac14}\|\partial_y^{-1}g_{v,0}\|_{L^\infty_y}\right),\label{3.23}\\
	\|\partial_y^2u_{0,\delta}\|_{L^2_y}&\leq C\nu^{-1}\left( \|g_{u,0}\|_{L^2_y}+\nu^{-\frac14}\|\frac{v_{0,\delta}}{1+y}\|_{L^\infty_y}+\nu^{\frac14}\|\rho_{0,\delta}\|_{L^\infty_y} \right)\nonumber\\
	&\leq C\nu^{-1}\left(\|g_{u,0}\|_{L^2_y}+\nu^{-\frac14}\|\mathcal{I}(g_{\rho,0})\|_{L^\infty_y}+\nu^{\frac14}\|\partial_y^{-1}g_{v,0}\|_{L^\infty_y}+{\nu}^{\frac54}\|g_{\rho,0}\|_{L^\infty_y}+\delta\nu^{-\frac14}\|\partial_y^{-1}g_{v,0}\|_{L^\infty_y}\right).\label{3.24}
\end{align}

Based on the uniform-in-$\delta$ estimates \eqref{3.15}, \eqref{3.16}, \eqref{3.18}-\eqref{3.21}, \eqref{3.22}-\eqref{3.24}, we can take the limit $
\delta\rightarrow 0^+$.  It is straightforward to see that the limit function $(\rho_0,u_0,v_0)=\lim_{\delta\rightarrow 0^+}(\rho_{0,\delta},u_{0,\delta},v_{0,\delta})$  solves linear system \eqref{3.1}, and belongs in $\left(H^1(\mathbb{R}_+)\cap L^p(\mathbb{R}_+)\right)\times \left( \dot{H}^2(\mathbb{R}_+)\cap L^\infty(\mathbb{R}_+)\right)\times \left( \dot{H}^2(\mathbb{R}_+)\cap L^p_{loc}(\mathbb{R}_+)\right)$. Moreover, it satisfies the bounds \eqref{3.6}-\eqref{3.12}.  Furthermore, taking $\delta\rightarrow 0^+$ in \eqref{3.21} gives
 $$\|v_0\mathbf{1}_{0\leq y\leq N}\|_{L^p_y}\leq \|\mathcal{I}(g_{\rho,0})\|_{L^p_y}.
 $$
Thus taking $N\rightarrow \infty$ we deduce that $v_0$ belongs to $L^p(\mathbb{R}_+)$ and satisfies \eqref{3.5}. 

Finally,  by integrating $\eqref{3.14}_1$ over $\mathbb{R}_+$ and using boundary condition $v_{0,\delta}|_{y=0}=0$, we obtain
\begin{align}
	\int_{0}^\infty \rho_{0,\delta}(y)d y=\frac{1}{\delta}\int_{0}^\infty g_{\delta,0}(y)d y=0.\label{3.25}
\end{align}
The zero mass condition \eqref{3.13} is obtained by taking $\delta\rightarrow 0^+$ in \eqref{3.25}. The proof of Theorem \ref{thmo} is complete. 
\end{proof}

\section{Compressible Orr-Sommerfeld equation}\label{S3}
\subsection{Derivation of an Orr-Sommerfeld approximation}\label{S3.1}
In Section \ref{S3}--\ref{S7}, we will solve the linearized Navier-Stokes system \eqref{4.1} with no-slip boundary condition \eqref{4.1-1}, The main difficulty arises from the stretching term $v\p_YU_s$, and strong coupling among density, vorticity, and divergence field. For this, we introduce the following quasi-compressible approximate system
\begin{align}\label{4.3}
\begin{cases}
i\alpha U_s\varrho+i\alpha\mathfrak{u}+\partial_Y\mathfrak{v}=0,\\
-\sqrt{\nu}\Delta_\alpha(U_s\varrho+\mathfrak{u})+i\alpha U_s\mathfrak{u}+\mathfrak{v}\partial_YU_s+i\alpha m^{-2}\varrho=f_{\mathfrak{u}},\\
-\sqrt{\nu}\Delta_\alpha\mathfrak{v}+i\alpha U_s\mathfrak{v}+m^{-2}\partial_Y\varrho=f_{\mathfrak{v}}.
\end{cases}
\end{align}
We stress that the inhomogeneity appears only in momentum equations.  For simplicity, we denote \eqref{4.3} by $\mathcal{L}_Q(\varrho, \mfu, \mfv)=(0, f_{\mfu}, f_{\mfv})$. 

The key feature of \eqref{4.3} is that  associated fluid quantities can be decoupled. Due to the continuity equation $\eqref{4.3}_1$, we can define an {\it effective stream function} $\phi$ which satisfies
\begin{align}\label{4.4}
\mathfrak{u}=\partial_Y\phi-U_s\varrho,~~~~ \mathfrak{v}=-i\alpha\phi.
\end{align}
Then by  \eqref{4.3}$_2$, we can express $\varrho$ in terms of $\phi$ as follows
\begin{align}\label{4.5}
m^{-2}\varrho=-A^{-1}(Y)\left(i\eps\Delta_\alpha\partial_Y\phi+U_s\partial_Y\phi-\phi\partial_YU_s+{i}{\alpha^{-1}}f_{\mathfrak{u}}\right).
\end{align}
where $A(Y)=1-m^2U_s^2$. Note that if the boundary layer profile is uniformly subsonic, that is, $m\in (0,1)$, then the function $A$ has a strictly positive lower bound: $A(Y)\geq 1-m^2>0$.

Substituting \eqref{4.5} into the third equation of \eqref{4.3}, we can derive the following equation for $\phi$ that can be viewed as the Orr-Sommerfeld equation in the compressible setting:
\begin{align}\label{4.6}
\text{OS}_{\text{CNS}}(\phi)\eqdef i\eps\Lambda(\Delta_\alpha\phi)+U_s\Lambda(\phi)-\phi\partial_Y(A^{-1}\partial_YU_s)=-f_{\mathfrak{v}}-\frac{i}{\alpha}\partial_Y(A^{-1}f_{\mathfrak{u}}),~~~Y>0.
\end{align}
Here,  we define the modified vorticity operator $\Lambda(\cdot)$ as
\begin{align}\label{vo} \Lambda(\phi)\eqdef\partial_Y(A^{-1}\partial_Y\phi)-\alpha^2\phi,
\end{align}
and denote
\begin{align}\nonumber
	\eps=1/\hat{n}.
\end{align}
Note that if Mach number $m=0$, then $\Lambda=\Delta_\alpha$, and  ${\rm OS_{CNS}}$ reduces to the classical Orr-Sommerfeld operator for the incompressible fluid. Upon obtaining solutions $\phi$ to \eqref{4.6}, the fluid variables $(\varrho, \mathfrak{u}, \mathfrak{v})$ can be subsequently determined by substituting $\phi$ into \eqref{4.4} and \eqref{4.5}.

In next two subsections, we will solve the Orr-Sommerfeld equation
\begin{align}\label{4.122-1}
	\begin{cases}
		{\rm{OS}_{CNS}}(\phi)=f,~Y>0,\\
		\phi|_{Y=0}=0,
	\end{cases}
\end{align}
with a general source term $f$. For convenience, we firstly study the following symmetrized version of \eqref{4.122-1}:
\begin{align}\label{4.8}
	\begin{cases}
		\widetilde{\text{OS}}_{\text{CNS}}(\phi)\eqdef {i}\eps\Delta_\alpha\Lambda(\phi)+U_s\Lambda(\phi)-\phi\partial_Y(A^{-1}\partial_YU_s)= f,\quad Y>0,\\
		\phi|_{Y=0}=0.
	\end{cases}
\end{align}
In the  low and middfle frequency regimes when $\hat{n}\lesssim \nu^{-\frac34}$, the equation \eqref{4.8} has a similar structure as the classical Orr-Sommerfeld equation asscciated with the incompressible Navier-Stokes system. In the stationary setting, this equation has been studied by G\'erard-Verat and Maekawa \cite{GM19}. In this section, we will extend their approach to address the compressible Orr-Sommerfeld \eqref{4.8}. For high frequencies when $\hat{n}\gg \nu^{-\frac34}$, we solve \eqref{4.122-1} by the energy method.

Henceforth, we denote the $L^2$-inner product on $\mathbb{R}_+$ with respect to $dY$ by $\langle\cdot, \cdot\rangle$ , and the associated $L^2$ norm by $\|\cdot\|_{L^2}$.

\subsection{Some estimates on Rayleigh and Airy solutions}\label{S3.2}

Note that $\eps=L/{n}$ which is small if $L\ll1$. Consequently, the first step in the analysis is to study the following Rayleigh equation for compressible fluid. This equation is derived from equation \eqref{4.8} by neglecting the high-order diffusion term.
\begin{align}\label{4.9}
\begin{cases}
\text{Ray}(\varphi)\eqdef U_s\Lambda(\varphi)-\varphi\partial_Y(A^{-1}\partial_YU_s)=h,\quad Y>0,\\
\varphi(Y)|_{Y=0}=0.
\end{cases}
\end{align}
%In order to get good estimates of the above equation for $0<\alpha\le 1$, we rewrite the Rayleigh operator in the following form
%\begin{align}\label{4.10}
%Ray(\varphi)&=U_s\partial_Y\left\{A^{-1}\partial_Y\left(\frac{\varphi}{U_s}\cdot U_s\right)\right\}-\partial_Y(A^{-1}\partial_YU_s)\varphi-\alpha^2U_s\varphi\notag\\
%&=U_s\partial_Y\left\{A^{-1}\partial_Y\left(\frac{\varphi}{U_s}\right)U_s+A^{-1}\frac{\varphi}{U_s}\partial_YU_s\right\}-\partial_Y(A^{-1}\partial_YU_s)\varphi-\alpha^2U_s\varphi\notag\\
%&=U_s^2\partial_Y\left\{A^{-1}\partial_Y\left(\frac{\varphi}{U_s}\right)\right\}+2U_sA^{-1}\partial_Y\left(\frac{\varphi}{U_s}\right)\partial_YU_s-\alpha^2U_s\varphi\notag\\
%&=\partial_Y\left\{A^{-1}U_s^2\partial_Y\left(\frac{\varphi}{U_s}\right)\right\}-\alpha^2U_s\varphi,
%\end{align}
%and consider the system as follows
%\begin{align}\label{4.11}
%\begin{cases}
%\partial_Y\left\{A^{-1}U_s^2\partial_Y(\frac{\varphi}{U_s})\right\}-\alpha^2U_s\varphi=h,\quad Y>0,\\
%\varphi(Y)|_{Y=0}=0.
%\end{cases}
%\end{align}
\begin{lemma}
\label{P4.1} Suppose that $m\in (0,1)$. Let $h/U_s \in L^2(\mathbb{R}_+)$.  Then there exists a unique solution $\varphi\in H^2(\mathbb{R}_+)\cap H_0^1(\mathbb{R}_+)$ to \eqref{4.9}. Moreover, the following two statements hold.
\begin{itemize}
\item[{\rm(i)}] If $\alpha\ge1$,
\begin{align}\label{4.12}
\|\partial_Y\varphi\|_{L^2}+\alpha\|\varphi\|_{L^2}&\le C\min\left\{\|\frac{Yh}{U_s}\|_{L^2},\frac{1}{\alpha}\|\frac{h}{U_s}\|_{L^2}\right\},\\
\|\Delta_\alpha\varphi\|_{L^2}&\le C\min\left\{\|\frac{Yh}{U_s}\|_{L^2},\frac{1}{\alpha}\|\frac{h}{U_s}\|_{L^2}\right\}+C\|\frac{h}{U_s}\|_{L^2}.\label{4.13}
\end{align}
\item[{\rm(ii)}] If $0\le\alpha\le 1$ and $(1+Y)\partial_Y^{-1}h\in L^2(\mathbb{R}_+)$,
\begin{align}\label{4.14}
\alpha\|\varphi\|_{L^2}&\le C\alpha\|(1+Y)\partial_Y^{-1}h\|_{L^2}+\frac{C}{\alpha^{\frac12}}|\int_0^\infty h(Y)\;dY|,\\
\|\partial_Y\varphi\|_{L^2}&\le C\left(\|(1+Y)\partial_Y^{-1}h\|_{L^2}+\|h\|_{L^2}\right)+\frac{C}{\alpha}|\int_0^\infty h(Y)\;dY|,\label{4.15}\\
\|\Delta_\alpha\varphi\|_{L^2}&\le C\left(\|(1+Y)\partial_Y^{-1}h\|_{L^2}+\|\frac{h}{U_s}\|_{L^2}\right)+\frac{C}{\alpha}|\int_0^\infty h(Y)\;dY|.\label{4.16}
\end{align}
\end{itemize}
\end{lemma}
\begin{proof}
	If $m\in (0,1)$,  $A(Y)>0$. Consequently, the Rayleigh equation \eqref{4.9} has the similar structure as its counterpart for the incompressible flow studied in \cite[Proposition 5.1]{GM19}. The proof follows a similar approach and we omit the details for brevity.
\end{proof}

To capture the viscous effects, we solve the following compressible Airy equation:
\begin{align}\label{4.51}
\begin{cases}
\widetilde{\text{Airy}}(\psi)\eqdef i\epsilon\Lambda(\psi)+U_s\psi=h,\quad Y>0,\\
\psi(Y)|_{Y=0}=0,
\end{cases}
\end{align}
where the operator $\Lambda(\cdot)$ is defined in \eqref{vo}.
\begin{lemma}
\label{P4.2}
Suppose that $m\in (0,1)$ and $\alpha\eps^{\frac13}\lesssim 1$. Let $h\in L^2(\mathbb{R}_+)$. There exists a unique solution $\psi\in H^2(\mathbb{R}_+)\cap H_0^1(\mathbb{R}_+)$ to the system \eqref{4.51} satisfying
\begin{align}\label{4.52}
\|U_s\psi\|_{L^2}+\epsilon^{\frac16}\|\sqrt{U_s}\psi\|_{L^2}+\epsilon^{\frac13}\|\psi\|_{L^2}+\epsilon^{\frac23}(\|\partial_Y\psi\|_{L^2}+\alpha\|\psi\|_{L^2})+\epsilon\|\Delta_\alpha\psi\|_{L^2}\le C\|h\|_{L^2}.
\end{align}
If $h$ takes form of  $\partial_Yg$ or $g/Y$, 
\begin{align}\label{4.53}
\eps^{\frac13}\|U_s\psi\|_{L^2}+\epsilon^{\frac12}\|\sqrt{U_s}\psi\|_{L^2}+\epsilon^{\frac23}\|\psi\|_{L^2}+\epsilon(\|\partial_Y\psi\|_{L^2}+\alpha\|\psi\|_{L^2})\le C\|g\|_{L^2}.
\end{align}
Moreover, if $h\in H^1(\mathbb{R}_+)$, $(1+Y)^2h\in L^2(\mathbb{R}_+)$ and $(1+Y)^2\partial_Yh\in L^2(\mathbb{R}_+)$,  the solution satisfies following weighted estimates
\begin{align}\label{Ai1}
\|U_sY\psi\|_{L^2}+\epsilon^{\frac13}\|Y\psi\|_{L^2}+\epsilon^{\frac23}\|Y\partial_Y\psi\|_{L^2}+\epsilon\|Y\partial_Y^2\psi\|_{L^2}
&\le C\|Yh\|_{L^2}+C\epsilon^{\frac13}\|h\|_{L^2},\\
\|U_sY\partial_Y\psi\|_{L^2}+\epsilon^{\frac23}\|Y\partial_Y^2\psi\|_{L^2}+\eps\|Y\p_Y^3\psi\|_{L^2}&\le  C\|Y\partial_Yh\|_{L^2}+C\|h\|_{L^2},\label{Ai3}\\
\|U_sY^2\partial_Y\psi\|_{L^2}+\eps^{\frac23}\|Y^2\partial_Y^2\psi\|_{L^2}+\epsilon\|Y^2\partial_Y^3\psi\|_{L^2}&\le C\|Y^2\partial_Yh\|_{L^2}+C\|Yh\|_{L^2}+C\epsilon^{\frac13}\|h\|_{L^2}.\label{Ai4}
\end{align}
If in addition $h=g/Y$, 
\begin{align}\label{Ai5}
\|U_sY\psi\|_{L^2}+\epsilon^{\frac13}\|Y\psi\|_{L^2}+\epsilon^{\frac23}\|Y\partial_Y\psi\|_{L^2}+\epsilon\|Y\partial_Y^2\psi\|_{L^2}
&\le C\|g\|_{L^2},\\
\|U_sY^2\partial_Y\psi\|_{L^2}+\eps^{\frac13}\|Y^2\p_Y\psi\|_{L^2}+\eps^{\frac{2}{3}}\|Y^2\p_Y^2\psi\|_{L^2}+\epsilon\|Y^2\partial_Y^3\psi\|_{L^2}
&\le C\|Y\partial_Yg\|_{L^2}+C\|g\|_{L^2}.\label{Ai6}
\end{align}
\end{lemma}
\begin{remark}\label{rmk3.1}
	Note that 
	$$	
	\begin{aligned}
		\widetilde{\rm  OS}_{\rm CNS}(\psi)&=\Delta_{\alpha}\widetilde{\rm Airy}(\psi)+U_s\Lambda(\psi)-\Delta_\alpha(U_s\psi)-\partial_{Y}(A^{-1}\partial_YU_s)\psi\nonumber\\
		&=\Delta_{\alpha}\widetilde{\rm Airy}(\psi)+U_s(A^{-1}-1)\partial_Y^2\psi+\text{ lower order terms}.
		\end{aligned}
	$$
If $m>0$, then $A^{-1}-1\sim U_s^2$, and the error introduced by $\widetilde{\rm Airy}$ approximation involves the second order derivative (co-normal one) of $\phi$ . Consequently, it is necessary to derive estimates  \eqref{Ai3}-\eqref{Ai6} on conormal derivatives of $\psi$.
\end{remark}
\begin{proof}
If $m\in (0,1)$, then  $\Lambda$ is an elliptic operator in divergence form. Therefore, the estimates \eqref{4.52} and \eqref{4.53} can be derived by the same arguments as those presented in \cite[Proposition 6.1]{GM19}. We only show  the co-normal estimates \eqref{Ai1}-\eqref{Ai6}. Note that $Y\psi$ satisfies
\begin{align*}
\widetilde{\text{Airy}}(Y\psi)=Yh+2i\epsilon A^{-1}\partial_Y\psi+i\epsilon\partial_Y(A^{-1})\psi.
\end{align*}
Applying \eqref{4.52} to $Y\psi$, then using 
$\|Y\p_Y\psi\|_{L^2}\leq \|\p_Y(Y\psi)\|_{L^2}+\|\psi\|_{L^2}$ and $\|Y\p^2_Y\psi\|_{L^2}\leq\|\Delta_\alpha(Y\psi)\|_{L^2}+2\|\p_Y\psi\|_{L^2}+\eps^{-\frac23}\|Y\psi\|_{L^2}$, we obtain
\begin{align}
&\|U_sY\psi\|_{L^2}+\epsilon^{\frac13}\|Y\psi\|_{L^2}+\epsilon^{\frac23}\|Y\partial_Y\psi\|_{L^2}+\epsilon\|Y\partial_Y^2\psi\|_{L^2}\nonumber\\
&\qquad\lesssim \|U_sY\psi\|_{L^2}+\eps^{\frac13}\|Y\psi\|_{L^2}+\eps^{\frac23}\|\partial_Y(Y\psi)\|_{L^2}+\eps\|\Delta_\alpha(Y\psi)\|_{L^2}+\eps^{\frac{2}{3}}\|\psi\|_{L^2}+\eps\|\partial_Y\psi\|_{L^2}\nonumber\\
&\qquad\lesssim \|Yh\|_{L^2}+\eps^{\frac23}\|\psi\|_{L^2}+\eps\|\partial_Y\psi\|_{L^2}\label{Ai7-1}.
\end{align}
For general $L^2$ function $h$, we input the bounds \eqref{4.52} into \eqref{Ai7-1} to get \eqref{Ai1}. If $h=g/Y$, then we input the bounds \eqref{4.53} into \eqref{Ai7-1} to get \eqref{Ai5}.\\ 
%Next, observe that
%\begin{align*}
%mAiry(Y^2\psi)=Y^2h+2i\epsilon (A^{-1}+Y\partial_YA^{-1})\psi+4i\epsilon YA^{-1}\partial_Y\psi.
%\end{align*}
%Similar derivation gives
%\begin{align*}
%\|U_sY^2\psi\|_{L^2}+\epsilon^{\frac13}\|Y^2\psi\|_{L^2}+\epsilon^{\frac23}\|Y^2\partial_Y\psi\|_{L^2}+\epsilon\|Y^2\partial_Y^2\psi\|_{L^2}
%&\lesssim \|Y^2h\|_{L^2}+\|Yh\|_{L^2}+\epsilon^{\frac13}\|h\|_{L^2},\\
%\|U_sY^2\psi\|_{L^2}+\epsilon^{\frac13}\|Y^2\psi\|_{L^2}+\epsilon^{\frac23}\|Y^2\partial_Y\psi\|_{L^2}+\epsilon\|Y^2\partial_Y^2\psi\|_{L^2}
%&\lesssim \|(1+Y)g\|_{L^2},\quad \mbox{if}\,\, h=g/Y.
%\end{align*}
%Similarly, observe that
%\begin{align*}
%mAiry(Y^2\psi)=Y^2h+2i\epsilon YA^{-1}\partial_Y\psi+2i\epsilon\partial_Y(A^{-1}Y\psi).
%\end{align*}
%Therefore, we obtain from \eqref{4.52} and \eqref{Ai7} that
%\begin{align}\label{Ai9}
%\epsilon^{\frac13}\|Y^2\psi\|_{L^2}+\epsilon^{\frac23}\|Y^2\partial_Y\psi\|_{L^2}
%&\lesssim \|Y^2h\|_{L^2}+\epsilon\|Y\partial_Y\psi\|_{L^2}+\epsilon\|\psi\|_{L^2}\notag\\
%&\lesssim \|Y^2h\|_{L^2}+\epsilon^{\frac13}\|Yh\|_{L^2}+\epsilon^{\frac23}\|h\|_{L^2}.
%\end{align}

Next we show higher order estimates \eqref{Ai3}. Note that $Y\p_Y\psi$ satisfies
\begin{align*}
\widetilde{\text{Airy}}(Y\partial_Y\psi)&=Y\partial_Yh+i\epsilon\left[- Y\partial_Y(\partial_YA^{-1}\partial_Y\psi)+\partial_Y(A^{-1}\partial_Y\psi)+ A^{-1}\partial_Y^2\psi\right]-Y\partial_YU_s\psi,
\end{align*}
and boundary condition $Y\partial_Y\psi|_{Y=0}=0.$ Then it follows from \eqref{4.52} that
\begin{align}
\|U_sY\partial_Y\psi\|_{L^2}+\epsilon^{\frac23}\|Y\partial_Y^2\psi\|_{L^2}+\eps\|Y\p_Y^3\psi\|_{L^2}&\lesssim \|Y\partial_Yh\|_{L^2}+\epsilon\|\partial_Y^2\psi\|_{L^2}+\eps^{\frac23}\|\partial_Y\psi\|_{L^2}+\|U_s\psi\|_{L^2}\notag\\
&\lesssim \|Y\partial_Yh\|_{L^2}+\|h\|_{L^2},\nonumber
\end{align}
which is \eqref{Ai3}. The estimates \eqref{Ai4} and \eqref{Ai6} can be derived  similarly. The proof of Lemma \ref{P4.2} is complete.
%It is direct to check
%\begin{align*}
%mAiry(Y^2\partial_Y\psi)=Y^2\partial_Yh-i\epsilon Y^2\partial_Y(\partial_YA^{-1}\partial_Y\psi)+2i\epsilon\partial_Y(A^{-1}Y\partial_Y\psi)+2i\epsilon A^{-1}Y\partial_Y^2\psi-Y^2\partial_YU_s\psi.
%\end{align*}
%Thus we have
%\begin{align}\label{Ai10-1}
%&\|U_sY^2\partial_Y\psi\|_{L^2}+\epsilon^{\frac13}\|Y^2\partial_Y\psi\|_{L^2}+\epsilon^{\frac23}\|Y^2\partial_Y^2\psi\|_{L^2}+\epsilon\|Y^2\partial_Y^3\psi\|_{L^2}\notag\\
%&\quad\lesssim \|Y^2\partial_Yh\|_{L^2}+\|Yh\|_{L^2}+\epsilon^{\frac13}\|h\|_{L^2}
%\end{align}
%and if in addition $h=g/Y$,
%\begin{align}\label{Ai11}
%&\|U_sY^2\partial_Y\psi\|_{L^2}+\epsilon^{\frac13}\|Y^2\partial_Y\psi\|_{L^2}+\epsilon^{\frac23}\|Y^2\partial_Y^2\psi\|_{L^2}+\epsilon\|Y^2\partial_Y^3\psi\|_{L^2}\notag\\
%&\quad\lesssim \|Y\partial_Yg\|_{L^2}+\|g\|_{L^2}.
%\end{align}
%This finishes the proof of Lemma \ref{P4.2}.
\end{proof}

 Finally, we introduce the properties of solutions to following classical Airy equation with the Neumann boundary condition. 
\begin{align}\label{4.67}
\begin{cases}
\text{Airy}(\xi)\eqdef i\epsilon\Delta_\alpha\xi+U_s\xi=\partial_Yg, \quad Y>0,\\
\partial_Y\xi|_{Y=0}=0.
\end{cases}
\end{align}
\begin{lemma}[Proposition 6.3, \cite{GM19}]
\label{P4.3}
Let $g\in H_0^1(\mathbb{R}_+)$. There exists a unique solution $\xi\in H^2(\mathbb{R}_+)$ to system \eqref{4.67} satisfying
\begin{align}\label{4.68}
\epsilon^{\frac12}\|\sqrt{U_s}\xi\|_{L^2}+\epsilon^{\frac23}\|\xi\|_{L^2}+\epsilon(\|\partial_Y\xi\|_{L^2}+\alpha\|\xi\|_{L^2})\le C\|g\|_{L^2},
\end{align}
and
\begin{align}\label{4.69}
\|U_s\xi\|_{L^2}\le C\epsilon^{-\frac13}\|g\|_{L^2}+C\epsilon^{-\frac23}\|U_sg\|_{L^2}.
\end{align}
If $(1+Y)^2g\in H^1(\mathbb{R}_+)$ in addition, 
\begin{align}\label{4.70}
\|Y\xi\|_{L^2}&\le C\epsilon^{-\frac13}\|g\|_{L^2}+C\epsilon^{-\frac23}\|Yg\|_{L^2},\\
\|U_sY^2\xi\|_{L^2}&\le C\epsilon^{-\frac13}\|Y^2g\|_{L^2}+C\epsilon^{-\frac23}\|U_sY^2g\|_{L^2}+C\|Yg\|_{L^2}+C\epsilon^{\frac13}\|g\|_{L^2},\label{4.71}
\end{align}
and
\begin{align}\label{4.72}
|\int_0^\infty U_s\xi\;dY|\le C\alpha^2\left(\epsilon^{\frac16}\|Yg\|_{L^2}+\epsilon^{\frac12}\|g\|_{L^2}\right).
\end{align}
Moreover, for $0< \alpha\le 1$,  
\begin{align}\label{4.73}
\|(1+Y)\partial_Y^{-1}(U_s\xi)\|_{L^2}\le C\|(1+Y)^2g\|_{L^2}.
\end{align}
\end{lemma}

\subsection{Rayleigh-Airy iteration}\label{S3.3}
In this subsection, we use the Rayleigh-Airy iteration method to solve the symmetrized compressible Orr-Sommerfeld equation \eqref{4.8} in the regime $\hat{n}\le C\nu^{-\frac34}$ where $C$ is a given positive constant. Let $\varphi^{1}\in H^2(\mathbb{R}_+)\cap H_0^1(\mathbb{R}_+)$ be the solution to the compressible Rayleigh equation $\text{Ray}(\varphi^{1})=f$ with $\varphi^1|_{Y=0}=0$. Then it holds that
\begin{align}
\widetilde{\text{OS}}_{\text{CNS}}(\varphi^{1})=f+i\epsilon\Delta_\alpha\Lambda(\varphi^{1}).\label{4.83-1}
\end{align}
To eliminate the error term induced by $\varphi^1$, we rewrite $\widetilde{\text{OS}}_{\text{CNS}}$ as
\begin{align}%
\widetilde{\text{OS}}_{\text{CNS}}(\phi)&=\Delta_\alpha(i\epsilon\Lambda(\phi)+U_s\phi)+U_s\Lambda(\phi)-\Delta_\alpha(U_s\phi)-\partial_Y(A^{-1}\partial_YU_s)\phi\nonumber\\
&=\Delta_\alpha \widetilde{\text{Airy}}(\phi)-\partial_Y\left((1+A^{-1})\partial_YU_s\phi+(1-A^{-1})U_s\partial_Y\phi\right),\nonumber
\end{align}
and introduce $\psi^{1}\in H^2(\mathbb{R}_+)\cap H_0^1(\mathbb{R}_+)$ solving the following compressible Airy equation
\begin{align}\label{4.83}
\begin{cases}
\widetilde{\text{Airy}}(\psi^{1})=-i\epsilon\Lambda(\varphi^{1})=-i\epsilon\left(\frac{f}{U_s}+\frac{\partial_Y(A^{-1}\partial_YU_s)\varphi^1}{U_s}\right),\quad Y>0,\\
\psi^1|_{Y=0}=0.
\end{cases}
\end{align}
From \eqref{4.83-1} and \eqref{4.83}, we obtain
\begin{align}\label{4.83-2}
\widetilde{\text{OS}}_{\text{CNS}}(\varphi^1+\psi^1)=-\partial_Y\left((1+A^{-1})\partial_YU_s\psi^1+(1-A^{-1})U_s\partial_Y\psi^1\right).
\end{align}
However, the error term on the right hand side of \eqref{4.83-2} does not vanish at $Y=0$. So it is not a good source term for the Rayleigh equation. To overcome this difficulty, we follow the approach in \cite{GM19} by rewriting $\widetilde{\text{OS}}_{\text{CNS}}$ as
\begin{align}
\widetilde{\text{OS}}_{\text{CNS}}(\phi)&=(i\epsilon\Delta_\alpha+U_s)\Lambda(\phi)-\partial_Y(A^{-1}\partial_YU_s)\phi\notag\\
&=(i\epsilon\Delta_\alpha+U_s)\left(\Lambda(\phi)-\frac{\partial_Y(A^{-1}\partial_YU_s)\phi}{U_s}\right)+i\epsilon\Delta_\alpha\left(\frac{\partial_Y(A^{-1}\partial_YU_s)\phi}{U_s}\right),\notag\\
&=\text{Airy}\left(\frac{1}{U_s}\text{Ray}(\phi)\right)+i\epsilon\Delta_\alpha\left(\frac{\partial_Y(A^{-1}\partial_YU_s)\phi}{U_s}\right),\nonumber
\end{align}
where  Airy$(\cdot)$ is the classical Airy operator defined in \eqref{4.67}.

Based on this decomposition, we introduce $\varphi^{2}$ as the solution to the compressible Rayleigh equation:
\begin{align}\nonumber
	\begin{cases}
		\text{Ray}(\varphi^{2})=U_s\xi^1,\quad Y>0,\\
		\varphi^{2}|_{Y=0}=0,
	\end{cases}
\end{align}
where $\xi^1$ solves the following classical Airy equation with Neumann boundary condition:
\begin{align}\nonumber
\begin{cases}
\text{Airy}(\xi^1)=\partial_Y\left((1+A^{-1})\partial_YU_s\psi^1+(1-A^{-1})U_s\partial_Y\psi^1\right),\quad Y>0,\\
\partial_Y\xi^1|_{Y=0}=0.
\end{cases}
\end{align}
It is straightforward to check that
\begin{align*}
\widetilde{\text{OS}}_{\text{CNS}}(\varphi^2)=\partial_Y[(1+A^{-1})\partial_YU_s\psi^1+(1-A^{-1})U_s\partial_Y\psi^1]
+i\epsilon\Delta_\alpha\left(\frac{\partial_Y(A^{-1}\partial_YU_s)\varphi^{2}}{U_s}\right).
\end{align*}
Combining $\varphi^1, \psi^1, \varphi^2$ together, we obtain
\begin{align*}
\widetilde{\text{OS}}_{\text{CNS}}(\varphi^1+\psi^1+\varphi^2)=f+i\epsilon\Delta_\alpha\left(\frac{\partial_Y(A^{-1}\partial_YU_s)\varphi^2}{U_s}\right).
\end{align*}
Then we can introduce $\psi^2$ which solves the compressible Airy equations with source term $-i\epsilon\left(\frac{\partial_Y(A^{-1}\partial_YU_s)\varphi^2}{U_s}\right)$. Thus, $\varphi^2+\psi^2$ yields the corrector at the second step.

Now we define the following iteration scheme. Let $\varphi^k(k\ge2)$ be the Rayleigh solution satisfying
\begin{align}\label{4.87}
\begin{cases}
\text{Ray}(\varphi^k)=U_s\xi^{k-1},\quad Y>0,\\
\varphi^k|_{Y=0}=0.
\end{cases}
\end{align}
Note that $\varphi^k$ introduces an error  $\widetilde{\text{OS}}_{\text{CNS}}(\varphi^k)=i\epsilon\left(\frac{\partial_Y(A^{-1}\partial_YU_s)\varphi^k}{U_s}\right)$. Then we set $\psi^{k}$ as the solution to the following compressible Airy equation:
\begin{align}\label{4.88}
\begin{cases}
\widetilde{\text{Airy}}(\psi^{k})=-i\epsilon\left(\frac{\partial_Y(A^{-1}\partial_YU_s)\varphi^{k}}{U_s}\right),\quad Y>0,\\
\psi^{k}|_{Y=0}=0.
\end{cases}
\end{align}
The error induced by $\psi^{k}$ is $\widetilde{\text{OS}}_{\text{CNS}}(\psi^k)=-\partial_Y((1+A^{-1})\partial_YU_s\psi^k+(1-A^{-1})U_s\partial_Y\psi^k)$. Finally, we introduce $\xi^{k}$ solving
\begin{align}\label{4.89}
\begin{cases}
\text{Airy}(\xi^{k})=\partial_Y\left((1+A^{-1})\partial_YU_s\psi^k+(1-A^{-1})U_s\partial_Y\psi^k\right),\quad Y>0,\\
\partial_Y\xi^{k}|_{Y=0}=0.
\end{cases}
\end{align}

Now we set the $k$-th order approximation
 $\Psi_k=\varphi^{1}+\psi^1+\sum\limits_{j=2}^{k}(\varphi^{j}+\psi^{j})$, which satisfies
\begin{align}
\widetilde{\text{OS}}_{\text{CNS}}(\Psi_k)=f-\partial_Y\left((1+A^{-1})\partial_YU_s\psi^k+(1-A^{-1})U_s\partial_Y\psi^k\right).\nonumber
\end{align}
In constrast to the incompressible Orr-Sommerfeld equation studied in \cite{GM19}, the error term introduced by the $k$-th approximation depends on higher-order derivatives of Airy solutions due to the compressibility. Fortunately, these terms are accompanied with weights that vanish on the boundary, and can be controlled using co-normal bounds of Airy solutions.

Note that the series $\Psi=\varphi^{1}+\psi^1+\sum\limits_{k=2}^{\infty}(\varphi^{k}+\psi^{k})$ defines a solution to  \eqref{4.8}. In the following proposition, we will prove the convergence of this series in the regime $\alpha\eps^{\frac13}\lesssim1$, which corresponds to $\hat{n}\lesssim \nu^{-\frac34}$. For convenience, we introduce $\psi^0$ satisfying
\begin{align}\label{4.89-1}
 \widetilde{\text{Airy}}(\psi^0)=-i\epsilon\frac{f}{U_s},~~\psi^0|_{Y=0}=0.
 \end{align}
\begin{proposition}{\bf (Convergence)}
\label{P4.4}
Suppose that $m\in (0,1)$, and $f/U_s\in L^2(\mathbb{R}_+)$. If $\eps\in (0,1)$ is sufficiently small and $\alpha\eps^{\frac13}\lesssim 1$, the symmetrized Orr-Sommerfeld equation \eqref{4.8} admits a solution $\phi\in H^4(\mathbb{R}_+)\cap H_0^1(\mathbb{R}_+)$ . Moreover, the following two statements hold
%\begin{align}
%	\|\partial_Y(\phi-\varphi^1-\psi^0)\|_{L^2}+\alpha\|\phi-\varphi^1-\psi^0\|_{L^2}
%	\le& C\epsilon^{-\frac23}(\|U_s\psi^{0}\|_{L^2}+\|U_sY^2\partial_Y\psi^0\|_{L^2})\nonumber\\
%	&+C\epsilon^{-\frac13}\|\psi^{0}\|_{L^2}+C\epsilon^{\frac13}\|\partial_Y\varphi^1\|_{L^2}, \label{4.90}
%\end{align}
%and
%\begin{align}
%\|\Delta_\alpha(\phi-\varphi^0-\psi^1)\|_{L^2}\leq  C\epsilon^{-\frac23}(\|\psi^{0}\|_{L^2}+\|U_sY^2\partial_Y\psi^0\|_{L^2})+C\|\partial_Y\varphi^1\|_{L^2}.\label{4.91}
%\end{align}

\begin{itemize}
	\item[{\rm(i)}] If $\alpha\ge 1$, 
\begin{align}\label{4.90}
&\|\partial_Y(\phi-\varphi^1-\psi^0)\|_{L^2}+\alpha\|\phi-\varphi^1-\psi^0\|_{L^2}\notag\\
&\qquad\le C\epsilon^{-\frac23}\left(\|U_s\psi^{0}\|_{L^2}+\epsilon^{\frac13}\|\psi^{0}\|_{L^2}+\|U_sY^2\partial_Y\psi^0\|_{L^2}\right)+C\epsilon^{\frac13}\|\partial_Y\varphi^1\|_{L^2},\\
&\|\Delta_\alpha(\phi-\varphi^1-\psi^0)\|_{L^2}\notag\\
&\qquad\le C\epsilon^{-\frac23}\left(\|\psi^{0}\|_{L^2}+\|U_sY^2\partial_Y\psi^0\|_{L^2}\right)+C\|\partial_Y\varphi^1\|_{L^2}.\label{4.91}
\end{align}
\item[{\rm(ii)}] If $\alpha\in (0,1)$, 
\begin{align}
&\|\partial_Y(\phi-\varphi^1-\psi^0)\|_{L^2}\nonumber\\
&\qquad\le C\epsilon^{-\frac23}\left(\|U_s\psi^{0}\|_{L^2}+\epsilon^{\frac13}\|\psi^{0}\|_{L^2}+\|U_sY^2\partial_Y\psi^0\|_{L^2}\right)+C\epsilon^{\frac13}\|\partial_Y\varphi^1\|_{L^2},\label{4.92}\\
&\alpha\|(\phi-\varphi^1-\psi^0)\|_{L^2}\nonumber\\
&\qquad\le C\alpha\left(\|\psi^{0}\|_{L^2}+\|U_sY^2\partial_Y\psi^0\|_{L^2}\right)+C\alpha\epsilon^{\frac23}\|\partial_Y\varphi^1\|_{L^2},\label{4.94}\\
&\|\Delta_\alpha(\phi-\varphi^1-\psi^0)\|_{L^2}\nonumber\\
&\qquad\le C\epsilon^{-\frac23}\left(\|\psi^{0}\|_{L^2}+\|U_sY^2\partial_Y\psi^{0}\|_{L^2}\right)
+C\|\partial_Y\varphi^1\|_{L^2}.\label{4.93}
\end{align}
\end{itemize}
\end{proposition}

\begin{proof}
We divide the proof into two cases.

\underline{\it Case 1: $\alpha \geq 1$}. Recall that $\varphi^k$ solves \eqref{4.87}.  Applying \eqref{4.12} to $\varphi^k$, we obtain
\begin{align}
\|\partial_Y\varphi^{k+1}\|_{L^2}+\alpha\|\varphi^{k+1}\|_{L^2}\le C\|Y\xi^{k}\|_{L^2},~ k\geq 1.\label{4.93-1}
\end{align}
Since $\xi^k$ solves \eqref{4.89}, then using Airy estimate \eqref{4.70} with $$g=(1+A^{-1})\partial_YU_s\psi^k+(1-A^{-1})U_s\partial_Y\psi^k,$$ $|1-A^{-1}|\leq Cm^2U_s^2$ and $|U_s/Y|\leq C$,  we can get
\begin{align}
\|Y\xi^{k}\|_{L^2}&\leq C \epsilon^{-\frac23}(\|Y\partial_YU_s\psi^{k}\|_{L^2}+\|YU_s^3\partial_Y\psi^k\|_{L^2})
+C\epsilon^{-\frac13}(\|\psi^{k}\|_{L^2}+\|U_s^3\partial_Y\psi^k\|_{L^2})\nonumber\\
&\le C\epsilon^{-\frac23}(\|U_s\psi^{k}\|_{L^2}+\|U_sY^2\partial_Y\psi^k\|_{L^2})+C\epsilon^{-\frac13}\|\psi^{k}\|_{L^2}.\label{4.93-2}
\end{align}
To estimate $\psi^k$, we set $h^{k}=\frac{\partial_Y(A^{-1}\partial_YU_s)\varphi^{k}}{U_s}$. Since $\psi^k$ solves \eqref{4.88}, and
\begin{align}\label{4.93-4}
	\|Y^2\p_Yh^k\|_{L^2}+\|Yh^k\|_{L^2}+\|h^k\|_{L^2}\leq C\|\p_Y\varphi^k\|_{L^2},
\end{align}
then from \eqref{4.52} and \eqref{Ai4} we can deduce that
\begin{align}\label{4.93-3}
\|U_sY^2\partial_Y\psi^k\|_{L^2}+\|U_s\psi^{k}\|_{L^2}+\epsilon^{\frac13}\|\psi^{k}\|_{L^2} &\leq C\eps(\|Y^2\p_Yh^k\|_{L^2}+\|Yh^k\|_{L^2}+\|h^k\|_{L^2})\nonumber\\
&\leq C \epsilon\|\partial_Y\varphi^{k}\|_{L^2}.
\end{align}
Substituting \eqref{4.93-3} and \eqref{4.93-2} into \eqref{4.93-1}, we get
\begin{align*}
\|\partial_Y\varphi^{k+1}\|_{L^2}+\alpha\|\varphi^{k+1}\|_{L^2}\le C\epsilon^{\frac13}\|\partial_Y\varphi^{k}\|_{L^2}.
\end{align*}

Taking $\eps\in (0,1)$ sufficiently small, we deduce that the series $\sum_{k=2}^{\infty}\varphi^{k}$  converges in $H^1(\mathbb{R}_+)$. Moreover, we have
\begin{align}\label{4.95}
\sum\limits_{k=2}^{\infty}\|\partial_Y\varphi^{k}\|_{L^2}+\sum\limits_{k=2}^{\infty}\alpha\|\varphi^{k}\|_{L^2}&\le C\|\partial_Y\varphi^2\|_{L^2}+\alpha\|\varphi^2\|_{L^2}\notag\\
&\le C\epsilon^{-\frac23}(\|U_s\psi^{1}\|_{L^2}+\|U_sY^2\partial_Y\psi^1\|_{L^2})+C\epsilon^{-\frac13}\|\psi^{1}\|_{L^2}.
\end{align}
From \eqref{4.13}, \eqref{4.68},  \eqref{4.93-3} and \eqref{4.95}, we have
\begin{align}\label{4.96}
\sum\limits_{k=2}^{\infty}\|\Delta_\alpha\varphi^k\|_{L^2}&\le C\left(\sum_{k=1}^{\infty}\|\xi^k\|_{L^2}\right)
\leq C\eps^{-\frac23}\left(\sum_{k=1}^\infty\|\psi^k\|_{L^2}+\|U_sY^2\p_Y\psi^k\|_{L^2}\right)\nonumber\\
&
\leq C\epsilon^{-\frac23}(\|\psi^{1}\|_{L^2}+\|U_sY^2\partial_Y\psi^1\|_{L^2})+C\left(\sum_{k=2}^\infty\|\p_Y\psi^k\|_{L^2}\right)\nonumber\\
&\leq C\epsilon^{-\frac23}(\|\psi^{1}\|_{L^2}+\|U_sY^2\partial_Y\psi^1\|_{L^2}).
\end{align}

For convergence of $\psi^k$, by using \eqref{4.52} and \eqref{4.93-4}, we obtain
\begin{align}\label{4.98}
\sum\limits_{k=2}^{\infty}\|\partial_Y\psi^{k}\|_{L^2}+\alpha\sum\limits_{k=2}^{\infty}\|\psi^{k}\|_{L^2}+\sum\limits_{k=2}^{\infty}\eps^{\frac13}\|\Delta_\alpha\psi^{k}\|_{L^2}
&\le C\epsilon^{\frac13}\left(\sum_{k=2}^\infty \|h^k\|_{L^2}\right)\leq C\epsilon^{\frac13}\left(\sum_{k=2}^\infty \|\partial_Y\phi^k\|_{L^2}\right)\nonumber\\
&\leq C\epsilon^{-\frac13}(\|U_s\psi^{1}\|_{L^2}+\|U_sY^2\partial_Y\psi^1\|_{L^2})+C\|\psi^{1}\|_{L^2}.
\end{align}
Combining \eqref{4.95}, \eqref{4.96} and \eqref{4.98} together, we obtain

\begin{align}
	\sum_{k=2}^\infty\|(\p_Y\varphi^k,\alpha\varphi^k)\|_{L^2}+\sum_{k=2}^\infty\|(\p_Y\psi^k,\alpha\psi^k)\|_{L^2}&\leq C\epsilon^{-\frac 23}(\|U_s\psi^{1}\|_{L^2}+\|U_sY^2\partial_Y\psi^1\|_{L^2})+C\epsilon^{-\frac 13}\|\psi^{1}\|_{L^2},\label{4.99-1}\\
	\sum_{k=2}^\infty\|\Delta_\alpha\varphi^k\|_{L^2}+\sum_{k=2}^\infty\|\Delta_\alpha\psi^k\|_{L^2}&\leq C\epsilon^{-\frac 23}(\|\psi^{1}\|_{L^2}+\|U_sY^2\partial_Y\psi^1\|_{L^2}).\label{4.99-2}
\end{align}

Finally, to bound $\psi^{1}$,  we decompose $\psi^1=\psi^0+\psi^{1, 1}$ with $\psi^0$ and $\psi^{1, 1}$ satisfy compressible Airy equations with source term $-i\epsilon\frac{f}{U_s}$ and 
$-i\epsilon\frac{\partial_Y(A^{-1}\partial_YU_s)\varphi^1}{U_s}$ respectively. Applying \eqref{4.52}, \eqref{Ai4} and \eqref{4.93-4} to $\psi^{1, 1}$,  we can obtain
\begin{align}\label{4.99-3}
\epsilon^{-\frac23}(\|U_s\psi^{1, 1}\|_{L^2}+\|U_sY^2\partial_Y\psi^{1, 1}\|_{L^2})+\epsilon^{-\frac13}\|\psi^{1, 1}\|_{L^2}+\|\partial_Y\psi^{1, 1}\|_{L^2}+\eps^{\frac13}\|\Delta_\alpha\psi^{1, 1}\|_{L^2}\le C\epsilon^{\frac13}\|\partial_Y\varphi^1\|_{L^2}.
\end{align}
Substituting it into \eqref{4.99-1} and \eqref{4.99-2}, we obtain \eqref{4.90} and \eqref{4.91}. The proof of Proposition \ref{P4.4} when $\alpha\ge 1$ is complete.\\

\underline{\it Case 2: $\alpha \in(0,1)$.} Recall that $\varphi^k$ solves \eqref{4.87}. Applying \eqref{4.15} to $\varphi^k$ yields that
\begin{align}
\|\partial_Y\varphi^{k+1}\|_{L^2}\le C\|(1+Y)\partial_Y^{-1}(U_s\xi^{k})\|_{L^2}+C\|U_s\xi^k\|_{L^2}+\frac{C}{\alpha}|\int_0^\infty U_s\xi^{k}\;dY|.\label{4.100-4}
\end{align}
Now we bound the right hand side of \eqref{4.100-4}. Recall that $\xi_k$ solves \eqref{4.89}. We use bounds \eqref{4.69}, \eqref{4.72} and \eqref{4.73} in Lemma \ref{P4.3} to obtain
\begin{align}
|\int_0^\infty U_s\xi^{k}\;dY|&\le C\alpha^2\left(\epsilon^{\frac16}\|YU_s^\prime\psi^{k}\|_{L^2}+\epsilon^{\frac16}\|YU_s^3\partial_Y\psi^k\|_{L^2}+\epsilon^{\frac12}\|U_s^\prime\psi^{k}\|_{L^2}+\epsilon^{\frac12}\|U_s^3\partial_Y\psi^k\|_{L^2}
\right)\nonumber\\
&\le C\alpha^2\left(\|\psi^k\|_{L^2}+\|U_sY^2\p_Y\psi^k\|_{L^2}\right),\label{4.100-2}\\
	\|U_s\xi^k\|_{L^2}&\leq C\eps^{-\frac{1}{3}}\|\psi^k\|_{L^2}+C\eps^{-\frac{2}{3}}\left(\|U_s\psi^k\|_{L^2}+\|U_sY^2\p_Y\psi^k\|_{L^2}\right),\label{4.100-1}
\end{align}
and
\begin{align}
\|(1+Y)\partial_Y^{-1}(U_s\xi^{k})\|_{L^2}&\leq C\|(1+Y)^2\partial_YU_s\psi^k\|_{L^2}+C\|(1+Y)^2U_s^3\partial_Y\psi^k\|_{L^2}\nonumber\\
&\le C\|\psi^k\|_{L^2}+C\|U_sY^2\partial_Y\psi^k\|_{L^2}.\label{4.100-3}
%\|U_s\xi^{k}\|_{L^2}+\|U_sY^2\xi^{k}\|_{L^2}\notag\\
%&\lesssim \epsilon^{-\frac13}\|U_s^\prime\psi^{k}\|_{L^2}+\epsilon^{-\frac23}\|YU_s^\prime\psi^{k}\|_{L^2}+\epsilon^{-\frac13}\|Y^2U_s^\prime\psi^k\|_{L^2}+\epsilon^{-\frac23}\|U_sY^2U_s^\prime\psi^{k}\|_{L^2}\notag\\
%&\quad+\|YU_s^\prime\psi^{k}\|_{L^2}+\epsilon^{\frac13}\|U_s^\prime\psi^{k}\|_{L^2}+\epsilon^{-\frac23}\|U_sY^2\partial_Y\psi^{k}\|_{L^2}\notag\\
%&\le C\epsilon^{-\frac23}(\|U_s\psi^{k}\|_{L^2}+\|U_sY^2\partial_Y\psi^k\|_{L^2})+C\epsilon^{-\frac13}\|\psi^{k}\|_{L^2}.
\end{align}
Substituting \eqref{4.100-2}, \eqref{4.100-1}, and \eqref{4.100-3} into \eqref{4.100-4}, we deduce that
\begin{align}\label{4.100-6}
	\|\partial_Y\varphi^{k+1}\|_{L^2}\leq C\eps^{-\frac13}\|\psi^k\|_{L^2}+C\eps^{-\frac23}\|U_s\psi^k\|_{L^2}+C\eps^{-\frac23}\|U_sY^2\psi^k\|_{L^2}.
\end{align}

Recall that $\psi^k$ solves \eqref{4.88}. Then similar to \eqref{4.93-3}, we obtain
\begin{align}
\eps^{\frac13}\|\psi^k\|_{L^2}+\|U_s\psi^k\|_{L^2}+\|U_sY^2\partial_Y\psi^k\|_{L^2}\leq C\eps\|\p_Y\varphi^k\|_{L^2}.\label{4.100-5}
\end{align}
Substituting \eqref{4.100-5} into \eqref{4.100-6}, we deduce that
\begin{align}\label{4.100}
\|\partial_Y\varphi^{k+1}\|_{L^2}&\le C\eps^{-\frac13}\|\psi^k\|_{L^2}+C\eps^{-\frac23}(\|U_s\psi^k\|_{L^2}+\|U_sY^2\partial_Y\psi^k\|_{L^2})\nonumber\\
&\leq C\epsilon^{\frac13}\|\partial_Y\varphi^{k}\|_{L^2}.
\end{align}
Similarly, we can use \eqref{4.14}, \eqref{4.16}, \eqref{4.68}, \eqref{4.100-2}-\eqref{4.100-5} to  obtain
\begin{align}\label{4.101-1}
	\alpha\|\varphi^{k+1}\|&\leq C\alpha\|(1+Y)\partial_Y^{-1}(U_s\xi^{k})\|_{L^2}+\frac{C}{\alpha^{\frac12}}|\int_0^\infty U_s\xi^{k}\;dY|\nonumber\\
	&\leq C\alpha \left( \|\psi^k\|_{L^2}+\|U_sY^2\p_Y\psi^k\|_{L^2}   \right)\nonumber\\
	&\leq 
	C\alpha \eps^{\frac23}\|\p_Y\varphi^k\|_{L^2},
\end{align}
and
\begin{align}\label{4.101}
\|\Delta_\alpha\varphi^{k+1}\|_{L^2}&\le C\|(1+Y)\partial_Y^{-1}(U_s\xi^{k})\|_{L^2}+C\|\xi^k\|_{L^2}+\frac{C}{\alpha}|\int_0^\infty U_s\xi^{k}\;dY|\nonumber\\
&\leq C\eps^{-\frac23}\left( \|\psi^k\|_{L^2}+\|U_sY^2\p_Y\psi^k\|_{L^2} \right)\nonumber\\
&\leq C\|\p_Y\varphi^k\|_{L^2}.
\end{align}
From \eqref{4.100}, \eqref{4.101-1}, and \eqref{4.101}, we derive the following estimates:
\begin{equation}
\begin{aligned}\label{4.102}
\sum\limits_{k=2}^{\infty}\|\partial_Y\varphi^{k}\|_{L^2}&\le C\|\partial_Y\varphi^2\|_{L^2}\le C\eps^{
-\frac13}\|\psi^1\|_{L^2}+C\eps^{-\frac23}\left(\|U_s\psi^1\|_{L^2}+\|U_sY^2\p_Y\psi^1\|_{L^2}\right),\\
\sum\limits_{k=2}^{\infty}\alpha\|\varphi^{k}\|_{L^2}&\le C\alpha\eps^{\frac23}\|\partial_Y\varphi^2\|_{L^2}+\alpha\|\varphi^2\|_{L^2}\le C\alpha\left(\|\psi^1\|_{L^2}+\|U_sY^2\p_Y\psi^1\|_{L^2}\right),\\
\sum\limits_{k=2}^{\infty}\|\Delta_\alpha\varphi^{k}\|_{L^2}&\leq C\|\p_Y\varphi^2\|_{L^2}+\|\Delta_\alpha\varphi^2\|_{L^2}\leq C\epsilon^{-\frac23}(\|\psi^{1}\|_{L^2}+\|U_sY^2\partial_Y\psi^{1}\|_{L^2}).
\end{aligned}
\end{equation}
Since $\psi^k$ solves \eqref{4.88}, then from \eqref{4.52} and \eqref{4.93-4}, we obtain
\begin{equation}
\begin{aligned}\label{4.104} \sum\limits_{k=2}^{\infty}\|\partial_Y\psi^{k}\|_{L^2}&\leq C\eps^{\frac13} \sum\limits_{k=2}^{\infty}\|\partial_Y\varphi^{k}\|_{L^2} \leq  C\|\psi^1\|_{L^2}+C\eps^{-\frac13}\left(\|U_s\psi^1\|_{L^2}+\|U_sY^2\p_Y\psi^1\|_{L^2}\right),\\
\sum\limits_{k=2}^{\infty}\alpha\|\psi^{k}\|_{L^2}&\leq C\alpha\eps^{\frac23} \sum\limits_{k=2}^{\infty}\|\partial_Y\varphi^{k}\|_{L^2} \leq  C\alpha\eps^{\frac13}\|\psi^1\|_{L^2}+C\alpha\left(\|U_s\psi^1\|_{L^2}+\|U_sY^2\p_Y\psi^1\|_{L^2}\right),\\
\sum\limits_{k=2}^{\infty}\|\Delta_\alpha\psi^{k}\|_{L^2}&\leq C \sum\limits_{k=2}^{\infty}\|\partial_Y\varphi^{k}\|_{L^2} \leq  C\eps^{-\frac13}\|\psi^1\|_{L^2}+C\eps^{-\frac23}\left(\|U_s\psi^1\|_{L^2}+\|U_sY^2\p_Y\psi^1\|_{L^2}\right).
\end{aligned}
\end{equation}
Finally, $\psi^1$ can be decomposed into $\psi^0+\psi^{1, 1}$ as in Case 1, where $\psi^{1,1}$ can be estimated the same  as \eqref{4.99-3}. The proof of Proposition \ref{P4.4} is complete.
\end{proof}

Now we can bound $\varphi^1$ and $\psi^0$ in terms of $f$ and then obtain the following Corollary on the solvability of symmtrized Orr-Sommerfeld equation \eqref{4.8}.
\begin{corollary}
\label{C4.5}
Under the same assumption as in Proposition \ref{P4.4}, the solution $\phi\in H^4(\mathbb{R}_+)\cap H_0^1(\mathbb{R}_+)$ to \eqref{4.8} satisfies:
\begin{itemize}
\item[{\rm(i)}]  If $\alpha\ge 1$, then it  holds
\begin{align}\label{4.106}
\|\partial_Y\phi\|_{L^2}+\alpha\|\phi\|_{L^2}&\le C\min\left\{\|\frac{Yf}{U_s}\|_{L^2}, \frac{1}{\alpha}\|\frac{f}{U_s}\|_{L^2}\right\}+C\eps^{\frac13}\left(\|\frac{Y^2\p_Yf}{U_s}\|_{L^2}+\|\frac{Yf}{U_s}\|_{L^2}\right),\\
\|\Delta_\alpha\phi\|_{L^2}&\le C\min\left\{\|\frac{Yf}{U_s}\|_{L^2}, \frac{1}{\alpha}\|\frac{f}{U_s}\|_{L^2}\right\}+C\|\frac{f}{U_s}\|_{L^2}+C\eps^{\frac13}\left(\|\frac{Y^2\p_Yf}{U_s}\|_{L^2}+\|\frac{Yf}{U_s}\|_{L^2}\right). \label{4.107}
\end{align}
\item[{\rm(ii)}] If $0<\alpha\le 1$, then it holds
\begin{align}\label{4.108}
\|\partial_Y\phi\|_{L^2}&\le C\|(1+Y)^2f\|_{L^2}+\frac{C}{\alpha}|\int_0^\infty f(Y)\;dY|+C\eps^{\frac13}\left(\|\frac{Y^2\p_Yf}{U_s}\|_{L^2}+\|\frac{Yf}{U_s}\|_{L^2}\right), \\
\alpha\|\phi\|_{L^2}&\leq C\alpha \|(1+Y)^2f\|_{L^2}+\frac{C}{\alpha^{\frac12}}|\int_0^\infty f(Y)\;dY|+C\eps^{\frac13}\left(\|\frac{Y^2\p_Yf}{U_s}\|_{L^2}+\|\frac{Yf}{U_s}\|_{L^2}\right),\label{4.110}\\
\|\Delta_\alpha\phi\|_{L^2}&\le C\left(\|(1+Y)^2f\|_{L^2}+\|\frac{f}{U_s}\|_{L^2}\right)+\frac{C}{\alpha}|\int_0^\infty f(Y)\;dY|+C\eps^{\frac13}\left(\|\frac{Y^2\p_Yf}{U_s}\|_{L^2}+\|\frac{Yf}{U_s}\|_{L^2}\right).\label{4.109}
\end{align}
\end{itemize}
\end{corollary}
\begin{proof}
	Since $\psi^0$ solves \eqref{4.89-1}, by applying \eqref{4.52}, \eqref{4.53} and \eqref{Ai6} to $\psi^0$, we obtain
	
% Indeed, from Proposition \ref{P4.4}, we first have
%\begin{align*}
%&\|\partial_Y\phi\|_{L^2}+\alpha\|\phi\|_{L^2}\\
%&\quad\lesssim \|\partial_Y\varphi^1\|_{L^2}+\alpha\|\varphi^1\|_{L^2}+\|\partial_Y\psi^0\|_{L^2}+\alpha\|\psi^0\|_{L^2}\\
%&\qquad+\epsilon^{-\frac23}(\|U_s\psi^{0}\|_{L^2}+\|U_sY^2\partial_Y\psi^0\|_{L^2})+\epsilon^{-\frac13}\|\psi^{0}\|_{L^2}.
%\end{align*}
%On the other hand, we get from Proposition \ref{P4.1} and \ref{P4.2} that
\begin{align*}
\|U_s\psi^0\|_{L^2}+\epsilon^{\frac13}\|\psi^0\|_{L^2}+\epsilon^{\frac23}\|\partial_Y\psi^0\|_{L^2}&\le C\eps^{\frac{2}{3}}\min\left\{ \|\frac{Yf}{U_s}\|_{L^2},\eps^{\frac13}\|\frac{f}{U_s}\|_{L^2}
\right\},\\
\|\Delta_\alpha\psi^0\|_{L^2}&\leq C\|\frac{f}{U_s}\|_{L^2},\\
\|U_sY^2\partial_Y\psi^0\|_{L^2}&\le C\epsilon \|Y\partial_Y(\frac{Yf}{U_s})\|_{L^2}+C\epsilon\|\frac{Yf}{U_s}\|_{L^2}\\
&\leq C\eps\|\frac{Y^2\partial_Yf}{U_s}\|_{L^2}+C\eps\|\frac{Yf}{U_s}\|_{L^2}.
\end{align*}
Plugging them back into \eqref{4.90}-\eqref{4.93} and applying estimates \eqref{4.12}-\eqref{4.16} to $\varphi^1$, we obtain \eqref{4.106}-\eqref{4.109}. 
%Thus we conclude with the fact $\alpha\epsilon^{\frac13}\le 1$ that
%\begin{align*}
%\|\partial_Y\phi\|_{L^2}+\alpha\|\phi\|_{L^2}&\le C\min\left\{\|\frac{Yf}{U_s}\|_{L^2}, \frac{1}{\alpha}\|\frac{f}{U_s}\|_{L^2}\right\}
%+\epsilon^{\frac13}\|Y\partial_Y(\frac{Yf}{U_s})\|_{L^2}+\epsilon^{\frac13}\|\frac{Yf}{U_s}\|_{L^2}\\
%\|\Delta_\alpha\phi\|_{L^2}&\le C\min\left\{\|\frac{Yf}{U_s}\|_{L^2}, \frac{1}{\alpha}\|\frac{f}{U_s}\|_{L^2}\right\}+\|\frac{f}{U_s}\|_{L^2}+\epsilon^{\frac13}\|Y\partial_Y(\frac{Yf}{U_s})\|_{L^2}+\epsilon^{\frac13}\|\frac{Yf}{U_s}\|_{L^2}.
%\end{align*}
%This finishes the proof of case (i). Case (ii) follows the similar procedure and the details are omitted here.
\end{proof}
\begin{remark}\label{R4.1}
Same as the incompressble case  discussed in \cite{GM19}, the solution obtained  to symmetrized Orr-Sommerfeld equation \eqref{4.8} satisfies one extra boundary condition
\begin{align}\label{4.110-1}
	\p_Y\Lambda(\phi)|_{Y=0}=0.
\end{align}
Note that $\psi^1$ solves $\widetilde{{\rm Airy}}(\psi^1)=-i\epsilon\Lambda(\varphi^1)$ with $\psi^1|_{Y=0}=0$. Then differentiating the equation and taking trace at $Y=0$, we get
\begin{align*}
i\epsilon\partial_Y\Lambda(\psi^1+\varphi^1)|_{Y=0}=-\partial_Y(U_s\psi^1)|_{Y=0}=0.
\end{align*}
From the iteration scheme \eqref{4.87}-\eqref{4.89}, we observe that
\begin{align*}
&i\eps\Lambda(\varphi^k)=i\eps\frac{\varphi^k\partial_Y(A^{-1}\partial_YU_s)}{U_s}+i\eps\xi^{k-1},\\
&i\epsilon\Lambda(\psi^k)=-i\epsilon\frac{\varphi^k\partial_Y(A^{-1}\partial_YU_s)}{U_s}-U_s\psi^k.
\end{align*}
By differentiating  above two equations and adding them together,  taking trace at $Y=0$ yields
\begin{align*}
i\epsilon\partial_Y\Lambda(\psi^k+\varphi^k)|_{Y=0}=-\partial_Y(U_s\psi^k)|_{Y=0}+i\epsilon\partial_Y\xi^{k-1}|_{Y=0}=0.
\end{align*}
Here we have used boundary conditions $\psi^k|_{Y=0}=0$ and $\partial_Y\xi^{k-1}|_{Y=0}=0$. Thus the boundary condition \eqref{4.110-1} is preserved at each step in the iteration. 
\end{remark}

Next we derive bounds on the solution to \eqref{4.8} if the inhomogeneous source term is in the form of derivative.

\begin{proposition}
\label{P4.6}
Let $m\in (0,1)$ and $\alpha\epsilon^{\frac13}\lesssim 1$. If $f=\p_Yg$ or $f=Y^{-1}g$ or $f=i\alpha g$ with $g\in L^2(\mathbb{R}_+)$, the solution $\phi$ to \eqref{4.8} satifies
\begin{align}\label{4.113}
\|\partial_Y\phi\|_{L^2}+\alpha\|\phi\|_{L^2}&\leq C\epsilon^{-\frac13}\|g\|_{L^2},\\
\|\Delta_\alpha\phi\|_{L^2}&\leq C\epsilon^{-\frac23}\|g\|_{L^2}.\label{4.114}
\end{align}
\end{proposition}
\begin{proof}
We look for the solution in the form $\phi=\phi_0+\phi_1$, where $\phi_0$ and $\phi_1$ are  solutions to following equations respectively: 
	\begin{align}\label{4.111}
		\begin{cases}
			i\eps\Delta_\alpha\Lambda(\phi_0)+\partial_Y(A^{-1}U_s\partial_Y\phi_0)-\alpha^2 U_s\phi_0=f,\\
			\phi_0|_{Y=0}=\partial_Y\phi_0|_{Y=0}=0,
		\end{cases}
	\end{align}
and
\begin{align}\label{4.112}
	\begin{cases}
		\widetilde{\text{OS}}_{\text{CNS}}(\phi_1)=\partial_Y(A^{-1}\partial_YU_s\phi_0),\\
		\phi_1|_{Y=0}=0.
	\end{cases}
\end{align}

\underline{Estimate on $\phi_0$}. The existence of solution $\phi_0$ can be shown by the same argument as in \cite[Proposition 7.12]{GM19}. Thus we only focus on the a priori estimate. Taking inner product of \eqref{4.111} with ${\phi}_0$ and then integrating by parts, we obtain
\begin{align}\label{4.115}
-i\epsilon\langle\Delta_\alpha\Lambda(\phi_0), \phi_0\rangle+\int_0^\infty A^{-1}|\sqrt{U_s}\partial_Y\phi_0|^2\;dY+\alpha^2\|\sqrt{U_s}\phi_0\|_{L^2}^2=\langle f, \phi_0\rangle.
\end{align}
If $f=\p_Yg$, or $f=Y^{-1}g$ or $f=i\alpha g$,  integration by parts or Hardy's inequality yields
\begin{align}\label{4.115-1}
	|\langle f, \phi_0\rangle|\leq \|g\|_{L^2}\|(\p_Y\phi_0,\alpha\phi_0)\|_{L^2}.
\end{align}

For the first terms on the left hand side of \eqref{4.115}, integrating by parts yields
\begin{align*}
-\langle\Delta_\alpha\Lambda(\phi_0), \phi_0\rangle&=-\langle \Lambda(\phi_0), \partial_Y^2\phi_0\rangle+\alpha^2\langle\Lambda(\phi_0), \phi_0\rangle\\
&=-\langle A^{-1}\p_Y^2\phi_0, \partial_Y^2\phi_0\rangle-\langle\partial_Y(A^{-1})\p_Y\phi_0, \partial_Y^2\phi_0\rangle+\alpha^2\langle\phi_0, \partial_Y^2\phi_0\rangle+\alpha^2\langle\Lambda(\phi_0), \phi_0\rangle\\
&=-\langle A^{-1}\partial_Y^2\phi_0, \partial_Y^2\phi_0\rangle
-\alpha^2\langle(A^{-1}+1)\p_Y\phi_0, \partial_Y\phi_0\rangle-\alpha^4\|\phi_0\|_{L^2}-\langle\partial_YA^{-1}\partial_Y^2\phi_0, \partial_Y\phi_0\rangle.
\end{align*}
Then by taking the real and imaginary parts of \eqref{4.115} respectively and using \eqref{4.115-1}, we obtain
\begin{align}\label{4.116}
\|\sqrt{U_s}\partial_Y\phi_0\|_{L^2}^2+\alpha^2\|\sqrt{U_s}\phi_0\|_{L^2}^2&\le C\|g\|_{L^2}\|(\p_Y\phi_0,\alpha\phi_0)\|_{L^2}+C\epsilon\|\partial_Y\phi_0\|_{L^2}\|\partial_Y^2\phi_0\|_{L^2},\\
\epsilon\|\Delta_\alpha\phi_0\|_{L^2}^2&\le C\|g\|_{L^2}\|(\p_Y\phi_0,\alpha\phi_0)\|_{L^2}+C\epsilon\|\partial_Y\phi_0\|_{L^2}^2.\label{4.117}
\end{align}

Set
$
\mathcal{E}=\|\partial_Y\phi_0\|_{L^2}^2+\|\alpha\phi_0\|_{L^2}^2+\eps^{\frac23}\|\Delta_\alpha\phi_0\|_{L^2}^2.
$ From \eqref{4.116} and \eqref{4.117}, we obtain 
\begin{align}\label{4.118}
\|\sqrt{U_s}\partial_Y\phi_0\|_{L^2}^2+\alpha^2\|\sqrt{U_s}\phi_0\|_{L^2}^2&\le C\|g\|_{L^2}\mathcal{E}^{\frac12}+C\eps^{\frac23}\mathcal{E},\\
\epsilon^{\frac{2}{3}}\|\Delta_\alpha\phi_0\|_{L^2}^2&\le C\eps^{-\frac13}\|g\|_{L^2}\mathcal{E}^{\frac12}+C\eps^{\frac23}\mathcal{E}.\label{4.119}
\end{align}
Moreover, applying the interpolation inequality
$\|f\|_{L^2}^2\leq C\|\sqrt{U_s}f\|_{L^2}^{\frac43}\|\p_Yf\|_{L^2}^{\frac23}+C\|\sqrt{U_s}f\|_{L^2}^2$ (cf. \cite[Propostion 2.4]{GM19}) to $(\p_Y\phi_0,\alpha\phi_0)$, we get
\begin{align*}
\|(\p_Y\phi_0,\alpha\phi_0)\|_{L^2}^2&\leq C\|\sqrt{U_s}(\partial_Y\phi_0,\alpha\phi_0)\|_{L^2}^{\frac43}\|\Delta_\alpha\phi_0\|_{L^2}^{\frac23}+C\|\sqrt{U_s}(\partial_Y\phi_0,\alpha\phi_0)\|_{L^2}^2\\
&\leq C\eps^{-\frac{1}{3}}\|g\|_{L^2}\mathcal{E}^{\frac12}+C\eps^{\frac{1}{9}}\|g\|_{L^2}^{\frac13}\mathcal{E}^{\frac56}+\|g\|_{L^2}^{\frac23}\mathcal{E}^{\frac23}+C\eps^{\frac{4}{9}}\mathcal{E}.
\end{align*}
Here we have used \eqref{4.118} and \eqref{4.119} in the last inequality. Now combining it with \eqref{4.119} and using Young's inequality, we obtain
\begin{align}
\mathcal{E}&\le C\eps^{-\frac{1}{3}}\|g\|_{L^2}\mathcal{E}^{\frac12}+C\eps^{\frac{1}{9}}\|g\|_{L^2}^{\frac13}\mathcal{E}^{\frac56}+\|g\|_{L^2}^{\frac23}\mathcal{E}^{\frac23}+C\eps^{\frac{4}{9}}\mathcal{E}\leq C\eps^{-\frac23}\|g\|_{L^2}^2+o(1)\mathcal{E},\nonumber
\end{align}
which implies
\begin{align}
	\|\partial_Y\phi_0\|_{L^2}+\|\alpha\phi_0\|_{L^2}+\eps^{\frac13}\|\Delta_\alpha\phi_0\|_{L^2}\leq C\eps^{-\frac13}\|g\|_{L^2}.\label{4.120}
\end{align}
The estimate on $\phi_0$ is complete.

\underline{Estimate on $\phi_1$.} Consider the case of $\alpha\geq 1$. By Hardy inequality and boundary conditions $\phi_0|_{Y=0}=\p_Y\phi_0|_{Y=0}=0$, we have
\begin{align}
	\|\frac{Y\p_Y(A^{-1}\p_YU_s\phi_0)}{U_s}\|_{L^2}&\leq \|\frac{YA^{-1}\p_YU_s\p_Y\phi_0}{U_s}\|_{L^2}+\|\frac{Y\p_Y(A^{-1}\p_YU_s)\phi_0}{U_s}\|_{L^2}\leq C\|\p_Y\phi_0\|_{L^2},\nonumber\\
	\|\frac{\p_Y(A^{-1}\p_YU_s\phi_0)}{U_s}\|_{L^2}&\leq \|\frac{A^{-1}\p_YU_s\p_Y\phi_0}{U_s}\|_{L^2}+\|\frac{\p_Y(A^{-1}\p_YU_s)\phi_0}{U_s}\|_{L^2}\leq C\|\p_Y\phi_0\|_{L^2}+C\|\p_Y^2\phi_0\|_{L^2},\nonumber\\
	\|\frac{Y^2\p_Y^2(A^{-1}\p_YU_s\phi_0)}{U_s}\|_{L^2}&\leq C\|\p_Y\phi_0\|_{L^2}+C\|\p_Y^2\phi_0\|_{L^2}.\nonumber
\end{align}
Applying \eqref{4.106} and \eqref{4.107} to $\phi_1$ with $f=\partial_Y(A^{-1}\partial_YU_s\phi_0)$, we obtain
\begin{equation}
\begin{aligned}
	\|\p_Y\phi_1\|_{L^2}+\alpha\|\phi_1\|_{L^2}&\leq C\|\p_Y\phi_0\|_{L^2}+C\eps^{\frac13}\|\p_Y^2\phi_0\|_{L^2}\leq C\eps^{-\frac13}\|g\|_{L^2},\\
		\|\Delta_\alpha\phi_1\|_{L^2}&\leq C\|\p_Y\phi_0\|_{L^2}+C\|\p_Y^2\phi_0\|_{L^2}\leq C\eps^{-\frac23}\|g\|_{L^2}.\label{4.121}
\end{aligned}
\end{equation}
Here we have used bounds \eqref{4.120}. The case of $\alpha\in (0,1)$ can be treated
 in the same way by noting that $\|(1+Y)^2\p_Y(A^{-1}\p_YU_s\phi_0)\|_{L^2}\leq C\|\p_Y\phi_0\|_{L^2}$ and
$\int_0^\infty \p_Y(A^{-1}\p_YU_s\phi_0) dY=0. 
$
Finally, combining  \eqref{4.120} with \eqref{4.121} we get \eqref{4.113} and \eqref{4.114}. The proof of Proposition \ref{P4.6} is complete.
\end{proof}

Now we can solve the original Orr-Sommerfeld equation \eqref{4.122-1}. Let $\tilde{\phi}$ be the solution to the symmetrized Orr-Sommerfeld equation
$\widetilde{\text{OS}}_{\text{CNS}}(\tilde{\phi})=f.
$
\begin{corollary}
\label{C4.9}
Let $\alpha\eps^{\frac13}\lesssim1$. The Orr-Sommerfeld equation ${\rm OS_{CNS}}(\phi)=f$ admits a solution in the form of $\phi=\tilde{\phi}+\phi_r$, such that  $\phi_r\in H^4(\mathbb{R}_+)\cap H_0^1(\mathbb{R}_+)$, and
\begin{align}
	\|\p_Y\phi_r\|_{L^2}+\alpha \|\phi_r\|_{L^2}&\leq C\eps^{\frac23}\|\p_Y\tilde{\phi}\|_{H^1},\label{4.122-2}\\
	\|\Delta_\alpha\phi_r\|_{L^2}&\leq C\eps^{\frac13}\|\p_Y\tilde{\phi}\|_{H^1}.\label{4.122-3}
\end{align}
\end{corollary}
\begin{proof}
	We observe
${\rm OS_{CNS}}(\tilde{\phi})={\rm \widetilde{ OS}_{CNS}}(\tilde{\phi})+i\eps[\Lambda,\Delta_{\alpha}]\tilde{\phi}$, where the commutator 
$$[\Lambda,\Delta_{\alpha}]\tilde{\phi}=\Lambda \Delta_\alpha (\tilde{\phi})-\Delta_\alpha \Lambda(\tilde{\phi})=-\p_Y\left( 2\p_Y(A^{-1})\p_Y^2\tilde{\phi}+\p_Y^2(A^{-1})\p_Y\tilde{\phi}\right).
$$
To eliminate the error induced by $\tilde{\phi}$, we inductively define $\phi^k (k\geq 1)$ as solutions to
$$\left\{ \begin{aligned}
	&{\rm\widetilde{OS}_{CNS}}(\phi^k)=-i\eps[\Lambda,\Delta_\alpha]\phi^{k-1},\\
	&\phi^k|_{Y=0}=0,
\end{aligned}\right.
$$
with $\phi^0=\tilde{\phi}.$ Let $\phi_r=\sum_{k=1}^\infty\phi^k$.
 It is straightforward to see that $\phi\eqdef\tilde{\phi}+\phi_r$ is the solution that we desire. Now we show the convergence of series. By applying \eqref{4.113} and \eqref{4.114} to $\phi^k$ with $$g=i\eps\left( 2\p_Y(A^{-1})\p_Y^2{\phi}^{k-1}+\p_Y^2(A^{-1})\p_Y{\phi}^{k-1}\right),$$ we obtain
\begin{align}
	\|(\p_Y\phi^k,\alpha\phi^k)\|_{L^2}&\leq C\eps^{\frac23}\left(\|\p_Y\phi^{k-1}\|_{L^2}+\|\p_Y^2\phi^{k-1}\|_{L^2}\right),\nonumber\\
	\|\Delta_\alpha\phi^k\|_{L^2}&\leq C\eps^{\frac13}\left(\|\p_Y\phi^{k-1}\|_{L^2}+\|\p_Y^2\phi^{k-1}\|_{L^2}\right).\nonumber
\end{align}
Therefore, taking $\eps\in (0,1)$ suffciently small, the series $\phi^k$ converges in $H^2(\mathbb{R}_+)$. Moreover, we have
\begin{align}
	\sum_{k=1}^\infty\|(\p_Y\phi^k,\alpha\phi^k)\|_{L^2}+\eps^{\frac13}\sum_{k=1}^\infty\|\Delta_\alpha\phi^k\|_{L^2}\leq C\eps^{\frac23}\|\p_Y\tilde{\phi}\|_{H^1}.\nonumber
\end{align}
Then \eqref{4.122-2} and \eqref{4.122-3} follow. The proof is complete.
\end{proof}
\subsection{Solvability at high frequencies}
In this subsection, we establish the solvability of Orr-Sommerfeld equation \eqref{4.122-1} for high frequencies when  $\alpha \eps^{\frac13}\gg1$ (corresponding to $\hat{n}\gg \nu^{-\frac34}$). In this regime, the diffusion dominates so that we can solve it by a direct energy method.
\begin{proposition}\label{prop4.8}
Let $m\in (0,1)$, and $f\in L^{2}(\mathbb{R}_+)$. There exists $\delta_1\in (0,1)$, such that if $\alpha\eps^{\frac13}\geq \delta_1^{-1}$, the Orr-Sommerfeld equation \eqref{4.122-1} admits a solution $\phi\in H^4(\mathbb{R}_+)\cap H^1_0(\mathbb{R}_+)$, which satisfies the following bounds
	\begin{align}
		\|(\partial_Y\phi, \alpha\phi)\|_{L^2}&\le\frac{C}{\eps\alpha^3}\|f\|_{L^2},\label{4.239}\\
		\|(\partial_Y^2-\alpha^2)\phi\|_{L^2}&\le\frac{C}{\eps\alpha^2}\|f\|_{L^2},\label{4.240}\\
		\|\sqrt{U_s}(\p_Y\phi,\alpha\phi)\|_{L^2}&\le\frac{C}{\eps^{\frac12}\alpha^{2}}\|f\|_{L^2},\label{4.241}\\
		\|\partial_Y(\partial_Y^2-\alpha^2)\phi\|_{L^2}&\le \frac{C}{\eps \alpha}\|f\|_{L^2}.\label{4.343}
	\end{align}
\end{proposition}
\begin{proof}
	To determine the solution, we supplement \eqref{4.122-1} with one extra boundary condition $\Delta_\alpha \phi|_{Y=0}=0$. The existence of solutions is classical, therefore we only focus on a priori estimates \eqref{4.239}-\eqref{4.343}.	Taking inner product of \eqref{4.122-1} with $\bar{\phi}$ and integrating by parts, we obtain
	\begin{align}\label{4.242}
		&i\epsilon\int_0^\infty A^{-1}\left[|\partial_Y^2\phi|^2+\alpha^2(1+A^{-1})|\partial_Y\phi|^2+\alpha^4|\phi|^2\right]\;dY\notag\\
		&\quad-\int_0^\infty U_s\left(A^{-1}|\partial_Y\phi|^2+\alpha^2|\phi|^2\right)\;dY-\int_0^\infty \partial_Y(A^{-1}\partial_YU_s)|\phi|^2\;dY\notag\\
		&\qquad\qquad=\int_0^\infty A^{-1}\partial_YU_s\bar{\phi}\partial_Y\phi dY-i\epsilon\int_0^\infty\partial_YA^{-1}\partial_Y^2\phi \overline{\partial_Y\phi}dY+\int_0^\infty f \bar{\phi}d Y.
	\end{align}
By Cauchy-Schwarz, 
$$
\text{\rm The R.H.S of \eqref{4.242} }\leq\eps\|\p_Y^2\phi\|_{L^2}\|\p_Y\phi\|_{L^2}+\|\p_Y\phi\|_{L^2}\|\phi\|_{L^2}+\|f\|_{L^2}\|\phi\|_{L^2}.
$$
Then by taking imaginary part of \eqref{4.242} and using Young's inequality, we obtain
	\begin{align}\label{4.243}
		\epsilon\|(\partial_Y^2\phi, \alpha\partial_Y\phi, \alpha^2\phi)\|_{L^2}^2&\le\eps\|\p_Y^2\phi\|_{L^2}\|\p_Y\phi\|_{L^2}+\|\p_Y\phi\|_{L^2}\|\phi\|_{L^2}+\|f\|_{L^2}\|\phi\|_{L^2}\nonumber\\
		&\leq \frac{\eps}{2}\|(\partial_Y^2\phi, \alpha\partial_Y\phi, \alpha^2\phi)\|_{L^2}^2+C\eps\left( \frac{1}{\alpha^2}\|\alpha\p_Y\phi\|_{L^2}^2+\frac{1}{(\alpha \eps^{\frac13})^6}\|\alpha^2\phi\|_{L^2}^2\right)+\frac{C}{\eps\alpha^4}\|f\|_{L^2}^2.
	\end{align}
By taking $\delta_1\in (0,1)$ suitably small such that $\frac{C}{\alpha^2}\leq C\delta_1^2\eps^{\frac23}<\frac14$ and $\frac{C}{(\alpha\eps^{\frac13})^6}\leq C\delta_1^6<\frac14$, we can absorb the first two terms on the R.H.S of \eqref{4.243} by the L.H.S.. Thus the bounds \eqref{4.239} and \eqref{4.240} follow from \eqref{4.243}. Moreover, the real part of \eqref{4.242} gives
	\begin{align*}
		\|\sqrt{U_s}(\partial_Y\phi, \alpha\phi)\|_{L^2}^2&\leq C\|\partial_Y\phi\|_{L^2}\|\phi\|_{L^2}+C\epsilon\|\partial_Y^2\phi\|_{L^2}\|\partial_Y\phi\|_{L^2}+C\|f\|_{L^2}\|\phi\|_{L^2}.
	\end{align*}
Using \eqref{4.239} and \eqref{4.240} we get
$$\|\sqrt{U_s}(\partial_Y\phi, \alpha\phi)\|_{L^2}^2\leq C\left(  \frac{1}{\eps^2\alpha^7}+\frac{1}{\eps\alpha^5}+\frac{1}{\eps\alpha^4}\right)\|f\|_{L^2}^2\leq \frac{C}{\eps\alpha^4}\|f\|_{L^2}^2,
$$
where we have used $\alpha\eps^{\frac13}\gtrsim1$. Estimate \eqref{4.241} is proved.

Finally, we turn to \eqref{4.343}. Note that
	\begin{align*}
		\Lambda(\phi)=A^{-1}\Delta_\alpha\phi+(A^{-1}-1)\alpha^2\phi+\partial_YA^{-1}\partial_Y\phi.
	\end{align*}
	Then taking inner product of \eqref{4.122-1} with $\Delta_\alpha\bar{\phi}$ and integrating by parts, we obtain
\begin{equation}
	\begin{aligned}
		&-i\epsilon\int_0^\infty A^{-1}|\partial_Y\Delta_\alpha\phi|^2+\alpha^2|\Delta_\alpha\phi|^2\;dY+\|\sqrt{A^{-1}U_s}\Delta_\alpha\phi\|_{L^2}^2\\
		&\quad\qquad-\int_0^\infty (A^{-1}-1)U_s\alpha^2(|\partial_Y\phi|^2+\alpha^2 |\phi|^2)\;dY\\
		&\quad\quad=-\int_0^\infty U_s\partial_YA^{-1}\partial_Y\phi\overline{\Delta_\alpha\phi} dY
	-\alpha^2\int_0^\infty\partial_Y[U_s(A^{-1}-1)] \phi\overline{ \partial_Y\phi }dY+\int_0^\infty f \overline{\Delta_\alpha\phi} dY.\label{4.244}
	\end{aligned}
\end{equation}
By using $|\partial_Y[U_s(A^{-1}-1)]|\leq CU_s$, and bounds \eqref{4.239}-\eqref{4.241}, we obtain
$$
\begin{aligned}
\text{\rm R.H.S of \eqref{4.244}}&\leq C\|\p_Y\phi\|_{L^2}\|\Delta_\alpha\phi\|_{L^2}+C\alpha\|\sqrt{U_s}(\p_Y\phi,\alpha\phi)\|_{L^2}^2+C\|f\|_{L^2}\|\Delta_\alpha\phi\|_{L^2}\nonumber\\
&\leq C\left(\frac{1}{\eps^2\alpha^5}+\frac{1}{\eps\alpha^3}+\frac{1}{\eps\alpha^2}\right)\|f\|_{L^2}^2\leq \frac{C}{\eps\alpha^2}\|f\|_{L^2}^2,
\end{aligned}
$$
where we have used $\alpha\eps^{\frac13}\gtrsim 1$. Then by taking imaginary part of \eqref{4.244}, we obtain \eqref{4.343}.	The proof of Proposition \ref{prop4.8} is complete.
\end{proof}
\section{Inhomogeneous quasi-compressible system}\label{S4}

In this section, we recover the solution $(\varrho,\mathfrak{u},\mathfrak{v})$ to the quasi-compressible system \eqref{4.3}. Let $\phi\in H^2(\mathbb{R}_+)$ be a solution to the compressible Orr-Sommerfeld equation \eqref{4.6}. Note that $\phi$ does not necessarily vanish on the boundary. Recall the density $\varrho$ and the velocity field $(\mathfrak{u},\mathfrak{v})$ defined in \eqref{4.5} and \eqref{4.4} resspectively. For simplicity, we set the divergence field $\mathcal{D}_q\eqdef \div_{\alpha}(\mathfrak{u},\mathfrak{v})$. In this section, the symbol $\langle f,g\rangle$ represents the standard $L^2(\mathbb{R}_+)$ inner product with respect to $Y$ variable.
 
\begin{proposition}\label{P8.1}
Let $(f_{\mathfrak{u}},f_{\mathfrak{v}})\in L^2(\mathbb{R}_+)^2$ and $\p_Yf_{\mathfrak{v}}\in L^2(\mathbb{R}_+)$. Suppose that $m\in (0,1)$, and $\alpha\eps^{\frac13}\lesssim 1$. Then $(\varrho,\mathfrak{u},\mathfrak{v})\in H^2(\mathbb{R}_+)^3$ is a solution to the system \eqref{4.3}. Moreover, it satisfies
\begin{align}
	\|\p_Y^k\mathfrak{v}\|_{L^2}&\lesssim \alpha\|\p_Y^k\phi\|_{L^2},~k=0,1,2,\label{8.0-1}\\
	m^{-2}\|\varrho\|_{L^2}&\lesssim \eps(1+\alpha^2)\|\Delta_\alpha\phi\|_{L^2}+\alpha^2\|\phi\|_{L^2}+\|U_s'\phi\|_{L^2}+\left(1+\frac{1}{\alpha} \right)\|(f_{\mathfrak{u}},f_{\mathfrak{v}})\|_{L^2},\label{8.0-2}\\
	m^{-2}\|(\p_Y\varrho,\Delta_\alpha\varrho)\|_{L^2}&\lesssim  \alpha^2\left( \eps\|\Delta_\alpha\phi\|_{L^2}+\|\phi\|_{L^2}  \right)+(1+\alpha)\|(f_{\mathfrak{u}},f_{\mathfrak{v}})\|_{L^2}+\|\p_Yf_\mathfrak{v}\|_{L^2},\label{8.0-3}\\
	\|\mathfrak{u}\|_{L^2}&\lesssim \|\p_Y\phi\|_{L^2}+C\|U_s'\phi\|_{L^2}+\eps\|\Delta_\alpha\phi\|_{L^2}+\frac{1}{\alpha}\|(f_{\mathfrak{u}},f_{\mathfrak{v}})\|_{L^2},\label{8.0-4}\\
	\|\p_Y\mathfrak{u}\|_{L^2}&\lesssim \|\Delta_\alpha\phi\|_{L^2}+\alpha^2\|\phi\|_{L^2}+(1+\alpha)\|U_s'\phi\|_{L^2}+\left(1+\frac{1}{\alpha}\right)\|(f_{\mathfrak{u}},f_{\mathfrak{v}})\|_{L^2},\label{8.0-5}\\
\|\Delta_\alpha\mathfrak{u}\|_{L^2}&\lesssim (1+\alpha^2)\|\Delta_{\alpha}\phi\|_{L^2}+\frac{1}{\eps}\left(\|\p_Y\phi\|_{L^2}+\alpha^2\|\phi\|_{L^2}+\|U_s'\phi\|_{L^2}\right)\nonumber\\
&\qquad+	\frac{1}{\eps}\left(1+\frac{1}{\alpha}\right)\|(f_{\mathfrak{u}},f_{\mathfrak{v}})\|_{L^2}+\|\p_Yf_{\mathfrak{v}}\|_{L^2},\label{8.0-6}\\
	\|(\partial_Y\mathcal{D}_q,\alpha \mathcal{D}_q)\|_{L^2}&\lesssim \sqrt{\nu}(1+\alpha)\|\Delta_\alpha\phi\|_{L^2}+\alpha^2\|\phi\|_{L^2}+\alpha\|U_s'\phi\|_{L^2}+(1+\alpha)\|(f_{\mathfrak{u}},f_{\mathfrak{v}})\|_{L^2}.\label{8.0-7}
\end{align}
\end{proposition}

\begin{proof}
	It is straightforward to check that $(\varrho,\mathfrak{u},\mathfrak{v})$ is a solution to \eqref{4.3}. Therefore, we only focus on the estimates \eqref{8.0-1}-\eqref{8.0-7}. The proof involves
	four steps.

	\underline{\it Step 1. Estimates on $\mathfrak{v}$}. From $\mathfrak{v}=-i\alpha\phi$ we can directly get
	\begin{align}
		\|\p_Y^k\mathfrak{v}\|_{L^2}\leq \alpha\|\p_Y^k\phi\|_{L^2},\nonumber
	\end{align}
 which is \eqref{8.0-1}.
 
	\underline{\it Step 2. Estimates on $\varrho$.}
	Unfortunately, directly estimating $\|\varrho\|_{L^2}$ from \eqref{4.5} is not good enough for  later use. The key point to use cancellation in quasi-compressible system \eqref{4.3} to bound $\|\p_Y\varrho\|_{L^2}$ and $\|U_s\varrho\|_{L^2}$. Then $\|\varrho\|_{L^2}$ is recoverd by interpolation.
	
	 Differentiating \eqref{4.5} yields
	\begin{align}
		m^{-2}\p_Y\varrho&=-i\eps\p_Y\left(A^{-1}\p_Y\Delta_\alpha\phi\right)-\p_Y\left[A^{-1}\left(U_s\p_Y\phi-\phi \p_YU_s\right)\right]-i\alpha^{-1}\p_Y\left(A^{-1}f_{\mathfrak{u}}\right)\nonumber\\
		&=-{\rm OS_{CNS}}(\phi)-i\eps\alpha^2\Delta_\alpha\phi-\alpha^2U_s\phi-i\alpha^{-1}\p_Y\left(A^{-1}f_{\mathfrak{u}}\right)\nonumber\\
		&=-i\eps\alpha^2\Delta_\alpha\phi-\alpha^2U_s\phi+f_{\mathfrak{v}},\nonumber
	\end{align}
where the equation \eqref{4.6} has been used in the last equality. Then we obtain
\begin{align}
	m^{-2}\|\p_Y\varrho\|_{L^2}\leq \alpha^2\left(\eps\|\Delta_\alpha\phi\|_{L^2}+\|\phi\|_{L^2}\right)+\|f_{\mathfrak{v}}\|_{L^2}.\label{8.2}
\end{align}

Next we derive the weighted estimate $\|U_s\varrho\|_{L^2}$. We rewrite momentum equations of \eqref{4.3} in terms of $\phi$ by substituting \eqref{4.4} into $\eqref{4.3}_2$ and $\eqref{4.3}_3$:
\begin{equation}\label{8.2-1}
\begin{aligned}
	-\sqrt{\nu}\Delta_\alpha\p_Y\phi-i\alpha\phi \p_YU_s+i\alpha U_s\p_Y\phi+i\alpha m^{-2}A(Y)\varrho&=f_{\mathfrak{u}},\\
	i\alpha\sqrt{\nu}\Delta_\alpha\phi+\alpha^2U_s\phi+m^{-2}\p_Y\varrho&=f_{\mathfrak{v}}.
\end{aligned}
\end{equation}
Taking inner product with $m^{-2}U_s^2\overline{i\alpha\varrho}$ and $m^{-2}U_s^2\p_Y\bar{\varrho}$ respectively, we get
\begin{align}
	m^{-4}\|U_s(\p_Y\varrho,\alpha{A}^{\frac12}\varrho)\|_{L^2}^2=&\underbrace{\langle \sqrt{\nu}\Delta_\alpha\p_Y\phi,m^{-2}U_s^2\overline{i\alpha\varrho} \rangle+\langle -i\alpha\sqrt{\nu}\Delta_\alpha\phi,m^{-2}U_s^2\p_Y\bar{\varrho} \rangle}_{I_1}\nonumber\\
	&+\underbrace{\langle -i\alpha U_s\p_Y\phi,m^{-2}U_s^2\overline{i\alpha\varrho} \rangle+\langle -\alpha^2U_s\phi,m^{-2}U_s^2\p_Y\bar{\varrho} \rangle}_{I_2}\nonumber\\
	&+\underbrace{\langle i\alpha\phi U_s',m^{-2}U_s^2\overline{i\alpha\varrho}  \rangle}_{I_3}+\underbrace{\langle f_{\mathfrak{u}},m^{-2}U_s^2\overline{i\alpha\varrho} \rangle+\langle f_{\mathfrak{v}},m^{-2}U_s^2\p_Y\bar{\varrho} \rangle}_{I_4}.\label{8.3}
\end{align}

Now we estimate $I_1-I_4$ term by term. For $I_1$, the weight $U_s$ allows for integrating by parts:
\begin{align}
|I_1|=2\left|\langle \sqrt{\nu}\Delta_\alpha\phi, m^{-2} U_sU_s'i\alpha \rho \rangle\right|\leq C\sqrt{\nu}\|\Delta_\alpha\phi\|_{L^2}\|m^{-2}U_s\alpha\varrho\|_{L^2}.\nonumber
\end{align}
Similarly, we have
\begin{align}
	|I_2|=3\left|\langle  i\alpha \phi,U_s^2U_s'i\alpha\varrho   \rangle\right|\leq C\alpha\|U_s'\phi\|_{L^2}\|m^{-2}U_s\alpha\varrho\|_{L^2}.\nonumber
\end{align}
By Cauchy-Schwarz, we obtain
\begin{align}
	|I_3|\leq C\alpha\|U_s'\phi\|_{L^2}\|m^{-2}U_s\alpha\varrho\|_{L^2},\nonumber
\end{align}
and
\begin{align}
	|I_4|\leq Cm^{-2
	}\|(f_{\mathfrak{u}},f_{\mathfrak{v}})\|_{L^2}\|U_s(\p_Y\varrho,\alpha\rho)\|_{L^2}.\nonumber
\end{align}
Substituting these bounds into \eqref{8.3} and noting that $A^{-1}\sim 1$ when $m\in (0,1)$,  we get
\begin{align}
	m^{-2}\|U_s(\p_Y\varrho,\alpha\varrho)\|_{L^2}\leq C\sqrt{\nu}\|\Delta_\alpha\phi\|_{L^2}+C\alpha\|U_s'\phi\|_{L^2}+C\|(f_{\mathfrak{u}},f_{\mathfrak{v}})\|_{L^2}. \label{8.4}
\end{align}
The bound $\|\varrho\|_{L^2}$ can be recovered by the interpolation
\begin{align}%\label{8.5}
	m^{-2}\|\varrho\|_{L^2}&\leq m^{-2}\left(
	\|\p_Y\varrho\|_{L^2}+\|U_s\varrho\|_{L^2}\right)\nonumber\\
	&\leq C\eps(1+\alpha^2)\|\Delta_\alpha\phi\|_{L^2}+C\alpha^2\|\phi\|_{L^2}+C\|U_s'\phi\|_{L^2}+C\left(1+\frac{1}{\alpha}\right)\|(f_{\mathfrak{u}},f_{\mathfrak{v}})\|_{L^2},\label{8.4-2}
\end{align}
which is \eqref{8.0-2}. Moreover, from \eqref{8.2-1} we have
\begin{align}
	m^{-2}\Delta_\alpha\varrho=-\alpha^2U_s^2\varrho-2\alpha^2U_s'\phi+i\alpha f_{\mathfrak{u}}+\p_Yf_{\mathfrak{v}},
	\label{8.6}
\end{align}
which implies
\begin{align}%\label{8.7}
	m^{-2}\|\Delta_\alpha\rho\|_{L^2}&\leq C\alpha^2\|U_s\varrho\|_{L^2}+C\alpha^2\|U_s'\phi\|_{L^2}+C\alpha\|f_{\mathfrak{u}}\|_{L^2}+C\|\p_Yf_\mathfrak{v}\|_{L^2}\nonumber\\
	&\leq C\alpha^2\left( \eps\|\Delta_\alpha\phi\|_{L^2}+\|\phi\|_{L^2}  \right)+C\alpha\|(f_{\mathfrak{u}},f_{\mathfrak{v}})\|_{L^2}+C\|\p_Yf_\mathfrak{v}\|_{L^2}.\nonumber
\end{align}
Here we have used \eqref{8.4} in the last inequality. Combining it with \eqref{8.2} yields \eqref{8.0-3}. 

\underline{\it Step 3. Estimate on $\mathfrak{u}$}. Recall from \eqref{4.4} that $\mathfrak{u}=\p_Y\phi-U_s\varrho$. Then it holds that
\begin{align}%\label{8.8}
	\|\mathfrak{u}\|_{L^2}&\leq \|\p_Y\phi\|_{L^2}+\|U_s\varrho\|_{L^2}\nonumber\\
	&\leq C\eps\|\Delta_\alpha\phi\|_{L^2}+C\|\p_Y\phi\|_{L^2}+C\|U_s'\phi\|_{L^2}+\frac{C}{\alpha}\|(f_{\mathfrak{u}},f_{\mathfrak{v}})\|_{L^2},\nonumber
\end{align}
which implies \eqref{8.0-4}. Similarly, using \eqref{8.4} and \eqref{8.4-2}, we have
\begin{align}
	&\|\p_Y\mathfrak{u}\|_{L^2}\leq \|\p_Y^2\phi\|_{L^2}+\|U_s\p_Y\varrho\|_{L^2}+\|\varrho\|_{L^2}\nonumber\\
	&\quad\leq C\left(1+\sqrt{\nu}+\eps(1+\alpha^2)\right)\|\Delta_\alpha\phi\|_{L^2}+C\alpha^2\|\phi\|_{L^2}+(1+\alpha)\|U_s'\phi\|_{L^2}+C\left(1+\frac{1}{\alpha}\right)\|(f_{\mathfrak{u}},f_{\mathfrak{v}})\|_{L^2}\nonumber\nonumber\\
	&\quad\leq C\|\Delta_\alpha\phi\|_{L^2}+C\alpha^2\|\phi\|_{L^2}+(1+\alpha)\|U_s'\phi\|_{L^2}+C\left(1+\frac{1}{\alpha}\right)\|(f_{\mathfrak{u}},f_{\mathfrak{v}})\|_{L^2},\nonumber %\label{8.9}
\end{align}
where $\alpha\eps^{\frac13}\lesssim 1$ has been used in the last inequality. Thus \eqref{8.0-5} follows. Moreover, from $\eqref{4.3}_2$ we have
\begin{align}
\|\Delta_\alpha \mathfrak{u}\|_{L^2}\leq& C\|\Delta_\alpha\varrho\|_{L^2}+C\|\varrho\|_{H^1}+\frac{C}{\eps}\left(\|\mathfrak{u}\|_{L^2}+\|U_s'\phi\|_{L^2}+m^{-2}\|\varrho\|_{L^2}\right)+\frac{C}{\sqrt{\nu}}\|f_{\mathfrak{u}}\|_{L^2}\nonumber\\
\leq& C(1+\alpha^2)\|\Delta_{\alpha}\phi\|_{L^2}+\frac{C}{\eps}\left(\|\p_Y\phi\|_{L^2}+\alpha^2\|\phi\|_{L^2}+\|U_s'\phi\|_{L^2}\right)\nonumber\\
&+\frac{C}{\eps}\left(1+\frac{1}{\alpha}\right)\|(f_{\mathfrak{u}},f_{\mathfrak{v}})\|_{L^2}+\|\p_Yf_{\mathfrak{v}}\|_{L^2},\nonumber %\label{8.10}
\end{align}
which is \eqref{8.0-6}. Thus we have completed the estimates on $\mathfrak{u}$.

\underline{\it Step 4. Estimate on $\div_\alpha(\mathfrak{u},\mathfrak{v})$.} Recall the notation $\mathcal{D}_q\eqdef\div_\alpha(\mathfrak{u},\mathfrak{v})$. From $\eqref{4.3}_1$, we have $\mathcal{D}_q=-i\alpha U_s\varrho$. Then it holds that
\begin{align}
	\|(\partial_Y\mathcal{D}_q,\alpha \mathcal{D}_q)\|_{L^2}\leq \alpha\|U_s(\partial_Y\varrho,\alpha\varrho)\|_{L^2}+\alpha\|U_s'\varrho\|_{L^2}.\label{8.6-1}
\end{align}

It sufficies to bound $\|U_s'\rho\|_{L^2}$. We divide the discussion into two cases. If $\alpha \in (0,1)$, we use $\|U_s'\varrho\|_{L^2}\leq C\|\varrho\|_{L^2}$, and bounds \eqref{8.4}, \eqref{8.4-2} to obtain 
\begin{align}
\|(\partial_Y\mathcal{D}_q,\alpha \mathcal{D}_q)\|_{L^2}&\lesssim \alpha\|U_s(\partial_Y\varrho,\alpha\varrho)\|_{L^2}+\alpha\|\varrho\|_{L^2}\nonumber\\
&\lesssim \sqrt{\nu}\|\Delta_\alpha \phi\|_{L^2}+\alpha^2\|\phi\|_{L^2}+\alpha\|U_s'\phi\|_{L^2}+\|(f_{\mathfrak{u}},f_\mathfrak{v})\|_{L^2}.\nonumber
\end{align}
The estimate of $\|U_s'\varrho\|_{L^2}$ for the case $\alpha\geq 1$ is more subtle. According to the structural condition \eqref{2.1}, we can find $Y_0>0$ and $\eta>0$, such that
$$U_s'(Y)\geq \frac12~{\rm in }~[0,2Y_0],~{\rm and }~ U_s(Y)\geq \eta~{\rm for }~Y\geq Y_0.
$$ 
We introduce a cut-off function $\chi(Y)\in [0,1]$,  with $\chi(Y)\equiv 1$ when $Y\in [0,Y_0]$ and $\chi(Y)\equiv 0$ when $Y\geq 2Y_0.$
Then we have
$$\begin{aligned}
	\|U_s'\varrho\|_{L^2}^2&=\int_0^\infty |U_s'|^2|\varrho|^2\chi dY+\int_0^\infty|U_s'|^2|\varrho|^2(1-\chi) dY\nonumber\\
	&= \int_0^{2Y_0}U_s\left(U_s'\chi\right)'|\varrho|^2 dY+ {\rm Re}\left(\int_0^{2Y_0}U_sU_s'\chi\varrho \partial_Y\bar{\varrho} dY\right)+\int_{Y_0}^\infty \frac{|U_s'|^2(1-\chi)}{U_s^2}|U_s\varrho|^2dY\nonumber\\
	&\leq C\|U_s\varrho\|_{L^2}\|\varrho\|_{H^1}+C\|U_s\varrho\|_{L^2}^2\leq C\eps^2\alpha^2\|\Delta_\alpha\phi\|_{L^2}^2+C\alpha^2\|\phi\|_{L^2}^2+C\|(f_{\mathfrak{u}},f_{\mathfrak{v}})\|_{L^2}^2,
\end{aligned}
$$
where we have used \eqref{8.2}, \eqref{8.4} and \eqref{8.4-2} in the last inequality. Substituting this into \eqref{8.6-1} and using \eqref{8.4} again, we can obtain \eqref{8.0-7}. Therefore,  the Proposition \ref{P8.1} is complete.
\end{proof}
Now we construct the inhomogeneous quasi-compressible solution with the boundary condition $\mathfrak{v}|_{Y=0}=0.$
\begin{corollary}
\label{P4.9}
Under the same assumption as Proposition \ref{P8.1},  
the system \eqref{4.3} with boundary condition $\mathfrak{v}|_{Y=0}=0$ admits a solution $(\varrho, \mathfrak{u}, \mathfrak{v})\in H^2(\mathbb{R}_+)^3$ satisfying the following estimates:
\begin{align}
\alpha\|\mathfrak{u}, \mathfrak{v}\|_{L^2}+\| \partial_Y\mathfrak{v}\|_{L^2}&\le C\eps^{-\frac13}\|(f_{\mathfrak{u}},f_{\mathfrak{v}})\|_{L^2},\label{4.127}\\
m^{-2}\|\varrho\|_{H^1}&\le C\eps^{-\frac13}\left(1+\frac{1}{\alpha}\right)\|(f_{\mathfrak{u}},f_{\mathfrak{v}})\|_{L^2},\label{4.126}\\
\alpha\|\partial_Y\mathfrak{u}\|_{L^2}+\|\partial_Y^2\mathfrak{v}\|_{L^2}&\le C\eps^{-\frac23}\|(f_{\mathfrak{u}},f_{\mathfrak{v}})\|_{L^2}, \label{4.128}\\
m^{-2}\|\Delta_\alpha\varrho\|_{L^2}&\le C\eps^{-\frac13}\|(f_{\mathfrak{u}},f_{\mathfrak{v}})\|_{L^2}+C\|\partial_Yf_{\mathfrak{v}}\|_{L^2}, \label{4.130}\\
\|\partial_Y^2\mathfrak{u}\|_{L^2}&\le C\eps^{-\frac43}\left(1+\frac1\alpha\right)\|(f_{\mathfrak{u}},f_{\mathfrak{v}})\|_{L^2}+C\|\partial_Yf_{\mathfrak{v}}\|_{L^2},\label{4.131}\\
	\|(\partial_Y\mathcal{D}_q,\alpha \mathcal{D}_q)\|_{L^2}&\leq C\eps^{-\frac13}\|(f_{\mathfrak{u}},f_{\mathfrak{v}})\|_{L^2}.\label{4.131-1}
\end{align}
\end{corollary}
\begin{proof}
Recall Corollary \ref{C4.9}. We construct the solution to the Orr-Sommerfeld equation \eqref{4.6} in the following form: $ \phi=\tilde{\phi}+\phi_r$, where $\tilde{\phi}$ is the solution to symmetrized Orr-Sommerfeld equation ${\rm\widetilde{ OS}_{CNS}}(\tilde{\phi})=-f_{\mathfrak{v}}-\frac{i}{\alpha}\partial_Y(A^{-1}f_{\mathfrak{u}}),$ and $\phi_r$ is the remainder. Then $(\varrho,\mathfrak{u},\mathfrak{v})$ in Proposition \ref{P8.1} defines a solution to \eqref{4.3} with $\mathfrak{v}|_{Y=0}=0$. Now we prove the estimates \eqref{4.127}-\eqref{4.131}.

Firstly,  applying Proposition \ref{P4.6} to $\tilde{\phi}$ yields
\begin{align}
	\|\p_Y\tilde{\phi}\|_{L^2}+\alpha\|\tilde{\phi}\|_{L^2}&\leq \frac{C}{\eps^{\frac13}\alpha}\|(f_{\mathfrak{u}},f_{\mathfrak{v}})\|_{L^2},\nonumber\\
	\|\Delta_\alpha\tilde{\phi}\|_{L^2}&\leq \frac{C}{\eps^{\frac23}\alpha}\|(f_{\mathfrak{u}},f_{\mathfrak{v}})\|_{L^2}.\nonumber
\end{align}
Then by using bounds of $\phi_r$ in Corollary \ref{C4.9}, we can obtain
\begin{align}
	\|\p_Y\phi\|_{L^2}+\alpha\|\phi\|_{L^2}&\leq C\|\p_Y\tilde{\phi}\|_{L^2}+C\alpha\|\tilde{\phi}\|_{L^2}+C\eps^{\frac23}\|\Delta_{\alpha}\tilde{\phi}\|_{L^2}\leq \frac{C}{\eps^{\frac13}\alpha}\|(f_{\mathfrak{u}},f_{\mathfrak{v}})\|_{L^2}\label{4.132},\\
	\|\Delta_\alpha\phi\|_{L^2}&\leq C\|\Delta_\alpha\tilde{\phi}\|_{L^2}+C\eps^{\frac13}\|\p_Y\tilde{\phi}\|_{L^2}+C\eps^{\frac13}\alpha^2\|\phi\|_{L^2}\nonumber\\
	&\leq C\left(\frac{1}{\eps^{\frac{2}{3}}\alpha}+\frac{1}{\alpha}+\frac{1}{\eps^{\frac13}}\right)\|(f_{\mathfrak{u}},f_{\mathfrak{v}})\|_{L^2}\leq \frac{C}{\eps^{\frac23}\alpha}\|(f_{\mathfrak{u}},f_{\mathfrak{v}})\|_{L^2},\label{4.133}\\
	\|U_s'\phi\|_{L^2}&\leq C\|\p_Y\phi\|_{L^2}\leq \frac{C}{\eps^{\frac13}\alpha}\|(f_{\mathfrak{u}},f_{\mathfrak{v}})\|_{L^2}.\label{4.134}
\end{align}
Substituting these bounds into Proposition \ref{P8.1}, we can prove \eqref{4.127}-\eqref{4.131-1}: we only show \eqref{4.131} since other estimates can be obtained similarly. From \eqref{8.0-6} and \eqref{4.132}-\eqref{4.134}, we obtain
\begin{align}
	\|\p_Y^2\mathfrak{u}\|_{L^2}&\leq C\left(\frac{1+\alpha^2}{\eps^{\frac23}\alpha}+\frac{1}{\eps^{\frac43}}\left(1+\frac1\alpha\right)+\frac{1}{\eps}\left(1+\frac1\alpha\right)\right)\|(f_{\mathfrak{u}},f_{\mathfrak{v}})\|_{L^2}+\|\p_Yf_{\mathfrak{v}}\|_{L^2}\nonumber\\
	&\leq \frac{C}{\eps^{\frac43}}\left(1+\frac1\alpha\right)\|(f_{\mathfrak{u}},f_{\mathfrak{v}})\|_{L^2}+\|\p_Yf_{\mathfrak{v}}\|_{L^2},\nonumber
\end{align}
where $\alpha\eps^{\frac13}\lesssim 1$ has been used in the last inequality. The proof of Corollary \ref{P4.9} is complete.	
\end{proof}
Finally, we solve the quasi-compressible system in the high-frequency regime.
\begin{proposition}\label{prop4.2}
	Let $(f_{\mathfrak{u}},f_{\mathfrak{v}})\in L^2(\mathbb{R}_+)^2$ and $\partial_Y f_{\mathfrak{v}}\in L^2(\mathbb{R}_+)$. Assume that $m\in (0,1)$, and $\alpha\geq 1$. There exists a solution $(\varrho,\mathfrak{u},\mathfrak{v})\in H^2(\mathbb{R}_+)^3$ to the system \eqref{4.3} with boundary condition $\mathfrak{v}|_{Y=0}=0$. Moreover, $(\varrho,\mathfrak{u},\mathfrak{v})$ satisfies
\begin{align}
	\|\p_Y^k\mathfrak{v}\|_{L^2}&\lesssim \alpha\|\p_Y^k\phi\|_{L^2},~k=0,1,2,\label{8.1-1}\\
	m^{-2}\|\varrho\|_{L^2}&\lesssim \eps\|\partial_Y\Delta_\alpha\phi\|_{L^2}+\|\phi\|_{H^1}+\alpha^{-1}\|f_{\mathfrak{u}}\|_{L^2},\label{8.1-3}\\
	m^{-2}\|(\partial_Y\varrho,\Delta_\alpha\varrho)\|_{L^2}&\lesssim\eps\alpha^2\|\Delta_\alpha\phi\|_{L^2}+\alpha^2\|\phi\|_{L^2}+\alpha\|(f_{\mathfrak{u}},f_{\mathfrak{v}})\|_{L^2}+\|\partial_Yf_{\mathfrak{v}}\|_{L^2},\label{8.1-2}\\
	\|\mathfrak{u}\|_{L^2}&\lesssim \eps\|\Delta_\alpha\phi\|_{L^2}+\|\phi\|_{H^1}+{\alpha}^{-1}\|(f_{\mathfrak{u}},f_{\mathfrak{v}})\|_{L^2},\label{8.1-4}\\
	\|\p_Y\mathfrak{u}\|_{L^2}&\lesssim \eps\|\partial_Y\Delta_\alpha\phi\|_{L^2}+\|\Delta_\alpha\phi\|_{L^2}+\|\partial_Y\phi\|_{L^2}+\alpha^2\|\phi\|_{L^2}+\|(f_{\mathfrak{u}},f_{\mathfrak{v}})\|_{L^2},\label{8.1-5}\\
	\|\Delta_\alpha\mathfrak{u}\|_{L^2}&\lesssim \|\partial_Y\Delta_\alpha\phi\|_{L^2}+\eps\alpha^2\|\Delta_\alpha\phi\|_{L^2}+\|\partial_Y\phi\|_{L^2}+\alpha^2\|\phi\|_{L^2}\nonumber\\
	&\qquad+\alpha\|(f_{\mathfrak{u}},f_{\mathfrak{v}})\|_{L^2}+\|\partial_Yf_{\mathfrak{v}}\|_{L^2}\label{8.1-6}\\
	\|(\partial_Y\mathcal{D}_q,\alpha\mathcal{D}_q)\|_{L^2}&\lesssim\eps\alpha\|\partial_Y\Delta_\alpha\phi\|_{L^2}+\eps\alpha^2\|\Delta_\alpha\phi\|_{L^2}+\alpha\|\partial_Y\phi\|_{L^2}+\alpha^2\|\phi\|_{L^2}+\alpha\|(f_{\mathfrak{u}},f_{\mathfrak{v}})\|_{L^2}.\label{8.1-7}
\end{align}
\end{proposition}
\begin{proof}
	The estimate \eqref{8.1-1} is obvious. For $\varrho,$ we use \eqref{4.5} to get
	\begin{align}
		m^{-2}\|\varrho\|_{L^2}\lesssim \eps\|\partial_Y\Delta_\alpha\phi\|_{L^2}+\|\phi\|_{H^1}+\alpha^{-1}\|f_{\mathfrak{u}}\|_{L^2}\nonumber
	\end{align}
which is \eqref{8.1-3}. The estimate $\|\partial_Y\varrho\|_{L^2}$ is obtained in \eqref{8.2}. Moreover, from \eqref{8.4} we get
\begin{align}
m^{-2}\|U_s(\partial_Y\varrho,\alpha\varrho)\|_{L^2}\lesssim \sqrt{\nu}\|\Delta_\alpha\phi\|_{L^2}+\alpha\|\phi\|_{L^2}+\|(f_u,f_v)\|_{L^2}.\label{8.4-1}
\end{align}
Combining this with \eqref{8.6}, we obtain
\begin{align}
	m^{-2}\|\Delta_\alpha\varrho\|_{L^2}&\lesssim \alpha^2\|U_s\varrho\|_{L^2}+\alpha^2\|\phi\|_{L^2}+\alpha\|f_{\mathfrak{u}}\|_{L^2}+\|\partial_Yf_{\mathfrak{v}}\|_{L^2}\nonumber\\
	&\lesssim \eps\alpha^2\|\Delta_{\alpha}\phi\|_{L^2}+\alpha^2\|\phi\|_{L^2}+\alpha\|(f_{\mathfrak{u}},f_{\mathfrak{v}})\|_{L^2}+\|\partial_Yf_{\mathfrak{v}}\|_{L^2},\nonumber
\end{align}
which is \eqref{8.1-2}.
 The estimates \eqref{8.1-4}-\eqref{8.1-6} on $\mathfrak{u}$ directly follow from the equality $\mathfrak{u}=\partial_Y\phi-U_s\varrho$, and bounds \eqref{8.1-3}, \eqref{8.1-2}, and \eqref{8.4-1} of the density $\varrho$. The proof of Proposition \ref{prop4.2} is complete.
\end{proof}

\section{Homogeneous quasi-compressible solutions} \label{S5}
In this section, we construct non-trivial solutions to the following homogeneous quasi-compressible system
\begin{align}\label{4.151}
\begin{cases}
\mathcal{L}_Q(\varrho_H, \mathfrak{u}_H, \mathfrak{v}_H)={\bf{0}},\\
\mathfrak{v}_H|_{Y=0}=0,\\
\mathfrak{u}_H|_{Y=0}\neq 0.
\end{cases}
\end{align}
These solutions will be used in Section \ref{S7} to construct boundary sublayer correctors to recover the no-slip boundary conditions in \eqref{4.1}. We divide the construction into three regimes: low frequencies $\hat{n}\ll \nu^{-\frac34}$, middle frequencies $\hat{n}\sim \nu^{-\frac34}$ and high frequencies $\hat{n}\gg \nu^{-\frac34}$.

\subsection{Low-frequency regime $\hat{n}\ll \nu^{-\frac34}$} \label{S5.1} In this regime, the structure of compressible Orr-Sommerfeld operator is similar to the classical one studied by G\'erard-Varet-Maekawa in \cite{GM19}. Thus the construction of homogeneous solution to the Orr-Sommerfeld equation \eqref{4.162-1} follows the same line as  \cite{GM19}. That is,  firstly we construct approximate slow and fast modes, then we use  bounds obtained in Section \ref{S4} to control the remainder. For completeness, we briefly introduce the construction as follows.
\subsubsection{Slow modes.}
We start from a homogenous solution $\varphi_{\rm Ray}$ to the compressible Rayleigh equation
 \begin{align}\label{4.154}
\begin{cases}
U_s\Lambda(\varphi_{\rm Ray})-\varphi_{\rm Ray}\partial_Y(A^{-1}\partial_YU_s)=0,\\
\varphi_{\rm Ray}(Y)|_{Y=0}=1.
\end{cases}
\end{align}
\begin{lemma}
\label{P4.10}
Let $m\in (0,1)$. There exists a solution $\varphi_{\rm Ray}\in H^1(\mathbb{R}_+)$ to \eqref{4.154} such that following statements hold.
\begin{itemize}
\item[{\rm(i)}] If $\alpha\in (0,1)$, there exists a positive constant $c_E>0$, such that $\varphi_{\rm Ray}$ has the following form:
 \begin{align}\label{4.155}
\varphi_{\rm Ray}=\varphi_{\rm Ray, 0}+\varphi_{\rm Ray, 1}+\varphi_{\rm Ray, 2},
\end{align}
where
\begin{align*}
&\varphi_{Ray, 0}=\frac{c_E}{\alpha}U_se^{-\alpha A_{\infty}^{\frac12}Y},\quad \varphi_{Ray, 1}|_{Y=0}=1,\\
&\|\partial_Y\varphi_{Ray, 1}\|_{L^2}+\|\varphi_{Ray, 1}\|_{L^2}\le C,\\
&\|\partial_Y\varphi_{Ray, 2}\|_{L^2}+\alpha\|\varphi_{Ray, 2}\|_{L^2}\le C\alpha^{\frac12},
\end{align*}
and $A_{\infty}\eqdef\lim\limits_{Y\to\infty}A(Y)=1-m^2$.
\item[{\rm(ii)}] If $\alpha\geq 1$, then $\varphi_{\rm Ray}$ has the following form:
 \begin{align}\label{4.161}
	\varphi_{\rm Ray}=\varphi_{\rm Ray,0}+{\varphi}_{\rm Ray,1},
\end{align}
where $\varphi_{\rm Ray,0}=e^{-\alpha Y},$ and
\begin{align}\label{4.162}
	\|\partial_Y{\varphi}_{\rm Ray,1}\|_{L^2}+\alpha\|\varphi_{\rm Ray,1}\|_{L^2}\le C\alpha^{-\frac12}.
\end{align}
\end{itemize}
 Moreover, if in addition $\frac{U_s^{\prime\prime}}{U_s}\in L^2(\mathbb{R}_+)$, then $\varphi_{\rm  Ray}$ belong to $H^2(\mathbb{R}_+)$. 
\end{lemma}
\begin{proof}
	The compressible Rayleigh equation \eqref{4.154} has a similar structure as the classical one for incompressible fluid, which has been studied in Proposition 5.4 (for $\alpha\in (0,1)$), and Proposition 5.6 (for $\alpha \geq 1$) in \cite{GM19}. 
	%We can construct $\varphi_{\rm Ray}$ in the same manner. 
	We omit the details  for brevity.
\end{proof}
By a perturbative argument,
the slow mode $\phi_{s}$ of the following compressible Orr-Sommerfeld equation:
\begin{equation}
	\left\{
\begin{aligned}
	&{\rm {OS}_{CNS}}(\phi_{s})=0,~Y>0,\label{4.162-1}\\
	&\phi_{s}|_{Y=0}=1,
\end{aligned}\right.
\end{equation}
can be constructed near $\varphi_{\rm Ray}$.
\begin{lemma}
\label{P4.12}
Let $m\in (0,1)$ and $\alpha\eps^{\frac13}\lesssim 1$. There exists a solution $\phi_{s}\in H^2(\mathbb{R}_+)$ to \eqref{4.162-1}, such that
\begin{itemize}
	\item[{\rm(i)}] if $\alpha\in (0,1)$, 
\begin{align}
&\|\partial_Y\phi_{s}\|_{L^2}+\alpha\|\phi_{s}\|_{L^2}+\|U_s'\phi_{s}\|_{L^2}\leq \frac{C}{\alpha},\label{4.163-1}\\
&\|\Delta_\alpha\phi_{s}\|_{L^2}\le C\left(\frac{1}{\alpha\eps^{\frac16}}+\frac{1}{\epsilon^{\frac13}}\right),\label{4.163-2}\\
&\p_Y\phi_{s}(0)=\frac{c_E}{\alpha}+O(1)\left(\frac{\eps^{\frac{1}{12}}}{\alpha}+\frac{1}{\eps^{\frac14}}\right)\label{4.163-3};
\end{align}
\item[{\rm(ii)}] if $\alpha\geq 1$, 
\begin{align}
	&\|\partial_Y\phi_{s}\|_{L^2}+\alpha\|\phi_{s}\|_{L^2}\leq C\alpha^{\frac12},\label{4.164-1}\\
	&\|\Delta_\alpha\phi_{s}\|_{L^2}\le C\eps^{-\frac13},\label{4.164-2}\\
	&\p_Y\phi_{s}(0)=-\alpha+O(1)\eps^{-\frac16}\label{4.164-3}.
\end{align}
\end{itemize}
\end{lemma}
\begin{proof}
	We firstly construct the solution $\tilde{\phi}_s$ to the symmetrized Orr-Sommerfeld equation 
	\begin{equation}
		\left\{
		\begin{aligned}
			&{\rm \widetilde{OS}_{CNS}}(\tilde{\phi}_s)=0,~Y>0,\\
			&\tilde{\phi}_s|_{Y=0}=1,\nonumber
		\end{aligned}
	\right.
	\end{equation}
then recover $\phi_s$ by Corrollary \ref{C4.9}. 

The solution takes the form $\tilde{\phi}_s=\varphi_{\rm Ray}+\psi+\xi$, where $\psi$ and $\xi$ satisfy the following two equations respectively:
\begin{equation}\label{4.165}
	\left\{
\begin{aligned}
	&{\rm \widetilde{Airy}}(\psi)=-i\eps \Lambda(\varphi_{\rm Ray})=-i\eps \frac{\p_Y(A^{-1}\p_YU_s)\varphi_{\rm Ray}}{U_s}:=-i\eps f_\psi,\\
	&\psi|_{Y=0}=0,
\end{aligned}
\right.
\end{equation}
 and 
\begin{equation}\label{4.166}
	\left\{
\begin{aligned}
	&{\rm \widetilde{OS}_{CNS}}(\xi)=\p_Y\left[(1+A^{-1})\p_YU_s\psi+(1-A^{-1})U_s\p_Y\psi)\right]:=f_{\xi},\\
&	\xi|_{Y=0}=0.
\end{aligned}\right.
	\end{equation}
To bound $\psi$ and $\xi$, we divide the discussion into two cases.

\underline{\it Case 1. $\alpha\in (0,1)$}. Recall the decomposition  \eqref{4.155} of $\varphi_{\rm Ray}$ . We further decompose $\psi=\psi_0+\psi_1$ where $\psi_0$ and $\psi_1$ satisfy $${\rm \widetilde{Airy}}(\psi_0)=-i\eps \frac{\p_Y(A^{-1}\p_YU_s)(\varphi_{\rm Ray,0}+\varphi_{\rm Ray,2})}{U_s}:=-i\eps f_{\psi,0},$$ and $${\rm \widetilde{Airy}}(\psi_1)=-i\eps \frac{\p_Y(A^{-1}\p_YU_s)\varphi_{\rm Ray,1}}{U_s}:=-i\eps f_{\psi,1}$$ respectively. We derive bounds of $\psi_0$ first. It is straightforward to check that for $k=0,1,2$,
\begin{equation}
\begin{aligned}
\|Y^kf_{\psi,0}\|_{L^2}&\leq C\|Y^k\p_Y(A^{-1}\p_YU_s)\|_{L^\infty}\|\frac{\varphi_{\rm Ray,0}}{U_s}\|_{L^\infty}+C\|\frac{Y^{k+1}\p_Y(A^{-1}\p_YU_s)}{U_s}\|_{L^\infty}\|\p_Y\varphi_{\rm Ray, 2}\|_{L^2}\nonumber\\
&\leq \frac{C}{\alpha},
\end{aligned}
\end{equation}
and for $k=1,2$, that
\begin{equation}
	\begin{aligned}
		\|Y^k\p_Yf_{\psi,0}\|_{L^2}\leq& C\left(\sum_{j=1}^3\|\frac{Y^{k+1}\p_Y^jU_s}{U_s^2}\|_{L^2}\right)\cdot\|\frac{\varphi_{\rm Ray,0}}{U_s}\|_{L^\infty}\nonumber\\
		&+C\left(\sum_{j=1}^3\|\frac{Y^{k+1}\p_Y^jU_s}{U_s^2}\|_{L^\infty}\right)\cdot\|\left(\p_Y\varphi_{\rm Ray,0}\|_{L^2}+\|\p_Y\varphi_{\rm Ray, 2}\|_{L^2}\right)\nonumber\\
		\leq& \frac{C}{\alpha}.
		\end{aligned}
\end{equation}
Then applying \eqref{4.52}, \eqref{Ai3} and \eqref{Ai4} to $\psi_1$ with $h=-i\eps f_{\psi,0}$, we obtain
\begin{equation}\label{4.167}
\begin{aligned}
	\|(\p_Y\psi_0,\alpha\psi_0)\|_{L^2}+\|Y\p_Y^2\psi_0\|_{L^2}+\|Y^2\p_Y^2\psi_0\|_{L^2}&\leq \frac{C\eps^{\frac13}}{\alpha},\\
	\|Y^2\p_Y^3\psi_0\|_{L^2}&\leq \frac{C}{\alpha}.
\end{aligned}
\end{equation}
For $\psi_1$, note that $f_{\psi,1}=Y^{-1}g_1$, where
$g_1=\frac{Y\p_Y(A^{-1}\p_YU_s)\varphi_{\rm Ray,1}}{U_s} \in L^2(\mathbb{R}_+)\cap H^1_{\rm loc}(\mathbb{R}_+)
$
and satisfies
\begin{align}
	\|g_1\|_{L^2}&\leq C\|\frac{Y\p_Y(A^{-1}\p_YU_s)}{U_s}\|_{L^2}\|\varphi_{\rm Ray,1}\|_{L^\infty}\leq C,\nonumber\\
	\|Y\p_Yg_1\|_{L^2}&\leq C\left(\sum_{j=1}^3\|\frac{Y^{2}\p_Y^jU_s}{U_s^2}\|_{L^\infty}\right)\cdot\left(\|\p_Y\varphi_{\rm Ray,1}\|_{L^2}+\|\varphi_{\rm Ray,1}\|_{L^2}\right)\leq C.\nonumber
\end{align}
Then applying \eqref{4.53}, \eqref{Ai5} and \eqref{Ai6} to $\psi_1$, we obtain
\begin{equation}
\begin{aligned}\label{4.167-1}
		\|(\p_Y\psi_1,\alpha\psi_1)\|_{L^2}+\|Y\p_Y^2\psi_1\|_{L^2}+\|Y^2\p_Y^2\psi_1\|_{L^2}+\|Y^2\p_Y^3\psi_1\|_{L^2}\leq C.
\end{aligned}
\end{equation}
Combining \eqref{4.167} and \eqref{4.167-1} yields
\begin{align}
\|(\p_Y\psi,\alpha\psi)\|_{L^2}+\|Y\p_Y^2\psi\|_{L^2}+\|Y^2\p_Y^2\psi\|_{L^2}&\leq C\left(\frac{\eps^{\frac13}}{\alpha}+1\right),\label{4.167-2}\\
\|Y^2\p_Y^3\psi\|_{L^2}&\leq \frac{C}{\alpha}.\label{4.167-3} 
\end{align}

Now we estimate $\xi$ which solves \eqref{4.166}. By using \eqref{4.167-2}, \eqref{4.167-3} and $|1-A^{-1}|\leq CU_s^2$, we  deduce 
\begin{equation}
	\begin{aligned}
		\|\frac{Yf_\xi}{U_s}\|_{L^2}+\|(1+Y)^2f_{\xi}\|_{L^2}\leq C\|\p_Y\psi\|_{L^2}+C\|Y^2\p_Y^2\psi\|_{L^2}\leq C\left(\frac{\eps^{\frac13}}{\alpha}+1\right),
	\end{aligned}\nonumber
\end{equation}
and
\begin{align}
	\|\frac{Y^2\p_Yf_\xi}{U_s}\|_{L^2}\leq C\|\p_Y\psi\|_{L^2}+C\|Y\p_Y^2\psi\|_{L^2}+C\|Y^2\p_Y^3\psi\|_{L^2}\leq \frac{C}{\alpha}.\nonumber
\end{align}
Plugging these estimates into \eqref{4.108} and \eqref{4.109}, we obtain
\begin{align}\label{4.167-4}
	\|\p_Y\xi\|_{L^2}+\alpha\|\xi\|_{L^2}&\leq C	\|(1+Y)^2f_{\xi}\|_{L^2}+C\eps^{\frac13}\left(\|\frac{Y}{U_s}f_\xi\|_{L^2}+\|\frac{Y^2\p_Yf_\xi}{U_s}\|_{L^2}\right)\nonumber\\
	&\leq C\left(\frac{\eps^{\frac13}}{\alpha}+1\right).
\end{align}
Combining bounds of $\varphi_{\rm Ray}$ in Lemma \ref{P4.10}, \eqref{4.167-2}, and \eqref{4.167-4} together, we deduce the following estimate
\begin{align}
		\|(\p_Y\tilde{\phi}_s,\alpha\tilde{\phi}_s)\|_{L^2}&\leq \|(\p_Y\varphi_{\rm Ray},\alpha\varphi_{\rm Ray})\|_{L^2}+\|(\p_Y\psi,\alpha\psi)\|_{L^2}+\|(\p_Y\xi,\alpha\xi)\|_{L^2}\nonumber\\
	&\leq \frac{C}{\alpha}+C\left(\frac{\eps^{\frac13}}{\alpha}+1\right)\leq \frac{C}{\alpha}.\label{4.167-6}
\end{align}

For the estimate on $\|\Delta_\alpha\tilde{\phi}_s\|_{L^2}$, note that the following boundary condition holds 
\begin{align}
	\p_Y\Lambda(\tilde{\phi}_s)|_{Y=0}=\p_Y\Lambda(\varphi_{\rm Ray})|_{Y=0}+\p_Y\Lambda(\psi)|_{Y=0}+\p_Y\Lambda(\xi)|_{Y=0}=0,\nonumber
\end{align}
because $i\eps\p_Y\Lambda(\varphi_{\rm Ray}+\psi)|_{Y=0}=-\p_Y(U_s\psi)|_{Y=0}=0$ and $\p_Y\Lambda(\xi)|_{Y=0}=0$, cf.  \eqref{4.110-1}. Then using the same argument as in \cite[Proposition 7.9]{GM19}, we can derive the following bounds on $\tilde{\phi}_{s,r}=\tilde{\phi}_s-\varphi_{\rm Ray, 0}:$
\begin{align}\label{4.167-5}
	\|\Delta_\alpha\tilde{\phi}_{s,r}\|_{L^2}\leq C\left(\frac{1}{\alpha\eps^{\frac16}}+\frac{1}{\eps^{\frac13}}\right).
\end{align}
Combinining \eqref{4.167-5} with the explicit formula of $\varphi_{\rm Ray}$ in Lemma \ref{P4.10} yields
\begin{align}
	\|\Delta_\alpha\tilde{\phi}_s\|_{L^2}&\leq \|\Delta_{\alpha}\varphi_{\rm Ray,0}\|_{L^2}+\|\Delta_{\alpha}\tilde{\phi}_{s,r}\|_{L^2}\nonumber\\
	&\leq \frac{C}{\alpha}+C\left(\frac{1}{\alpha\eps^{\frac16}}+\frac{1}{\eps^{\frac13}}\right)\le C\left(\frac{1}{\alpha\eps^{\frac16}}+\frac{1}{\eps^{\frac13}}\right)\label{4.167-7}.
\end{align}

Finally, we look for the original Orr-Sommerfeld solution in the form of $\phi_{s}=\tilde{\phi}_s+\phi_r$. By Corollary \ref{C4.9}, \eqref{4.167-6}, and \eqref{4.167-5}, we deduce that
\begin{align}
	\|\p_Y\phi_r\|_{L^2}+\alpha\|\phi_r\|_{L^2}&\leq C\eps^{\frac23}\|(\partial_{Y}\tilde{\phi}_s,\alpha\tilde{\phi}_s)\|_{L^2}+C\eps^{\frac23}\|\Delta_\alpha\tilde{\phi}_s\|_{L^2}\nonumber\\
	&\leq \frac{C\eps^{\frac23}}{\alpha}+C\eps^{\frac13}\left(\frac{\eps^{\frac16}}{\alpha}+1
	\right)\leq C\eps^{\frac13}\left(\frac{\eps^{\frac16}}{\alpha}+1
	\right),\label{4.167-8}
	\end{align}
and
\begin{align}
	\|\Delta_\alpha\phi_r\|_{L^2}
	&\leq C\eps^{\frac13}\|(\partial_{Y}\tilde{\phi}_s,\alpha\tilde{\phi}_s)\|_{L^2}+C\eps^{\frac13}\|\Delta_\alpha\tilde{\phi}_s\|_{L^2}\nonumber\\
	&\leq \frac{C\eps^{\frac13}}{\alpha}+C\left(\frac{\eps^{\frac16}}{\alpha}+1\right)\leq C\left(\frac{\eps^{\frac16}}{\alpha}+1\right).\label{4.167-9}
\end{align}
Combining these two bounds with \eqref{4.167-7} yields \eqref{4.163-2}. In addition,
\begin{equation}%
	\begin{aligned}
		\|U_s'\phi_{s}\|_{L^2}&\leq \|U_s'\varphi_{\rm Ray,1}\|_{L^2}+\|U_s'(\phi_{s}-\varphi_{\rm Ray,1})\|_{L^2}\\
		&\leq C\|U_s'\|_{L^2}\|\varphi_{\rm Ray, 1}\|_{L^\infty}+C\|YU_s'\|_{L^\infty}\|\p_Y(\phi_{s}-\varphi_{\rm Ray,1})\|_{L^2}\leq \frac{C}{\alpha}.\nonumber
	\end{aligned}
\end{equation}
Combining it with \eqref{4.167-6}, \eqref{4.167-7}, and \eqref{4.167-8}, we obtain \eqref{4.163-1}. Moreover, recall the decomposition $$\phi_{s}=\varphi_{\rm Ray,0}+\tilde{\phi}_{s,r}+\phi_r,$$
where $$\|\partial_Y\tilde{\phi}_{s,r}\|_{L^2}\leq \|\partial_Y\varphi_{\rm Ray,1}\|_{L^2}+\|\partial_Y\varphi_{\rm Ray ,2}\|_{L^2}+\|\partial_Y\psi\|_{L^2}+\|\partial_Y\xi\|_{L^2}\leq C\left(\frac{\eps^{\frac13}}{\alpha}+1\right).$$
Then by \eqref{4.167-5} and the interpolation, we get
$\|\p_Y\tilde{\phi}_{s,r}\|_{L^\infty}\leq C\left(\frac{\eps^{\frac1{12}}}{\alpha}+\frac{1}{\eps^{\frac14}}\right) $. Similarly, from \eqref{4.167-8} and \eqref{4.167-9}, we deduce
$\|\partial_Y\phi_r\|_{L^\infty}\leq C\eps^{\frac16}\left(\frac{\eps^{\frac16}}{\alpha}+1\right)$. Combining these two bounds with 
$\p_Y\varphi_{\rm Ray,0}(0)=\frac{c_E}{\alpha}$ yields \eqref{4.163-3}. The proof of Proposition \ref{P4.12} for $\alpha\in (0,1)$ is complete.

\underline{\it Case 2. $\alpha\geq 1$.} Note that $\psi$ solves \eqref{4.165}. We can write $f_{\psi}=Y^{-1}g_{\psi}$, where
  $g_{\psi}=\frac{Y\partial_Y(A^{-1}\partial_YU_s)\varphi_{\rm Ray}}{U_s}$. Recall \eqref{4.161} and \eqref{4.162} for the estimates on $\varphi_{\rm Ray}$. It is straightforward to compute
\begin{align}
	\|g_{\psi}\|_{L^2}&\leq \|\frac{Y\p_Y(A^{-1}\p_YU_s)}{U_s}\|_{L^\infty}\|\varphi_{\rm Ray}\|_{L^2}\leq C\alpha^{-\frac12},\nonumber\\
		\|Y\p_Yg_{\psi}\|_{L^2}&\leq C\left(\sum_{j=1}^3\|\frac{Y^2\p_Y^jU_s}{U_s^2}\|_{L^\infty}\right)\cdot\left(\|\varphi_{\rm Ray,0}\|_{L^2}+\|\varphi_{\rm Ray, 1}\|_{H^1}\right)\nonumber\\
		&\qquad+C\left(\sum_{j=1}^2\|\frac{Y\p_Y^jU_s}{U_s}\|_{L^\infty}\right)\cdot\|Y\p_Y\varphi_{\rm Ray,0}\|_{L^2}\leq C\alpha^{-\frac12}.\nonumber
\end{align}
 Plugging these bounds into \eqref{4.53}, \eqref{Ai5} and \eqref{Ai6}, we obtain
\begin{align}\label{4.164-4}
	\|(\p_Y\psi,\alpha\psi)\|_{L^2}+\|Y\p_Y^2\psi\|_{L^2}+\|Y^2\p_Y^2\psi\|_{L^2}+\|Y^2\p_Y^3\psi\|_{L^2} &\leq C\|g_{\psi}\|_{L^2}+C\|Y\p_Yg_{\psi}\|_{L^2}\nonumber\\
	&\leq C\alpha^{-\frac12}.
\end{align}

Next step is to estimate $\xi$ which solves \eqref{4.166}. It is straightforward to check
\begin{align}
	\|\frac{Yf_\xi}{U_s}\|_{L^2}+\|\frac{Y^2\p_Yf_{\xi}}{U_s}\|_{L^2}&\leq 	C\|\p_Y\psi\|_{L^2}+C\|Y\p_Y^2\psi\|_{L^2}+C\|Y^2\p_Y^2\psi\|_{L^2}+C\|Y^2\p_Y^3\psi\|_{L^2}\nonumber\\
	&\leq C\alpha^{-\frac12}.\nonumber
\end{align}
Applying \eqref{4.106} and \eqref{4.107} to $\xi$, we can obtain
\begin{align}\label{4.164-5}
	\|\p_Y\xi\|_{L^2}+\alpha\|\xi\|_{L^2}\leq C\alpha^{-\frac12}.
\end{align}
Combining \eqref{4.162}, \eqref{4.164-4}, \eqref{4.164-5} with equalities $
	\|\varphi_{\rm Ray,0}\|_{L^2}=\alpha^{-\frac12}$ and $ \|\p_Y\varphi_{\rm Ray,0}\|_{L^2}=\alpha^{\frac12},$ we get
	\begin{align}\label{4.164-6}
		\|\p_Y\tilde{\phi}_s\|_{L^2}+\alpha\|\tilde{\phi}_s\|_{L^2}\leq C\alpha^{\frac12}.
	\end{align}
Given $H^1$-bound \eqref{4.164-6}   and the additional boundary condition $\partial_Y\Lambda(\tilde{\phi}_s)$ as the Case 1, we can obtain 
\begin{align}\label{4.164-7}
	\|\Delta_\alpha\tilde{\phi}_s\|_{L^2}=\|\Delta_{\alpha}\tilde{\phi}_{s,r}\|_{L^2}\leq C\eps^{-\frac13}
\end{align}
by the same argument as in \cite[Proposition 7.10]{GM19}.  This completes $H^2$-estimate of $\tilde{\phi}_s$.

 Finally, we follow the same procedure as in Case 1 to construct the original Orr-Sommerfeld solution in the form of $\phi_s=\tilde{\phi}_s+\phi_r$. From Corollary \ref{C4.9}, bounds \eqref{4.164-6} and \eqref{4.164-7} of $\tilde{\phi}_s$, we get
\begin{align}
	\|\p_Y\phi_r\|_{L^2}+\alpha\|\phi_r\|_{L^2}&\leq C\eps^{\frac23}(\|\partial_{Y}\tilde{\phi}_s\|_{L^2}+\|\Delta_\alpha\tilde{\phi}_s\|_{L^2}+\alpha^2\|\tilde{\phi}_s\|_{L^2})\nonumber\\
	&\leq C\eps^{\frac23}\left(\alpha^{\frac32}+\eps^{-\frac13}\right)\leq C\alpha^{-\frac12},\label{4.164-8}
\end{align}
and
\begin{align}
	\|\Delta_\alpha\phi_r\|_{L^2}&\leq C\eps^{\frac13}(\|\partial_{Y}\tilde{\phi}_s\|_{L^2}+\|\Delta_\alpha\tilde{\phi}_s\|_{L^2}+\alpha^2\|\tilde{\phi}_s\|_{L^2})\nonumber\\
	 &\leq C\eps^{\frac13}\left(\alpha^{\frac32}+\eps^{-\frac13}\right)\leq C\eps^{-\frac13},\label{4.164-9}
\end{align}
where $\alpha\eps^{\frac13}\lesssim 1$ has been used in the last inequalities of \eqref{4.164-8} and \eqref{4.164-9}. Combining \eqref{4.164-6}-\eqref{4.164-9} together,
we obtain \eqref{4.164-1} and \eqref{4.164-2}. Estimate \eqref{4.164-3} follows from interpolation. The proof of Lemma \ref{P4.12} is complete.
\end{proof}
Now we recover the slow mode $(\varrho_{s},\mathfrak{u}_{s},\mathfrak{v}_{s})$ of the homogeneous
quasi-compressible system.
\begin{proposition}
\label{P4.14}
Let $m\in(0, 1)$ and $\alpha\eps^{\frac13}\lesssim 1$.  There exists a solution $(\varrho_{s},\mathfrak{u}_{s},\mathfrak{v}_{s})\in H^2(\mathbb{R}_+)^3$ to $$\mathcal{L}_{Q}(\varrho_{s},\mathfrak{u}_{s},\mathfrak{v}_{s})={\bf 0}.$$ Moreover, the following statements hold.
\begin{itemize}
	\item [{\rm(i)}] If $0<\alpha\le 1$,
\begin{align}\label{4.173}
m^{-2}\|\varrho_{s}\|_{L^2}+\|(\mathfrak{u}_{ s},\mathfrak{v}_{s})\|_{L^2}&\leq \frac{C}{\alpha},\\
m^{-2}\|\p_Y\varrho_{s}\|_{L^2}+m^{-2}\|\Delta_\alpha\varrho_{s}\|_{L^2}+\|\p_Y\mathfrak{v}_{s}\|_{L^2}&\leq C,\label{4.174}\\
\|\partial_Y\mathfrak{u}_{ s}\|_{L^2}+\|\Delta_\alpha\mathfrak{v}_{s}\|_{L^2}&\le C\left(\frac{1}{\alpha\eps^{\frac16}}+\frac{1}{\eps^{\frac13}}\right),\label{4.176}\\
\|\Delta_\alpha \mathfrak{u}_{s}\|_{L^2}&\le \frac{C}{\alpha\epsilon},\label{4.177}\\
\|\partial_Y\div_{\alpha}(\mathfrak{u}_s,\mathfrak{v}_s),\alpha\div_{\alpha}(\mathfrak{u}_s,\mathfrak{v}_s)\|_{L^2}&\leq C.\label{4.177-1}
\end{align}
On the boundary, it holds that
\begin{align}
	\mathfrak{v}_{s}(0)=-i\alpha, ~~\mathfrak{u}_{s}(0)=\frac{c_E}{\alpha}+O(1)\left( \frac{\eps^{\frac{1}{12}}}{\alpha}+\frac{1}{\eps^{\frac14}}\right).\label{4.178}
\end{align}
\item[{\rm(ii)}] If $\alpha\ge 1$,
\begin{align}\label{4.179}
\|(\mathfrak{u}_{s}, \mathfrak{v}_{s})\|_{L^2}&\le C\alpha^{\frac12},\\
m^{-2}\|(\varrho_{s},\partial_Y\varrho_{s},\Delta_{\alpha}\varrho_{s})\|_{L^2}+\|\partial_Y\mathfrak{v}_{s}\|_{L^2}&\le C\alpha^{\frac32},\label{4.180}\\
\|\p_Y \mathfrak{u}_{s}\|_{L^2}+\alpha^{-1}\|\Delta_\alpha \mathfrak{v}_{ s}\|_{L^2}&\le C\left(\alpha^{\frac32}+\frac{1}{\epsilon^{\frac13}}\right),\label{4.181}\\
\|\Delta_\alpha\mathfrak{u}_{s}\|_{L^2}&\le \frac{C\alpha^{\frac32}}{\eps},\label{4.182}\\
\|\partial_Y\div_{\alpha}(\mathfrak{u}_s,\mathfrak{v}_s),\alpha\div_{\alpha}(\mathfrak{u}_s,\mathfrak{v}_s)\|_{L^2}&\leq C\alpha^{\frac32}.\label{4.82-1}
\end{align}
On the boundary, it holds that
\begin{align}
	\mathfrak{v}_{s}(0)=-i\alpha, ~~\mathfrak{u}_{s}(0)=-\alpha+O(1)\eps^{-\frac16}.\label{4.183}
\end{align}
\end{itemize}
\end{proposition}
\begin{proof}
	Let $\phi_{s}$ be the slow mode constructed in Lemma \ref{P4.12}. Then it is straightforward to check $(\varrho_s,\mathfrak{u}_s,\mathfrak{v}_s)$  where
	\begin{align}
		\mathfrak{v}_{s}=-i\alpha\phi_{s},~\varrho_{s}=-m^2A^{-1}\left(i\eps\Delta_\alpha\p_Y\phi_{s}+U_s\p_Y\phi_{s}-\phi_{s}\p_YU_s\right),~\mathfrak{u}_{s}=\p_Y\phi_{s}-U_s\varrho_{s} \label{4.182-1}
	\end{align}
	is a homogeneous solution to quasi-compressible system $\mathcal{L}_{Q}(\varrho_{s},\mathfrak{u}_{s},\mathfrak{v}_{s})={\bf 0}$. The bounds \eqref{4.173}-\eqref{4.177} and \eqref{4.179}-\eqref{4.182} follow from \eqref{4.163-1}, \eqref{4.163-2}, \eqref{4.164-1}, \eqref{4.164-2} and Propositiom \ref{P8.1} with $f_{\mathfrak{u}}=f_{\mathfrak{v}}=0$.  Boundary conditions
	 \eqref{4.178} and \eqref{4.183} have been shown in \eqref{4.163-3} and \eqref{4.164-4} respectively, by noting that $\mathfrak{v}_s(0)=-i\alpha\phi_s(0)$ and $\mathfrak{u}_s(0)=\partial_Y\phi_s(0)$. The proof of Proposition \ref{P4.14} is complete. 
\end{proof}

\subsubsection{Fast mode}\label{S4.5.2}
To eliminate boundary values of $(\mathfrak{u}_s,\mathfrak{v}_s)$, we construct a boundary sublayer solution $(\varrho_{f,l}, \mathfrak{u}_{f,l}, \mathfrak{v}_{ f,l})$ to 
\begin{align}\label{4.198}
\begin{cases}
\mathcal{L}_Q(\varrho_{ f,l}, \mathfrak{u}_{f,l}, \mathfrak{v}_{ f,l})={\bf{0}},\\
v_{f,l}|_{Y=0}=-i\alpha.
\end{cases}
\end{align}
Unfortunately, if inserting the boundary layer ansatz into \eqref{4.198}, the leading order profile is not accurate for later use. For this, we firstly construct a boundary sublayer solution to the compressible Orr-Sommerfeld equation
\begin{equation}\label{4.198-1}
\left\{\begin{aligned}
&{\rm OS_{CNS}}(\phi_{f,l})=0,~Y>0,\\
&\phi_{f,l}|_{Y=0}=1,
\end{aligned}
\right.
\end{equation}
then recover $(\varrho_{f,l}, \mathfrak{u}_{f,l}, \mathfrak{v}_{f,l})$ by Proposition \ref{P8.1}.

Observe that
\begin{align}
{\rm OS_{CNS}}(\phi)=&A^{-1}\left( i\eps\Delta_{\alpha}^2\phi+U_s\Delta_\alpha\phi    \right)+\alpha^2(A^{-1}-1)\left(i\eps\Delta_\alpha\phi+U_s\phi\right)\nonumber\\
&+\p_Y(A^{-1})\left(i\eps\p_Y\Delta_{\alpha}\phi+U_s\p_Y\phi\right)-\p_Y(A^{-1}\p_YU_s)\phi.\label{4.200}
\end{align}
From \eqref{4.200}, one can see that in the low-frequency regime, ${\rm OS_{CNS}}$ is asymptotically equivalent to the classical Orr-Sommerfeld equation for incompressible model inside the boundary sublayer. Thus, it is natural to construct the fast mode that is closed to the incompressible one studied by Ger\'ard-Varet and Maekawa \cite{GM19}. For completeness, we briefly introduce the construction as follows.

Set $\delta=e^{-\frac{1}{6}\pi i}\eps^{\frac13}$ as the scale of sublayer. Let ${\rm Ai}(z)$ be the classical Airy function satisfying $\frac{d^2}{dz^2}{\rm Ai}-z{\rm Ai}=0$. Set $\psi_{0,l}(Y)={\rm Ai}(\delta^{-1}Y)$.
Then the approximate fast mode $\phi_{f,app,l}$ is constructed as follows.
%Let $C_0>0$ be a positive constant which will be determined later. Set
%\begin{align}
%	\phi_{f,app}(Y)=C_0\eps^{-\frac23}\int_Y^\infty e^{\alpha(Y-Y')}\int_{Y'}^\infty e^{\alpha(Y''-Y')}\psi_0(Y'')\dd Y'' dY'.\nonumber
%\end{align}
%Then we have
\begin{lemma}[Proposition 7.11, \cite{GM19}]\label{lm5.4}
	There exists $\kappa_0\in (0,1)$, such that if $\alpha\eps^{\frac13}\leq \kappa_0$, there exists a solution $\phi_{f,app,l}$ to the equation
	\begin{align}
		(\p_Y^2-\alpha^2)\phi_{f,app,l}=\psi_{0,l},\label{4.200-1}
	\end{align}
which satisfies the following bounds
\begin{align}
	\|Y^j\p_Y^k\phi_{f,app,l}\|_{L^2}\leq C\eps^{\frac16+\frac{j-k}{3}}, \label{4.200-2}
\end{align}
for $j,k=0,1,2,3.$ Moreover, the boundary conditions of $\phi_{f,app,l}$ satisfy
\begin{equation}
\begin{aligned}
	\phi_{f,app,l}(0)&=1,\\
	\p_Y\phi_{f,app,l}(0)&=\left(
e^{\frac{1}{6}\pi i}3^{-\frac23}\Gamma(\frac13)+O(\eps^{\frac13}\alpha)+O(\eps^{\frac13})	\right)\eps^{-\frac13}.\label{4.200-3}
\end{aligned}
\end{equation}
\end{lemma}

Now we construct the fast mode $\phi_{f,l}$ to the compressible Orr-Sommerfeld equation near $\phi_{f,app,l}$. Firstly, from \eqref{4.200}, we compute the error $E_{f,app,l}$ created by $\phi_{f,app,l}$:
\begin{align}
	E_{f,app,l}={\rm OS_{CNS}}(\phi_{f,app,l})=&A^{-1}\left[U_s(Y)-U_s'(0)Y\right]\psi_{0,l}+\alpha^2\left[ -i\eps \psi_{0,l}+(A^{-1}-1)U_s\phi_{f,app,l}  \right]\nonumber\\
	&+\p_Y(A^{-1})\left(i\eps\p_Y\psi_{0,l}+U_s\p_Y\phi_{f,app,l}\right)-\p_Y(A^{-1}\p_YU_s)\phi_{f,app,l}.\nonumber
\end{align}
Then we look for the solution $\phi_{f,l}=\phi_{f,app,l}+\phi_{f,r,l}$, where $\phi_{f,r,l}$ satisfies ${\rm OS_{CNS}}(\phi_{f,r,l})=-E_{f,app,l}$ and $\phi_{f,r,l}|_{Y=0}=0.$ It is straightforward to compute
\begin{align}
\|YE_{f,app,l}\|_{L^2}\leq& C\|Y^3\Delta_\alpha\phi_{f,app,l}\|_{L^2}+C\alpha^2\eps\|Y\Delta_\alpha\phi_{f,app,l}\|_{L^2}+C\|Y^3\phi_{f,app,l}\|_{L^2}\nonumber\\
&+C\eps\|Y\p_Y\Delta_\alpha\phi_{f,app,l}\|_{L^2}+C\|Y^2\p_Y\phi_{f,app,l}\|_{L^2}+C\|Y\phi_{f,app,l}\|_{L^2}\nonumber\\
\leq &C\eps^{\frac12}+\alpha^2\eps^{\frac56}.\nonumber
\end{align}
Applying Proposition \ref{P4.6} with $g=YE_{f,app,l}$ and Corollary \ref{C4.9} to $\phi_{f,r,l}$ yields
\begin{align}
	\|\p_Y\phi_{f,r,l}\|_{L^2}+\alpha\|\phi_{f,r,l}\|_{L^2}&\leq C\eps^{\frac16}+C\alpha^2\eps^{\frac12},\nonumber\\
	\|\Delta_\alpha\phi_{f,r,l}\|_{L^2}&\leq \frac{C}{\eps^{\frac13}}(\eps^{\frac16}+\alpha^2\eps^{\frac12}).\nonumber
\end{align}
By interpolation, $\|\p_Y\phi_{f,r,l}\|_{L^\infty}\leq C(1+\alpha^2\eps^{\frac13})\leq C(1+\alpha)$, where $\alpha\eps^{\frac13}\lesssim 1$ has been used. Finally, it holds that
$$\|U_s'\phi_{f,l}\|_{L^2}\leq C\|\phi_{f,app,l}\|_{L^2}+C\|YU_s'\|_{L^\infty}\|\p_Y\phi_{f,r,l}\|_{L^2}\leq C\eps^{\frac16}+C\alpha^2\eps^{\frac12}\leq C\eps^{-\frac16}.
$$

We summarize the above construction as follows.
\begin{lemma}
\label{P4.16}
Let $m\in (0, 1)$. If  $\alpha \eps^{\frac13}\leq \kappa_0$, where $\kappa_0$ is given in Lemma \ref{lm5.4}, there exists a solution $\phi_{f,l}$ to \eqref{4.198-1}, such that
\begin{align}\label{4.213}
\|U_s'\phi_{f,l}\|_{L^2}+\|(\partial_Y\phi_{f,l}, \alpha\phi_{f,l})\|_{L^2}&\leq C\epsilon^{-\frac16},\\
\|(\partial_Y^2-\alpha^2)\phi_{f,l}\|_{L^2}&\leq C\epsilon^{-\frac12}.\label{4.214}
\end{align}
Moreover, on the boundary 
\begin{align}
\partial_Y\phi_{f,l}(0)=\left(e^{\frac{1}{6}\pi i}3^{-\frac23}\Gamma(\frac13)+O(\epsilon^{\frac13}\alpha)+O(\epsilon^{\frac13})\right)\epsilon^{-\frac13}.\nonumber
\end{align}
\end{lemma}
Substituting \eqref{4.213} and \eqref{4.214} into bounds \eqref{8.0-1}-\eqref{8.0-7} in Proposition \ref{P8.1}, we recover the fast mode $(\varrho_{f,l}, \mathfrak{u}_{f,l}, \mathfrak{v}_{f,l})$ of the quasi-compressible system.
\begin{proposition}
\label{P4.17}
Let $m\in (0, 1)$. If $0<\hat{n}\le \kappa_0\nu^{-\frac34}$,  the homogeneous quasi-compressible system \eqref{4.198} admits a solution $(\varrho_{f,l}, \mathfrak{u}_{f,l}, \mathfrak{v}_{f,l})\in  H^2(\mathbb{R}_+)^3$ satisfying
\begin{align}
\|(\mathfrak{u}_{f,l},\mathfrak{v}_{f,l})\|_{L^2}&\leq C\eps^{-\frac16},\label{4.215}\\
\|\p_Y\mathfrak{u}_{f,l}\|_{L^2}+\alpha^{-1}\|\Delta_\alpha \mathfrak{v}_{f,l}\|_{L^2}&\leq C\eps^{-\frac12}\label{4.216},\\
m^{-2}\|(\varrho_{f,l},\p_Y\varrho_{f,l},\Delta_\alpha\varrho_{f,l})\|_{L^2}+\|\partial_Y\mathfrak{v}_{f,l}\|_{L^2}&\leq C(\alpha+1)\eps^{-\frac16},\label{4.217}\\
\|\Delta_\alpha \mathfrak{u}_{f,l}\|_{L^2}&\leq C(\alpha+1)\eps^{-\frac76},\label{4.218-1}\\
\|(\partial_Y\div_{\alpha}(\mathfrak{u}_{f,l},\mathfrak{v}_{f,l}),i\alpha\div_{\alpha}(\mathfrak{u}_{f,l},\mathfrak{v}_{f,l}))\|_{L^2}&\leq C(\alpha+1)\eps^{-\frac16}.\label{4.218-3}
\end{align}
In particular, on the boundary it holds that
\begin{align}
\mathfrak{v}_{f,l}(0)=-i\alpha,~~ \mathfrak{u}_{f,l}(0)=\left(e^{\frac{1}{6}\pi i}3^{-\frac23}\Gamma(\frac13)+O(\alpha\epsilon^{\frac13})+O(\epsilon^{\frac13})\right)\epsilon^{-\frac13}.\label{4.218-2}
\end{align}
\end{proposition}
\subsubsection{Homogeneous quasi-compressible solutions}
Now we set 
\begin{align}\label{H1}
	(\varrho_{H,l},\mathfrak{u}_{H,l},\mathfrak{v}_{H,l})=(\varrho_{f,l},\mathfrak{u}_{f,l},\mathfrak{v}_{f,l})-(\varrho_{s},\mathfrak{u}_{s},\mathfrak{v}_{s})
\end{align}
as the homogeneous quasi-compressible solution. Here $(\varrho_{f,l},\mathfrak{u}_{f,l},\mathfrak{v}_{f,l})$ and $(\varrho_{s},\mathfrak{u}_{s},\mathfrak{v}_{s})$  are the fast and slow modes given by Propositions \ref{P4.17} and  \ref{P4.14} respectively.
\begin{proposition}\label{H2}
	Let $m\in (0,1)$. For any $0<\hat{n}\leq \kappa_0\nu^{-\frac34}$, $	(\varrho_{H,l},\mathfrak{u}_{H,l},\mathfrak{v}_{H,l})\in H^2(\mathbb{R}_+)^3$ is a solution to the homogeneous quasi-compressible system \eqref{4.151}. Moreover, the following two statements hold
	\begin{itemize}
		\item[(i)] If $\alpha\in (0,1)$, 
		\begin{align}
			m^{-2}\|(\varrho_{H,l},\partial_Y\varrho_{H,l},\Delta_\alpha\varrho_{H,l})\|_{L^2}+\|\mathfrak{u}_{H,l}\|_{L^2}+\|\mathfrak{v}_{H,l}\|_{H^1}&\leq C\left(\frac1\alpha+\frac{1}{\eps^{\frac16}}\right)\label{H3},\\
			\|\partial_Y\mathfrak{u}_{H,l}\|_{L^2}+\|\Delta_\alpha \mathfrak{v}_{H,l}\|_{L^2}&\leq C\left(\frac{1}{\alpha \eps^{\frac16}}+\frac{1}{\eps^{\frac12}}\right),\label{H4}\\ 
			\|\Delta_\alpha\mathfrak{u}_{H,l}\|_{L^2}&\leq \frac{C}{\eps}\left(\frac1\alpha+\frac{1}{\eps^{\frac16}}\right),\label{H5}\\
			\|(\partial_Y\div_\alpha(\mathfrak{u}_{H,l},\mathfrak{v}_{H,l}),\alpha\div_\alpha(\mathfrak{u}_{H,l},\mathfrak{v}_{H,l})\|_{L^2}&\leq C{\eps^{-\frac16}}\label{H6}.
		\end{align}
	On the boundary, it holds that
	\begin{align}\label{H7}
	\mathfrak{v}_{H,l}(0)=0,~~{\rm and }~~	\left|\mathfrak{u}_{H,l}(0)\right|\geq C\left(\frac{1}{\alpha}+\frac1{\eps^{\frac13}}\right);
	\end{align}
\item[(ii)] If $\alpha>1$, 
\begin{align}
	\|(\mathfrak{u}_{H,l},\mathfrak{v}_{H,l})\|_{L^2}&\leq C\eps^{-\frac16}\label{H8},\\
	m^{-2}\|(\varrho_{H,l},\partial_Y\varrho_{H,l},\Delta_\alpha\varrho_{H,l})\|_{L^2}+\|\partial_Y \mathfrak{v}_{H,l}\|_{L^2}&\leq C\alpha\eps^{-\frac16},\label{H9}\\
		\|\partial_Y\mathfrak{u}_{H,l}\|_{L^2}+\alpha^{-1}\|\Delta_\alpha\mathfrak{v}_{H,l}\|_{L^2}&\leq C\eps^{-\frac12},\label{H11}\\
		\|\Delta_\alpha\mathfrak{u}_{H,l}\|_{L^2}&\leq C\alpha \eps^{-\frac76},
		\label{H12}\\
	\|(\partial_Y\div_\alpha(\mathfrak{u}_{H,l},\mathfrak{v}_{H,l}),\alpha\div_\alpha(\mathfrak{u}_{H,l},\mathfrak{v}_{H,l})\|_{L^2}&\leq C\alpha\eps^{-\frac16}.\label{H13}
\end{align}
	On the boundary, it holds that
\begin{align}\label{H14}
\mathfrak{v}_{H,l}(0)=0,~~{\rm and }~~	\left|\mathfrak{u}_{H,l}(0)\right|\geq C\eps^{-\frac13}.
\end{align}
	\end{itemize}
\end{proposition} 
\begin{proof}
	 Using bounds \eqref{4.173}-\eqref{4.177-1} for $\alpha\in (0,1)$, \eqref{4.179}-\eqref{4.82-1} for $\alpha\geq 1$, and \eqref{4.215}-\eqref{4.218-2}, then noting that $\alpha\eps^{\frac13}\lesssim 1$, we can get estimates \eqref{H3}-\eqref{H6}, and \eqref{H8}-\eqref{H13}. Next we estimate boundary values of $(\mathfrak{u}_{H,l},\mathfrak{v}_{H,l})$. From the construction \eqref{H1}, it is directly to see that $\mathfrak{v}_{H,l}(0)=0$. For $\alpha\in (0,1)$, using \eqref{4.178}, \eqref{4.218-2}, and $\alpha\eps^{\frac13}\leq \kappa_0$, we can find a positive constant $c_0$, such that
	 \begin{align}
	 	\left|\mathfrak{u}_{H,l}(0)\right|\geq c_0\left( \frac{1}{\alpha}+\frac{1}{\eps^{\frac13}}\right)-C(\kappa_0+\eps^{\frac1{12}})\left( \frac{1}{\alpha}+\frac{1}{\eps^{\frac13}}\right).\nonumber
	 \end{align}
 Then \eqref{H7} follows by taking $\kappa_0$ suitably small. Similarly, for $\alpha\geq 1$, we use \eqref{4.183} and \eqref{4.218-2} to get
 \begin{align}
 	\left|\mathfrak{u}_{H,l}(0)\right|\geq c_0\eps^{-\frac13}-C(\alpha\eps^{\frac13}+\eps^{\frac16})\eps^{-\frac13}\geq (c_0-C\kappa_0-C\eps^{\frac16})\eps^{-\frac13}.\nonumber
 \end{align}
Thus, \eqref{H14} follows by taking $\kappa_0$ suitably small. The proof of Proposition \ref{H2} is complete.
\end{proof}

\subsection{Middle-frequency regime $\hat{n}\sim \nu^{-\frac{3}{4}}$}\label{S5.2}
%In this subsection, we construct a solution $(\rho_{q}, u_{q}, v_{q})$ to the Quasi-compressible equations in middle frequency regime $\hat{n}\sim\nu^{-\frac34}$.
%\begin{align}\label{4.219}
%\begin{cases}
%\mathcal{L}_Q(\rho_{q}, u_{q}, v_{q})={\bf{0}},\\
%v_q|_{Y=0}=0,\quad |u_q||_{Y=0}>0.
%\end{cases}
%\end{align}
%Similar as \eqref{4.6}, there exists a function $\phi$ such that
%\begin{align}\label{4.220}
%\partial_Y\phi=u_q+U_s\rho_q,\quad -i\alpha\phi=v_q,\quad \phi|_{Y=0}=0,
%\end{align}
%and satisfy the homogeneous compressible Orr-Sommerfeld equations
%\begin{align}\label{4.221}
%\begin{cases}
%OS_{CNS}(\phi)=\frac{i}{\hat{n}}\Lambda(\Delta_\alpha\phi)+U_s\Lambda(\phi)-\phi\partial_Y(A^{-1}\partial_YU_s)=0,\\
%\phi|_{Y=0}=0,\quad |\partial_Y\phi||_{Y=0}>0.
%\end{cases}
%\end{align}
%Rewrite this equation in the following form
%\begin{align*}
%&A^{-1}\left[i\epsilon(\partial_Y^2-\alpha^2)^2\phi+U_s(\partial_Y^2-\alpha^2)\phi\right]+(A^{-1}-1)\alpha^2(i\epsilon\Delta_\alpha\phi+U_s\phi)\notag\\
%&\quad+\partial_YA^{-1}(i\epsilon\partial_Y\Delta_\alpha\phi+U_s\partial_Y\phi)-\phi\partial_Y(A^{-1}\partial_YU_s)=0.
%\end{align*}
In this regime, we can directly construct a fast mode satifying \eqref{4.151}, without relying on slow mode. We start from the compressible Orr-Sommerfeld equation \eqref{4.200}. We aim to construct a fast mode $\phi_{f,m}$
, such that the following leading order $$i\eps\Delta_\alpha^2\phi_{f,m}+Y\Delta_{\alpha}\phi_{f,m}\approx 0.$$
For this, we set the scale of boundary sublayer $\delta=\eps^{-\frac{1}{6}\pi i}\eps^{\frac13}$. Let ${\rm Ai}(z)$ be the classical Airy function
and \begin{equation}W(Y)=\frac{{\rm Ai}(\delta^{-1}(Y+Y_0))}{{\rm Ai}(\delta^{-1}Y_0)}, {\rm ~with~} Y_0=-i\eps\alpha^2.\label{4.222-9}
\end{equation}
Then we introduce $\phi_{f,app,m}$ solving
\begin{align}\label{4.222}
\begin{cases}
(\partial_Y^2-\alpha^2)\phi_{f,app,m}={W}(Y),\\
\phi_{f,app,m}|_{Y=0}=0.
\end{cases}
\end{align}
It is straightforward to compute the following error term induced by $\phi_{f,app,m}$:
\begin{align}
E_{f,app,m}\eqdef{\rm OS_{CNS}}(\phi_{f,app,m})=&A^{-1}\left(U_s-U_s'(0)Y\right)W+\alpha^2(A^{-1}-1)\left(i\eps W+U_s\phi_{f,app,m}\right)\nonumber\\
&+\p_Y(A^{-1})\left(i\eps \p_YW+U_s\p_Y\phi_{f,app,m}\right)-\p_Y(A^{-1}\p_YU_s)\phi_{f,app,m}.\label{4.222-3}
\end{align}

Let $C_1\leq \alpha \eps^{\frac13}\leq C_2$, which corresponds to $C_1\nu^{-\frac34}\leq \hat{n}\leq C_2\nu^{-\frac34}$). Here $C_1$ and $C_2$ are two given positive constants. Then bounds on $\phi_{f,app,m}$ can be summarized in the following.
\begin{lemma}\label{lm5.2}
	Let $k\geq 0$ be an integer. Then it holds that
	\begin{align}
		\|U_s^k\p_Y\phi_{f,app,m}\|_{L^2}+\alpha\|U_s^k\phi_{f,app,m}\|_{L^2}\leq C\eps^{\frac{1}{2}+\frac{k}{3}}\label{4.222-1}.
	\end{align}
Moreover, on the boundary we have
\begin{align}\label{4.222-2}
	|\p_Y\phi_{f,app,m}(0)|\geq C\eps^{\frac13}.
\end{align}
\end{lemma}
\begin{proof}
Since the zeros of function ${\rm Ai}(z)$ are on the negative real axis, cf. \cite{AS}, and
$\delta^{-1}Y_0=e^{-\frac{1}{3}\pi i}(\alpha\eps^{\frac13})^2\sim 1$ in the middle frequency regime, we obtain ${\rm Ai}(\delta^{-1}Y_0)\sim 1$. Thus the profile \eqref{4.222-9} is well-defined. According to the decay of Airy function ${\rm Ai}(z)\sim z^{-\frac14}e^{-\frac{2}{3}z^{\frac32}}$ for $|z|\gg1$ and $|\arg z|\leq \frac{5\pi}{6}$,  cf. \cite{AS}, we can find $\lambda_0>0$, such that 
\begin{align}\label{4.222-4}
	|W(Y)|\leq Ce^{-\lambda_0 \eps^{-\frac13}Y}.
\end{align}
Then taking inner product of \eqref{4.222} with  $U_s^{2k}\phi_{f,app}$ for $k\geq 1$, using  boundary layer structure \eqref{4.222-4}, and $\alpha\sim \eps^{-\frac13}$, we obtain
\begin{align}
	\|U_s^k\p_Y\phi_{f,app,m}\|_{L^2}^2+\alpha^2\|U_s^k\phi_{f,app,m}\|_{L^2}^2&\leq C\|U_s^k\phi_{f,app,m}\|_{L^2}\left(\|U_s^kW\|_{L^2}+\|U_s^{k-1}\p_Y\phi_{f,app,m}\|_{L^2}\right)\nonumber\\
	&\leq C\alpha^{-2}\|U_s^kW\|_{L^2}^2+C\alpha^{-2}\|U_s^{k-1}\p_Y\phi_{f,app,m}\|_{L^2}^2\nonumber\\
	&\leq C \eps^{\frac{2k}3+1 } +C\eps^{\frac23}\|U_s^{k-1}\p_Y\phi_{f,app,m}\|_{L^2}^2.\nonumber
\end{align}
Thus, the bounds \eqref{4.222-1} can be obtained by induction. 

Now we show \eqref{4.222-2}. From \eqref{4.222}, we have
\begin{align}
	\phi_{f,app,m}(Y)=-\int_0^Ye^{-\alpha(Y-Y')} d Y'\int_{Y'}^\infty e^{\alpha(Y'-Y'')}W(Y'')d Y''.\nonumber
\end{align}
Thus taking derivative and evaluating  it at $Y=0$, we obtain
\begin{align}
\partial_Y\phi_{f,app,m}(0)=-\int_0^\infty e^{-\alpha Y}W(Y)d Y.\label{4.22-5}
\end{align}
To estimate \eqref{4.22-5}, we follow the argument as Lemma B.3 in \cite{CWZ23}. Introduce the anti-derivative ${\rm Ai}_0$ of the classical Airy function: $${\rm Ai}_0(z)=\int_{e^{\frac{i}{6}\pi i}z}^\infty{\rm Ai}(z')dz'.$$ 
Then it holds that
\begin{align}
	\partial_Y\phi_{f,app,m}(0)&=-\int_0^\infty e^{-\alpha Y}\frac{{\rm Ai}\left(e^{\frac{1}{6}\pi i}\eps^{-\frac13}(Y+Y_0)\right)}{{\rm Ai}(e^{\frac{1}{6}\pi i}\eps^{-\frac13}Y_0)}d Y\nonumber\\
	&=-\int_0^\infty e^{-\alpha Y}\frac{{\rm Ai}_0'\left(\eps^{-\frac13}(Y+Y_0)\right)}{{\rm Ai}_0'(\eps^{-\frac13}Y_0)}d Y\nonumber\\
	&=\frac{\eps^{\frac13}{\rm Ai}_0(\eps^{-\frac13}Y_0)}{{\rm Ai}_0'(\eps^{-\frac13}Y_0)}\left(1-\alpha\int_0^\infty e^{-\alpha Y}\frac{{\rm Ai}_0\left(\eps^{-\frac13}(Y+Y_0)\right)}{{\rm Ai}_0(\eps^{-\frac13}Y_0)}\right)\label{4.22-6}.
\end{align}
Since ${\rm Im}(\eps^{-\frac13}Y_0)=-(\alpha\eps^{\frac13})^2\leq-C_1$ for $\alpha\eps^{\frac13}\geq C_1,$ we can use Lemmas 8.2 and 8.3  in \cite{CLWZ20} respectively to obtain $$\left|\frac{{\rm Ai}_0(\eps^{-\frac13}Y_0)}{{\rm Ai}_0'(\eps^{-\frac13}Y_0)}\right|\gtrsim 1,~~{\rm  and}~~
\left|\frac{{\rm Ai}_0\left(\eps^{-\frac13}(Y+Y_0)\right)}{{\rm Ai}_0(\eps^{-\frac13}Y_0)}\right|\leq e^{-\frac{1}{3}\eps^{-\frac13}Y}.$$ Substituting these two bounds into \eqref{4.22-6}, we deduce that
\begin{align}
	|\partial_Y\phi_{f,app,m}(0)|&\geq \eps^{\frac13}\left|\frac{{\rm Ai}_0(\eps^{-\frac13}Y_0)}{{\rm Ai}_0'(\eps^{-\frac13}Y_0)}\right|\left(1-\alpha\int_0^\infty e^{-\frac{1}{3}\eps^{-\frac13}Y-\alpha Y} dY\right)\nonumber\\
	&\geq \frac{\eps^{\frac13}}{1+3\alpha\eps^{\frac13}}\left|\frac{{\rm Ai}_0(\eps^{-\frac13}Y_0)}{{\rm Ai}_0'(\eps^{-\frac13}Y_0)}\right|\geq C\eps^{\frac13},\nonumber
\end{align}
where $C_1\leq \alpha\eps^{\frac13}\leq C_2$ has been used in the last inequality. Thus, \eqref{4.222-2} follows and the proof of Lemma \ref{lm5.2} is complete.
\end{proof}
Now we can construct fast mode which is closed to $\phi_{f,app,m}$.
\begin{lemma}
\label{P4.18}
Let $m\in (0, 1)$ and $ C_1\nu^{-\frac34}\leq \hat{n}\leq C_2\nu^{-\frac34}$, where positive constants $C_1$ and $C_2$ are given. There exists a solution $\phi_{f,m}\in H^4(\mathbb{R}_+)\cap H_0^1(\mathbb{R}_+)$ to the homogeneous Orr-Sommerfied equation ${\rm OS_{CNS}}(\phi_{f,m})=0$ satisfying 
\begin{align}
\|(\partial_Y\phi_{f,m}, \alpha\phi_{f,m})\|_{L^2}&\leq C\eps^{\frac12},\label{4.224}\\
\|(\partial_Y^2-\alpha^2)\phi_{f,m}\|_{L^2}&\leq C \eps^{\frac16}.\label{4.225}
\end{align}
In particular, on the boundary it holds that
\begin{align}
|\partial_Y\phi_{f,m}(0)|\ge C\eps^{\frac13}.\label{4.225-1}
\end{align}
\end{lemma}
\begin{proof}
We look for the solution in the form $\phi_{f,m}=\phi_{f,app,m}+\phi_{r,m}$, where the remainder $\phi_{r,m}$ solves
\begin{align*}
\begin{cases}
{\rm OS_{CNS}}(\phi_{r,m})=-E_{f,app,m},\\
{\phi}_{r,m}|_{Y=0}=0.
\end{cases}
\end{align*}
Here $E_{f,app,m}$ is the error term defined in \eqref{4.222-3}. By using \eqref{4.222-1}, we have 
\begin{align}
	\|E_{f,app,m}\|_{L^2}&\lesssim \|Y^2W\|_{L^2}+\alpha^2\eps\|Y^2W\|_{L^2}+\alpha^2\|U_s^3\phi_{f,app,m}\|_{L^2}\nonumber\\
	&\qquad+\eps\|\p_YW\|_{L^2}+\|U_s\p_Y\phi_{f,app,m}\|_{L^2}+\|\phi_{f,app,m}\|_{L^2}\nonumber\\
	&\leq C\eps^{\frac56}.\nonumber
\end{align}
Then applying Proposition \ref{P4.6} with $g=\frac{1}{i\alpha}E_{f,app,m}$ and Corollary \ref{C4.9} to $\phi_r$ , we can obtain
\begin{equation}
	\|(\p_Y\phi_{r,m},\alpha\phi_{r,m})\|_{L^2}+\eps^{\frac13}\|\Delta_\alpha\phi_{r,m}\|_{L^2}\leq \eps^{\frac56}\nonumber.
\end{equation}
Combining this with \eqref{4.222-1}, we get \eqref{4.224} and \eqref{4.225}. Moreover, by interpolation we have
\begin{align}\label{4.225-2}
	\|\partial_{Y}\phi_{r,m}\|_{L^\infty}\leq C\|\partial_{Y}\phi_{r,m}\|_{L^2}^{\frac12}\|\partial_{Y}^2\phi_{r,m}\|_{L^2}^{\frac12}\leq C\eps^{\frac23}.
\end{align}
Recall that $\phi_{f,m}=\phi_{f,app,m}+\phi_{r,m}$. Then the estimate \eqref{4.225-1} follows from \eqref{4.222-2} and \eqref{4.225-2}. The proof is complete.
%From Appendix B, we get
%\begin{align*}
%\|(\partial_Y\phi_{app}, \alpha\phi_{app})\|_{L^2}\lesssim\nu^{\frac38},\quad \|(U_s\partial_Y\phi_{app}, U_s\alpha\phi_{app})\|_{L^2}\lesssim\nu^{\frac58},\quad |\partial_Y\phi_{app}(0)|\ge C\nu^{\frac14}.
%\end{align*}
%Hence, we have
%\begin{align*}
%\|f\|_{L^2}\lesssim\nu^{\frac58}.
%\end{align*}
%On the other hand, we obtain from Corollary \ref{C4.9}
%\begin{align*}
%\|(\partial_Y\tilde{\phi}, \alpha\tilde{\phi})\|_{L^2}\lesssim\|f\|_{L^2}\lesssim\nu^{\frac58},\quad \|(\partial_Y^2-\alpha^2)\tilde{\phi}\|_{L^2}\lesssim\epsilon^{-\frac13}\|f\|_{L^2}\lesssim\nu^{\frac38}.
%\end{align*}
%Thus for sufficiently small $\nu$, it holds
%\begin{align*}
%\|\partial_Y\tilde{\phi}\|_{L^\infty}\lesssim\nu^{\frac12}\ll \nu^{\frac14}.
%\end{align*}
%This finishes the proof together with \eqref{B.4}.  
\end{proof}
Finally, Based on $\phi_{f,m}$, we can recover the homogeneous quasi-compressible solution $(\varrho_{H,m},\mathfrak{u}_{H,m},\mathfrak{v}_{H,m})$ at middle frequencies.
\begin{proposition}
\label{P4.19}
Let $m\in (0, 1)$. If  $C_1\nu^{-\frac34}\leq \hat{n}\leq C_2\nu^{-\frac34}$, then the quasi-compressible system \eqref{4.151} admits a solution $(\varrho_{H,m}, \mathfrak{u}_{H,m}, \mathfrak{v}_{H,m})\in  H^2(\mathbb{R}_+)^3$ satisfying
\begin{align}
\|(\mathfrak{u}_{H,m}, \mathfrak{v}_{H,m})\|_{L^2}&\leq C\eps^{\frac12},\label{4.226}\\
m^{-2}\|\varrho_{H,m}\|_{H^2}+\|(\p_Y\mathfrak{u}_{H,m}, \p_Y\mathfrak{v}_{H,m})\|_{L^2}&\leq C\eps^{\frac16},\label{4.226-1}\\
\|\Delta_\alpha \mathfrak{u}_{H,m}\|_{L^2}&\leq C\eps^{-\frac56}\label{4.226-2},\\
 \|\Delta_\alpha \mathfrak{v}_{H,m}\|_{L^2}&\leq C\eps^{-\frac16}\label{4.226-3},\\
 \|(\partial_Y\div_\alpha(\mathfrak{u}_{H,m},\mathfrak{v}_{H,m}),\alpha\div_\alpha(\mathfrak{u}_{H,m},\mathfrak{v}_{H,m}))\|_{L^2}&\leq C\eps^{\frac16}.\label{4.226-4}
\end{align}
Moreover, on the boundary 
\begin{align}\label{4.229}
\mathfrak{v}_{H,m}(0)=0,~~{\rm and }~~|\mathfrak{u}_{H,m}(0)|\ge C\eps^{\frac13}, ~
\end{align}
\end{proposition}
\begin{proof}
	From Lemma \ref{P4.18} we can obtain
	$$\|(\partial_Y\phi_{f,m}, \alpha\phi_{f,m})\|_{L^2}\leq C\eps^{\frac12},~
	\|(\partial_Y^2-\alpha^2)\phi_{f,m}\|_{L^2}\leq C \eps^{\frac16}
	,~\|U_s'\phi_{f,m}\|_{L^2}\leq C\eps^{\frac16}.$$ Subtituting these bounds into Proposition \ref{P8.1} and using $\alpha\sim \eps^{-\frac13}$, we can obtain estimates \eqref{4.226}-\eqref{4.226-4}. The bound \eqref{4.229} follows from \eqref{4.225-1}, by noting that $\mathfrak{v}_{H,m}(0)=-i\alpha\phi_{f,m}(0)=0$  and $\mathfrak{u}_{H,m}(0)=\partial_Y\phi_{f,m}(0)$. The proof of Proposition \ref{P4.19} is complete.
\end{proof}

\subsection{High frequency regime $\hat{n}\gg \nu^{-\frac{3}{4}}$}
	Homogeneous solutions constructed in Sections \ref{S5.1}  and \ref{S5.2} do not work in this regime. 
	Our idea is to construct a fast mode $\phi_{f,h}$ to the Orr-Sommerfeld equation around the leading order profile 
\begin{align}\label{hl}
	\phi_{f,app,h}(Y)=\frac{Ye^{-\alpha Y}}{2\alpha},
\end{align}
which satisfies
$$
\Delta_\alpha^2\phi_{f,app,h}(Y)=0,~{\rm and}~\phi_{f,app,h}(0)=0.
$$
%\begin{align}\label{4.230}
%\begin{cases}
%\mathcal{L}_Q(\rho_{q}, u_{q}, v_{q})=0,\\
%v_q|_{Y=0}=0,\quad u_q|_{Y=0}>0.
%\end{cases}
%\end{align}
%Similar as \eqref{4.6}, there exists a function $\phi$ such that
%\begin{align}\label{4.231}
%\partial_Y\phi=u_q+U_s\rho_q,\quad -i\alpha\phi=v_q,\quad \phi|_{Y=0}=0,
%\end{align}
%and satisfy the homogeneous compressible Orr-Sommerfeld equation
%\begin{align}\label{4.232}
%\begin{cases}
%OS_{CNS}(\phi)=\frac{i}{\hat{n}}\Lambda(\Delta_\alpha\phi)+U_s\Lambda(\phi)-\phi\partial_Y(A^{-1}\partial_YU_s)=0,\\
%\phi|_{Y=0}=0,\quad \partial_Y\phi|_{Y=0}>0.
%\end{cases}
%\end{align}
\begin{lemma}\label{P4.20}
Let $m\in (0, 1)$. There exists $\kappa_1\in (0, 1)$ such that if $\hat{n}\ge \kappa_1^{-1}\nu^{-\frac34}$, there exists a solution $\phi_{f,h}\in H^4(\mathbb{R}_+)\cap H_0^1(\mathbb{R}_+)$ to the homogeneous Orr-Sommerfeld equation ${\rm OS_{CNS}}(\phi)=0$. Moreover, the solution satisfies
\begin{align}
\|(\partial_Y\phi_{f,h}, \alpha\phi_{f,h})\|_{L^2}&\leq C\alpha^{-\frac32},\label{4.233}\\
\|\Delta_\alpha\phi_{f,h}\|_{L^2}&\leq C\alpha^{-\frac12},\label{4.234}\\
\|\sqrt{U_s}\phi_{f,h}\|_{L^2}&\leq C\alpha^{-3},\label{4.235}\\
\|\partial_Y\Delta_\alpha\phi_{f,h}\|_{L^2}&\leq C\alpha^{\frac12}.\label{4.333}
\end{align}
In particular, on the boundary 
\begin{align}
\left|\partial_Y\phi_{f,h}(0)\right|\ge\frac{1}{4\alpha}.\label{4.237-1}
\end{align}
\end{lemma}
%\begin{remark}
%The weighted estimate \eqref{4.235} and the third order derivative estimate \eqref{4.333} are essential for the boundary corrector, which is hard to obtain in low and middle frequency.
%\end{remark}
\begin{proof}
We look for the solution $\phi_{f,h}$ in the form
$\phi_{f,h}=\phi_{f,app,h}+\phi_{r,h}$ where $\phi_{f,app,h}$ is defined in \eqref{hl}, and $\phi_{r,h}$ is the remainder. It is straightforward to check that $\phi_{f,app,h}$ satisfies the bounds \eqref{4.233}-\eqref{4.333}.  Direct computation yields
\begin{align}
\phi_{f,app,h}|_{Y=0}=0,~\text{ and }~ \partial_Y\phi_{f,app,h}|_{Y=0}=\frac{1}{2\alpha}.\label{4.237}
\end{align}
Then $\phi_{r,h}$ satisfies
\begin{align}
\begin{cases}
{\rm OS_{CNS}}(\phi_{r,h})=-E_{f,app,h},\quad Y>0,\\
\phi_{r,h}|_{Y=0}=0.\nonumber
\end{cases}
\end{align}
Here $E_{f,app,h}$ is the error term induced by $\phi_{f,app,h}$:
\begin{align}%
E_{f,app,h}&=i\epsilon e^{-\alpha Y}\left[\alpha \partial_Y(A^{-1})+\alpha^2(A^{-1}-1)\right]-\partial_Y(A^{-1}\partial_YU_s)\phi_{f,app,h}   \notag\\
&\quad+U_s\left[(\partial_YA^{-1})\partial_Y\phi_{f,app,h}-A^{-1}e^{-\alpha Y}+(A^{-1}-1)\alpha^2\phi_{f,app,h}\right].\nonumber
\end{align}
Using $|A^{-1}(Y)-1|\leq CY^2$ and $|\p_Y A(Y)|\leq CY$, we obtain
\begin{align}
	\|E_{f,app,h}\|_{L^2}\leq& C\eps \left(\alpha^2\|Y^2e^{-\alpha Y}\|_{L^2}+\alpha\|Ye^{-\alpha Y}\|_{L^2}\right)+C\|\phi_{f,app,h}\|_{L^2}\nonumber\\
	&+C\|Y^2\partial_Y\phi_{f,app,h}\|_{L^2}+C\|Ye^{-\alpha Y}\|_{L^2}+C\alpha^2\|Y^3\phi_{f,app,h}\|_{L^2}\nonumber\\
	\leq& C\alpha^{-\frac32}(1+\eps\alpha)\leq C\alpha^{-\frac32},\nonumber
\end{align}
where $\eps\alpha=\sqrt{\nu}\leq 1$ has been used in the last inequality. 

Let $\delta_1$ be the positive number in Proposition \ref{prop4.8}. Then applying the bounds \eqref{4.239}-\eqref{4.343} to $\phi_{r,h}$ with $f=E_{f,app,h}$, we have, for $\kappa_1\in (0,\delta_1)$, that
\begin{equation}
\begin{aligned}\label{4.249}
\|(\partial_Y{\phi}_{r,h}, \alpha{\phi}_{r,h})\|_{L^2}&\leq \frac{C}{\eps\alpha ^3}\|E_{f,app,h}\|_{L^2}\leq \frac{C}{\eps \alpha^{\frac92}},\\
\|\Delta_{\alpha}\phi_{r,h}\|_{L^2}&\leq \frac{C}{\eps\alpha ^2}\|E_{f,app,h}\|_{L^2}\leq \frac{C}{\eps \alpha^{\frac72}},\\
\|\sqrt{U_s}\phi_{r,h}\|_{L^2}&\leq\frac{C}{\eps^{\frac12}\alpha^3}\|E_{f,app,h}\|_{L^2}\leq \frac{C}{\eps^{\frac12}\alpha^{\frac92}},\\
\|\partial_Y\Delta_\alpha\phi_{r,h}\|_{L^2}&\leq 
\frac{C}{\eps \alpha}\|E_{f,app,h}\|_{L^2}\leq \frac{C}{\eps \alpha^{\frac52}}.
\end{aligned}
\end{equation}
Then bounds in \eqref{4.233}-\eqref{4.333} follow from $\alpha\eps^{\frac13}\gtrsim 1$ and \eqref{4.249}. Moreover,  we obtain from Sobolev inequality that
\begin{align*}
\|\partial_Y\phi_{r,f}\|_{L^\infty}&\leq C\|\partial_Y{\phi}_{r,f}\|_{L^2}^{\frac12}\|\partial_Y^2\phi_{r,f}\|_{L^2}^{\frac12}\leq\frac{C}{\eps\alpha^4}\nonumber\\
&\leq \frac{C}{\alpha(\hat{n}\nu^{\frac34})^2}\leq \frac{C\kappa_1^2}{\alpha}\leq \frac{1}{4\alpha},
\end{align*}
provided that $\kappa_1$ is sufficiently small. Combining this with \eqref{4.237}, we obtain \eqref{4.237-1}. The proof of Lemma \ref{P4.20} is complete.
\end{proof}
Substituting \eqref{4.233}-\eqref{4.333} into bounds \eqref{8.1-1}-\eqref{8.1-7} in the high-frequency regime,  we can recover the homogeneous solution $(\varrho_{H,h},\mathfrak{u}_{H,h},\mathfrak{v}_{H,h})$ in the following
\begin{proposition}
\label{P4.22}
Let $m\in (0, 1)$ and $\kappa_1$ be the constant in Lemma \ref{P4.20}. If $\hat{n}\geq \kappa_1^{-1}\nu^{-\frac34}$, then the system \eqref{4.151} admits a solution $(\varrho_{H,h}, \mathfrak{u}_{H,h}, \mathfrak{v}_{H,h})\in  H^2(\mathbb{R}_+)^3$ satisfying
\begin{align}
\|(\mathfrak{u}_{H,h},\mathfrak{v}_{H,h})\|_{L^2}&\leq C\alpha^{-\frac32},\label{4.252}\\
\|(\partial_Y\mathfrak{u}_{H,h},\partial_Y\mathfrak{v}_{H,h})\|_{L^2}&\leq C\alpha^{-\frac12},\label{4.253}\\
\|(\Delta_\alpha\mathfrak{u}_{H,h},\Delta_\alpha\mathfrak{v}_{H,h})\|_{L^2}&\leq C\alpha^{\frac12},\label{4.254}\\
m^{-2}\|(\alpha\varrho_{H,h},\partial_Y\varrho_{H,h},\Delta_{\alpha}\varrho_{H,h})\|_{L^2}&\leq C(\alpha^{-\frac12}+\alpha^{\frac12}\nu^{\frac12})\label{4.254-1},\\
\|(\partial_Y\div_\alpha(\mathfrak{u}_{H,h},\mathfrak{v}_{H,h}),\alpha \div_\alpha(\mathfrak{u}_{H,h},\mathfrak{v}_{H,h}))\|_{L^2}&\leq C(\alpha^{-\frac12}+\alpha^{\frac12}\nu^{\frac12}).\label{4.254-2}
\end{align}
In particular, on the boundary 
\begin{align}\label{4.255}
\mathfrak{v}_{H,h}(0)=0,~\text{ and }~\mathfrak{u}_{H,h}(0)\ge\frac{1}{4\alpha}.
\end{align}
\end{proposition}

\section{Linearized system with only $v|_{Y=0}=0$}\label{S6}
In this section, we will estabilish solvability of the following linear Navier-Stokes equations 
\begin{align}\label{9.1}
	\begin{cases}
		i\alpha U_s\rho_{n,sl}+\div_\alpha(u_{n,sl},v_{n,sl})=f_{\rho, n},\\
		i\alpha U_su_{n,sl}+v_{n,sl}\partial_YU_s+i\alpha m^{-2}\rho_{n,sl}-\sqrt{\nu}\Delta_\alpha u_{n,sl}-\lambda\sqrt{\nu} i\alpha\div_\alpha(u_{n,sl},v_{n,sl})=f_{u, n},\\
		i\alpha U_s v_{n,sl}+m^{-2}\partial_Y\rho_{n,sl}-\sqrt{\nu}\Delta_\alpha v_{n,sl}-\lambda\sqrt{\nu}\partial_Y\div_\alpha(u_{n,sl},v_{n,sl})=f_{v, n},
	\end{cases}
\end{align}
with the boundary condition
\begin{align}
v_{n,sl}|_{Y=0}=0.\label{9.1-1}
\end{align}
Here $f_{\rho,n}\in H^1(\mathbb{R}_+)$ and $(f_{u,n},f_{v,n})\in L^2(\mathbb{R}_+)^2$ are inhomogeneous source terms. Note that we only impose the zero normal velocity condition \eqref{9.1-1} on the boundary. For simplicity, we denote  $\mathcal{D}_{sl}\eqdef i\alpha u_{n,sl}+\partial_Yv_{n,sl}$.

The main results in the section are the following
\begin{proposition}[Low and middle frequencies]
	\label{C5.3}
Let $m\in (0,1)$ and $\hat{n}\lesssim \nu^{-\frac34}$. Then for $f_{\rho,n}\in H^1(\mathbb{R}_+)$ and $(f_{u,n},f_{v,n})\in L^2(\mathbb{R}_+)^2$, there exists a solution $(\rho_{n,sl} ,u_{n,sl}. v_{n,sl})\in H^1(\mathbb{R}_+)\times H^2(\mathbb{R}_+)^2$ to \eqref{9.1}. Moreover, the following two statements hold.
\begin{itemize}
	\item[{\rm(i)}] If $\alpha\ge 1$,
	\begin{align}
		m^{-2}\|\rho_{n,sl}\|_{H^1}+\alpha\|(u_{n,sl}, v_{n,sl})\|_{L^2}&\le\frac{C}{\hat{n}^{\frac13}\sqrt{\nu}}(\|(f_{\rho,n},f_{u,n},f_{v,n})\|_{L^2}+\sqrt{\nu}\|(\partial_Y f_{\rho,n}, \alpha f_{\rho,n})\|_{L^2}),\label{5.60}\\
		\|\partial_Yv_{n,sl}\|_{L^2}+	\|(\partial_{Y}\mathcal{D}_{sl},\alpha\mathcal{D}_{sl})\|_{L^2}&\leq \frac{C}{\hat{n}^{\frac13}\sqrt{\nu}}(\|(f_{\rho,n},f_{u,n},f_{v,n})\|_{L^2}+\sqrt{\nu}\|(\partial_Y f_{\rho,n}, \alpha f_{\rho,n})\|_{L^2}), 
		\label{5.62}\\
		\|\partial_Yu_{n,sl}\|_{L^2}&\le \frac{C}{\sqrt{\nu}}(\|(f_{\rho,n},f_{u,n},f_{v,n})\|_{L^2}+\sqrt{\nu}\|(\partial_Y f_{\rho,n}, \alpha f_{\rho,n})\|_{L^2}),\label{5.61}\\
		\|\partial_Y^2u_{n,sl}\|_{L^2}+\hat{n}^{\frac12}\|\partial_Y^2v_{n,sl}\|&\le \frac{C\hat{n}^{\frac23}}{\sqrt{\nu}}(\|(f_{\rho,n},f_{u,n},f_{v,n})\|_{L^2}+\sqrt{\nu}\|(\partial_Y f_{\rho,n}, \alpha f_{\rho,n})\|_{L^2}).\label{5.63}
	\end{align}
On the boundary, it holds that
\begin{align}
	|u_{n,sl}(0)|\leq \frac{C}{\alpha^{\frac12}\hat{n}^{\frac16}\sqrt{\nu}}(\|(f_{\rho,n},f_{u,n},f_{v,n})\|_{L^2}+\sqrt{\nu}\|(\partial_Y f_{\rho,n}, \alpha f_{\rho,n})\|_{L^2}).\label{5.62-1}
\end{align}
\item[{\rm(ii)}] If $0<\alpha\le 1$,
	\begin{align}\label{5.64}
		m^{-2}\|\rho_{n,sl}\|_{H^1}+\|(u_{n,sl}, v_{n,sl})\|_{L^2}&
		\leq \frac{C\hat{n}^{\frac23}}{\alpha}(\|(f_{\rho,n},f_{u,n},f_{v,n})\|_{L^2}+\sqrt{\nu}\|\partial_Y f_{\rho,n}\|_{L^2}),\\
		\|\partial_Yv_{n,sl}\|_{L^2}+\|\partial_{Y}\mathcal{D}_{sl},\alpha\mathcal{D}_{sl}\|_{L^2}&\leq C\hat{n}^{\frac23}(\|(f_{\rho,n},f_{u,n},f_{v,n})\|_{L^2}+\sqrt{\nu}\|\partial_Y f_{\rho,n}\|_{L^2}),\label{5.65}\\
			\|\partial_Yu_{n,sl}\|_{L^2}&\leq \frac{C}{\sqrt{\nu}}(\|(f_{\rho,n},f_{u,n},f_{v,n})\|_{L^2}+\sqrt{\nu}\|\partial_Y f_{\rho,n}\|_{L^2}),\label{5.65-1}\\
		\|(\partial_Y^2u_{n,sl}, \partial_Y^2v_{n,sl})\|_{L^2}&\le \frac{C\hat{n}^{\frac23}}{\sqrt{\nu}}(\|(f_{\rho,n},f_{u,n},f_{v,n})\|_{L^2}+\sqrt{\nu}\|\partial_Y f_{\rho,n}\|_{L^2}).\label{5.66}
	\end{align}
On the boundary, it holds that
\begin{align}
	|u_{n,sl}(0)|\leq \frac{C}{\hat{n}^{\frac16}\sqrt{\nu}}(\|(f_{\rho,n},f_{u,n},f_{v,n})\|_{L^2}+\sqrt{\nu}\|\partial_Y f_{\rho,n}\|_{L^2}).\label{5.62-2}
\end{align}
\end{itemize}
\end{proposition}
\begin{proposition}[High frequencies]
	\label{P5.4}
	Let $m\in (0, 1)$. There exists a positive constant $\kappa_2\in (0, \kappa_1)$, where $\kappa_1$ is given in the Proposition \ref{P4.22}, such that if $\hat{n}\ge \kappa_2^{-1}\nu^{-\frac34}$,
	  then for any $f_{\rho,n}\in H^1(\mathbb{R}_+)$ and $(f_{u,n}, f_{v,n})\in L^2(\mathbb{R}_+)^2$, the system \eqref{9.1} admits a solution $(\rho_{n,sl}, u_{n,sl}, v_{n,sl})\in H^1(\mathbb{R}_+)\cap H^2(\mathbb{R}_+)^2$ satisfying
	\begin{align}\label{5.75}
		\|(\partial_Yu_{n,sl}, \alpha u_{n,sl})\|_{L^2}+\|(\partial_Yv_{n,sl}, \alpha v_{n,sl})\|_{L^2}&\le \frac{C}{\alpha\sqrt{\nu}}(\|(f_{\rho,n},f_{u,n},f_{v,n})\|_{L^2}+\sqrt{\nu}\|(\partial_Y f_{\rho,n}, \alpha f_{\rho,n})\|_{L^2}),\\
		m^{-2}\|(\partial_Y\rho_{n,sl}, \alpha\rho_{n,sl})\|_{L^2}&\le  \frac{C}{\nu^{\frac14}}(\|(f_{\rho,n},f_{u,n},f_{v,n})\|_{L^2}+\sqrt{\nu}\|(\partial_Y f_{\rho,n}, \alpha f_{\rho,n})\|_{L^2}),\label{5.76}\\
		\|\Delta_\alpha u_{n,sl}, \Delta_\alpha v_{n,sl})\|_{L^2}&\le \frac{C}{\nu^{\frac34}}(\|(f_{\rho,n},f_{u,n},f_{v,n})\|_{L^2}+\sqrt{\nu}\|(\partial_Y f_{\rho,n}, \alpha f_{\rho,n})\|_{L^2}).\label{5.77}
	\end{align}
Moreover,  it holds that
\begin{align}\label{5.78}
	|u_{n,sl}(0)|\leq \frac{C}{\alpha^{\frac32}\nu^{\frac12}}(\|(f_{\rho,n},f_{u,n},f_{v,n})\|_{L^2}+\sqrt{\nu}\|(\partial_Y f_{\rho,n}, \alpha f_{\rho,n})\|_{L^2}).
\end{align}
\end{proposition}

In the low and middle frequency regimes we solve the linear Navier-Stokes system \eqref{9.1} by using the quasi-compressible-Stokes iteration. In the subsection \ref{S6.1}, we introduce the Stokes regularizing system. The iteration scheme and its convergence will be shown in the subection \ref{S6.2}. Finally, we will show the solvability of \eqref{9.1} in high frequences in the subsection \ref{S6.3}. In this section, the symbol $\langle f,g\rangle$ represents the standard $L^2(\mathbb{R}_+)$ inner product with respect to $Y$ variable.
\subsection{Stokes regularization}\label{S6.1}
Note that the error terms created by the quasi-compressible system do not have enough regularity to iterate. To smooth out these singular errors,  we introduce the following Stokes regularizing system 
\begin{align}\label{5.1}
\begin{cases}
i\alpha U_s\udr+\mathrm{div}_\alpha(\udu, \udv)=q_\rho,\\
i\alpha U_s\udu+(i\alpha m^{-2}+\sqrt{\nu}\partial_Y^2U_s)\udr-\sqrt{\nu}\Delta_\alpha\udu-\lambda\sqrt{\nu} i\alpha\mathrm{div}_\alpha(\udu, \udv)=q_u,\\
i\alpha U_s\udv+m^{-2}\partial_Y\udr-\sqrt{\nu}\Delta_\alpha\udv-\lambda\sqrt{\nu}\partial_Y\mathrm{div}_\alpha(\udu, \udv)=q_{v},\\
\partial_Y\udu|_{Y=0}=\udv|_{Y=0}=0,
\end{cases}
\end{align}
with a given inhomogeneous source term $(q_\rho,q_u,q_v)$. Compared with original linear system \eqref{9.1}, we remove the stretching term $\udv\partial_YU_s$ in the second equation, and replace boundary condition on $\udu$ by $\partial_Y\udu|_{Y=0}=0$. For simplicity, we denoted \eqref{5.1} by $\mathcal{L}_S(\udr,\udu,\udv)=(q_\rho, q_u, q_v)$.
\begin{proposition}
\label{P5.1}
Let $m\in (0, 1)$, $q_\rho\in H^1(\mathbb{R}_+)$ and $q_u, q_v\in L^2(\mathbb{R}_+)$. The Stokes regularizing system \eqref{5.1} admits a unique solution $(\udr, \udu, \udv)\in H^1(\mathbb{R}_+)\times H^2(\mathbb{R}_+)^2$ satisfying
\begin{align}\label{5.2}
\|(\udu,\udv)\|_{L^2}&\le \frac{C\hat{n}^{\frac13}}{\alpha}(\|(q_\rho,q_u,q_v)\|_{L^2}+\sqrt{\nu}\|(\partial_Y q_\rho, \alpha q_\rho)\|_{L^2}),\\
m^{-2}\|(\partial_Y\udr, \alpha\udr)\|_{L^2}+\|\mathrm{div}_\alpha(\udu,\udv)\|_{L^2}&\le C\hat{n}^{\frac16}(\|(q_\rho,q_u,q_v)\|_{L^2}+\sqrt{\nu}\|(\partial_Y q_\rho, \alpha q_\rho)\|_{L^2}),\label{5.3}\\
\|\partial_Y\mathrm{div}_\alpha(\udu, \udv)\|_{L^2}&\le C(1+\alpha)\hat{n}^{\frac16}(\|(q_\rho,q_u,q_v)\|_{L^2}+\sqrt{\nu}\|(\partial_Y q_\rho, \alpha q_\rho)\|_{L^2}),\label{5.4}\\
\|(\partial_Y\udu, \partial_Y\udv)\|_{L^2}&\le \frac{C\hat{n}^{\frac23}}{\alpha}(\|(q_\rho,q_u,q_v)\|_{L^2}+\sqrt{\nu}\|(\partial_Y q_\rho, \alpha q_\rho)\|_{L^2}),\label{5.5}\\
\|(\Delta_\alpha\udu, \Delta_\alpha\udv)\|_{L^2}&\le C(\alpha+\nu^{-\frac12})\hat{n}^{\frac16}(\|(q_\rho,q_u,q_v)\|_{L^2}+\sqrt{\nu}\|(\partial_Y q_\rho, \alpha q_\rho)\|_{L^2}),\label{5.6}
\end{align}
for some positive constant $C$ independent of $\nu$, $m$ and $\alpha$.  Moreover, if $\alpha\in (0,1)$, it holds that
\begin{align}
	\|\partial_Y\udv\|_{L^2}\leq  C\hat{n}^{\frac13}(\|(q_\rho,q_u,q_v)\|_{L^2}+\sqrt{\nu}\|(\partial_Y q_\rho, \alpha q_\rho)\|_{L^2}).\label{5.6-1}
\end{align}
\end{proposition}
\begin{proof}
The existence of solution follows the argument in \cite[Proposition 3.9]{YZ23}. Thus, we focus on the a priori estimates for brevity. The proof is divided into four steps.

\underline{\it Step 1. Elliptic estimate.} Taking inner product of \eqref{5.1}$_2$ and \eqref{5.1}$_3$ with $(\udu, \udv)$ respectively, then integrating by parts, we derive
\begin{align}\label{5.7}
&\sqrt{\nu}\left(\|(\partial_Y\udu, \alpha\udu)\|_{L^2}^2+\|(\partial_Y\udv, \alpha\udv)\|_{L^2}^2\right)+\lambda\sqrt{\nu}\|\mathrm{div}_\alpha(\udu, \udv)\|_{L^2}^2+i\alpha \|\sqrt{U_s}(\udu,\udv)\|_{L^2}^2\notag\\
&\qquad=m^{-2}\int_0^\infty \udr\overline{\mathrm{div}_\alpha(\udu, \udv)}\;dY+\int_0^\infty (q_u-\sqrt{\nu}\partial_Y^2U_s\udr)\bar{\udu}+q_v\bar{\udv}\;dY.
\end{align} 
By Schwarz inequality, we get
\begin{align}\label{5.8}
\left|\int_0^\infty (q_u-\sqrt{\nu}\partial_Y^2U_s\udr)\bar{\udu}+q_v\bar{\udv}\;dY\right|\le \sqrt{\nu}\|\udr\|_{L^2}\|\udu\|_{L^2}+\|q_u\|_{L^2}\|\udu\|_{L^2}+\|q_v\|_{L^2}\|\udv\|_{L^2}.
\end{align}
For the first term on the right-hand side of \eqref{5.7}, we use the continuity equation \eqref{5.1}$_1$ to get
\begin{align}\label{5.9}
m^{-2}\int_0^\infty \udr\overline{\mathrm{div}_\alpha(\udu, \udv)}\;dY=m^{-2}i\alpha \|\sqrt{U_s}\udr\|_{L^2}+m^{-2}\int_0^\infty \udr\bar{q}_\rho dY.
\end{align}
%Then the real part implies
%\begin{align}\label{5.9}
%Re\left(m^{-2}\int_0^\infty \zeta\overline{\mathrm{div}_\alpha(\theta, \vartheta)}\;dY\right)\le \|m^{-1}\zeta\|_{L^2}\|m^{-1}q_0\|_{L^2}.
%\end{align}
Plugging \eqref{5.8} and \eqref{5.9} into \eqref{5.7}, then taking real and imaginary parts  respectively, we obtain
\begin{align}\label{5.10}
&\sqrt{\nu}\left(\|(\partial_Y\udu, \alpha\udu)\|_{L^2}^2+\|(\partial_Y\udv, \alpha\udv)\|_{L^2}^2\right)+\lambda\sqrt{\nu}\|\mathrm{div}_\alpha(\udu, \udv)\|_{L^2}^2\nonumber\\
&\qquad\leq C\|m^{-2}\udr\|_{L^2}\|q_\rho\|_{L^2}+C\|(\udu,\udv)\|_{L^2}\|(q_u,q_v)\|_{L^2}+C\sqrt{\nu}\|\udr\|_{L^2}\|\udu\|_{L^2},
\end{align}
and
\begin{align}\label{5.11}
\alpha\|(\sqrt{U_s}\udu, \sqrt{U_s}\udv)\|_{L^2}^2\le& \alpha m^{-2}\|\sqrt{U_s}\udr\|_{L^2}^2+C\sqrt{\nu}\|\udr\|_{L^2}\|\udu\|_{L^2}\nonumber\\
&+C\|m^{-2}\udr\|_{L^2}\|q_\rho\|_{L^2}+C\|(\udu,\udv)\|_{L^2}\|(q_u,q_v)\|_{L^2}.
\end{align}
Note that unlike the time-dependent problem studied in \cite{YZ23}, for the steady case,  the elliptic estimate \eqref{5.11} is not closed. We have to obtain bounds on $\udr$.

\underline{\it Step 2. Estimates on density}. Denote the vorticity and divergence by $\omega:=\partial_Y\udu-i\alpha\udv$ and  $\mathcal{D}:=\mathrm{div}_\alpha(\udu, \udv)$ respectively. Then momentum equations  $\eqref{5.1}_1$ and $\eqref{5.1}_2$ can be written as
%\begin{align}\label{5.12}
%\Delta_\alpha\theta=\partial_Y\mathcal{A}+i\alpha\mathcal{B},\quad \Delta_\alpha\vartheta=-i\alpha\mathcal{A}+\partial_Y\mathcal{B}.
%\end{align}
%Furthermore $\mathcal{A}|_{Y=0}=0$ and the second and the third equation of system becomes
\begin{align}\label{5.13}
i\alpha m^{-2}\udr&=\sqrt{\nu}\partial_Y\omega+(1+\lambda)\sqrt{\nu}i\alpha\mathcal{D}-i\alpha U_s\udu-\sqrt{\nu}\udr\partial_Y^2U_s+q_u,\\
m^{-2}\partial_Y\udr&=-\sqrt{\nu}i\alpha\omega+(1+\lambda)\sqrt{\nu}\partial_Y\mathcal{D}-i\alpha U_s\udv+q_v.\label{5.14}
\end{align}
Note that boundary conditions in \eqref{5.1} imply $\omega|_{Y=0}=0$.

Taking inner product of \eqref{5.13}-\eqref{5.14} with $\overline{i\alpha\udr}$ and $\overline{\partial_Y\udr}$ respectively and using the boundary condition $\omega|_{Y=0}=0$, we have
\begin{align}\label{5.15}
m^{-4}\|(\partial_Y\udr, \alpha\udr)\|_{L^2}^2&=\underbrace{(1+\lambda)\sqrt{\nu}\left(\langle i\alpha\mathcal{D}, m^{-2}\overline{i\alpha \udr}\rangle+\langle\partial_Y\mathcal{D}, m^{-2}\overline{\partial_Y\udr}\rangle\right)}_{J_1}\notag\\
&\quad+\underbrace{\langle -i\alpha U_s\udu-\sqrt{\nu}\udr\partial_Y^2U_s+q_u, m^{-2}\overline{i\alpha\udr}\rangle+\langle -i\alpha U_s\udv+q_v, m^{-2}\partial_Y\bar{\udr}\rangle}_{J_2}.
\end{align}
From the continuity equation $\eqref{5.1}_1$, we get
\begin{align*}
\mathcal{D}=q_\rho-i\alpha U_s\udr,\quad \partial_Y\mathcal{D}=\partial_Yq_\rho-i\alpha U_s\partial_Y\udr-i\alpha\partial_YU_s\udr.
\end{align*}
Substituting these into $J_1$, we get
\begin{align*}
J_1=(1+\lambda)\sqrt{\nu}\left(-i\alpha m^{-2}\|(\sqrt{U_s}\partial_Y\udr, \sqrt{U_s}\alpha\udr)\|_{L^2}^2+ m^{-2}\langle i\alpha q_\rho, \overline{i\alpha\udr}\rangle+\langle \partial_Yq_0-i\alpha\partial_YU_s\udr, \partial_Y\bar{\udr}\rangle\right)
\end{align*}
Thus, real part of \eqref{5.15} gives
\begin{align}\label{5.16}
m^{-4}\|(\partial_Y\udr, \alpha\udr)\|_{L^2}^2=&{\rm Re}J_2+(1+\lambda)\sqrt{\nu}{\rm Re}\langle i\alpha q_\rho, m^{-2}\overline{i\alpha\udr}\rangle\rangle\nonumber\\
&+(1+\lambda)\sqrt{\nu}
{\rm Re}\langle \partial_Yq_\rho-i\alpha\partial_YU_s\udr,m^{-2} \partial_Y\bar{\udr}\rangle.
\end{align}

Now we bound the right hand side of \eqref{5.16}. By Schwarz inequality, we obtain 
$$|J_2|\leq \left(\alpha\|U_s(\udu,\udv)\|_{L^2}+\|(q_u,q_v)\|_{L^2}\right)m^{-2}\|(\partial_Y\udr,\alpha\udr)\|_{L^2}+\frac{Cm^{-2}}{\hat{n}}\|\alpha\udr\|_{L^2}^2,
$$
$$
\begin{aligned}
\left|(1+\lambda)\sqrt{\nu}{\rm Re}\langle i\alpha q_0, m^{-2}\overline{i\alpha\udr}\rangle\rangle\right|
&\leq Cm^{-2}\sqrt{\nu}\|\alpha\udr\|_{L^2}\|\alpha q_0\|_{L^2}\\
\left|(1+\lambda)\sqrt{\nu}
{\rm Re}\langle \partial_Yq_\rho-i\alpha\partial_YU_s\udr, m^{-2}\partial_Y\bar{\udr}\rangle\right|&\leq Cm^{-2}\sqrt{\nu}\|\alpha\udr\|_{L^2}\|\partial_Y\udr\|_{L^2}+Cm^{-2}\sqrt{\nu}\|\partial_Yq_\rho\|_{L^2}\|\partial_Y\udr\|_{L^2}.
\end{aligned}
$$
Substituting these bounds back into \eqref{5.16} and using $0\leq U_s\leq 1$, we get
\begin{align}
	m^{-2}\|(\partial_Y\udr, \alpha\udr)\|_{L^2}\leq&
	\alpha\|\sqrt{U_s}(\udu,\udv)\|_{L^2}+C\sqrt{\nu}\|(\partial_Y q_\rho, \alpha q_\rho)\|_{L^2}\nonumber\\
	&+C\|(q_u, q_v)\|_{L^2}+C(\sqrt{\nu}+\frac{1}{\hat{n}})\|(\partial_Y\udr, \alpha\udr)\|_{L^2}.\label{rho}
\end{align}
Then substituting \eqref{5.11}  into the above estimate, using Cauchy-Schwarz and $\frac1{\hat{n}}\leq L$, we obtain, for any $\eta\in (0,1)$, that
\begin{align}
m^{-2}\|(\partial_Y\udr, \alpha\udr)\|_{L^2}
\leq &    \alpha m^{-1}\|\sqrt{U_s}\udr\|_{L^2}+C(\sqrt{\nu}+L+\eta)\|(\partial_Y\udr,\alpha\udr)\|_{L^2}^2+C_\eta \sqrt{\nu}\|\udu\|_{L^2}\nonumber\\
&+C\alpha^{\frac12}\|(\udu,\udv)\|_{L^2}^{\frac12}\|(q_u,q_v)\|_{L^2}^{\frac12}+C_\eta\|(q_\rho,q_u,q_v)\|_{L^2}+C\sqrt{\nu}\|(\partial_Yq_\rho,\alpha q_\rho)\|_{L^2},
\label{5.16-1}
\end{align}

 Since $m\in (0, 1)$ and $0\leq U_s\leq 1$, then we can take $\eta$, $L$ and $\nu$ sufficiently small, such that the first two terms on the right hand side of \eqref{5.16-1} can be absorbed by the left hand side. Thus, we deduce that
\begin{align}
\label{5.17}
m^{-2}\|(\partial_Y\udr, \alpha\udr)\|_{L^2}\le& C\alpha^{\frac12}\|(\udu, \udv)\|_{L^2}^{\frac12}\|(q_u, q_v)\|_{L^2}^{\frac12}+C\|(q_\rho, q_u, q_v)\|_{L^2}\nonumber\\
&+C\sqrt{\nu}\|\partial_Yq_\rho,\alpha q_\rho\|_{L^2}+C\sqrt{\nu}\|\udu\|_{L^2}.
\end{align}

\underline{\it Step 3. Bound on $\|(\udu,\udv)\|_{L^2}$}. Substituting \eqref{5.17} to \eqref{5.10} and \eqref{5.11}, then using $|\sqrt{U_s}\udr|\leq |\udr|$, we obtain
\begin{align}\label{5.18}
&\sqrt{\nu}\left(\|(\partial_Y\udu, \alpha\udu)\|_{L^2}^2+\|(\partial_Y\udv, \alpha\udv)\|_{L^2}^2\right)+\lambda\sqrt{\nu}\|\mathrm{div}_\alpha(\udu, \udv)\|_{L^2}^2+\alpha\|(\sqrt{U_s}\udu, \sqrt{U_s}\udv)\|_{L^2}^2\notag\\
&\quad\leq C\alpha m^{-2}\|\rho\|_{L^2}^2+C\|(\udu, \udv)\|_{L^2}\|(q_u, q_v)\|_{L^2}+\frac{C}{\alpha}(\|(q_\rho, q_u, q_v)\|_{L^2}^2+\nu\|(\partial_Y q_\rho, \alpha q_\rho)\|_{L^2}^2)+\frac{C\nu}{\alpha}\|\udu\|_{L^2}^2.\nonumber\\
&\quad\le C\|(\udu, \udv)\|_{L^2}\|(q_u, q_v)\|_{L^2}+\frac{C}{\alpha}(\|(q_\rho, q_u, q_v)\|_{L^2}^2+\nu\|(\partial_Y q_\rho, \alpha q_\rho)\|_{L^2}^2)+\frac{C\nu}{\alpha}\|\udu\|_{L^2}^2.
\end{align}
Let $E_1=\|(q_\rho,q_u,q_v)\|_{L^2}+\sqrt{\nu}\|(\partial_Y q_\rho, \alpha q_\rho)\|_{L^2}$. Then using interpolation inequality (cf. \cite[Proposition 2.4]{GM19})
$$\|f\|_{L^2}^2\leq C\|\sqrt{U}_sf\|_{L^2}^{\frac43}\|\partial_Yf\|_{L^2}^{\frac23}+C\|\sqrt{U}_sf\|_{L^2}^2,
$$
and \eqref{5.18}, one has
\begin{align*}
\|(\udu, \udv)\|_{L^2}^2&\le C\|(\sqrt{U_s}\udu, \sqrt{U_s}\udv)\|_{L^2}^{\frac43}\|(\partial_Y\udu, \partial_Y\udv)\|_{L^2}^{\frac23}+C\|(\sqrt{U_s}\udu, \sqrt{U_s}\udv)\|_{L^2}^2,\\
&\leq C(\frac{\nu}{\alpha^2}+\frac{\nu^{\frac56}}{\alpha^{\frac53}})\|(\udu,\udv)\|_{L^2}^2+C(\alpha^{-2}+\alpha^{-\frac{5}{3}}\nu^{-\frac16})E_1^2+C\alpha^{-\frac43}\nu^{\frac16}\|(\udu,\udv)\|_{L^2}^{\frac13}E_1^{\frac53}\\
&\quad+C(\alpha^{-1}\nu^{-\frac16}+C\alpha^{-\frac53}\nu^{\frac12})\|(\udu,\udv)\|_{L^2}^{\frac23}E_1^{\frac43}+C\alpha^{-\frac23}\nu^{-\frac16}\|(\udu,\udv)\|_{L^2}E_1\\
&\quad+C(\alpha^{-1}\nu^{\frac16}+C\alpha^{-\frac53}\nu^{\frac12})\|(\udu,\udv)\|_{L^2}^{\frac43} E_1^{\frac23}+C\alpha^{-\frac43}\nu^{\frac12}\|(\udu,\udv)\|_{L^2}^{\frac53}E_{1}^{\frac13}\\
&\leq (\eta+\frac{C\nu}{\alpha^2}+\frac{C\nu^{\frac56}}{\alpha^{\frac53}})\|(\udu,\udv)\|_{L^2}^2+ \frac{C_\eta\hat{n}^{\frac23}}{\alpha^2}E_1^2,
\end{align*}
where $\eta\in (0,1)$ can be arbitrarily small. Note that $\frac{\nu^{\frac56}}{\alpha^{\frac53}}=\frac{1}{\hat{n}^{\frac53}}\leq L^{\frac53}$ and $\frac{\nu}{\alpha^2}=\frac{1}{\hat{n}^2}\leq L^2$. Taking $L$ and $\eta$ sufficiently small so that the first term on the right hand side can be absorbed by the left hand side, we obtain
\begin{align}\label{5.19}
\|(\udu, \udv)\|_{L^2}\le \frac{C\hat{n}^{\frac13}}{\alpha}(\|(q_\rho,q_u,q_v)\|_{L^2}+\sqrt{\nu}\|(\partial_Y q_\rho, \alpha q_\rho)\|_{L^2}),
\end{align}
which is \eqref{5.2}.

\underline{\it Step 4. Higher order estimates.} Substituting \eqref{5.19} into \eqref{5.17} and \eqref{5.18} respectively yields
\begin{align}
m^{-2}\|(\partial_Y\udr, \alpha\udr)\|_{L^2}&\le C\hat{n}^{\frac16}(\|(q_\rho,q_u,q_v)\|_{L^2}+\sqrt{\nu}\|(\partial_Y q_\rho, \alpha q_\rho)\|_{L^2}),\nonumber\\
\|(\partial_Y\udu, \alpha\udu)\|_{L^2}+\|(\partial_Y\udv, \alpha\udv)\|_{L^2}&\leq \frac{C\hat{n}^{\frac16}}{\alpha^{\frac12}\nu^{\frac14}}(\|(q_\rho,q_u,q_v)\|_{L^2}+\sqrt{\nu}\|(\partial_Y q_\rho, \alpha q_\rho)\|_{L^2})\nonumber\\
&\leq \frac{C{\hat{n}}^{\frac23}}{\alpha}(\|(q_\rho,q_u,q_v)\|_{L^2}+\sqrt{\nu}\|(\partial_Y q_\rho, \alpha q_\rho)\|_{L^2}),\nonumber\\
\|\sqrt{U_s}(\udu,\udv)\|_{L^2}&\leq \frac{C\hat{n}^{\frac16}}{\alpha}(\|(q_\rho,q_u,q_v)\|_{L^2}+\sqrt{\nu}\|(\partial_Y q_\rho, \alpha q_\rho)\|_{L^2}).\nonumber
\end{align}
Moreover, by the continuity equation $\eqref{5.1}_1$, we have
\begin{align}
\|{\rm div}_{\alpha}(\udu,\udv)\|_{L^2}&\le \|q_\rho\|_{L^2}+\|\alpha U_s\udr\|_{L^2}\nonumber\\
&\le C\hat{n}^{\frac16}(\|(q_\rho,q_u,q_v)\|_{L^2}+\sqrt{\nu}\|(\partial_Y q_\rho, \alpha q_\rho)\|_{L^2}),\nonumber\\
\|\partial_Y{\rm div}_{\alpha}(\udu,\udv)\|_{L^2}&\le \|\partial_Yq_\rho\|_{L^2}+\|\alpha \partial_YU_s\udr\|_{L^2}+\|\alpha U_s\partial_Y\udr\|_{L^2}\notag\\
&\le C(1+\alpha)\hat{n}^{\frac16}(\|(q_\rho,q_u,q_v)\|_{L^2}+\sqrt{\nu}\|(\partial_Y q_\rho, \alpha q_\rho)\|_{L^2}).\nonumber
\end{align}
The $H^2$ estimate follows directly from $\eqref{5.1}_2$ and $\eqref{5.1}_3$
\begin{align}
\|(\Delta_\alpha\udu, \Delta_\alpha\udv)\|_{L^2}&\le \frac{C}{\sqrt{\nu}}\left(\alpha\|(U_s\udu, U_s\udv)\|_{L^2}+m^{-2}\|(\partial_Y\udr, \alpha\udr)\|_{L^2}+\|(q_u, q_v)\|_{L^2}\right)+C\|(\partial_Y\mathcal{D}, \alpha\mathcal{D})\|_{L^2}\notag\\
&\le C(\alpha+\nu^{-\frac12})\hat{n}^{\frac16}(\|(q_\rho,q_u,q_v)\|_{L^2}+\sqrt{\nu}\|(\partial_Y q_\rho, \alpha q_\rho)\|_{L^2}),\nonumber
\end{align}
which is \eqref{5.6}. Finally, for $\alpha\in (0,1)$, we have 
\begin{align*}
	\|\partial_Y\udv\|_{L^2}\le \alpha\|\udu\|_{L^2}+\|\mathrm{div}_{\alpha}(\udu, \udv)\|_{L^2}\le C\hat{n}^{\frac13}(\|(q_\rho,q_u,q_v)\|_{L^2}+\sqrt{\nu}\|(\partial_Y q_\rho, \alpha q_\rho)\|_{L^2}),
\end{align*}
where the bounds \eqref{5.2} and \eqref{5.6-1} have been used in the last inequality.  Therefore, \eqref{5.6-1} follows, and the proof of Proposition \ref{P5.1} is complete.
\end{proof}
%\begin{remark}\label{R5.1}
%For $0<\alpha\le 1$, the estimate of $\partial_Y\vartheta$ can be achieved from the divergence part
%which is essential for the convergence of the iteration.
%\end{remark}

\subsection{Quasi-compressible-Stokes iteration}\label{S6.2}
Now we introduce the quasi-compressible-Stokes iteration to solve \eqref{9.1} with the boundary condition \eqref{9.1-1}. Recall that $\mathcal{L}$ is the linearized Navier-Stokes operator associated to \eqref{9.1}, $\mathcal{L}_Q$ and $\mathcal{L}_S$  are the quasi-compressible and Stokes operators, which are defined in \eqref{4.3} and
\eqref{5.1} respectively. 

We start from a Stokes solution $(\udr^0, \udu^0, \udv^0)$ satisfying
\begin{align}\label{5.25}
\mathcal{L}_S(\udr^0, \udu^0, \udv^0)=(f_{\rho, n}, f_{u, n}, f_{v, n}).
\end{align}
A straightforward computation yields the following error term introduced by $(\udr^0, \udu^0, \udv^0)$:
\begin{align*}
 \mathcal{L}(\udr^0, \udu^0, \udv^0)-\mathcal{L}_S(\udr^0, \udu^0, \udv^0)=(0, \udv^0\partial_YU_s, 0).
\end{align*}
Note that $\udv^0\partial_YU_s\in H^2(\mathbb{R}_+)$. Then we can eliminate this error term by using a quasi-compressible solution $(\varrho^0, \mathfrak{u}^0, \mathfrak{v}^0)$:
\begin{align}\label{5.26}
\mathcal{L}_Q(\varrho^0, \mathfrak{u}^0, \mathfrak{v}^0)=(0, -\udv^0\partial_YU_s, 0).
\end{align}
Set $\mathbf{s}_0=(\udr^0, \udu^0, \udv^0)+(\varrho^0, \mathfrak{u}^0, \mathfrak{v}^0)$ as the approximate solution at zeroth step. It is straightforward to see that
\begin{align*}
\mathcal{L}(\mathbf{s}_0)=(f_{\rho, n}, f_{u, n}, f_{v, n})+\mathbf{e}_0,
\end{align*}
where the  error term
\begin{align*}
	\mathbf{e}_0&\eqdef \mathcal{L}(\varrho^0, \mathfrak{u}^0, \mathfrak{v}^0)-\mathcal{L}_Q(\varrho^0, \mathfrak{u}^0, \mathfrak{v}^0)\\
	&=(0, \sqrt{\nu}\Delta_\alpha(U_s\varrho^0)+\sqrt{\nu}U_s^{\prime\prime}\varrho^0-\lambda\sqrt{\nu}i\alpha\mathrm{div}(\mathfrak{u}^0, \mathfrak{v}^0), 
	-\lambda\sqrt{\nu}\partial_Y\mathrm{div}(\mathfrak{u}^0, \mathfrak{v}^0)).
\end{align*}
%To eliminate it, we can solve the Stokes system.
%\begin{align}\label{5.27}
%\mathcal{L}_S(\zeta^1, \theta^1, \vartheta^1)=-\vec{e}_1,
%\end{align}
%which leads to another error
%\begin{align*}
%\vec{e}_2\triangleq \mathcal{L}(\zeta^1, \theta^1, \vartheta^1)-L_S(\zeta^1, \theta^1, \vartheta^1)=(0, \vartheta^1\partial_YU_s, 0).
%\end{align*}
%Similarly, let $(\varrho^1, \mathfrak{u}^1, \mathfrak{v}^1)$ solves the following Quasi-compressible system
%\begin{align}\label{5.28}
%\mathcal{L}_Q(\varrho^1, \mathfrak{u}^1, \mathfrak{v}^1)=-\vec{e}_2,
%\end{align}
%and results in the error
%\begin{align*}
%\vec{e}_3\triangleq (0, \sqrt{\nu}\Delta_\alpha(U_s\varrho^1)+\sqrt{\nu}U_s^{\prime\prime}\varrho^1-\lambda\sqrt{\nu}i\alpha\mathrm{div}(\mathfrak{u}^1, \mathfrak{v}^1), -\lambda\sqrt{\nu}\partial_Y\mathrm{div}(\mathfrak{u}^1, \mathfrak{v}^1)).
%\end{align*}
%Set $\vec{s}_1=(\zeta^1, \theta^1, \vartheta^1)+(\varrho^1, \mathfrak{u}^1, \mathfrak{v}^1)$, and we have
%\begin{align}\label{5.29}
%\mathcal{L}(\vec{s}_0+\vec{s}_1)=(f_{0, n}, f_{1, n}, f_{2, n})+\vec{e}_3.
%\end{align}

Now we can iterate above process. For any integer $k\ge 1$, we define the $k$th-order corrector
\begin{align}
	\mathbf{s}_k=(\udr^k,\udu^k,\udv^k)+(\varrho^k,\mathfrak{u}^k,\mathfrak{v}^k), \nonumber
\end{align}
where $(\udr^k, \udu^k, \udv^k)$ solves the Stokes system
\begin{align}\label{5.30}
\mathcal{L}_S(\udr^k, \udu^k, \udv^k)=-\mathbf{e}_{k-1},
\end{align}
and $(\varrho^k, \mathfrak{u}^k, \mathfrak{v}^k)$ solves the quasi-compressible system
\begin{align}\label{5.31}
	\mathcal{L}_Q(\varrho^k, \mathfrak{u}^k, \mathfrak{v}^k)=(0,\udv^k\partial_YU_s,0).
\end{align}
Here the error function in \eqref{5.30} is given by
\begin{align*}
\mathbf{e}_{k-1}\triangleq (0, \sqrt{\nu}\Delta_\alpha(U_s\varrho^{k-1})+\sqrt{\nu}U_s^{\prime\prime}\varrho^{k-1}-\lambda\sqrt{\nu}i\alpha\mathrm{div}(\mathfrak{u}^{k-1}, \mathfrak{v}^{k-1}), -\lambda\sqrt{\nu}\partial_Y\mathrm{div}(\mathfrak{u}^{k-1}, \mathfrak{v}^{k-1})).
\end{align*}

%and yields an error
%\begin{align*}
%\vec{e}_{2k+1}&\triangleq\mathcal{L}(\varrho^k, \mathfrak{u}^k, \mathfrak{v}^k)-L_Q(\varrho^k, \mathfrak{u}^k, \mathfrak{v}^k)\\
%&=(0, \sqrt{\nu}\Delta_\alpha(U_s\varrho^{k})+\sqrt{\nu}U_s^{\prime\prime}\varrho^k-\lambda\sqrt{\nu}i\alpha\mathrm{div}(\mathfrak{u}^{k}, \mathfrak{v}^{k}), -\lambda\sqrt{\nu}\partial_Y\mathrm{div}(\mathfrak{u}^{k}, \mathfrak{v}^{k})).
%\end{align*}
For any $N\geq 1$, it is straightforward to check that
\begin{align*}
\mathcal{L}\left(\sum\limits_{k=0}^{N}\mathbf{s}_k\right)=(f_{0, n}, f_{1, n}, f_{2, n})+\mathbf{e}_{N},
\end{align*}
where the error term in $N$-th step is
$$\mathbf{e}_{N}\triangleq (0, \sqrt{\nu}\Delta_\alpha(U_s\varrho^{N})+\sqrt{\nu}U_s^{\prime\prime}\varrho^{N}-\lambda\sqrt{\nu}i\alpha\mathrm{div}(\mathfrak{u}^{N}, \mathfrak{v}^{N}), -\lambda\sqrt{\nu}\partial_Y\mathrm{div}(\mathfrak{u}^{N}, \mathfrak{v}^{N})).
$$
Therefore the series $\mathbf{s}=\sum_{k=0}^{\infty}\vec{s}_j$ formally defines a solution to \eqref{9.1}. Next proposition shows the convegence of this series. 
\begin{proposition}
\label{P5.2}
{\bf{(Convergence)}}
Let $m\in(0, 1)$.  There exists a consistant $L_0>0$, such that if $L\in (0,L_0)$, then for each $\hat{n}\lesssim \nu^{-\frac34}$, and for any $(f_{\rho, n},f_{u, n}, f_{v, n})\in H^1(\mathbb{R}_+)\times L^2(\mathbb{R}_+)^2$, there exists a solution $\mathbf{s}=(\rho_{n,sl} ,u_{n,sl}, v_{n,sl})\in H^1(\mathbb{R}_+)\times H^2(\mathbb{R}_+)^2$ to the linearized Navier-Stokes system \eqref{9.1} with the boundary condition \eqref{9.1-1}. Moreover, the solution satisfies following properties.
\begin{itemize}
\item[({\rm i})] If $\alpha\ge 1$, 
\begin{align}\label{5.32}
m^{-2}\|\rho_{n,sl}-\udr^0\|_{H^1}+\alpha\|(u_{n,sl}-\udu^0, v_{n,sl}-\udv^0)\|_{L^2}
&\le C\hat{n}^{\frac13}\|\udv^0\|_{L^2},\\
\|\partial_Yv_{n,sl}-\partial_Y\udv^0\|_{L^2}+\|\partial_Y\mathcal{D}_{sl}-\partial_Y\mathcal{D}^0\|_{L^2}+\alpha \|\mathcal{D}_{sl}-\mathcal{D}^0\|_{L^2}&\le C\hat{n}^{\frac13}\|\udv^0\|_{L^2},\label{5.33}\\
\|\partial_Yu_{n,sl}-\partial_Y\udu^0\|_{L^2}+\|\partial_Y^2v_{n,sl}-\partial_Y^2\udv^0\|_{L^2}&\le C\hat{n}^{\frac23}\|\udv^0\|_{L^2},\label{5.34}\\
\|\partial_Y^2u_{n,sl}-\partial_Y^2\udu^0\|_{L^2}&\leq C\hat{n}^{\frac43}\|\udv^0\|_{L^2}.\label{5.37-1}
\end{align}
\item[({\rm ii})] If $\alpha\in (0,1)$, 
\begin{align}\label{5.35}
m^{-2}\|(\rho_{n,sl}-\udr^0)\|_{H^1}+\|(u_{n,sl}-\udu^0, v_{n,sl}-\udv^0)\|_{L^2}
&\le \frac{C\hat{n}^{\frac13}}{\alpha}\|\partial_Y\udv^0\|_{L^2},\\
\|\partial_Yv_{n,sl}-\partial_Y\udv^0\|_{L^2}+\|\partial_Y\mathcal{D}_{sl}-\partial_Y\mathcal{D}^0\|_{L^2}+\alpha \|\mathcal{D}_{sl}-\mathcal{D}^0\|_{L^2}&\leq C\hat{n}^{\frac13}\|\partial_Y\udv^0\|_{L^2},\label{5.36-1}\\
\|\partial_Yu_{n,sl}-\partial_Y\udu^0\|_{L^2}+\|\partial_Y^2v_{n,sl}-\partial_Y^2\udv^0)\|_{L^2}&\le \frac{C\hat{n}^{\frac23}}{\alpha}\|\partial_Y\udv^0\|_{L^2},\label{5.36}\\
\|\partial_Y^2u_{n,sl}-\partial_Y^2\udu^0\|_{L^2} &\le \frac{C\hat{n}^{\frac43}}{\alpha}\|\partial_Y\udv^0\|_{L^2}.\label{5.37}
\end{align}
\end{itemize}
where  $\mathcal{D}^0\eqdef i\alpha \udu^0+\partial_Y\udv^0$.
\end{proposition}
\begin{remark}
	To prove the convergence of iteration, we only need the smallness of torus length $L$.
\end{remark}	
\begin{proof}
\underline{\it Case (i): $\alpha\ge 1$.} For simplicity, we denote $\mathcal{D}^k\eqdef \div_{\alpha}(\mathfrak{u}^k,\mathfrak{v}^k)= i\alpha\mathfrak{u}^k+\partial_Y\mathfrak{v}^k$. Since $(\udr^{k+1},\udu^{k+1},\udv^{k+1})$ solves \eqref{5.30}, we apply the bound \eqref{5.2} to $(\udu^{k+1},\udv^{k+1})$ to obtain
\begin{align}\label{5.38}
\|(\udu^{k+1}, \udv^{k+1})\|_{L^2}&\leq \frac{C\hat{n}^{\frac13}}{\alpha}\|\mathbf{e}_k\|_{L^2}\leq  \frac{C\hat{n}^{\frac13}\nu^{\frac12}}{\alpha}\left(\|\Delta_\alpha(U_s\varrho^{k})\|_{L^2}+\|\varrho^k\|_{L^2}+\|(\partial_Y\mathcal{D}_k,i\alpha\mathcal{D}_k)\|_{L^2}\right).
\end{align}
Since $(\varrho^{k},\mathfrak{u}^k,\mathfrak{v}^k)$ is the solution to \eqref{5.31}, then
 it follows from Corollary \ref{P4.9} that
\begin{align}\label{5.39}
\|\Delta_\alpha(U_s\varrho^{k})\|_{L^2}+\|\varrho^k\|_{L^2}+\|(\partial_Y\mathcal{D}_k,i\alpha \mathcal{D}_k)\|_{L^2}\le C\hat{n}^{\frac13}\|\udv^{k}\partial_YU_s\|_{L^2}.
\end{align}
Plugging \eqref{5.39} back into \eqref{5.38} and using $\alpha=\hat{n}\sqrt{\nu}$, we obtain 
\begin{align}\nonumber
\|(\udu^{k+1}, \udv^{k+1})\|_{L^2}\le \frac{Cn^{\frac23}\nu^{\frac12}}{\alpha}\|\udv^{k}\partial_YU_s\|_{L^2}\le \frac{C}{\hat{n}^{\frac13}}\|\udv^{k}\|_{L^2}\leq CL^{\frac13}\|\udv^{k}\|_{L^2}.
\end{align}
Therefore, for sufficiently small $L$, such that $CL^{\frac13}\le\frac12$, the series $\sum_{k=1}^{\infty}(\udu^{k}, \udv^{k})$  converges absolutely in $H^1(\mathbb{R}_+)$, and satisfies
\begin{align}
\sum\limits_{k=1}^{\infty}\|(\udu^{k}, \udv^{k})\|_{L^2}\le \frac{C}{\hat{n}^{\frac13}}\|\udv^0\|_{L^2}.\label{5.40-1}
\end{align}
Moreover, from \eqref{5.39} and \eqref{5.40-1} we obtain
\begin{align}
	\sum_{k=0}^\infty\|\mathbf{e}_k\|_{L^2}&\leq \nu^{\frac12}\left(\sum_{k=0}^\infty\|\Delta_\alpha(U_s\varrho^{k})\|_{L^2}+\|\varrho^k\|_{L^2}+\|(\partial_Y\mathcal{D}_k,i\alpha\mathcal{D}_k)\|_{L^2}\right)\nonumber\\
&\leq C\hat{n}^{\frac13}\nu^{\frac12}\left(\sum_{k=0}^\infty\|\udv^k\partial_YU_s\|_{L^2}\right)\leq  C\hat{n}^{\frac13}\nu^{\frac12}\|\udv^0\|_{L^2}\nonumber.
\end{align}
%which implies
%\begin{align}\label{5.41}
%\sum\limits_{k=0}^{\infty}\|(\theta^{k}, \vartheta^{k})\|_{L^2}\le C\|\vartheta^0\|_{L^2}.
%\end{align}
Then applying \eqref{5.3}-\eqref{5.6} in Proposition \ref{P5.1} to $(\udr^{k+1},\udu^{k+1},\udv^{k+1})$ yields
\begin{align}\label{5.42}
\sum\limits_{k=1}^{\infty}\left(m^{-2}\|(\partial_Y\udr^k, \alpha\udr^k)\|_{L^2}+\|\mathrm{div}_\alpha(\udu^k, \udv^k)\|_{L^2}\right)&\le C\hat{n}^{\frac16}\left(\sum\limits_{k=0}^{\infty}\|\mathbf{e}^k\|_{L^2}\right)\le C\hat{n}^{\frac12}\nu^{\frac12}\|\udv^0\|_{L^2},\\
\sum\limits_{k=1}^{\infty}\|\partial_Y\mathrm{div}_\alpha(\udu^k, \udv^k)\|_{L^2}&\leq C\alpha\hat{n}^{\frac16}\left(\sum\limits_{k=0}^{\infty}\|\mathbf{e}^k\|_{L^2}\right)\leq C\hat{n}^{\frac32}\nu\|\udv^0\|_{L^2},\label{5.44}\\
\sum\limits_{k=1}^{\infty}\|(\partial_Y\udu^k, \partial_Y\udv^k)\|_{L^2}&\le \frac{C\hat{n}^{\frac23}}{\alpha}\left(\sum\limits_{k=0}^{\infty}\|\mathbf{e}^k\|_{L^2}\right)\leq \|\udv^0\|_{L^2},\label{5.43}\\
\sum\limits_{k=1}^{\infty}\|(\Delta_\alpha\udu^k, \Delta_\alpha\udv^k)\|_{L^2}&\le C(\alpha+\nu^{-\frac12})\hat{n}^{\frac16}\left(\sum\limits_{k=0}^{\infty}\|\mathbf{e}^k\|_{L^2}\right)\leq C\hat{n}^{\frac12}\|\udv^0\|_{L^2}.\label{5.45}
\end{align}
As for the $H^2$-convergence of $\sum_{k=1}^{\infty}(\varrho^{k}, \mathfrak{u}^{k}, \mathfrak{v}^{k})$, we use Corollary \ref{P4.9} to obtain
\begin{align}
\sum_{k=1}^{\infty}\alpha\|(\mathfrak{u}^k, \mathfrak{v}^k)\|_{L^2}+ \|\partial_Y\mathfrak{v}^k\|_{L^2}&\leq  C\hat{n}^{\frac13}\left(\sum_{k=0}^{\infty}\|\udv^{k}\|_{L^2}\right)\leq C\hat{n}^{\frac13}\|\udv^0\|_{L^2},\label{5.46}\\
\sum_{k=1}^{\infty}m^{-2}\|\varrho^k\|_{H^1}&\leq C\hat{n}^{\frac13}\left(\sum_{k=0}^{\infty}\|\udv^{k}\|_{L^2}\right)\leq C\hat{n}^{\frac13}\|\udv^0\|_{L^2}\label{5.47}\\
\sum_{k=1}^{\infty}\|\partial_Y\frak{u}^k\|_{L^2}+\|\partial_Y^2\mathfrak{v}^k\|_{L^2}&\leq C\hat{n}^{\frac23}\left(\sum_{k=0}^{\infty}\|\udv^k\|_{L^2}\right)\le C\hat{n}^{\frac23}\|\udv^0\|_{L^2},\label{5.48}\\
\sum_{k=1}^{\infty}\|\partial_Y^2\mathfrak{u}^k\|_{L^2}&\leq C\hat{n}^{\frac43}\left(\sum_{k=0}^{\infty}\|\udv^k\|_{L^2}\right)\le C\hat{n}^{\frac43}\|\udv^0\|_{L^2},\label{5.48-2}\\
\sum_{k=1}^{\infty}\|(\partial_Y\mathcal{D}_k,\alpha\mathcal{D}_k)\|_{L^2}&\leq C\hat{n}^{\frac13}\left(\sum_{k=0}^{\infty}\|\udv^k\|_{L^2}\right)\le C\hat{n}^{\frac13}\|\udv^0\|_{L^2}.\label{5.48-3}
\end{align}
Combining \eqref{5.40-1} with \eqref{5.42}-\eqref{5.48-3}, we obtain \eqref{5.32}-\eqref{5.37-1}. 

\underline{Case (ii): $\alpha\in (0,1)$}. Applying the bound \eqref{5.6-1} to $\udv^{k+1}$, we get
\begin{align}\label{5.49}
\|\partial_Y\udv^{k+1}\|_{L^2}\le C\hat{n}^{\frac13}\|\mathbf{e}_k\|_{L^2}.
\end{align}
Similar to \eqref{5.38} and \eqref{5.39}, we use the bounds \eqref{4.127}-\eqref{4.131-1} in Corollary \ref{P4.9} to deduce that
\begin{align}
\|\mathbf{e}_k\|_{L^2}&\leq \sqrt{\nu}\left(\|\Delta_\alpha(U_s\varrho^{k})\|_{L^2}+\|\varrho^k\|_{L^2}+\|(\partial_Y\mathcal{D}_k, \alpha\mathcal{D}_k)\|_{L^2}\right)\notag\\
&\le \frac{C\hat{n}^{\frac13} \nu^{\frac12}}{\alpha}\|\udv^{k}\partial_YU_s\|_{L^2}\leq \frac{C\hat{n}^{\frac13}\nu^{\frac12}}{\alpha}\|\partial_Y\udv^{k}\|_{L^2}.\label{5.50}
\end{align}
Here we have used the Hardy's inequality: $\|\udv^{k}\partial_YU_s\|_{L^2}\leq \|Y\partial_YU_s\|_{L^\infty}\|Y^{-1}\udv^k\|_{L^2}\leq C\|\partial_Y\udv^k\|_{L^2}$.
Substituting \eqref{5.50} into \eqref{5.49}  and using $\alpha=\hat{n}\sqrt{\nu}$, we obtain
\begin{align}
\|\partial_Y\udv^{k+1}\|_{L^2}\leq \frac{C\hat{n}^{\frac23}{\nu}^{\frac12}}{\alpha}\|\partial_Y\udv^{k}\|_{L^2}\le \frac{C}{\hat{n}^{\frac13}}\|\partial_Y\udv^{k}\|_{L^2}\leq CL^{\frac13}\|\partial_Y\udv^{k}\|_{L^2}.\nonumber
\end{align}
Therefore, for sufficiently small $L$, the series $\sum_{k=0}^{\infty}\partial_Y\udv^{k}$  converges 
absolutely in $L^2(\mathbb{R}_+)$, and 
\begin{align}
	\sum_{k=0}^\infty \|\partial_Y\udv^k\|_{L^2}\leq C\|\partial_Y\udv^0\|_{L^2}.\nonumber
\end{align}
Since $(\varrho^k,\mathfrak{u}^k,\mathfrak{v}^k)$ solves \eqref{5.31}, then applying the bounds in \eqref{4.127}-\eqref{4.131-1} in Corollary \ref{P4.9} to $(\varrho^k,\mathfrak{u}^k,\mathfrak{v}^k)$ again, we get
\begin{align}
	\sum_{k=1}^{\infty}\|(\mathfrak{u}^k, \mathfrak{v}^k)\|_{L^2}+ m^{-2}\|\varrho^k\|_{H^1}&\leq  \frac{C\hat{n}^{\frac13}}{\alpha}\left(\sum_{k=0}^{\infty}\|\udv^{k}\partial_YU_s\|_{L^2}\right)\nonumber\\
	&\leq \frac{C\hat{n}^{\frac13}}{\alpha}\left(\sum_{k=0}^{\infty}\|\partial_Y\udv^{k}\|_{L^2}\right)\leq \frac{C\hat{n}^{\frac13}}{\alpha}\|\partial_Y\udv^0\|_{L^2},\nonumber\\
	\sum_{k=1}^{\infty}\|\partial_Y\mathfrak{v}^k\|_{L^2}+\|(\partial_Y\mathcal{D}_k,\alpha\mathcal{D}_k)\|_{L^2}&\leq C\hat{n}^{\frac13}\left(\sum_{k=0}^{\infty}\|\udv^{k}\partial_YU_s\|_{L^2}\right)\nonumber\\
	&\leq C\hat{n}^{\frac13}\left(\sum_{k=0}^{\infty}\|\partial_Y\udv^{k}\|_{L^2}\right)\leq  C\hat{n}^{\frac13}\|\partial_Y\udv^0\|_{L^2},\nonumber\\
	\sum_{k=1}^{\infty}\alpha\|\partial_Y\frak{u}^k\|_{L^2}+\|\partial_Y^2\mathfrak{v}^k\|_{L^2}&\leq C\hat{n}^{\frac23}\left(\sum_{k=0}^{\infty}\|\udv^{k}\partial_YU_s\|_{L^2}\right)\nonumber\\
	&\leq C\hat{n}^{\frac23}\left(\sum_{k=0}^{\infty}\|\partial_Y\udv^{k}\|_{L^2}\right)\leq  C\hat{n}^{\frac23}\|\partial_Y\udv^0\|_{L^2},\nonumber\\
	\sum_{k=1}^{\infty}\|\partial_Y^2\mathfrak{u}^k\|_{L^2}&\leq \frac{C\hat{n}^{\frac43}}{\alpha}\left(\sum_{k=0}^{\infty}\|\udv^{k}\partial_YU_s\|_{L^2}\right)\nonumber\\
	&\leq \frac{C\hat{n}^{\frac43}}{\alpha}\left(\sum_{k=0}^{\infty}\|\partial_Y\udv^{k}\|_{L^2}\right)\leq  \frac{C\hat{n}^{\frac43}}{\alpha}\|\partial_Y\udv^0\|_{L^2}.\nonumber
\end{align}
 The remaining estimates on $(\udr^k,\udu^k,\udv^k)$ is similar to Case (i). We omit the details for brevity. The proof of Proposition \ref{P5.2} is complete.
\end{proof}

We are now ready to prove Proposition \ref{C5.3}.
\bigbreak
{\bf Proof of Proposition \ref{C5.3}:} We only prove \eqref{5.66}.  Note that $(\udr^0,\udu^0,\udv^0)$ solves \eqref{5.25}. Then by using \eqref{5.6}, \eqref{5.6-1}, \eqref{5.36}, and \eqref{5.37}, we obtain
\begin{align}
	\|(\partial_Y^2u_{n,sl},\partial_Y^2v_{n,sl})\|_{L^2}&\leq \|(\partial_Y^2\udu^0,\partial_Y^2\udv^0)\|_{L^2}+\|(\partial_Y^2u_{n,sl}-\partial_Y^2\udu^0,\partial_Y^2v_{n,sl}-\partial_Y^2\udv^0)\|_{L^2}\nonumber\\
	&\leq  \frac{C\hat{n}^{\frac16}}{\sqrt{\nu}}(\|(f_{\rho,n},f_{u,n},f_{v,n})\|_{L^2}+\sqrt{\nu}\|\partial_Y f_{\rho,n}\|_{L^2})+\frac{C\hat{n}^{\frac43}}{\alpha}\|\partial_Y\udv^0\|_{L^2}\nonumber\\
	&\leq C\left(\frac{\hat{n}^{\frac16}}{\sqrt{\nu}}+\frac{\hat{n}^{\frac53}}{\alpha}\right)(\|(f_{\rho,n},f_{u,n},f_{v,n})\|_{L^2}+\sqrt{\nu}\|\partial_Y f_{\rho,n}\|_{L^2})\nonumber\\
	&\leq \frac{C\hat{n}^{\frac23}}{\sqrt{\nu}}(\|(f_{\rho,n},f_{u,n},f_{v,n})\|_{L^2}+\sqrt{\nu}\|\partial_Y f_{\rho,n}\|_{L^2}),\nonumber
\end{align}
which is \eqref{5.66}. The other inequalities can be shown similarly. Boundary values \eqref{5.62-1} and \eqref{5.62-2} follow from the Sobolev inequality. The proof of Proposition \ref{C5.3} is complete. \qed

\subsection{Solvability in high frequencies}\label{S6.3}
In this subsection we will solve \eqref{9.1} for $\hat{n}\geq \kappa_2^{-1}\nu^{-\frac34}$, where $\kappa_2\in (0,1)$ will be determined later. In this regime, the stretching term $v_{n,sl}\partial_YU_s$ can be controlled by the tangential diffusion. Thus, \eqref{9.1} can be solved by the energy method.
\bigbreak
{\bf Proof of Proposition \ref{P5.4}:} We supplement \eqref{9.1} with an extra boundary condition $\partial_Yu_{n,sl}|_{Y=0}=0$. Taking the inner product of momentum equations $\eqref{9.1}_2$ and $\eqref{9.1}_3$ with $ \bar{u}_{n,sl}$, and $\bar{v}_{n,sl}$ respectively, then integrating by parts, we obtain
\begin{align}\label{5.80}
&\sqrt{\nu}\left(\|(\partial_Yu_{n,sl}, \alpha u_{n,sl})\|_{L^2}^2+\|(\partial_Yv_{n,sl}, \alpha v_{n,sl})\|_{L^2}^2\right)+\lambda\sqrt{\nu}\|\mathcal{D}_{sl}\|_{L^2}^2+i\alpha \|(\sqrt{U_s}u_{n,sl}, \sqrt{U_s}v_{n,sl})\|_{L^2}^2\notag\\
&\quad=-\langle v_{n,sl}\partial_YU_s, u_{n,sl}\rangle+m^{-2}\langle\rho_{n,sl}, \mathrm{div}_\alpha(u_{n,sl}, v_{n,sl})\rangle+\langle f_{u,n}+\sqrt{\nu}\rho_{n,sl}\partial_Y^2U_s, u_{n,sl}\rangle+\langle f_{v,n}, v_{n,sl}\rangle.
\end{align}
Compared with \eqref{5.7}, the only difference in \eqref{5.80} is the first term on the right hand side, which is estimated as 
\begin{align*}
|\langle v_{n,sl}\partial_YU_s, u_{n,sl}\rangle|\le \sqrt{\nu}\|\alpha u_{n,sl}\|_{L^2}\|\alpha v_{n,sl}\|_{L^2}\frac{\|\partial_YU_s\|_{L^\infty}}{\hat{n}^2\nu^{\frac32}}.
\end{align*}
By taking $\kappa_2\in (0,1)$ small enough, such that $\frac{\|\partial_YU_s\|_{L^\infty}}{\hat{n}^2\nu^{\frac32}}\leq \kappa_2^2\|\partial_YU_s\|_{L^\infty}\leq \frac12$, we have
\begin{align*}
|\langle v_{n,sl}\partial_YU_s, u_{n,sl}\rangle|\le \frac{\sqrt{\nu}}{2}\left(\|(\alpha u_{n,sl}, \partial_Yu_{n,sl})\|_{L^2}^2+\|(\alpha v_{n,sl}, \partial_Yv_{n,sl})\|_{L^2}^2\right),
\end{align*}
which can be absorbed by the left hand side of \eqref{5.80}. The remaining terms can be treated 
by the same argument as the Stokes system. Recall $E_1=\|(f_{\rho,n},f_{u,n},f_{v,n})\|_{L^2}+\sqrt{\nu}\|(\partial_Yf_{\rho,n},\alpha f_{\rho,n})\|_{L^2}$. Then similar to \eqref{5.18}, we can deduce 
\begin{align}
&\sqrt{\nu}\left(\|(\partial_Yu_{n,sl}, \alpha u_{n,sl})\|_{L^2}^2+\|(\partial_Yv_{n,sl}, \alpha v_{n,sl})\|_{L^2}^2\right)\nonumber\\
&\qquad\leq C\|(u_{n,sl},v_{n,sl})\|_{L^2}E_1+\frac{C}{\alpha}E_1^2+\frac{C\nu}{\alpha}\|u_{n,sl}\|_{L^2}^2\nonumber\\
&\qquad\leq (\eta+\frac{C\sqrt{\nu}}{\alpha^2})\sqrt{\nu}\|\alpha(u_{n,sl},v_{n,sl})\|_{L^2}^2+C_\eta(\alpha^{-2}\nu^{-\frac12}+\alpha^{-1})E_{1}^2.
\nonumber
\end{align}
By taking $\eta\in (0,1)$ sufficiently small, we arrive at
\begin{align}\label{5.85}
&\|(\partial_Yu_{n,sl}, \alpha u_{n,sl})\|_{L^2}+\|(\partial_Yv_{n,sl}, \alpha v_{n,sl})\|_{L^2}\leq  \frac{C}{\alpha\sqrt{\nu}}E_1,
\end{align}
which is \eqref{5.75}. The estimate \eqref{5.78} follows from the Sobolev inequality. Moreover, similar to \eqref{5.17}, we can deduce that
\begin{align}
m^{-2}\|(\partial_Y\rho_{n,sl}, \alpha\rho_{n,sl})\|_{L^2}&\leq C\alpha^{\frac12}\|(u_{n,sl},v_{n,sl})\|_{L^2}^{\frac12}E_1^{\frac12}+CE_1+C\sqrt{\nu}\|u_{n,sl}\|_{L^2}\nonumber\\ 
&\leq \frac{C}{\nu^{\frac14}}E_1,\nonumber
\end{align}
where we have used \eqref{5.85} in the last inequality.  Thus the estimate on density field \eqref{5.76}   is obtained. 

Finally, $H^2$-estimate of $(u_{n,sl},v_{n,sl})$ can be obtained by using the equation $\eqref{9.1}_{2}$ and $\eqref{9.1}_3$ as follows.
\begin{align}
&\|(\partial_Y^2u_{n,sl},\alpha^2u_{n,sl})\|_{L^2}+\|(\partial_Y^2v_{n,sl},\alpha^2v_{n,sl})\|_{L^2}\nonumber\\
&\qquad\leq\frac{C\alpha}{\sqrt{\nu}}\|(u_{n,sl}, v_{n,sl})\|_{L^2}+C\alpha\|(\partial_Yu_{n,sl},\partial_Y v_{n,sl})\|_{L^2}\nonumber\\
&\qquad\qquad+\frac{Cm^{-2}}{\sqrt{\nu}}\|(\partial_Y\rho_{n,sl},\alpha\rho_{n,sl})\|_{L^2}+\frac{C}{\sqrt{\nu}}\|(f_{u,n},f_{v,n})\|_{L^2}\nonumber\\
&\qquad\leq C\left(\frac{1}{\alpha\nu}+\frac{1}{\nu^{\frac34}}\right)E_1\leq \frac{C}{\nu^{\frac34}}E_1,\nonumber
\end{align}
where in the last inequality we have used $\alpha=\hat{n}\sqrt{\nu}\gtrsim \nu^{-\frac14}$ in the regime $\hat{n}\gtrsim \nu^{-\frac34}$. Thus \eqref{5.77} follows and the Proof of Proposition \ref{P5.4} is complete.
\qed

\section{Construction of boundary layer corrections}\label{S7}
The solutions obtained in Section \ref{S6} only satisfy $v_{n,sl}|_{Y=0}=0$. To recover the full no-slip boundary conditions \eqref{1.7-1}, we need to construct boundary correctors satisfying
\begin{align}\label{6.2}
\begin{cases}
\mathcal{L}(\rho_{n,b}, u_{n,b}, v_{n,b})={\bf{0}},\\
v_{n,b}|_{Y=0}=0,\quad u_{n,b}|_{Y=0}\neq 0.
\end{cases}
\end{align}
These correctors will be constructed around homogeneous solutions to the quasi-compressible system obtained in Section \ref{S5}. 
\subsection{Boundary corrector at low frequencies}
Recall Proposition \ref{H2} for the definition of homogeneous quasi-compressible solution $(\varrho_{H,l},\mathfrak{u}_{H,l},\mathfrak{v}_{H,l})$ in the regime $0<\hat{n}\leq \kappa_0\nu^{-\frac34}$. We look for boundary layer corrections in following form:
$$(\rho_{n,b},u_{n,b},v_{n,b})=(\varrho_{H,l},\mathfrak{u}_{H,l},\mathfrak{v}_{H,l})+(\rho_{n,b,r},u_{n,b,r},v_{n,b,r}).$$
Here $(\rho_{n,b,r},u_{n,b,r},v_{n,b,r})$ is the remainder satisfying
$$
\left\{
\begin{aligned}
	&\mathcal{L}(\rho_{n,b,r},u_{n,b,r},v_{n,b,r})=(0,f_{u,n,r},f_{v,n,r}),\\
	&v_{n,b,r}|_{Y=0}=0,
\end{aligned}
\right.
$$
where
$$
\begin{aligned}
	f_{u,n,r}&=-\sqrt{\nu}\Delta_\alpha(U_s\varrho_{H,l})-\sqrt{\nu}U_s''\varrho_{H,l}+\lambda\sqrt{\nu}i\alpha\div_{\alpha}(\mathfrak{u}_{H,l},\mathfrak{v}_{H,l}),\\
	f_{v,n,r}&=\lambda\sqrt{\nu}\partial_Y\div_\alpha(\mathfrak{u}_{H,l},\mathfrak{v}_{H,l}).
\end{aligned}
$$
Using the bounds \eqref{H3}-\eqref{H6} for $\alpha\in (0,1)$, and \eqref{H8}-\eqref{H13} for $\alpha>1$,  we obtain
\begin{equation}\nonumber
\begin{aligned}
	\|(f_{u,n,r},f_{v,n,r})\|_{L^2}&\leq C\sqrt{\nu}\left(\|\varrho_{H,l}\|_{H^1}+\|\Delta_\alpha\varrho_{H,l}\|_{L^2}+\|\left(\partial_Y\div_\alpha(\mathfrak{u}_{H,l},\mathfrak{v}_{H,l}),\alpha\div_\alpha(\mathfrak{u}_{H,l},\mathfrak{v}_{H,l})\right)\|_{L^2}\right)\\
	&\leq\left\{\begin{aligned} &C\sqrt{\nu}\left(\frac1\alpha+\frac1{\eps^{\frac16}}\right),&&{\rm if}~\alpha\in (0,1),\\
 &\frac{C\alpha\sqrt{\nu}}{\eps^{\frac16}},&&{\rm if}~\alpha>1.
	\end{aligned}\right.
\end{aligned}
\end{equation}
Then using bounds \eqref{5.60}-\eqref{5.62-1} in Proposition \ref{C5.3} on $(\rho_{n,b,r},u_{n,b,r},v_{n,b,r})$, we obtain,
for $\alpha>1$, that
\begin{align}
	m^{-2}\|\rho_{n,b,r}\|_{H^1}+\alpha\|(u_{n,b,r},v_{n,b,r})\|_{L^2}+\|\partial_Yv_{n,b,r}\|_{L^2}&\leq \frac{C\eps^{\frac13}}{\sqrt{\nu}}\|(f_{u,n,r},f_{v,n,r})\|_{L^2}\leq C\alpha\eps^{\frac16},\nonumber\\
	\|\partial_Yu_{n,b,r}\|_{L^2}&\leq \frac{C}{\sqrt{\nu}}\|(f_{u,n,r},f_{v,n,r})\|_{L^2}\leq C\alpha\eps^{-\frac16},\nonumber\\
	\|\partial_Y^2u_{n,b,r}\|_{L^2}+\eps^{-\frac12}\|\partial_Y^2v_{n,b,r}\|_{L^2}&\leq \frac{C}{\eps^{\frac23}\sqrt{\nu}}\|(f_{1,n,r},f_{2,n,r})\|_{L^2}\leq C\alpha\eps^{-\frac56},\nonumber
\end{align}
and
$$\left|u_{n,b,r}(0)\right|\leq \frac{C\eps^{\frac16}}{\alpha^{\frac12}\sqrt{\nu}}\\|(f_{u,n,r},f_{v,n,r})\|_{L^2}\leq C\alpha^{\frac12}\leq C\eps^{-\frac16}.
$$
For $\alpha\in (0,1)$, we use \eqref{5.64}-\eqref{5.66} to obtain
\begin{align}
	m^{-2}\|\rho_{n,b,r}\|_{H^1}+\|(u_{n,b,r},v_{n,b,r})\|_{L^2}&\leq \frac{C\hat{n}^{\frac23}}{\alpha}\|(f_{u,n,r},f_{v,n,r})\|_{L^2}\leq C\eps^{\frac13}\left(\frac{1}{\alpha}+\frac{1}{\eps^{\frac16}}\right),\nonumber\\
	\|\partial_Yv_{n,b,r}\|_{L^2}&\leq C\hat{n}^{\frac23}\|(f_{u,n,r},f_{v,n,r})\|_{L^2}\leq C\eps^{\frac13}\left(\frac{1}{\alpha}+\frac{1}{\eps^{\frac16}}\right),\nonumber\\
	\|\partial_Yu_{n,b,r}\|_{L^2}&\leq \frac{C}{\sqrt{\nu}}\|(f_{u,n,r},f_{v,n,r})\|_{L^2}\leq C\left(\frac{1}{\alpha}+\frac{1}{\eps^{\frac16}}\right),\nonumber\\
	\|\partial_Y^2u_{n,b,r}\|_{L^2}+\|\partial_Y^2v_{n,b,r}\|_{L^2}&\leq \frac{C\hat{n}^{\frac23}}{\sqrt{\nu}}\|(f_{u,n,r},f_{v,n,r})\|_{L^2}\leq C\eps^{-\frac23}\left(\frac{1}{\alpha}+\frac{1}{\eps^{\frac16}}\right),\nonumber
\end{align}
and 
$$\left|u_{n,b,r}(0)\right|\leq C\eps^{\frac16}\left(\frac{1}{\alpha}+\frac{1}{\eps^{\frac16}}\right).$$
Combining these bounds on the  remainder $(\rho_{n,b,r},u_{n,b,r},v_{n,b,r})$ and the leading order profile $(\varrho_{H,l},\mathfrak{u}_{H,l},\mathfrak{v}_{H,l})$ constructed in Proposition \ref{H2}, we have the following proposition.

\begin{proposition}\label{B1}
	Let $m\in (0,1)$. Then for any $0<\hat{n}\leq \kappa_0\nu^{-\frac34}$, there exists a homogeneous solution $	(\rho_{n,b},u_{n,b},{v}_{n,b})\in H^1(\mathbb{R}_+)\times H^2(\mathbb{R}_+)^2$ to the linearized Navier-Stokes system \eqref{6.2}. Moreover, the following statements hold
	\begin{itemize}
		\item[{\rm(i)}] If $\alpha\in (0,1)$, it holds that
		\begin{align}
			m^{-2}\|\rho_{n,b}\|_{H^1}+\|u_{n,b}\|_{L^2}+\|v_{n,b}\|_{H^1}&\leq C\left(\frac1\alpha+\frac{1}{\eps^{\frac16}}\right)\label{B1.1},\\
			\|\partial_Yu_{n,b}\|_{L^2}&\leq C\left(\frac{1}{\alpha \eps^{\frac16}}+\frac{1}{\eps^{\frac12}}\right),\label{B1.2}\\ 
			\|(\Delta_\alpha{u}_{n,b},\Delta_\alpha v_{n,b})\|_{L^2}&\leq \frac{C}{\eps}\left(\frac1\alpha+\frac{1}{\eps^{\frac16}}\right)\label{B1.3}.
		\end{align}
		On the boundary, it holds that
		\begin{align}\label{B1.4}
			\left|u_{n,b}(0)\right|\geq C\left(\frac{1}{\alpha}+\frac1{\eps^{\frac13}}\right),~~v_{n,b}(0)=0.
		\end{align}
		\item[{\rm(ii)}] If $\alpha>1$, it holds that
		\begin{align}
			\|(u_{n,b},v_{n,b})\|_{L^2}&\leq C\eps^{-\frac16}\label{B2.1},\\
			m^{-2}\|\rho_{n,b}\|_{H^1}+\|\partial_Y v_{n,b}\|_{L^2}&\leq C\alpha\eps^{-\frac16},\label{B2.2}\\ 
			\|\partial_Yu_{n,b}\|_{L^2}+\alpha^{-1}\|\Delta_\alpha v_{n,b}\|_{L^2}&\leq C\eps^{-\frac12},\label{B2.4}\\
			\|\Delta_\alpha{u}_{n,b}\|_{L^2}&\leq C\alpha \eps^{-\frac76}.
			\label{B2.5}
		\end{align}
		On the boundary, it holds that
		\begin{align}\label{B2.6}
			\left|{u}_{n,b}(0)\right|\geq C\eps^{-\frac13},~~\text{and}~~v_{n,b}(0)=0.
		\end{align}
	\end{itemize}
\end{proposition} 

\subsection{Boundary corrector at middle frequencies} Recall Proposition \ref{P4.19} for the definition of homogeneous quasi-compressible solution $(\varrho_{H,m},\mathfrak{u}_{H,m},\mathfrak{v}_{H,m})$ in the regime $\hat{n}\sim \nu^{-\frac34}$, which corresponds to $\alpha\sim \eps^{-\frac13}$. The boundary corrector to the Navier-Stokes system \eqref{6.2} is in the form of
$$(\rho_{n,b},u_{n,b},v_{n,b})=(\varrho_{H,m},\mathfrak{u}_{H,m},\mathfrak{v}_{H,m})+(\rho_{n,b,r},u_{n,b,r},v_{n,b,r}),$$
where the remainder $(\rho_{n,b,r},u_{n,b,r},v_{n,b,r})$ satisifies
$$
\left\{
\begin{aligned}
	&\mathcal{L}(\rho_{n,b,r},u_{n,b,r},v_{n,b,r})=(0,f_{u,n,r},f_{v,n,r}),\\
	&v_{n,b,r}|_{Y=0}=0.
\end{aligned}
\right.
$$
Here the error terms $f_{u,n,r}$ and $f_{v,n,r}$ are given by
$$
\begin{aligned}
f_{u,n,r}&=-\sqrt{\nu}\Delta_\alpha(U_s\varrho_{H,m})-\sqrt{\nu}U_s''\varrho_{H,m}+\lambda\sqrt{\nu}i\alpha\div_{\alpha}(\mathfrak{u}_{H,m},\mathfrak{v}_{H,m}),\\
f_{v,n,r}&=\lambda\sqrt{\nu}\partial_Y\div_\alpha(\mathfrak{u}_{H,m},\mathfrak{v}_{H,m}).
\end{aligned}
$$
From the bounds \eqref{4.226}-\eqref{4.226-4} in Proposition \ref{P4.19}, we deduce that
\begin{align}
	\|(f_{u,n,r},f_{v,n,r})\|_{L^2}&\leq C\sqrt{\nu}\left(\|\varrho_{H,m}\|_{H^2}+\|\left(\partial_Y\div_\alpha(\mathfrak{u}_{H,m},\mathfrak{v}_{H,m}),\alpha\div_\alpha(\mathfrak{u}_{H,m},\mathfrak{v}_{H,m})\right)\|_{L^2}\right)\nonumber\\
	&\leq C\sqrt{\nu}\eps^{\frac16}.\nonumber
\end{align}
Then applying \eqref{5.60}-\eqref{5.63} in Proposition \ref{C5.3} to $(\rho_{n,b,r},u_{n,b,r},v_{n,b,r})$, we obtain
\begin{align}
	m^{-2}\|\rho_{n,b,r}\|_{H^1}+\eps^{-\frac13}\|(u_{n,b,r},v_{n,b,r})\|_{L^2}+\|\partial_Yv_{n,b,r}\|_{L^2}&\leq \frac{C\eps^{\frac13}}{\sqrt{\nu}}\|(f_{u,n,r},f_{v,n,r})\|_{L^2}\leq C\eps^{\frac12},\nonumber\\
	\|\partial_Yu_{n,b,r}\|_{L^2}&\leq \frac{C}{\sqrt{\nu}}\|(f_{u,n,r},f_{v,n,r})\|_{L^2}\leq C\eps^{\frac16},\nonumber\\
	\|\partial_Y^2u_{n,b,r}\|_{L^2}+\eps^{-\frac12}\|\partial_Y^2v_{n,b,r}\|_{L^2}&\leq \frac{C}{\eps^{\frac23}\sqrt{\nu}}\|(f_{u,n,r},f_{v,n,r})\|_{L^2}\leq C\eps^{-\frac12},\nonumber
\end{align}
and
$$\left|u_{n,b,r}(0)\right|\leq C\|u_{n,b,r}\|_{L^2}^{\frac12}\|\partial_Yu_{n,b,r}\|_{L^2}^{\frac12}\leq C\eps^{\frac{5}{12}+\frac{1}{12}}\leq C\eps^{\frac12}.
$$
Combining these bounds on the remainder $(\rho_{n,b,r},u_{n,b,r},v_{n,b,r})$ with estimates \eqref{4.226}-\eqref{4.229} on the leading order profile $(\varrho_{H,m},\mathfrak{u}_{H,m},\mathfrak{v}_{H,m})$ in Proposition \ref{P4.19}, we deduce the following bounds of boundary layer correctors $(\rho_{n,b},u_{n,b},v_{n,b})$.
\begin{proposition}\label{B2}
	Let $m\in (0,1)$. If $\hat{n}\sim \nu^{-\frac34}$, there exists a homogeneous solution $(\rho_{n,b},u_{n,b},{v}_{n,b})\in H^1(\mathbb{R}_+)\times H^2(\mathbb{R}_+)^2$ to the linearized Navier-Stokes system \eqref{6.2} satisfying
		\begin{align}
		\|(u_{n,b},v_{n,b})\|_{L^2}&\leq C\eps^{\frac12}\label{B3.1},\\
			m^{-2}\|\rho_{n,b}\|_{H^1}+\|(\partial_Yu_{n,b},\partial_Yv_{n,b})\|_{L^2}&\leq C\eps^{\frac16},\label{B3.2}\\ 
			\|\Delta_\alpha{u}_{n,b}\|_{L^2}&\leq C\eps^{-\frac56}\label{B3.3},\\
			\|\Delta_\alpha{v}_{n,b}\|_{L^2}&\leq C\eps^{-\frac16}\label{B3.4}.
		\end{align}
	Moreover, on the boundary
		\begin{align}\label{B3.5}
			\left|u_{n,b}(0)\right|\geq C\eps^{\frac13},~~\text{and}~~v_{n,b}(0)=0.
		\end{align}
\end{proposition} 

\subsection{Boundary corrector at high frequencies}
Recall Proposition \ref{P4.22} for the definition of homogeneous quasi-compressible solution $(\varrho_{H,h},\mathfrak{u}_{H,h},\mathfrak{v}_{H,h})$ in the regime $\hat{n}\geq \kappa_1^{-1}\nu^{-\frac34}$. Let $\hat{\kappa}_0=\min\{\kappa_1,\kappa_2\}$, where $\kappa_2$ is given in Proposition \ref{P5.4}.  Similar to the discussions in  low and middle frequency regimes, we look for the boundary corrector in the following form
$$(\rho_{n,b},u_{n,b},v_{n,b})=(\varrho_{H,h},\mathfrak{u}_{H,h},\mathfrak{v}_{H,h})+(\rho_{n,b,r},u_{n,b,r},v_{n,b,r}).$$
Here the remainder $(\rho_{n,b,r},u_{n,b,r},v_{n,b,r})$ satisfies
$$
\left\{
\begin{aligned}
	&\mathcal{L}(\rho_{n,b,r},u_{n,b,r},v_{n,b,r})=(0,f_{u,n,r},f_{v,n,r}),\\
	&v_{n,b,r}|_{Y=0}=0.
\end{aligned}
\right.
$$
where the error terms $f_{u,n,r}$ and $f_{v,n,r}$ are given by
$$
\begin{aligned}
	f_{u,n,r}&=-\sqrt{\nu}\Delta_\alpha(U_s\varrho_{H,h})-\sqrt{\nu}U_s''\varrho_{H,h}+\lambda\sqrt{\nu}i\alpha\div_{\alpha}(\mathfrak{u}_{H,h},\mathfrak{v}_{H,h}),\\
	f_{v,n,r}&=\lambda\sqrt{\nu}\partial_Y\div_\alpha(\mathfrak{u}_{H,h},\mathfrak{v}_{H,h}).
\end{aligned}
$$
From the bounds \eqref{4.252}-\eqref{4.254-2} in Proposition \ref{P4.22}, we obtain
\begin{align}
	\|(f_{u,n,r},f_{v,n,r})\|_{L^2}&\leq C\sqrt{\nu}\left(\|\varrho_{H,h}\|_{H^1}+\|\Delta_\alpha\varrho_{H,h}\|_{L^2}+\|\left(\partial_Y\div_\alpha(\mathfrak{u}_{H,h},\mathfrak{v}_{H,h}),\alpha\div_\alpha(\mathfrak{u}_{H,h},\mathfrak{v}_{H,h})\right)\|_{L^2}\right)\nonumber\\
	&\leq C\sqrt{\nu}\left( \alpha^{\frac12}\nu^{\frac12}+\alpha^{-\frac12}  \right).\nonumber
\end{align}
Then applying \eqref{5.75}-\eqref{5.78} in Proposition \ref{P5.4} to $(\rho_{n,b,r},u_{n,b,r},v_{n,b,r})$, we deduce that
\begin{align}
\|(\partial_Yu_{n,b,r},\alpha u_{n,b,r})\|_{L^2}+\|(\partial_Yv_{n,b,r},\alpha v_{n,b,r})\|_{L^2}&\leq \frac{C}{\alpha\sqrt{\nu}}\|(f_{u,n,r},f_{v,n,r})\|_{L^2}\leq C(\alpha^{-\frac12}\nu^{\frac12}+\alpha^{-\frac32}),\nonumber\\
	m^{-2}\|(\partial_Y\rho_{n,b,r},\alpha\rho_{n,b,r})\|_{L^2}&\leq \frac{C}{\nu^{\frac14}}\|(f_{u,n,r},f_{v,n,r})\|_{L^2}\leq C\nu^{\frac14}(\alpha^{\frac12}\nu^{\frac12}+\alpha^{-\frac12}),\nonumber\\
	\|(\Delta_\alpha u_{n,b,r},\Delta_\alpha v_{n,b,r})\|_{L^2}&\leq \frac{C}{\nu^{\frac34}}\|(f_{u,n,r},f_{v,n,r})\|_{L^2}\leq C\nu^{-\frac14}(\alpha^{\frac12}\nu^{\frac12}+\alpha^{-\frac12}),\nonumber
\end{align}
and
$$\left|u_{n,b,r}(0)\right|\leq C\|u_{n,b,r}\|_{L^2}^{\frac12}\|\partial_Yu_{n,b,r}\|_{L^2}^{\frac12}\leq C\alpha^{-1}(\nu^{\frac12}+\alpha^{-1})\ll \alpha^{-1}.
$$
Combining these bounds on the remainder $(\rho_{n,b,r},u_{n,b,r},v_{n,b,r})$ with the estimates \eqref{4.252}-\eqref{4.255} on the leading order profile $(\varrho_{H,h},\mathfrak{u}_{H,h},\mathfrak{v}_{H,h})$ in Proposition \ref{P4.22}, we obtain the following properties of $(\rho_{n,b},u_{n,b},v_{n,b})$ at high frequencies.
\begin{proposition}\label{B3}
	Let $m\in (0,1)$. There exists $\hat{\kappa}_0\in (0,1)$, such that if $\hat{n}\geq \hat{\kappa_0}^{-1}\nu^{-\frac34}$,  the linearized Navier-Stokes system \eqref{6.2} admits a homogeneous solution $(\rho_{n,b},u_{n,b},{v}_{n,b})\in H^1(\mathbb{R}_+)\times H^2(\mathbb{R}_+)^2$ satisfying
	\begin{align}
		\|(u_{n,b},v_{n,b})\|_{L^2}&\leq C\alpha^{-\frac32}\label{B4.1},\\
	\|(\partial_Yu_{n,b},\partial_Yv_{n,b})\|_{L^2}&\leq C\alpha^{-\frac12},\label{B4.2}\\ 
		\|(\Delta_\alpha{u}_{n,b},\Delta_\alpha{v}_{n,b})\|_{L^2}&\leq C\alpha^{\frac12}\label{B4.3},\\
		m^{-2}\|(\partial_Y\rho_{n,b},\alpha\rho_{n,b})\|_{L^2}&\leq C(\alpha^{-\frac12}+\alpha^{\frac12}\nu^{\frac12})\label{B4.4}.
	\end{align}
	Moreover, on the boundary
	\begin{align}\label{B4.5}
		\left|u_{n,b}(0)\right|\geq \frac{C}{\alpha},~~\text{and}~~v_{n,b}(0)=0.
	\end{align}
\end{proposition} 

\subsection{Linear stability for non-zero modes}\label{S7.4}
In this subsection, we construct the solution $(\rho_n,u_n,v_n)$ to the linearized  steady Navier-Stokes system \eqref{4.1} with no-slip boundary conditions \eqref{4.1-1}, thereby demonstrating the linear stability of non-zero modes. The solutions are constructed in the following form
\begin{align}\label{7.4-1}
	\left(\rho_n,u_n,v_n\right)=\left( \rho_{n,sl},u_{n,sl},v_{n,sl}   \right)-\frac{u_{n,sl}(0)}{u_{n,b}(0)}(\rho_{n,b},u_{n,b},v_{n,b}).
\end{align}
Here $(\rho_{n,sl},u_{n,sl},v_{n,sl})$ is the  solution with the boundary condition $v_{n,sl}|_{Y=0}=0$ only, which has been obtained in Proposition \ref{C5.3} for low and middle frequencies, and in Proposition \ref{P5.4} for high freqeuncies. The boundary layer corrector $(\rho_{n,b},u_{n,b},v_{n,b})$ has been obtained in Propositions \ref{B1}, \ref{B2} and \ref{B3} for low, middle, and high frequencies respectively. Note that from \eqref{B1.4}, \eqref{B2.6}, \eqref{B3.5} and \eqref{B4.5}, we have $|u_{n,b}(0)|\neq 0$.

It is straightforward to see that $	\left(\rho_n,u_n,v_n\right)$ satisfies the Navier-Stokes equations \eqref{4.1} with no-slip boundary conditions \eqref{4.1-1}. Now we summarize their bounds in the following proposition.

\begin{proposition}\label{B4}
	Let $m\in (0,1)$. For any $f_{\rho,n}\in H^1(\mathbb{R}_+)$ and $(f_{u,n},f_{v,n})\in L^2(\mathbb{R}_+)^2$,  there exists a solution $\left(\rho_n,u_n,v_n\right)\in H^1(\mathbb{R}_+)\times \left(H^2(\mathbb{R}_+)\cap H^1_0(\mathbb{R}_+)\right)^2$ to the linearized  steady Navier-Stokes system \eqref{4.1} with no-slip boundary condition \eqref{4.1-1}. Moreover, the following statements hold.
	\begin{itemize}
		\item[{\rm(i)}] In the low-frequency regime $0<\hat{n}\leq \kappa_0\nu^{-\frac34}$, where $\kappa_0$ is defined in Proposition \ref{B1}, it holds
		\begin{align}
			&\|(m^{-2}\rho_n,u_n,v_n)\|_{H^1}+\alpha\|(m^{-2}\rho_n,u_n,v_n)\|_{L^2}\nonumber\\
			&\qquad\leq \frac{C}{\nu^{\frac12}}(\|(f_{\rho,n},f_{u,n},f_{v,n})\|_{L^2}+\sqrt{\nu}\|(\partial_Y f_{\rho,n}, \alpha f_{\rho,n})\|_{L^2})\label{7.4-2},\\
			&\|(\Delta_\alpha u_n,\Delta_\alpha v_n)\|_{L^2}\nonumber\\
			&\qquad\leq \frac{C}{\nu^{\frac98}}(\|(f_{\rho,n},f_{u,n},f_{v,n})\|_{L^2}+\sqrt{\nu}\|(\partial_Y f_{\rho,n}, \alpha f_{\rho,n})\|_{L^2})\label{7.4-4}.
		\end{align}
	\item[{\rm(ii)}] In the high-frequency regime $\hat{n}\geq \hat{\kappa}_0^{-1}\nu^{-\frac34}$, where $\hat{\kappa}_0$ is defined in Proposition \ref{B3}, it holds
		\begin{align}
		&\|\partial_Y(m^{-2}\rho_n,u_n,v_n)\|_{L^2}+\alpha\|(m^{-2}\rho_n,u_n,v_n)\|_{L^2}\nonumber\\
		&\qquad\leq \frac{C}{\nu^{\frac14}}(\|(f_{\rho,n},f_{u,n},f_{v,n})\|_{L^2}+\sqrt{\nu}\|(\partial_Y f_{\rho,n}, \alpha f_{\rho,n})\|_{L^2})\label{7.4-6},\\
		&\|(\Delta_\alpha u_n,\Delta_\alpha v_n)\|_{L^2}\nonumber\\
		&\qquad\leq \frac{C}{\nu^{\frac34}}(\|(f_{\rho,n},f_{u,n},f_{v,n})\|_{L^2}+\sqrt{\nu}\|(\partial_Y f_{\rho,n}, \alpha f_{\rho,n})\|_{L^2})\label{7.4-7}.
	\end{align}
\item[(iii)] In the middle-frequency regime $\kappa_0\nu^{-\frac34}\leq \hat{n}\leq \hat{\kappa}_0^{-1}\nu^{-\frac34}$, it holds
	\begin{align}
	&\|\partial_Y(m^{-2}\rho_n,u_n,v_n)\|_{L^2}+\alpha\|(m^{-2}\rho_n,u_n,v_n)\|_{L^2}\nonumber\\
	&\qquad\leq \frac{C}{\nu^{\frac12}}(\|(f_{\rho,n},f_{u,n},f_{v,n})\|_{L^2}+\sqrt{\nu}\|(\partial_Y f_{\rho,n}, \alpha f_{\rho,n})\|_{L^2})\label{7.4-8}\\
	&\|(\Delta_\alpha u_n,\Delta_\alpha v_n)\|_{L^2}\nonumber\\
	&\qquad\leq \frac{C}{\nu^{\frac98}}(\|(f_{\rho,n},f_{u,n},f_{v,n})\|_{L^2}+\sqrt{\nu}\|(\partial_Y f_{\rho,n}, \alpha f_{\rho,n})\|_{L^2})\label{7.4-10}.
\end{align}
	\end{itemize}
\end{proposition}
\begin{proof}
	We only prove (ii). Other cases can be proved similarly. Note that $\alpha\gtrsim \nu^{-\frac14}$ in the high-frequency regime. By \eqref{5.78} and \eqref{B4.5}, we have
	$$|\frac{u_{n,sl}(0)}{u_{n,b}(0)}|\leq \frac{C}{\alpha^{\frac12}\nu^{\frac12}}(\|(f_{\rho,n},f_{u,n},f_{v,n})\|_{L^2}+\sqrt{\nu}\|(\partial_Y f_{\rho,n}, \alpha f_{\rho,n})\|_{L^2}).
	$$
	Substituting bounds \eqref{5.75}, \eqref{5.76},  \eqref{B4.1}, \eqref{B4.2} and \eqref{B4.4} into \eqref{7.4-1}, we obtain
	$$\begin{aligned}
&\|\partial_Y(m^{-2}\rho_n,u_n,v_n)\|_{L^2}+\alpha\|(m^{-2}\rho_n,u_n,v_n)\|_{L^2}\\
&\qquad\leq \|\partial_Y(m^{-2}\rho_{n,sl},u_{n,sl},v_{n,sl})\|_{L^2}+|\frac{u_{n,sl}(0)}{u_{n,b}(0)}|\|\partial_Y(m^{-2}\rho_{n,b},u_{n,b},v_{n,b})\|_{L^2}\\
&\qquad\qquad+\alpha\|(m^{-2}\rho_{n,sl},u_{n,sl},v_{n,sl})\|_{L^2}+\alpha|\frac{u_{n,sl}(0)}{u_{n,b}(0)}|\|(m^{-2}\rho_{n,b},u_{n,b},v_{n,b})\|_{L^2}\\
	&\qquad\leq C\left[\frac{1}{\alpha\nu^{\frac12}}+\frac{1}{\nu^{\frac14}}+\frac{1}{\alpha^{\frac12}\nu^{\frac12}}\left(\frac{1}{\alpha^{\frac12}}+\alpha^{\frac12}\nu^{\frac12}\right)\right](\|(f_{\rho,n},f_{u,n},f_{v,n})\|_{L^2}+\sqrt{\nu}\|(\partial_Y f_{\rho,n}, \alpha f_{\rho,n})\|_{L^2})\\
	&\qquad \leq \frac{C}{\nu^{\frac14}}(\|(f_{\rho,n},f_{u,n},f_{v,n})\|_{L^2}+\sqrt{\nu}\|(\partial_Y f_{\rho,n}, \alpha f_{\rho,n})\|_{L^2}),
	\end{aligned}
	$$
	which is \eqref{7.4-6}.  \eqref{7.4-7} can be proved similarly and we omit the details for brevity.
	%Other estimates can be obtained similarly. We omit the detail for brevity. 
	The proof of Proposition \ref{B4} is complete.
	
%Similarly, we have
%		$$\begin{aligned}
%		&\|m^{-2}(\partial_Y\rho_n,\alpha\rho_n)\|+\|(\partial_Yu_n,\alpha u_n)\|_{L^2}+\|(\partial_Yv_n,\alpha v_n)\|_{L^2}\\
%		&\qquad\leq C\left[\frac{1}{\alpha\nu^{\frac12}}+\frac{1}{\nu^{\frac14}}+\frac{1}{\alpha^{\frac12}\nu^{\frac12}}\left(\frac{1}{\alpha^{\frac12}}+\alpha^{\frac12}\nu^{\frac12}\right)\right]\|(f_{0,n},\alpha f_{0,n},\partial_Yf_{0,n},f_{1,n},f_{2,n})\|_{L^2}\\
%		&\qquad\leq \frac{C}{\nu^{\frac14}}\|(f_{0,n},\alpha f_{0,n},\partial_Yf_{0,n},f_{1,n},f_{2,n})\|_{L^2},
%	\end{aligned}
%	$$
%which gives \eqref{7.4-6}, and 
%	$$\begin{aligned}
%	\|(\Delta_\alpha u_n,\Delta_\alpha v_n)\|_{L^2}&\leq C\left(\frac{1}{\nu^{\frac34}}+\frac{}{}\right)
%	&\leq C\|(f_{0,n},\alpha f_{0,n},\partial_Yf_{0,n},f_{1,n},f_{2,n})\|_{L^2},
%\end{aligned}
%$$
\end{proof}
Rescaling the bounds given in \eqref{7.4-2}-\eqref{7.4-10} to the original variables $(x,y)$, using Parserval identity, and noting the relationship
$(f_{\rho,\neq},f_{u,\neq},f_{v,\neq})=\nu^{\frac12}(g_{\rho,\neq},g_{u,\neq}, g_{v,\neq }),
$ (cf. \eqref{gs}),  we have the following linear stability result for non-zero modes.
 
\begin{theorem}\label{thmno}
Let the Mach number $m\in (0,1)$. There exists a positive constant $L_0>0$, such that for any $L\in (0,L_0)$ and $0<\nu\ll 1$, if  
$g_{\rho,\neq}\in H^1(\Omega)$, and $(g_{u,\neq},g_{v,\neq})\in L^2(\Omega)^2$, there exists a unique solution to the linear Navier-Stokes equations \eqref{1.6}. Moreover, the non-zero mode $(\rho_{\neq},u_{\neq },v_{\neq })$ of the solution satisfies the following bounds 
\begin{align}
	\|(m^{-2}\rho_{\neq},u_{\neq},v_{\neq})\|_{L^2(\Omega)}&\leq C\left(\|(g_{\rho,\neq},g_{u,\neq},g_{v,\neq})\|_{L^2(\Omega)}+\nu\|\nabla_{x,y}g_{\rho,\neq}\|_{L^2(\Omega)}\right),\label{7.4-11}\\
	\|(m^{-2}\nabla_{x,y}\rho_{\neq},\nabla_{x,y} u_{\neq},\nabla_{x,y} v_{\neq})\|_{L^2(\Omega)}&\leq \frac{C}{\nu^{\frac12}}\left(\|(g_{\rho,\neq},g_{u,\neq},g_{v,\neq})\|_{L^2(\Omega)}+\nu\|\nabla_{x,y}g_{\rho,\neq}\|_{L^2(\Omega)}\right),\label{7.4-12}\\
	\|(\nabla^2_{x,y} u_{\neq},\nabla^2_{x,y} v_{\neq})\|_{L^2(\Omega)}&\leq \frac{C}{\nu^{\frac{13}8}}\left(\|(g_{\rho,\neq},g_{u,\neq},g_{v,\neq})\|_{L^2(\Omega)}+\nu\|\nabla_{x,y}g_{\rho,\neq}\|_{L^2(\Omega)}\right).\label{7.4-13}
\end{align}
\end{theorem}
\begin{remark}
	Compared with the linear stability of the incompressible Navier-Stokes system obtained in \cite{GM19}, the order of $\nu$ given on the right hand side of \eqref{7.4-11} and \eqref{7.4-12} are the same, and they are optimal. However, for the compressible Navier-Stokes equation, we need to bound higher-order derivatives to close the estimates. For this, the optimal order of $\nu$ in second order derivative estimate should be $\nu^{-\frac32}$ rather than $\nu^{-\frac{13}{8}}$ given in \eqref{7.4-13}. On the other hand, such an optimal order estimate is hard to obtain because of complexity of the system in the compressible setting.
\end{remark}

\section{Nonlinear stability}\label{S8}
\subsection{Modified linear system and low order estimate}
In this section, we solve the nonlinear equations \eqref{1.5}. Inspired by \cite{KN19}, we fix any velocity field $(\tilde{u},\tilde{v})$ and consider the following modified linear system 
\begin{align}\label{7.1}
	\begin{cases}
		U_s\partial_x\rho+\div_{x, y}(u, v)+\partial_x(\tilde{u}\rho)+\partial_y(\tilde{v}\rho)=g_\rho,\\
		U_s\partial_xu+v\partial_yU_s+m^{-2}\partial_x\rho-\nu\Delta_{x, y} u-\lambda{\nu}\partial_x\div_{x, y}(u, v)+{\nu}\rho\partial_y^2U_s=g_u,\\
		U_s\partial_xv+m^{-2}\partial_y\rho-{\nu}\Delta_{x, y} v-\lambda{\nu}\partial_y\div_{x, y}(u, v)=g_v,\\
		u|_{y=0}=v|_{y=0}=0,
	\end{cases}
\end{align}
Here $(g_\rho, g_u, g_v)$ are given inhomogeneous source terms. For later use, we define three sets of norms:
\begin{align}
	[[(g_\rho,g_u,g_v)]]_{1}\eqdef&\nu^{-1}\|\partial_y^{-1}g_{u,0})\|_{L^1_y}+\nu^{-\frac12}\|\mathcal{I}(g_{\rho,0})\|_{L^\infty_y}+\nu^{-\frac12}\|\partial_y^{-1}g_{u,0}\|_{L^2_y}+\|\mathcal{I}(g_{\rho,0})\|_{L^1_y}\nonumber\\
	&+\|\mathcal{I}(g_{\rho,0})\|_{L^2_{y}}+\|\partial_y^{-1}g_{v,0}\|_{L^2_y}+\|\partial_y^{-1}g_{v,0}\|_{L^\infty_y}+\nu^{\frac12}\|(g_{\rho,0},g_{v,0})\|_{L^2_y}+\nu^{\frac58}\|g_{u,0}\|_{L^2_y}
	\nonumber\\
	&+\nu\|g_{\rho,0}\|_{L^\infty_y}+\nu^{\frac{3}{2}}\|\partial_yg_{\rho,0}\|_{L^2_y}+\|(g_{\rho,\neq},g_{u,\neq},g_{v,\neq})\|_{L^2(\Omega)}+\nu\|\nabla_{x,y}g_{\rho,\neq}\|_{L^2(\Omega)},\label{g}\\
	[[(g_\rho,g_u,g_v)]]_{2}\eqdef& [[(g_\rho,g_u,g_v)]]_{1}+\nu^{\frac{13}{8}}\|\nabla_{x,y}(g_u,g_v)\|_{L^2(\Omega)}+\nu^{\frac{21}{8}}\|\nabla^2_{x,y}g_{\rho}\|_{L^2(\Omega)}, 	\label{gn}\\
	[[(g_\rho,g_u,g_v)]]_{3}\eqdef& [[(g_\rho,g_u,g_v)]]_{2}+\nu^{\frac{21}{8}}\|\nabla_{x,y}^2(g_u,g_v)\|_{L^2(\Omega)}+\nu^{\frac{29}{8}}\|\nabla^3_{x,y}g_{\rho}\|_{L^2(\Omega)}.	\label{gnh}
\end{align}

Recall the solution space $\mathfrak{X}$ and associated norm $\|\cdot\|_{\mathcal{X}}$ defined in \eqref{SS} and \eqref{S1} respectively. For convenience, we also introduce the low-order norm
\begin{align}
	\|(\rho,u,v)\|_{1}=&\|v_0\|_{L^1_y}+\|(m^{-2}\rho_0,u_0,v_0)\|_{L^\infty_y}+\|(m^{-2}\rho_0,v_0)\|_{L^2_y}+\|(m^{-2}\rho_{\neq},u_{\neq},v_{\neq})\|_{L^2(\Omega)}\nonumber\\
	&+\nu^{\frac12}\|\nabla_{x,y}(m^{-2}\rho_,u,v)\|_{L^2(\Omega)}+\nu^{\frac{13}{8}}\|\nabla_{x,y}^2(u,v)\|_{L^2(\Omega)},\label{el}
\end{align}
which contains derivatives up to the second order, and the-middle order norm 
\begin{align}
	\|(\rho,u,v)\|_{2}=&\|(\rho,u,v)\|_{1}+\nu^{\frac{13}{8}}\|\nabla^2_{x,y}\rho\|_{L^2(\Omega)}+\nu^{\frac{21}{8}}\|\nabla^3_{x,y}(u,v)\|_{L^2(\Omega)},\label{elm}
\end{align}
which contains derivatives up to the third order. In particular, for any $\beta>0$, by $L^\infty$-inequality \eqref{lw} and an interpotaion
$\|\partial_{x}^{1+\beta}f\|_{L^2(\Omega)}\leq C\|\partial_xf\|_{L^2(\Omega)}^{1-\beta}\|\partial_x^2f\|_{L^2(\Omega)}^{\beta}$, we obtain
\begin{align}\label{Lw1}
	&\|(m^{-2}\rho,u,v)\|_{L^\infty(\Omega)}\nonumber\\
	&\quad\leq C\|(m^{-2}\rho_0,u_0,v_0)\|_{L^\infty_y}+C_\beta\|\nabla_{x,y} (m^{-2}\rho_{\neq},u_{\neq},v_{\neq})\|_{L^2(\Omega)}^{1-\frac{\beta}{2}}\|\nabla_{x,y}^2(m^{-2}\rho_{\neq},u_{\neq},v_{\neq})\|_{L^2(\Omega)}^{\frac{\beta}{2}}\nonumber\\
	&\quad\leq C_\beta\nu^{-\frac12-\frac{9}{16}\beta}\|(\rho,u,v)\|_{2}.
\end{align}
Similar, it holds that
\begin{align}\label{Lw2}
	&\|(\nabla_{x,y}u,\nabla_{x,y}v)\|_{L^\infty(\Omega)}\nonumber\\
	&\quad\leq C\|(\partial_yu_0,\partial_yv_0)\|_{L^\infty_y}+ C_\beta\|(\nabla_{x,y}^2 u,\nabla_{x,y}^2v)\|_{L^2(\Omega)}^{1-\frac{\beta}{2}}\|(\nabla_{x,y}^3u,\nabla_{x,y}^3v)\|_{L^2(\Omega)}^{\frac{\beta}{2}}\nonumber\\
	&\quad\leq C\|(\nabla_{x,y}u,\nabla_{x,y})\|_{L^2(\Omega)}^{\frac12}\|(\nabla_{x,y}^2u,\nabla^2_{x,y}v)\|_{L^2(\Omega)}^{\frac12}+C_\beta\|(\nabla_{x,y}^2 u,\nabla_{x,y}^2v)\|_{L^2(\Omega)}^{1-\frac{\beta}{2}}\|(\nabla_{x,y}^3u,\nabla_{x,y}^3v)\|_{L^2(\Omega)}^{\frac{\beta}{2}}\nonumber\\
	&\quad\leq C_\beta\nu^{-\frac{13}{8}-\frac{1}{2}\beta}\|(\rho,u,v)\|_{2},
\end{align}
and
\begin{align}\label{Lw3}
	&m^{-2}\|\nabla_{x,y}\rho\|_{L^\infty(\Omega)}+\nu\|(\nabla_{x,y}^2u,\nabla_{x,y}^2v)\|_{L^\infty(\Omega)}\nonumber\\
	&\quad\leq C\|m^{-2}\nabla_{x,y}\rho\|_{L^2(\Omega)}^{\frac12}\|m^{-2}\nabla_{x,y}^2\rho\|_{L^2(\Omega)}^{\frac12}+\nu\|(\nabla_{x,y}^2u,\nabla_{x,y}^2v)\|_{L^2(\Omega)}^{\frac12}\|(\nabla_{x,y}^3u,\nabla_{x,y}^3v)\|_{L^2(\Omega)}^{\frac12}\nonumber\\ &\qquad+C_\beta\|m^{-2}\nabla_{x,y}^2\rho\|_{L^2(\Omega)}^{1-\frac{\beta}{2}}\|m^{-2}\nabla_{x,y}^3\rho\|_{L^2(\Omega)}^{\frac{\beta}{2}}+C_\beta\nu\|(\nabla_{x,y}^3 u,\nabla_{x,y}^3v)\|_{L^2(\Omega)}^{1-\frac{\beta}{2}}\|(\nabla_{x,y}^4u,\nabla_{x,y}^4v)\|_{L^2(\Omega)}^{\frac{\beta}{2}}\nonumber\\
	&\quad\leq C_\beta\nu^{-\frac{13}{8}-\frac{1}{2}\beta}\|(\rho,u,v)\|_{\mathfrak{X}}.
\end{align}
For the external force, since $\mathcal{P}_0F_{\rm ext,1}=\mathcal{P}_0F_{\rm ext,2}=0$, we have the following $L^\infty$-bound:
\begin{align}
	\|(F_{\rm ext,1},F_{\rm ext,2})\|_{L^\infty(\Omega)}&\leq C_\beta\|(\nabla_{x,y}F_{\rm ext,1},\nabla_{x,y}F_{\rm ext,2})\|_{L^2(\Omega)}^{1-\frac{\beta}{2}}\|(\nabla_{x,y}^2F_{\rm ext,1},\nabla_{x,y}^2F_{\rm ext,2})\|_{L^2(\Omega)}^{\frac{\beta}{2}}\nonumber\\
	& \leq C_\beta\nu^{-\frac{13}{8}-\frac{\beta}{2}}\|(F_{\rm ext,1},F_{\rm ext,2})\|_{w}.\label{Lw4}
\end{align}
 Moreover, the given the velocity field $(\tilde{u},\tilde{v})$ in \eqref{7.1} is in the space 
\begin{align}
	\mathfrak{X}_*\eqdef\bigg\{&(\tilde{u},\tilde{v})\mid \tilde{v}_0\in L^1(\mathbb{R}_+)\cap H^4(\mathbb{R}_+), \tilde{u}_0\in L^{\infty}(\mathbb{R}_+)\cap \dot{H^4}(\mathbb{R}_+),(\tilde{u}_{\neq},\tilde{v}_{\neq})\in H^4(\mathbb{T}_L\times \mathbb{R}_+)^2,\nonumber\\
	&\qquad\qquad\tilde{u}|_{y=0}=\tilde{v}|_{y=0}=\div_{x,y}(\tilde{u},\tilde{v})|_{y=0}=0, \|(\tilde{u},\tilde{v})\|_{\mathfrak{X}_*}<\infty\bigg\},\label{X*}
\end{align}
where the norm $\|(\tilde{u},\tilde{v})\|_{\mathfrak{X}_*}$ is defined by
\begin{align}
	\|(\tilde{u},\tilde{v})\|_{\mathfrak{X}_*}=&\|\tilde{v}_0\|_{L^1_y}+\|\tilde{v}_0\|_{L^2_y}+\|(\tilde{u}_0,\tilde{v}_0)\|_{L^\infty_y}+\|(\tilde{u}_{\neq},\tilde{v}_{\neq})\|_{L^2(\Omega)}+\nu^{\frac12}\|\nabla_{x,y}(\tilde{u},\tilde{v})\|_{L^2(\Omega)}\nonumber\\
	&+\nu^{\frac{13}{8}}\|\nabla_{x,y}^2(\tilde{u},\tilde{v})\|_{L^2(\Omega)}+\nu^{\frac{21}{8}}\|\nabla_{x,y}^3(\tilde{u},\tilde{v})\|_{L^2(\Omega)}+\nu^{\frac{29}{8}}\|\nabla_{x,y}^4(\tilde{u},\tilde{v})\|_{L^2(\Omega)}.\label{el1}
\end{align}
In \eqref{X*}, we introduce an additional boundary condition $\div_{x,y}(\tilde{u},\tilde{v})|_{y=0}=0$. This identity is crucial to obtain higher order estimates, cf.  \eqref{HE10-1} and \eqref{HE10-2}. We will show that it can be preserved in the nonlinear iteration later.

For later use, we define the middle-order norm of the given velocity fields $(\tilde{u},\tilde{v})$.
\begin{align}
	\|(\tilde{u},\tilde{v})\|_{2,*}=&\|\tilde{v}_0\|_{L^1_y}+\|\tilde{v}_0\|_{L^2_y}+\|(\tilde{u}_0,\tilde{v}_0)\|_{L^\infty_y}+\|(\tilde{u}_{\neq},\tilde{v}_{\neq})\|_{L^2(\Omega)}\nonumber\\
	&+\nu^{\frac12}\|\nabla_{x,y}(\tilde{u},\tilde{v})\|_{L^2(\Omega)}+\nu^{\frac{13}{8}}\|\nabla_{x,y}^2(\tilde{u},\tilde{v})\|_{L^2(\Omega)}+\nu^{\frac{21}{8}}\|\nabla_{x,y}^3(\tilde{u},\tilde{v})\|_{L^2(\Omega)}.\label{el2}
\end{align}

From Theorems \ref{thmo} and \ref{thmno}, we have the following bounds of low-order norm of the solution.
\begin{proposition}\label{propl}
	Let the Mach number $m\in (0,1)$. For any $L\in (0,L_0)$, where $L_0$ is given in Theorem \ref{thmno} and $0<\nu\ll 1$, if $\|(\tilde{u},\tilde{v})\|_{2,*}<\infty$, and
	$[[(g_\rho,g_u,g_v)]]_{1}<\infty$, the solution $(\rho,u,v)$ to \eqref{7.1} satisfies the following bound
	\begin{align}
		\|(\rho,u,v)\|_{1}\leq C [[(g_\rho,g_u,g_v)]]_{1}+C\nu^{-(\frac{9}{8})^-}\|(\tilde{u},\tilde{v})\|_{2,*}\|(\rho,u,v)\|_{2}.\label{low}
	\end{align} 
\end{proposition}
\begin{proof}
	We rewrite \eqref{7.1} into
	\begin{align}
		\begin{cases}
			U_s\partial_x\rho+\div_{x, y}(u, v)=g_\rho-\partial_x(\tilde{u}\rho)-\partial_y(\tilde{v}\rho)\eqdef\tilde{g}_\rho,\\
			U_s\partial_xu+v\partial_yU_s+m^{-2}\partial_x\rho-{\nu}\Delta_{x, y} u-\lambda{\nu}\partial_x\div_{x, y}(u, v)+\nu\rho\partial_y^2U_s=g_u,\\
			U_s\partial_xv+m^{-2}\partial_y\rho-{\nu}\Delta_{x, y} v-\lambda{\nu}\partial_y\div_{x, y}(u, v)=g_v,\\
			u|_{y=0}=v|_{y=0}=0.\nonumber
		\end{cases}
	\end{align}
Now we estimate the zero and non-zero modes of the solution separately.

\underline{\it Estimates on zero modes:} Recall \eqref{elm} and \eqref{el2} for the definition of $\|(\rho,u,v)\|_2$ and  $\|(\tilde{u},\tilde{v})\|_{2,*}$ respectively. Note that $\mathcal{I}(\tilde{g}_{\rho,0})=\mathcal{I}(g_{\rho,0})-\mathcal{P}_0(\tilde{v}\rho)$, and $\tilde{g}_{\rho,0}=g_{\rho,0}-\mathcal{P}_0\left(\partial_y(\tilde{v}\rho)\right)$. Then from \eqref{A9.1}-\eqref{A9.3} we obtain
\begin{align}
	\|\mathcal{I}(\tilde{g}_{\rho,0})\|_{L^1_y}&\leq \|\mathcal{I}(g_{\rho,0})\|_{L^1_y}+\|\mathcal{P}_0(\tilde{v}\rho)\|_{L^1_y}\nonumber\\
	&\leq \|\mathcal{I}(g_{\rho,0})\|_{L^1_y}+C\|\tilde{v}\|_{L^2(\Omega)}\|\rho\|_{L^2(\Omega)}\nonumber\\
	&\leq \|\mathcal{I}(g_{\rho,0})\|_{L^1_y}+C\|(\tilde{u},\tilde{v})\|_{2,*}\|(\rho,u,v)\|_{2},\label{L1}\\
		\|\mathcal{I}(\tilde{g}_{\rho,0})\|_{L^2_y}&\leq \|\mathcal{I}(g_{\rho,0})\|_{L^2_y}+\|\mathcal{P}_0(\tilde{v}\rho)\|_{L^2_y}&\nonumber\\
		&\leq \|\mathcal{I}(g_{\rho,0})\|_{L^2_y}+C\|\tilde{v}_0\|_{L^\infty_y}\|\rho_0\|_{L^2_y}+C\|\tilde{v}_{\neq}\|_{L^2(\Omega)}^{\frac12}\|\partial_y\tilde{v}_{\neq}\|_{L^2(\Omega)}^{\frac12}\|\rho_{\neq}\|_{L^2(\Omega)}\nonumber\\
		&\leq  \|\mathcal{I}(g_{\rho,0})\|_{L^2_y}+C\nu^{-\frac14}\|(\tilde{u},\tilde{v})\|_{2,*}\|(\rho,u,v)\|_{2} \label{L2},\\
		\|\mathcal{I}(\tilde{g}_{\rho,0})\|_{L^\infty_y}&\leq \|\mathcal{I}(g_{\rho,0})\|_{L^\infty_y}+\|\mathcal{P}_0(\tilde{v}\rho)\|_{L^\infty_y}\nonumber\\
		&\leq \|\mathcal{I}(g_{\rho,0})\|_{L^\infty_y}+C\|\tilde{v}_0\|_{L^\infty_y}\|\rho_0\|_{L^\infty_y}+C\|\tilde{v}_{\neq}\|_{L^2(\Omega)}^{\frac12}\|\partial_y\tilde{v}_{\neq}\|_{L^2(\Omega)}^{\frac12}\|\rho_{\neq}\|_{L^2(\Omega)}^{\frac12}\|\partial_y\rho_{\neq}\|_{L^2(\Omega)}^{\frac12}\nonumber\\
		&\leq \|\mathcal{I}(g_{\rho,0})\|_{L^\infty_y} +C\nu^{-\frac12}\|(\tilde{u},\tilde{v})\|_{2,*}\|(\rho,u,v)\|_{2}\label{L3},\\
		\|\tilde{g}_{\rho,0}\|_{L^2_y}&\leq \|g_{\rho,0}\|_{L^2_y}+\|\mathcal{P}_0(\rho\partial_y\tilde{v})\|_{L^2_{y}}+\|\mathcal{P}_0(\tilde{v}\partial_y\rho)\|_{L^2_{y}}\nonumber\\
		&\leq \|g_{\rho,0}\|_{L^2_y}+C\|\rho_0\|_{L^\infty_y}\|\partial_y\tilde{v}_0\|_{L^2_y}+C\|\tilde{v}_0\|_{L^\infty_y}\|\partial_y\rho_0\|_{L^2_y}\nonumber\\
		&\quad+C\|\tilde{v}_{\neq}\|_{L^2(\Omega)}^{\frac12}\|\partial_y\tilde{v}_{\neq}\|_{L^2(\Omega)}^{\frac12}\|\partial_y\rho_{\neq}\|_{L^2(\Omega)}+C\|\rho_{\neq}\|_{L^2(\Omega)}^{\frac12}\|\partial_y\rho_{\neq}\|_{L^2(\Omega)}^{\frac12}\|\partial_y\tilde{v}_{\neq}\|_{L^2(\Omega)}\nonumber\\
		&\leq \|g_{\rho,0}\|_{L^2_y}+C\nu^{-\frac34}\|(\tilde{u},\tilde{v})\|_{2,*}\|(\rho,u,v)\|_{2},\label{L4}\\
		\|\tilde{g}_{\rho,0}\|_{L^\infty_y}
		&\leq \|g_{\rho,0}\|_{L^\infty_y}+\|\mathcal{P}_0(\rho\partial_y\tilde{v})\|_{L^\infty_y}+\|\mathcal{P}_0(\tilde{v}\partial_y\rho)\|_{L^\infty_y}\nonumber\\
		&\leq \|g_{\rho,0}\|_{L^\infty_y} +C\|\rho_0\|_{L^\infty_y}\|\partial_y\tilde{v}_0\|_{L^2_y}^{\frac12}\|\partial_y^2\tilde{v}_0\|_{L^2_y}^{\frac12}+C\|\tilde{v}_0\|_{L^\infty_y}\|\partial_y\rho_0\|_{L^2_y}^{\frac12}\|\partial_y^2\rho_0\|_{L^2_y}^{\frac12}\nonumber\\
		&\quad+C\|\tilde{v}_{\neq}\|_{L^2(\Omega)}^{\frac12}\|\partial_y\tilde{v}_{\neq}\|_{L^2(\Omega)}^{\frac12}\|\partial_y\rho_{\neq}\|_{L^2(\Omega)}^{\frac12}\|\partial_y^2\rho_{\neq}\|_{L^2(\Omega)}^{\frac12}\nonumber\\
		&\quad+C\|\rho_{\neq}\|_{L^2(\Omega)}^{\frac12}\|\partial_y\rho_{\neq}\|_{L^2(\Omega)}^{\frac12}\|\partial_y\tilde{v}_{\neq}\|_{L^2(\Omega)}^{\frac12}\|\partial_y\tilde{v}_{\neq}\|_{L^2(\Omega)}^{\frac12}\nonumber\\
		&\leq  \|g_{\rho,0}\|_{L^\infty_y}+C\nu^{-\frac{21}{16}}\|(\tilde{u},\tilde{v})\|_{2,*}
		\|(\rho,u,v)\|_{2},\label{L5}\\
			\|\partial_y\tilde{g}_{\rho,0}\|_{L^2}&\leq \|\partial_yg_{\rho,0}\|_{L^2_y}+\|\mathcal{P}_0(\rho\partial_y^2\tilde{v})\|_{L^2_{y}}+\|\mathcal{P}_0(\partial_y\tilde{v}\partial_y\rho)\|_{L^2_y}+\|\mathcal{P}_0(\tilde{v}\partial_y^2\rho)\|_{L^2_{y}}\nonumber\\
			&\leq \|\partial_yg_{\rho,0}\|_{L^2_y}+C\|\rho_0\|_{L^\infty_y}\|\partial_y^2\tilde{v}_0\|_{L^2_y}+C\|\tilde{v}_0\|_{L^\infty_y}\|\partial_y^2\rho_0\|_{L^2_y}\nonumber\\
			&\quad+C\|\partial_y\tilde{v}_0\|_{L^2_y}^{\frac12}\|\partial_y^2\tilde{v}_0\|_{L^2_y}^{\frac12}\|\partial_y\rho_0\|_{L^2_y}+C\|\tilde{v}_{\neq}\|_{L^2(\Omega)}^{\frac12}\|\partial_y\tilde{v}_{\neq}\|_{L^2(\Omega)}^{\frac12}\|\partial_y^2\rho_{\neq}\|_{L^2(\Omega)}\nonumber\\
			&\quad +C\|\rho_{\neq}\|_{L^2(\Omega)}^{\frac12}\|\partial_y\rho_{\neq}\|_{L^2(\Omega)}^{\frac12}\|\partial_y^2\tilde{v}_{\neq}\|_{L^2(\Omega)}+C\|\partial_y\tilde{v}_{\neq}\|_{L^2_{x,y}}^{\frac12}\|\partial_y^2\tilde{v}_{\neq}\|_{L^2(\Omega)}^{\frac12}\|\partial_y\rho_{\neq}\|_{L^2(\Omega)}\nonumber\\
			&\leq \|\partial_yg_{\rho,0}\|_{L^2_y}+C\nu^{-\frac{15}{8}}\|(\tilde{u},\tilde{v})\|_{2,*}\|(\rho,u,v)\|_{2} .\label{L6}
\end{align}
Sustituting \eqref{L1}-\eqref{L6} into \eqref{3.5}-\eqref{3.12} in Theorem \ref{thmo}, we obtain
\begin{align}
	\|v_0\|_{L^p_y}&\le C\|\mathcal{I}[g_{\rho, 0}]\|_{L^p_y}+C\nu^{-\frac{1}{2}}\|(\tilde{u},\tilde{v})\|_{2,*}\|(\rho,u,v)\|_{2},~p=1,2,~{\rm and }~ \infty,\label{L7}\\ 
	\|\partial_yv_0\|_{L^2_y}&\le C\|g_{\rho, 0}\|_{L^2_y}+C\nu^{-\frac{3}{4}}\|(\tilde{u},\tilde{v})\|_{2,*}\|(\rho,u,v)\|_{2},\label{L7-1}\\
	\|\partial_y^2v_0\|_{L^2_y}&\le C\|\partial_yg_{\rho, 0}\|_{L^2_y}+C\nu^{-\frac{15}{8}}\|(\tilde{u},\tilde{v})\|_{2,*}\|(\rho,u,v)\|_{2},\label{L8}\\
	m^{-2}\|\rho_0\|_{L^2_y}&\le C\nu\|g_{\rho, 0}\|_{L^2_y}+C\|\partial_y^{-1}g_{v, 0}\|_{L^2_y}+C\nu^{\frac{1}{4}}\|(\tilde{u},\tilde{v})\|_{2,*}\|(\rho,u,v)\|_{2},\label{L9}\\ 
	m^{-2}\|\rho_0\|_{L^\infty_y}&\le C\nu\|g_{\rho, 0}\|_{L^\infty_y}+C\|\partial_y^{-1}g_{v, 0}\|_{L^\infty_y}+C\nu^{-\frac{5}{16}}\|(\tilde{u},\tilde{v})\|_{2,*}\|(\rho,u,v)\|_{2},\label{L9-1}\\
	m^{-2}\|\partial_y\rho_0\|_{L^2_y}&\le C{\nu}\|\partial_yg_{\rho, 0}\|_{L^2_y}+C\|g_{v, 0}\|_{L^2_y}+C\nu^{-\frac{7}{8}}\|(\tilde{u},\tilde{v})\|_{2,*}\|(\rho,u,v)\|_{2},\label{L10}\\
	\|u_0\|_{L^\infty_y}&\le C\nu^{-\frac12}\|\mathcal{I}[g_{\rho, 0}]\|_{L^\infty_y}+C{\nu}\|g_{\rho, 0}\|_{L^\infty_y}+C\nu^{-1}\|\partial_y^{-1}g_{u, 0}\|_{L^1_y}\nonumber\\
	&\qquad+C\|\partial_y^{-1}g_{v, 0}\|_{L^\infty_y}+C\nu^{-1}\|(\tilde{u},\tilde{v})\|_{2,*}\|(\rho,u,v)\|_{2},\label{L11}\\
	\|\partial_yu_0\|_{L^2_y}&\le C\nu^{-\frac34}\|\mathcal{I}[g_{\rho, 0}]\|_{L^\infty_y}+C{\nu}^{\frac34}\|g_{\rho, 0}\|_{L^\infty_y}+C\nu^{-1}\|\partial_y^{-1}g_{u, 0}\|_{L^2_y}\nonumber\\
	&\qquad+C\nu^{-\frac14}\|\partial_y^{-1}g_{v, 0}\|_{L^\infty_y}+C\nu^{-\frac54}\|(\tilde{u},\tilde{v})\|_{2,*}\|(\rho,u,v)\|_{2},\label{L12}\\
	\|\partial_y^2u_0\|_{L^2_y}&\le C\nu^{-\frac54}\|\mathcal{I}[g_{\rho, 0}]\|_{L^\infty_y}+C{\nu}^{\frac14}\|g_{\rho, 0}\|_{L^\infty_y}+C\nu^{-1}\|g_{u, 0}\|_{L^2_y}\nonumber\\
	&\qquad+C\nu^{-\frac34}\|\partial_y^{-1}g_{v, 0}\|_{L^\infty_y}+C\nu^{-\frac{7}4}\|(\tilde{u},\tilde{v})\|_{2,*}\|(\rho,u,v)\|_{2}\label{L13}.
\end{align}
\underline{\it Estimate on non-zero modes:} Note that $\tilde{g}_{\rho,\neq}=g_{\rho,\neq}-\mathcal{P}_{\neq}(\tilde{u}\partial_x\rho)-\mathcal{P}_{\neq}(\tilde{v}\partial_y\rho)-\mathcal{P}_{\neq}(\rho\partial_x\tilde{u})-\mathcal{P}_{\neq}(\rho\partial_y\tilde{v})$. Using \eqref{A9.5} and $L^\infty$-bounds \eqref{Lw1}, we obtain
\begin{align}
	\|\tilde{g}_{\rho,\neq}\|_{L^2(\Omega)}&\leq C\|g_{\rho,\neq}\|_{L^2(\Omega)}+\|(\tilde{u}\partial_x\rho)_{\neq}\|_{L^2(\Omega)}+\|(\tilde{v}\partial_y\rho)_{\neq}\|_{L^2(\Omega)}+\|(\rho\partial_x\tilde{u})_{\neq}\|_{L^2(\Omega)}+\|(\rho\partial_y\tilde{v})_{\neq}\|_{L^2(\Omega)}\nonumber\\
	&\leq C\|g_{\rho,\neq}\|_{L^2(\Omega)}+C\|(\tilde{u},\tilde{v})\|_{L^\infty(\Omega)}\|\nabla_{x,y}\rho\|_{L^2(\Omega)}+\|\rho\|_{L^\infty(\Omega)}\|(\nabla_{x,y}\tilde{u},\nabla_{x,y}\tilde{v})\|_{L^2(\Omega)}\nonumber\\
	&\leq C\|g_{\rho,\neq}\|_{L^2(\Omega)}+ C\nu^{-1^{-}}\|(\tilde{u},\tilde{v})\|_{2,*}\|(\rho,u,v)\|_{2}.\label{L14}
\end{align}
Let $\partial=\partial_x$ or $\partial_y$. Then it holds that
$$\partial \tilde{g}_\rho=\partial g_\rho-\tilde{u}\partial\partial_x\rho-\tilde{v}\partial\partial_y\rho-\partial\tilde{u}\partial_x\rho-\partial\tilde{v}\partial_y\rho-\partial\rho(\partial_x\tilde{u}+\partial_v\tilde{v})-\rho\partial(\partial_x\tilde{u}+\partial_v\tilde{v}).$$
Similar to \eqref{L14}, we use \eqref{A9.5}, and $W^{1,\infty}$ bounds \eqref{Lw1} and \eqref{Lw2} to obtain 
\begin{align}
	\|\nabla_{x,y}\tilde{g}_{\rho,\neq}\|_{L^2(\Omega)}
	&\leq\|\nabla_{x,y}g_{\rho,\neq}\|_{L^2(\Omega)}+C\|(\tilde{u},\tilde{v})\|_{L^\infty(\Omega)}\|\nabla^2_{x,y}\rho\|_{L^2(\Omega)}\nonumber\\
	&\qquad+C\|\rho\|_{L^\infty(\Omega)}\|\nabla_{x,y}^2(\tilde{u},\tilde{v})\|_{L^2(\Omega)}+C\|\nabla_{x,y}(\tilde{u},\tilde{v})\|_{L^\infty(\Omega)}\|\nabla_{x,y}\rho\|_{L^2(\Omega)}\nonumber\\
	&\leq C\|\nabla_{x,y}g_{\rho,\neq}\|_{L^2(\Omega)}+C\nu^{-(\frac{17}{8})^-}\|(\tilde{u},\tilde{v})\|_{2,*}\|(\rho,u,v)\|_{2}.\label{L15}
\end{align}
Then from \eqref{7.4-11}, \eqref{7.4-12} and \eqref{7.4-13} in Theorem \ref{thmno}, we have
\begin{align}
	&\nu^{\frac{13}{8}}\|\nabla_{x,y}^2(u_{\neq},v_{\neq})\|_{L^2(\Omega)}+\nu^{\frac12}\|\nabla_{x,y}(m^{-2}\rho_{\neq},u_{\neq},v_{\neq})\|_{L^2(\Omega)}+\|(m^{-2}\rho_{\neq},u_{\neq},v_{\neq})\|_{L^2(\Omega)}\nonumber\\
	&\quad\qquad\leq C\|(g_{\rho,\neq},g_{u,\neq},g_{v,\neq})\|_{L^2(\Omega)}+\nu\|\nabla_{x,y}g_{\rho,\neq}\|_{L^2(\Omega)}+C\nu^{-(\frac{9}{8})^{-}}\|(\tilde{u},\tilde{v})\|_{2,*}\|(\rho,u,v)\|_{2}\label{L16}.
\end{align}
Combining \eqref{L7}-\eqref{L13} for the estimates on the zero mode and \eqref{L16} for the estimates on non-zero modes, and then comparing the order of $\nu$,  we obtain the desired estimate \eqref{low}. The proof of Proposition \ref{propl} is complete.
\end{proof}

\subsection{Estimates on high order derivatives}
In this subsection, we denote the standard inner product in $L^2(\Omega)$ by $\langle f,g\rangle=\iint_{\Omega}f g dx dy.$ Next proposition gives the higher order estimates on the solution to \eqref{7.1}.

%\begin{align}\label{7.3}
%&\|\mathcal{P}_{\neq}\rho\|_{H^1}+\|(\mathcal{P}_{\neq}u, \mathcal{P}_{\neq}v)\|_{H^2}\notag\\
%&\quad\lesssim \nu^{-\frac54}(\|\mathcal{P}_{\neq}f_0\|_{H^1}+\|(\mathcal{P}_{\neq}f_{1}, \mathcal{P}_{\neq}f_{2})\|_{L^2}+\|
%\tilde{u}\partial_X\rho\|_{H^1}+\|\tilde{v}\partial_Y\rho\|_{H^1}+\|\rho(\partial_X\tilde{u}+\partial_Y\tilde{v})\|_{H^1}).
%\end{align}

\begin{proposition}\label{L7.1}
Let $m\in (0, 1)$. For any given velocity field  $(\tilde{u},\tilde{v})\in \mathfrak{X}_*$, and any inhomogeneous source term $(g_\rho,g_u,g_v)\in H^3(\Omega)\times H^2(\Omega)^2$, where $g_\rho\in H^3(\Omega)$ satisfies $g_\rho|_{y=0}=0$ and $\iint_{\Omega}g_\rho(x,y)d x dy=0$, the solution $(\rho, u, v)$ to the linear Navier-Stokes system \eqref{7.1} satisfies the zero mass condition:
$$\iint_{\Omega}\rho(x,y)dx dy=0,
$$
and the following bounds:
\begin{align}\label{HE0}
	&{\nu}\|\nabla^3_{x,y}(u, v)\|_{L^2(\Omega)}+m^{-2}\|\nabla^2_{x,y}\rho\|_{L^2(\Omega)}\nonumber\\
	&\qquad\leq C\nu^{-(\frac{31}{16})^-}(1+\nu^{-\frac{5}{16}}\|(\tilde{u},\tilde{v})\|_{\mathfrak{X}_*}^{\frac12})\|(\tilde{u},\tilde{v})\|_{\mathfrak{X}_*}^{\frac12}\|(\rho,u,v)\|_2\nonumber\\ 
	&\qquad\qquad+C\nu^{-\frac{13}{8}}\|(\rho,u,v)\|_1 +C\nu\|\nabla_{x,y}^2g_{\rho}\|_{L^2(\Omega)}+C\|(g_{u},g_{v})\|_{H^1(\Omega)},
\end{align}
and
\begin{align}\label{HF0}
	&{\nu}\|\nabla^4_{x,y}(u, v)\|_{L^2(\Omega)}+m^{-2}\|\nabla^3_{x,y}\rho\|_{L^2(\Omega)}\nonumber\\
	&\qquad\leq C\nu^{-(\frac{47}{16})^-}(1+\nu^{-\frac{5}{16}}\|(\tilde{u},\tilde{v})\|_{\mathfrak{X}_*}^{\frac12})\|(\tilde{u},\tilde{v})\|_{\mathfrak{X}_*}^{\frac12}\|(\rho,u,v)\|_{\mathfrak{X}}\nonumber\\ 
	&\qquad\qquad+C\nu^{-\frac{21}{8}}\|(\rho,u,v)\|_2 +C\nu\|\nabla_{x,y}^3g_{\rho}\|_{L^2(\Omega)}+C\|(g_{u},g_{v})\|_{H^2(\Omega)}.
\end{align}
\end{proposition}
We divide the proof of Proposition \ref{L7.1} into several steps. Firstly, the following lemma gives bounds on the tangential derivatives.
\begin{lemma} \label{lmHF1} Let  $(g_\rho,g_u,g_v)\in H^3(\Omega)\times H^2(\Omega)^2$. The solution $(\rho,u,v)\in \mathfrak{X}$ to the linear system \eqref{7.1} satisfies
	\begin{align}
		&\nu\|\nabla_{x,y}(\partial_x^2u,\partial_x^2v)\|_{L^2(\Omega)}+\nu\|\partial_x^2\div_{x,y}(u,v)\|_{L^2(\Omega)}\nonumber\\
		&\qquad\leq o(1)m^{-2}\|\nabla_{x,y}^2\rho\|_{L^2(\Omega)} +C\nu^{-(\frac{31}{16})^-}\|(\tilde{u},\tilde{v})\|_{\mathfrak{X}_*}^{\frac12}\|(\rho,u,v)\|_{2}\nonumber\\
		&\qquad\qquad+C\nu^{-\frac{11}{8}}\|(\rho,u,v)\|_1+C\nu\|\nabla_{x,y}^2g_\rho\|_{L^2(\Omega)}+C\|\nabla_{x,y}(g_u,g_v)\|_{L^2(\Omega)},\label{HF1}
	\end{align}
and 
\begin{align}
	&\nu\|\nabla_{x,y}(\partial_x^3u,\partial_x^3v)\|_{L^2(\Omega)}+\nu\|\partial_x^3\div_{x,y}(u,v)\|_{L^2(\Omega)}\nonumber\\
	&\qquad\leq o(1)m^{-2}\|\nabla_{x,y}^3\rho\|_{L^2(\Omega)} +C\nu^{-(\frac{47}{16})^-}\|(\tilde{u},\tilde{v})\|_{\mathfrak{X}_*}^{\frac12}\|(\rho,u,v)\|_{\mathfrak{X}}\nonumber\\
	&\qquad\qquad+C\nu^{-\frac{19}{8}}\|(\rho,u,v)\|_2+C\nu\|\nabla_{x,y}^3g_\rho\|_{L^2(\Omega)}+C\|\nabla_{x,y}^2(g_u,g_v)\|_{L^2(\Omega)}.\label{HF2}
\end{align}
Here $o(1)$ represents a positive number that can be chosen arbitrarily small.
\end{lemma}
\begin{proof} Let $k=2$ or $3$. Applying $\partial_x^k$ to $ \eqref{7.1}_2$ and $\eqref{7.1}_3$, then taking inner product of the resulting equations with $\nu\partial_x^ku$ and $\nu\partial_x^kv$ respectively, we obtain
\begin{equation}
\begin{aligned}
&{\nu}^2\|\nabla_{x,y}(\partial_x^ku, \partial_x^kv)\|_{L^2_{x,y}}^2+\lambda{\nu}^2\|\partial_x^k\div_{x, y}(u, v)\|_{L^2_{x,y}}^2\\
&\quad=
-\langle U_s\partial_x^{k+1}u, \nu\partial_x^ku\rangle-\langle U_s\partial_x^{k+1}v, \nu\partial_x^kv\rangle-m^{-2}\langle \partial_x^{k+1}\rho, \nu\partial_x^ku\rangle-m^{-2}\langle\partial_y\partial_x^k\rho, \nu\partial_x^kv\rangle\\
&\qquad-\langle\partial_x^kv\partial_yU_s, \nu\partial_x^ku\rangle-\langle {\nu}\partial_x^k\rho\partial_y^2U_s, \nu\partial_x^ku\rangle
+\langle\partial_x^kg_{u}, \nu\partial_x^ku\rangle+\langle\partial_x^kg_{v}, \nu\partial_x^kv\rangle=: \sum\limits_{i=1}^{8}K_i.\label{HE1}
\end{aligned}
\end{equation}

By integration by parts, we have
\begin{align}\label{HE2}
K_1=K_2=0,
\end{align}
and
\begin{align}\label{HE3}
|K_7|+|K_8|&=|\langle\partial_x^{k-1}g_{u}, \nu\partial_x^{k+1}u\rangle|+|\langle\partial_x^{k-1}g_{v}, \nu\partial_x^{k+1}v\rangle|\nonumber\\
&\le\frac{\nu^2}{2}\|\nabla_{x,y}(\partial_x^ku, \partial_x^kv)\|_{L^2(\Omega)}^2+C\|\nabla_{x,y}^{k-1}(g_{u}, g_{v})\|_{L^2(\Omega)}^2.
\end{align}

Next, we estimate $K_3$ and $K_4$. Consider $k=2$ first. Integrating by parts again and using $\eqref{7.1}_1$, and also the boundary conditions  $\tilde{u}|_{y=0}=\tilde{v}|_{y=0}=0$, we obtain
\begin{align}\label{HE4}
K_3+K_4&=m^{-2}\nu\langle \partial_x^2\rho, \partial_x^2\div_{x, y}(u, v)\rangle\notag\\
&=m^{-2}\nu\langle \partial_x^2\rho, -U_s\partial_x^3\rho-\partial_x^3(\tilde{u}\rho)-\partial_x^2\partial_y(\tilde{v}\rho)+\partial_x^2g_{\rho}\rangle\notag\\
&=m^{-2}\nu\langle \partial_x^2\rho, -U_s\partial_x^3\rho-\tilde{u}\partial_x^3\rho-\tilde{v}\partial_x^2\partial_y\rho\rangle\notag\\
&\quad+m^{-2}\nu\langle \partial_x^2\rho, -3\partial_x\tilde{u}\partial_x^2\rho-2\partial_x\tilde{v}\partial_{xy}\rho-\partial_y\tilde{v}\partial_x^2\rho\rangle\notag\\
&\quad+m^{-2}\nu\langle \partial_x^2\rho, -3\partial_x^2\tilde{u}\partial_x\rho-2\partial_{xy}\tilde{v}\partial_x\rho-\partial_x^2\tilde{v}\partial_{y}\rho\rangle\notag\\
&\quad+ m^{-2}\nu\langle \partial_x^2\rho, -\rho\div_{x, y}(\tilde{u},\tilde{v})\rangle+m^{-2}\langle \partial_x^2\rho, \partial_x^2g_{\rho}\rangle\nonumber\\
&=\frac12m^{-2}\nu\langle \partial_x^2\rho, \partial_x^2\rho\div_{x, y}(\tilde{u},\tilde{v})\rangle\notag\\
&\quad-m^{-2}\langle \partial_x^2\rho, 3\partial_x\tilde{u}\partial_x^2\rho+2\partial_x\tilde{v}\partial_{xy}\rho+\partial_y\tilde{v}\partial_x^2\rho\rangle\notag\\
&\quad-m^{-2}\nu\langle \partial_x^2\rho, 3\partial_x^2\tilde{u}\partial_x\rho+\partial_x^2\tilde{v}\partial_{y}\rho+2\partial_{xy}\tilde{v}\partial_x\rho\rangle\notag\\
&\quad-m^{-2}\nu\langle \partial_x^2\rho, \rho\div_{x, y}(\tilde{u},\tilde{v})\rangle-m^{-2}\langle\partial_x^2\rho, \partial_x^2g_{\rho}\rangle\nonumber\\
&\leq Cm^{-2}\nu\|\nabla_{x,y}^2\rho\|_{L^2(\Omega)}^2\|\nabla_{x,y}(\tilde{u},\tilde{v})\|_{L^\infty(\Omega)}+Cm^{-2}\nu\|\nabla^2_{x,y}\rho\|_{L^2(\Omega)}\|\nabla_{x,y}\rho\|_{L^2(\Omega)}\|\nabla_{x,y}^2(\tilde{u},\tilde{v})\|_{L^\infty(\Omega)}\nonumber\\
&\qquad+Cm^{-2}\nu\|\nabla_{x,y}^2\rho\|_{L^2(\Omega)}\|\rho\|_{L^\infty(\Omega)}\|\nabla_{x,y}^2(\tilde{u},\tilde{v})\|_{L^2(\Omega)}+Cm^{-2}\nu\|\nabla^2_{x,y}\rho\|_{L^2(\Omega)}\|\nabla^2_{x,y}g_\rho\|_{L^2(\Omega)}\nonumber\\
&\leq o(1)m^{-4}\|\nabla_{x,y}^2\rho\|_{L^2(\Omega)}^2+ C\nu^{-(\frac{31}{8})^-}\|(\tilde{u},\tilde{v})\|_{\mathfrak{X}_*}\|(\rho,u,v)\|_{2}^2+C\nu^2\|\nabla^2_{x,y}g_\rho\|_{L^2(\Omega)}.
\end{align}
Here we have used $W^{2,\infty}$ bounds \eqref{Lw1}-\eqref{Lw3} in the last inequality. The case  $k=3$ is similar to \eqref{HE4}. Therefore, we can obtain
\begin{align}
	|K_3+K_4|&\leq Cm^{-2}\nu\iint_{\Omega}|\nabla^3_{x,y}\rho|\left(\sum_{j=1}^3|\nabla_{x,y}^j(\tilde{u},\tilde{v})||\nabla^{4-j}_{x,y}\rho|+|\nabla^3_{x,y}g_\rho|\right)d x dy\nonumber\\
	&\leq  Cm^{-2}\nu\|\nabla_{x,y}^3\rho\|_{L^2(\Omega)}^2\|\nabla_{x,y}(\tilde{u},\tilde{v})\|_{L^\infty(\Omega)}+Cm^{-2}\nu\|\nabla^3_{x,y}\rho\|_{L^2(\Omega)}\|\nabla_{x,y}^2\rho\|_{L^2(\Omega)}\|\nabla_{x,y}^2(\tilde{u},\tilde{v})\|_{L^\infty(\Omega)}\nonumber\\
	&\qquad+Cm^{-2}\nu\|\nabla_{x,y}^3\rho\|_{L^2(\Omega)}\|\nabla_{x,y}\rho\|_{L^\infty(\Omega)}\|\nabla_{x,y}^3(\tilde{u},\tilde{v})\|_{L^2(\Omega)}+Cm^{-2}\nu\|\nabla^3_{x,y}\rho\|_{L^2(\Omega)}\|\nabla_{x,y}^3g_\rho\|_{L^2(\Omega)}\nonumber\label{HE4-1}\\
	&\leq o(1)m^{-4}\|\nabla_{x,y}^3\rho\|_{L^2(\Omega)}^2+ C\nu^{-(\frac{47}{8})^-}\|(\tilde{u},\tilde{v})\|_{\mathfrak{X}_*}\|(\rho,u,v)\|_{\mathfrak{X}}^2+C\nu^2\|\nabla_{x,y}^3g_\rho\|_{L^2(\Omega)}^2,
\end{align}
where  the $W^{2,\infty}$ bounds \eqref{Lw1}-\eqref{Lw3} has been used in the last inequality in \eqref{HE4-1}.

By Cauchy-Schwarz, the remaining two terms in \eqref{HE1} can be bounded as follows.
\begin{align}\label{HE5}
|K_5|+|K_6|&\leq C\nu\|\nabla_{x,y}^k\rho\|_{L^2(\Omega)}\|\nabla_{x,y}^ku\|_{L^2(\Omega)}+C\nu^{\frac12}\|\nabla^k_{x,y}u\|_{L^2(\Omega)}\|\nabla^k_{x,y}v\|_{L^2(\Omega)}\nonumber\\
&\leq o(1)\|\nabla_{x,y}^k\rho\|_{L^2(\Omega)}^2+ C\nu^{\frac12}\|\nabla^k_{x,y}(u,v)\|_{L^2(\Omega)}^2\nonumber\\
&\leq o(1)\|\nabla_{x,y}^k\rho\|_{L^2(\Omega)}^2+C\nu^{-2k+\frac{5}{4}}\|(\rho,u,v)\|_{k-1}^2.
\end{align}
Substituting \eqref{HE2}, \eqref{HE3}, \eqref{HE4} and \eqref{HE5} into \eqref{HE1}, we get \eqref{HF1}. The estimates \eqref{HF2} follows from \eqref{HE1}, \eqref{HE5} with $k=3$, and \eqref{HE4-1}. The proof of Lemma \ref{lmHF1} is complete.
\end{proof}
Next we bound $\Delta_{x,y}\rho$ and full derivatives of $\div_{x,y}(u,v)$ in terms of $\partial_x^k\div_{x,y}(u,v)$.
\begin{lemma}\label{lmHF2}
	If $g_{\rho}|_{y=0}=0$ and $(\tilde{u},\tilde{v})\in \mathfrak{X}_*$, we have
	\begin{align}\label{HE0-1}
		&m^{-2}\|\Delta_{x, y}\rho\|_{L^2(\Omega)}+{\nu}\|\nabla^2_{x,y}\div_{x, y}(u, v)\|_{L^2(\Omega)}\notag\\
		&\quad \leq o(1)\|\nabla_{x,y}^2\rho\|_{L^2(\Omega)}+C\nu\|\partial_x^2\div_{x, y}(u, v)\|_{L^2(\Omega)}+C\nu\|\nabla_{x,y}^2g_{\rho}\|_{L^2(\Omega)}+C\|\nabla_{x,y}(g_{u},g_{v})\|_{L^2(\Omega)}\notag\\
		&\quad\quad +C\nu^{-(\frac{31}{16})^-}(1+\nu^{-\frac{5}{16}}\|(\tilde{u},\tilde{v})\|_{\mathfrak{X}_*}^{\frac12})\|(\tilde{u},\tilde{v})\|_{\mathfrak{X}_*}^{\frac12}\|(\rho,u,v)\|_2 +C\nu^{-\frac{13}{8}}\|(\rho,u,v)\|_1,
	\end{align}
and
\begin{align}\label{HE0-2}
	&m^{-2}\|\Delta_{x, y}\nabla_{x,y}\rho\|_{L^2(\Omega)}+{\nu}\|\nabla^3_{x,y}\div_{x, y}(u, v)\|_{L^2(\Omega)}\notag\\
	&\quad\leq o(1)\|\nabla_{x,y}^3\rho\|_{L^2(\Omega)}+C\nu\|\partial_x^3 \div_{x, y}(u,v)\|_{L^2(\Omega)}+C\nu\|\nabla_{x,y}^3g_\rho\|_{L^2(\Omega)}+C\|\nabla_{x,y}^2(g_u,g_v)\|_{L^2(\Omega)}\nonumber\\
	&\quad\quad+C\nu^{-(\frac{47}{16})^-}(1+\nu^{-\frac{5}{16}}\|(\tilde{u},\tilde{v})\|_{\mathfrak{X}_*}^{\frac{1}{2}})\|(\tilde{u},\tilde{v})\|_{\mathfrak{X}_*}^{\frac{1}{2}}\|(\rho,u,v)\|_{\mathfrak{X}}+C\nu^{-\frac{21}{8}}\|(\rho,u,v)\|_2.
\end{align}
Here $o(1)$ represents a positive number that can be chosen arbitrarily small.
\end{lemma}

\begin{proof}
Taking  $\partial_x$ and $\partial_y$ to $\eqref{7.1}_2$ and $\eqref{7.1}_3$ respectively, and adding them together, we obtain
\begin{align}\label{HE7}
	m^{-2}\Delta_{x,y}\rho=&{\nu}(1+\lambda)\Delta_{x,y}\div_{x, y}(u,v)-U_s\partial_x\div_{x, y}(u,v)\nonumber\\
	&-2\partial_xv\partial_yU_s-\nu\partial_x\rho \partial_y^2U_s+\div_{x, y}(g_u,g_v).
\end{align}
Note that the first term on the right hand side of \eqref{HE7} contains $\partial_y^2\div_{x,y}(u,v)$, which is hard to be controlled directly.  To eliminate it, we apply ${\nu}(1+\lambda)\partial_y^2$ to $\eqref{7.1}_1$ and add the resulted equation to \eqref{HE7}. Then it holds that
\begin{align}\label{HE8}
m^{-2}\Delta_{x, y}\rho&=-(1+\lambda){\nu}\left[(U_s+\tilde{u})\partial_{x}\partial_y^2\rho+\tilde{v}\partial_y^3\rho\right]\notag\\
&\quad+ (1+\lambda){\nu}\partial_x^2\div_{x, y}(u, v)-U_s\partial_x\div_{x, y}(u,v)-2\partial_xv\partial_yU_s\notag\\
&\quad-(1+\lambda){\nu}(2\partial_yU_s\partial_{xy}\rho+2\partial_y\tilde{u}\partial_{xy}\rho+3\partial_y\tilde{v}\partial_y^2\rho+\partial_x\tilde{u}\partial_y^2\rho)\notag\\
&\quad-(1+\lambda){\nu}(\partial_y^2U_s\partial_{x}\rho+\partial_y^2\tilde{u}\partial_{x}\rho+3\partial_y^2\tilde{v}\partial_y\rho+2\partial_{xy}\tilde{u}\partial_y\rho)\notag\\
&\quad-(1+\lambda){\nu}\rho\partial_y^2\div_{x, y}(\tilde{u},\tilde{v})-{\nu}\partial_x\rho\partial_y^2U_s +(1+\lambda){\nu}\partial_y^2g_{\rho}+\div_{x, y}(g_u,g_v),\notag\\
&=-(1+\lambda){\nu}\left[(U_s+\tilde{u})\partial_{x}\partial_y^2\rho+\tilde{v}\partial_y^3\rho\right]+O(1)\nu|\partial_x^2\div_{x,y}(u,v)|\notag\\
&\quad+O(1)\nu^{-\frac{1}{2}}(|\nabla_{x,y}(u,v)|+\nu^{\frac12}|\nabla^2_{x,y}(u,v)|)+O(1)(|\nabla_{x,y}\rho|+\nu^{\frac12}|\nabla^2_{x,y}\rho|)\notag\\
&\quad+O(1)\nu\bigg( \sum_{j=0}^2|\nabla_{x,y}^j\rho||\nabla_{x,y}^{3-j}(\tilde{u},\tilde{v})|\bigg)+O(1)\nu|\partial_{y}^2g_\rho|+O(1)|\div_{x,y}(g_u,g_v)|.
\end{align}

Taking inner product of \eqref{HE8} with $m^{-2}\Delta_{x,y}\rho$, using Cauchy-Schwarz and $W^{2,\infty}$ bounds \eqref{Lw1}-\eqref{Lw3} of $(\tilde{u},\tilde{v})$, we obtain
\begin{align}\label{HE9}
m^{-4}\|\Delta_{x, y}\rho\|_{L^2(\Omega)}^2&\leq -m^{-2}(1+\lambda){\nu}\langle (U_s+\tilde{u})\partial_{x}\partial_y^2\rho+\tilde{v}\partial_y^3\rho, \partial_x^2\rho+\partial_y^2\rho\rangle\notag\\
&\quad +C\nu^2\|\partial_x^2\div_{x, y}(u, v)\|_{L^2(\Omega)}^2+o(1)\|\nabla_{x,y}^2\rho\|_{L^2(\Omega)}^2+C\nu^{-\frac{13}{4}}\|(\rho,u,v)\|_{1}^2\nonumber\\
&\quad+C\nu^{-(\frac92)^-}\|(\tilde{u},\tilde{v})\|_{\mathfrak{X}_*}^2\|(\rho,u,v)\|_2^2 +C\nu^2\|\nabla^2_{x,y}g_{\rho}\|_{L^2(\Omega)}^2+C\|\nabla_{x,y}(g_{u}, g_{v})\|_{L^2(\Omega)}^2.
\end{align}

It remains to control the first term. We separate it into two parts. Firstly, integrating by parts and using boundary conditions $\tilde{u}|_{y=0}=\tilde{v}|_{y=0}=0$, we obtain
\begin{align}
(1+\lambda){\nu}|\langle (U_s+\tilde{u})\partial_{x}\partial_y^2\rho+\tilde{v}\partial_y^3\rho, \partial_y^2\rho\rangle|
&=\frac12(1+\lambda){\nu}\left|\langle \div_{x, y}(\tilde{u},\tilde{v})\partial_y^2\rho, \partial_y^2\rho\rangle\right|\notag\\
&\quad\leq C\nu\|\nabla_{x,y}(\tilde{u}, \tilde{v})\|_{L^\infty(\Omega)}\|\nabla^2_{x,y}\rho\|_{L^2(\Omega)}^2\nonumber\\
&\quad\leq C\nu^{-(\frac{31}{8})^-}\|(\tilde{u},\tilde{v})\|_{\mathfrak{X}_*}\|(\rho,u,v)\|_2^2.\nonumber
\end{align}
For the second part, we integrate by parts again to compute
\begin{align*}
\langle (U_s+\tilde{u})\partial_{ x}\partial_y^2\rho, \partial_x^2\rho\rangle&=-\langle(\partial_yU_s+\partial_y\tilde{u})\partial_{xy}\rho, \partial_x^2\rho\rangle-\langle( U_s+\tilde{u})\partial_{xy}\rho, \partial_{x}^2\partial_y\rho\rangle\\
&=-\langle(\partial_yU_s+\partial_y\tilde{u})\partial_{xy}\rho, \partial_x^2\rho\rangle+\frac{1}{2}\langle \partial_x\tilde{u}\partial_{xy}\rho,\partial_{xy}\rho\rangle.
\end{align*}
Thus, it holds
\begin{align} 
|(1+\lambda){\nu}\langle (U_s+\tilde{u})\partial_{x}\partial_y^2\rho, \partial_x^2\rho\rangle|&\leq C \sqrt{\nu}(1+\sqrt{\nu}\|\nabla_{x,y}(\tilde{u},\tilde{v})\|_{L^\infty(\Omega)})\|\nabla^2_{x,y}\rho\|_{L^2(\Omega)}^2.\nonumber\\
&\leq o(1)\|\nabla_{x,y}^2\rho\|_{L^2(\Omega)}^2+C\nu^{-(\frac{31}{8})^-}\|(\tilde{u},\tilde{v})\|_{\mathfrak{X}_*}\|(\rho,u,v)\|_2^2.\nonumber
\end{align}
Similarly, we have
\begin{align*}
\langle \tilde{v}\partial_{y}^3\rho, \partial_x^2\rho\rangle&=-\langle\partial_y\tilde{v}\partial_{y}^2\rho, \partial_x^2\rho\rangle-\langle \tilde{v}\partial_{y}^2\rho, \partial_{x}^2\partial_y\rho\rangle\\
&=-\langle\partial_y\tilde{v}\partial_{y}^2\rho, \partial_x^2\rho\rangle+\langle\partial_x\tilde{v}\partial_{y}^2\rho, \partial_{xy}\rho\rangle
+\langle \tilde{v}\partial_{x}\partial_y^2\rho, \partial_{xy}\rho\rangle\\
&=-\langle\partial_y\tilde{v}\partial_{y}^2\rho, \partial_x^2\rho\rangle+\langle\partial_x\tilde{v}\partial_{y}^2\rho, \partial_{xy}\rho\rangle
-\frac12\langle \partial_y\tilde{v}\partial_{xy}\rho, \partial_{xy}\rho\rangle,
\end{align*}
which implies
\begin{align}
|(1+\lambda){\nu}\langle \tilde{v}\partial_{y}^3\rho, \partial_x^2\rho\rangle|&\leq C {\nu}\|\nabla_{x,y}(\tilde{u},\tilde{v})\|_{L^\infty(\Omega)}\|\nabla^2_{x,y}\rho\|_{L^2(\Omega)}^2.\nonumber\\
&\leq C\nu^{-(\frac{31}{8})^-}\|(\tilde{u},\tilde{v})\|_{\mathfrak{X}_*}\|(\rho,u,v)\|_2^2.\nonumber
\end{align}
Plugging these bounds back into \eqref{HE9}, we obtain
\begin{align}\label{HE10}
&m^{-2}\|\Delta_{x, y}\rho\|_{L^2(\Omega)}\nonumber\\
&\quad\leq o(1)\|\nabla_{x,y}^2\rho\|_{L^2(\Omega)}+C\nu\|\partial_x^2\div_{x, y}(u, v)\|_{L^2(\Omega)}+C\nu\|\nabla_{x,y}^2g_{\rho}\|_{L^2(\Omega)}+C\|\nabla_{x,y}(g_{u},g_{v})\|_{L^2(\Omega)}\notag\\
&\quad\quad +C\nu^{-(\frac{31}{16})^-}(1+\nu^{-\frac{5}{16}}\|(\tilde{u},\tilde{v})\|_{\mathfrak{X}_*}^{\frac12})\|(\tilde{u},\tilde{v})\|_{\mathfrak{X}_*}^{\frac12}\|(\rho,u,v)\|_2 +C\nu^{-\frac{13}{8}}\|(\rho,u,v)\|_1.
\end{align}

Next we estimate $\|\Delta_{x,y}\nabla_{x,y}\rho\|_{L^2(\Omega)}$. Let $\partial=\partial_x$ or $\partial_y$. Similar to \eqref{HE8}, taking derivative in \eqref{HE7},  we have
\begin{align}\label{HE8-1}
	m^{-2}\Delta_{x, y}\partial\rho&=-(1+\lambda){\nu}\left[(U_s+\tilde{u})\partial_{x}\partial_y^2\partial\rho+\tilde{v}\partial_y^3\partial\rho\right]+ O(1)\nu|\partial_x^2\partial\div_{x,y}(u,v)|\notag\\
	&\quad+O(1)\left[\sum_{j=1}^3\nu^{-\frac{3-j}{2}}\left(|\nabla_{x,y}^j(u,v)|+\nu^{\frac12}|\nabla^j_{x,y}\rho|\right)\right]+O(1)\nu\left(\sum_{j=0}^3|\nabla_{x,y}^j\rho||\nabla_{x,y}^{4-j}(\tilde{u},\tilde{v})|\right)\nonumber\\
	&\quad+O(1)\nu|\nabla_{x,y}^3g_\rho|+O(1)|\nabla_{x,y}^2(g_u,g_v)|.
\end{align}
Applying the same approach to \eqref{HE10},  we take inner product of \eqref{HE8-1} with $m^{-2}\Delta_{x,y}\partial\rho$ to obtain
\begin{align}\label{HE8-2}
		&m^{-2}\|\Delta_{x, y}\partial\rho\|_{L^2(\Omega)}\nonumber\\
		&\quad\leq o(1)\|\nabla_{x,y}^3\rho\|_{L^2(\Omega)}+C\nu\|\partial_x^2\partial \div_{x, y}(u,v)\|_{L^2(\Omega)}+C\nu\|\nabla_{x,y}^3g_\rho\|_{L^2(\Omega)}+C\|\nabla_{x,y}^2(g_u,g_v)\|_{L^2(\Omega)}\nonumber\\
		&\quad\quad+C\nu^{-(\frac{47}{16})^-}(1+\nu^{-\frac{5}{16}}\|(\tilde{u},\tilde{v})\|_{\mathfrak{X}_*}^{\frac{1}{2}})\|(\tilde{u},\tilde{v})\|_{\mathfrak{X}_*}^{\frac{1}{2}}\|(\rho,u,v)\|_{\mathfrak{X}}+C\nu^{-\frac{21}{8}}\|(\rho,u,v)\|_2.
\end{align}

With the bounds on $\Delta_{x,y}\rho,$ we can control
 $\|\nabla_{x,y}^k\div_{x,y}(u,v)\|_{L^2(\Omega)}$ by elliptic estimates. From \eqref{HE7} and \eqref{HE10}, we obtain
\begin{align}\label{HE8-3}
	&\nu\|\Delta_{x,y}\div_{x,y}(u,v)\|_{L^2(\Omega)}\nonumber\\&\quad\leq  Cm^{-2}\|\Delta_{x,y}\rho\|_{L^2(\Omega)}+C\left(\sum_{j=1}^2\nu^{\frac{j-2}{2}}\|\nabla_{x,y}^j(u,v)\|_{L^2(\Omega)}\right)+C\|\nabla_{x,y}\rho\|_{L^2(\Omega)}+C\|\nabla_{x,y}(g_u,g_v)\|_{L^2(\Omega)}\nonumber\\
	&\quad\leq o(1)\|\nabla_{x,y}^2\rho\|_{L^2(\Omega)}+ C\nu\|\partial_x^2\div_{x, y}(u, v)\|_{L^2(\Omega)}+C\nu\|\nabla_{x,y}^2g_{\rho}\|_{L^2(\Omega)}+C\|\nabla_{x,y}(g_{u},g_{v})\|_{L^2(\Omega)}\notag\\
	&\qquad\quad +C\nu^{-(\frac{31}{16})^-}(1+\nu^{-\frac{5}{16}}\|(\tilde{u},\tilde{v})\|_{\mathfrak{X}_*}^{\frac12})\|(\tilde{u},\tilde{v})\|_{\mathfrak{X}_*}^{\frac12}\|(\rho,u,v)\|_2 +C\nu^{-\frac{13}{8}}\|(\rho,u,v)\|_1.
\end{align}
Then differentiating \eqref{HE7} and using the same argument as \eqref{HE8-3}, we get
\begin{align}\label{HE8-4}
	&\nu\|\Delta_{x,y}\partial\div_{x,y}(u,v)\|_{L^2(\Omega)}\nonumber\\&\quad\leq  Cm^{-2}\|\Delta_{x,y}\partial\rho\|_{L^2(\Omega)}+\|\nabla_{x,y}^2(g_u,g_v)\|_{L^2(\Omega)}+C\left(\sum_{j=1}^3\nu^{\frac{j-3}{2}}\|\nabla_{x,y}^j(u,v)\|_{L^2(\Omega)}+\sum_{j=1}^2\nu^{\frac{j-2}{2}}\|\nabla_{x,y}^j\rho\|_{L^2(\Omega)}\right)\nonumber\\
	&\quad\leq o(1)\|\nabla_{x,y}^3\rho\|_{L^2(\Omega)}+C\nu\|\partial_x^2\partial \div_{x, y}(u,v)\|_{L^2(\Omega)}+C\nu\|\nabla_{x,y}^3g_\rho\|_{L^2(\Omega)}+C\|\nabla_{x,y}^2(g_u,g_v)\|_{L^2(\Omega)}\nonumber\\
	&\qquad\quad+C\nu^{-(\frac{47}{16})^-}(1+\nu^{-\frac{5}{16}}\|(\tilde{u},\tilde{v})\|_{\mathfrak{X}_*}^{\frac{1}{2}})\|(\tilde{u},\tilde{v})\|_{\mathfrak{X}_*}^{\frac{1}{2}}\|(\rho,u,v)\|_{\mathfrak{X}}+C\nu^{-\frac{21}{8}}\|(\rho,u,v)\|_2.
\end{align}

From $\eqref{7.1}_1$, we can derive the boundary value of $\div_{x,y}(u,v)$:
$$
\begin{aligned}\div_{x,y}(u,v)|_{y=0}&=g_\rho|_{y=0}-\rho\div_{x,y}(\tilde{u},\tilde{v})|_{y=0}-\tilde{u}\partial_{x}\rho|_{y=0}-\tilde{v}\partial_y\rho|_{y=0}-U_s\partial_x\rho|_{y=0}\\
&=0,
\end{aligned}
$$
because $g_{\rho}|_{y=0}=0$ and $(\tilde{u},\tilde{v})\in \mathfrak{X}_*$. Then by the following
  $H^2$-regularity  estimate of the Laplacian operator:
$$\|\nabla_{x,y}^2f\|_{L^2(\Omega)}\leq C\|\Delta_{x,y}f\|_{L^2(\Omega)}, ~\text{for any } f\in H^2(\Omega)\cap H^1_{0}(\Omega),
$$  
 we get
\begin{align}\label{HE10-1}
&\nu\|\nabla^2_{x,y}\div_{x, y}(u, v)\|_{L^2(\Omega)}\leq C\nu\|\Delta_{x,y}\div(x,y)(u,v)\|_{L^2(\Omega)}\nonumber\\
&\qquad\leq o(1)\|\nabla_{x,y}^2\rho\|_{L^2(\Omega)}+C\nu\|\partial_x^2\div_{x, y}(u, v)\|_{L^2(\Omega)}+C\nu\|\nabla_{x,y}^2g_{\rho}\|_{L^2(\Omega)}+C\|\nabla_{x,y}(g_{u},g_{v})\|_{L^2(\Omega)}\notag\\
&\qquad\quad +C\nu^{-(\frac{31}{16})^-}(1+\nu^{-\frac{5}{16}}\|(\tilde{u},\tilde{v})\|_{\mathfrak{X}_*}^{\frac12})\|(\tilde{u},\tilde{v})\|_{\mathfrak{X}_*}^{\frac12}\|(\rho,u,v)\|_2 +C\nu^{-\frac{13}{8}}\|(\rho,u,v)\|_1.
\end{align}
Here we have used the bound \eqref{HE8-3} in the last inequality. Similarly, using \eqref{HE8-4} with $\partial=\partial_x$, we deduce
\begin{align}
	&\nu\|\nabla^2_{x,y}\partial_x\div_{x, y}(u, v)\|_{L^2(\Omega)}\leq C\nu\|\Delta_{x,y}\partial_x\div(x,y)(u,v)\|_{L^2(\Omega)}\nonumber\\
	& \qquad\leq o(1)\|\nabla_{x,y}^3\rho\|_{L^2(\Omega)}+C\nu\|\partial_x^3 \div_{x, y}(u,v)\|_{L^2(\Omega)}+C\nu\|\nabla_{x,y}^3g_\rho\|_{L^2(\Omega)}+C\|\nabla_{x,y}^2(g_u,g_v)\|_{L^2(\Omega)}\nonumber\\
	&\qquad\quad+C\nu^{-(\frac{47}{16})^-}(1+\nu^{-\frac{5}{16}}\|(\tilde{u},\tilde{v})\|_{\mathfrak{X}_*}^{\frac{1}{2}})\|(\tilde{u},\tilde{v})\|_{\mathfrak{X}_*}^{\frac{1}{2}}\|(\rho,u,v)\|_{\mathfrak{X}}+C\nu^{-\frac{21}{8}}\|(\rho,u,v)\|_2.\label{HE10-2}
\end{align}
Moreover, by $H^3$-estimate of Laplacian: $\|\nabla_{x,y}^3f\|_{L^2(\Omega)}\leq C\|\Delta_{x,y}f\|_{H^1(\Omega)}+C\|f\|_{L^2(\Omega)}$ (cf.  \cite[P323]{E}), and the bound \eqref{HE8-4} we can obtain
\begin{align}
&	\nu\|\nabla^2_{x,y}\partial_y\div_{x, y}(u, v)\|_{L^2(\Omega)}\leq C\nu\|\Delta_{x,y}\div_{x,y}(u,v)\|_{H^1(\Omega)}+C\nu\|\div_{x,y}(u,v)\|_{L^2(\Omega)}\nonumber\\
&\quad \leq o(1)\|\nabla_{x,y}^3\rho\|_{L^2(\Omega)}+C\nu\|\partial_x^2\nabla_{x,y} \div_{x, y}(u,v)\|_{L^2(\Omega)}+C\nu\|\nabla_{x,y}^3g_\rho\|_{L^2(\Omega)}+C\|\nabla_{x,y}^2(g_u,g_v)\|_{L^2(\Omega)}\nonumber\\
&\qquad\quad+C\nu^{-(\frac{47}{16})^-}(1+\nu^{-\frac{5}{16}}\|(\tilde{u},\tilde{v})\|_{\mathfrak{X}_*}^{\frac{1}{2}})\|(\tilde{u},\tilde{v})\|_{\mathfrak{X}_*}^{\frac{1}{2}}\|(\rho,u,v)\|_{\mathfrak{X}}+C\nu^{-\frac{21}{8}}\|(\rho,u,v)\|_2.\label{HE10-4}
\end{align}
Note that the second term on the right hand side of \eqref{HE10-4}  which contains at least two order $x$-derivative of $\div_{x,y}(u,v)$, has been bounded in \eqref{HE10-2}.  Therefore, we arrive at:
\begin{align}
&\nu\|\nabla^3_{x,y}\div_{x, y}(u, v)\|_{L^2(\Omega)}\nonumber\\
&\qquad\leq o(1)\|\nabla_{x,y}^3\rho\|_{L^2(\Omega)}+ C\nu\|\partial_x^3 \div_{x, y}(u,v)\|_{L^2(\Omega)}+C\nu\|\nabla_{x,y}^3g_\rho\|_{L^2(\Omega)}+C\|\nabla_{x,y}^2(g_u,g_v)\|_{L^2(\Omega)}\nonumber\\
&\qquad\quad+C\nu^{-(\frac{47}{16})^-}(1+\nu^{-\frac{5}{16}}\|(\tilde{u},\tilde{v})\|_{\mathfrak{X}_*}^{\frac{1}{2}})\|(\tilde{u},\tilde{v})\|_{\mathfrak{X}_*}^{\frac{1}{2}}\|(\rho,u,v)\|_{\mathfrak{X}}+C\nu^{-\frac{21}{8}}\|(\rho,u,v)\|_2.
\label{HE10-3}
\end{align}

Combining \eqref{HE10} with \eqref{HE10-1}, we deduce \eqref{HE0-1}, while the higher order estimate \eqref{HE0-2} can be obtain by combining \eqref{HE8-2} and \eqref{HE10-3}. The proof of Lemma \ref{lmHF2} is complete.
\end{proof}

In next lemma, we bound the full derivatives of $\rho$, $u$, and $v$ in terms of $\div_{x,y}(u,v)$. 
\begin{lemma}\label{lmHF3}
	Let $k=2$ or $3$. If $(g_u,g_v)\in H^{k-1}(\Omega)^2$, $g_\rho|_{y=0}=0$ and $\iint_{\Omega}g_\rho(x,y)dx dy=0$, we have
\begin{align}\label{HF3}
	&{\nu}\|\nabla_{x,y}^{k+1}(u,v)\|_{L^2(\Omega)}+m^{-2}\|\nabla^{k}_{x,y}\rho\|_{L^2(\Omega)}\nonumber\\
	&\qquad\leq C{\nu}\|\nabla^k_{x,y}\div_{x,y}(u,v)\|_{L^2(\Omega)}+C\nu^{-k+\frac38}\|(\rho,u,v)\|_{k-1}+C\|(g_u,g_v)\|_{H^{k-1}(\Omega)}.
\end{align}
\end{lemma}
\begin{proof}
 We rewrite momentum equations \eqref{7.1} as 
\begin{equation}
	\left\{
\begin{aligned}
	-{\nu}\Delta_{x,y}u+m^{-2}\partial_x\rho&=g_u-U_s\partial_xu-v\partial_yU_s\\
	&\quad+\lambda{\nu}\partial_x\div_{x, y}(u,v)-{\nu}\rho\partial_y^2U_s\eqdef w_u,\\
	-{\nu}\Delta_{x,y}v+m^{-2}\partial_y\rho&=g_v-U_s\partial_xv+\lambda{\nu}\partial_y\div_{x, y}(u,v)\eqdef w_v.\label{St}
\end{aligned}
\right.
\end{equation}
Note that for $k=2$ or $3$, 
\begin{align}
	\|(w_u,w_v)\|_{H^{k-1}(\Omega)}\leq& C\|(g_u,g_v)\|_{H^{k-1}(\Omega)}+C\nu\|\nabla_{x,y}\div_{x,y}(u,v)\|_{H^{k-1}(\Omega)}\nonumber\\
	&+C\|U_s\partial_x(u, v)\|_{H^{k-1}(\Omega)}+C\nu\|\rho\partial_y^2U_s\|_{H^{k-1}(\Omega)}+C\|v\partial_yU_s\|_{H^{k-1}(\Omega)}\nonumber\\
	\leq& C\nu\|\nabla_{x,y}^k\div_{x,y}(u,v)\|_{L^2(\Omega)}+C\sum_{j=1}^k\|\partial_y^{k-j}U_s\|_{L^\infty}\|\nabla^j_{x,y}(u,v)\|_{L^2(\Omega)}\nonumber\\
	&+C\sum_{j=0}^{k-1}\nu\|\partial_y^{k+1-j}U_s\|_{L^\infty}\|\nabla^j_{x,y}\rho\|_{L^2(\Omega)} +C\sum_{j=0}^{k-1}\|\partial_y^{k-j}U_s\|_{L^\infty}\|\nabla_{x,y}^jv\|_{L^2(\Omega)}\nonumber\\
	&+C\|y\partial_y^{k}U_s\|_{L^\infty}\|vy^{-1}\|_{L^2(\Omega)}+C\|(g_u,g_v)\|_{H^{k-1}(\Omega)}\nonumber\\
\leq& C\nu\|\nabla_{x,y}^k\div_{x,y}(u,v)\|_{L^2(\Omega)}+C\|(g_u,g_v)\|_{H^{k-1}(\Omega)}\nonumber\\ &+C\sum_{j=1}^k\nu^{\frac{j-k}{2}}\left(\|\nabla_{x,y}^j(u,v)\|_{L^2(\Omega)}+\|\nabla_{x,y}^{j-1}\rho\|_{L^2(\Omega)}\right)\nonumber\\
\leq& C\nu\|\nabla_{x,y}^k\div_{x,y}(u,v)\|_{L^2(\Omega)}+ C\nu^{-k+\frac38}\|(\rho,u,v)\|_{k-1}+C\|(g_u,g_v)\|_{H^{k-1}(\Omega)}.\nonumber
\end{align}

Then applying Lemma \ref{lem9.1} to $(\rho,u,v)$, we obtain, for $k=2$ or $3$, that
\begin{align}
	&{\nu}\|\nabla_{x,y}^{k+1}(u,v)\|_{L^2(\Omega)}+m^{-2}\|\nabla^{k}_{x,y}\rho\|_{L^2(\Omega)}\nonumber\\
	&\quad\leq C{\nu}\|\div_{x,y}(u,v)\|_{H^k(\Omega)}+Cm^{-2}\|\rho\|_{H^{k-1}(\Omega)}+C\nu\|\nabla_{x,y}(u,v)\|_{H^{k-1}(\Omega)}+C\|(w_u,w_v)\|_{H^{k-1}(\Omega)}\nonumber\\
	&\quad\leq C{\nu}\|\nabla^k_{x,y}\div_{x,y}(u,v)\|_{L^2(\Omega)}+C\nu^{-k+\frac38}\|(\rho,u,v)\|_{k-1}+C\|(g_u,g_v)\|_{H^{k-1}(\Omega)},\nonumber
\end{align}
which is \eqref{HF3}. The proof of Lemma \ref{lmHF3} is complete.
\end{proof}
\bigbreak
{\bf Proof of Proposition \ref{L7.1}:} Subsitituting the bound \eqref{HF1} of $\|\partial_x^2\div_{x, y}(u,v)\|_{L^2(\Omega)}$ into \eqref{HE0-1}, we get the following bounds on the full derivatives of $\div_{x,y}(u,v)$:
\begin{align}
{\nu}\|\nabla^2_{x,y}\div_{x, y}(u, v)\|_{L^2(\Omega)}\leq& o(1)\|\nabla_{x,y}^2\rho\|_{L^2(\Omega)}+C\nu\|\nabla_{x,y}^2g_{\rho}\|_{L^2(\Omega)}+C\|\nabla_{x,y}(g_{u},g_{v})\|_{L^2(\Omega)}\notag\\
 &+C\nu^{-(\frac{31}{16})^-}(1+\nu^{-\frac{5}{16}}\|(\tilde{u},\tilde{v})\|_{\mathfrak{X}_*}^{\frac12})\|(\tilde{u},\tilde{v})\|_{\mathfrak{X}_*}^{\frac12}\|(\rho,u,v)\|_2 +C\nu^{-\frac{13}{8}}\|(\rho,u,v)\|_1.\nonumber
\end{align}
Then plugging the above bound into \eqref{HF3}, we can deduce \eqref{HE0}. The higher order estimate \eqref{HF0} can be obtained similarly. The proof of Proposition \ref{L7.1} is complete.

\begin{corollary}\label{corl}
	Let $m\in (0, 1)$. Suppose that $L\in (0,L_0)$, where $L_0$ is given in Theorem \ref{thmno}. For any given velocity field $(\tilde{u},\tilde{v})\in \mathfrak{X}_*$ satisfying $\|(\tilde{u},\tilde{v})\|_{\mathfrak{X}_*}\leq \nu^{\frac{9}{8}+}$, and any inhomogeneous source term $(g_\rho,g_u,g_v)\in H^2(\Omega)\times H^1(\Omega)^2$,  satisfying $g_\rho|_{y=0}=0$, $\iint_{\Omega}g_\rho(x,y)d x dy=0,$ and $[[g_\rho,g_u,g_v]]_{2}<\infty$, the system \eqref{7.1} admits a unique solution $(\rho, u, v)$ that satisfies the zero-mass condition
	\begin{align}
		\iint_{\Omega}\rho(x,y)d x dy=0,\label{zm}
	\end{align} 
the boundary condition
\begin{align}
	\div_{x,y}(u,v)|_{y=0}=0,\label{db}
\end{align}
and the bound
	\begin{align}
	\|(\rho,u,v)\|_{2}\leq C[[(g_\rho,g_u,g_v)]]_{2}.\label{LS}
	\end{align}
Moreover, if $g_\rho\in H^3(\Omega)$ and $(g_u,g_v)\in H^2(\Omega)^2$, we have
\begin{align}
	\|(\rho,u,v)\|_{\mathfrak{X}}\leq C[[(g_\rho,g_u,g_v)]]_{3}.\label{LH}
\end{align}
\end{corollary}
\begin{proof}
	Combining the low order estimate \eqref{low} in Proposition \ref{propl} with the high order estimate \eqref{HE0}, we obtain
		\begin{align}
		\|(\rho,u,v)\|_{2}\leq C[[g_\rho,g_u,g_v]]_{2}+C\left(\nu^{-(\frac{5}{16})-}\|(\tilde{u},\tilde{v})\|_{\mathfrak{X}_*}^{\frac12}+\nu^{-(\frac{9}{8})-}\|(\tilde{u},\tilde{v})\|_{\mathfrak{X}_*}\right)\|(\rho,u,v)\|_{2}.\nonumber
	\end{align}
Thus, by taking $(\tilde{u},\tilde{v})$ such that $\|(\tilde{u},\tilde{v})\|_{2}\lesssim \nu^{\frac{9}{8}+}$, we can absorb the second term on the right hand side by the left hand side. Therefore, the estimate \eqref{LS} follows. The bound \eqref{LH} can be obtained similarly. The proof of Corollary \ref{corl} is complete.
\end{proof}
\subsection{Proof of Theorem \ref{T1.1}}\label{S8.3}
To construct the solution to  nonlinear system \eqref{1.5}, we introduce the following iteration scheme
\begin{align*}
\begin{cases}
U_s\partial_x\rho^{i+1}+\div_{x,y}(u^{i+1},v^{i+1})+\partial_x(u^i\rho^{i+1})+\partial_y(v^i\rho^{i+1})=0,\\
U_s\partial_xu^{i+1}+v^{i+1}\partial_yU_s+m^{-2}\partial_x\rho^{i+1}-\nu\Delta_{x,y} u^{i+1}-\lambda\nu\partial_x\div_{x,y}(u^{i+1},v^{i+1})+\nu\rho^{i+1}\partial_y^2U_s={N}_u(\rho^i,u^i,v^i),\\
U_s\partial_xv^{i+1}+m^{-2}\partial_y\rho^{i+1}-{\nu}\Delta_{x,y} v^{i+1}-\lambda\nu\partial_y\div_{x, y}(u^{i+1}, v^{i+1})={N}_v(\rho^i,u^i,v^i),\\
u^{i+1}|_{y=0}=v^{i+1}|_{y=0}=0,\\
(\rho^0,u^0,v^0)=(0,0,0),
\end{cases}
\end{align*}
where nonlinear operators $N_u$ and $N_v$ are defined in \eqref{N1} and \eqref{N2} respectively. Now we show the convegence of  iteration by the following two steps:

\underline{\it Step 1. Uniform bounds of $(\rho^i,u^i,v^i)$ in $\|\cdot\|_{\mathfrak{X}}$ norm.} From Corollary \ref{corl}, we have 
\begin{align}
	\|(\rho^{i+1},u^{i+1},v^{i+1})\|_{\mathfrak{X}}\leq C[[(0,{N}_u(\rho^i,u^i,v^i),{N}_v(\rho^i,u^i,v^i))]]_3,\label{P1}
\end{align}
where the norm $[[~\cdot~]]_3$ is defined in \eqref{gnh}. Now we bound the right hand side of \eqref{P1}. Firstly we control the zero-mode part in the norm $[[~\cdot~]]_3$. Recall \eqref{N1} and \eqref{N2} for the definition of nonlinear term  $(N_u,N_v)$. Integrating \eqref{N1} and using $\mathcal{P}_0F_{\rm ext,1}=0$, we obtain
\begin{align}
	\partial_y^{-1}N_{u,0}(\rho,u,v)=-\mathcal{P}_0[(1+\rho)uv]-\partial_y^{-1}[\partial_yU_s\mathcal{P}_0(\rho v)]+\partial_y^{-1}\mathcal{P}_0(\rho F_{{\rm ext},1}).\label{Po0}
\end{align}
For the last two terms on the right hand side of \eqref{Po0}, we can obtain the following pointwise estimate:
\begin{align}
\left| \partial_y^{-1}[\partial_yU_s\mathcal{P}_0(\rho v)]\right|&\leq C\nu^{-\frac12}\|\mathcal{P}_0(\rho v)\|_{L^\infty_y}\int_y^\infty\left(1+\left|\frac{y'}{\sqrt{\nu}}\right|\right)^{-s} dy'\nonumber\\
&\leq C\|\mathcal{P}_0(\rho v)\|_{L^\infty_y}(1+Y)^{-s+1},\label{Po1}
\end{align}
where $Y=y/\sqrt{\nu}$. Similarly, it holds that
\begin{align}
	\left| \partial_y^{-1}\mathcal{P}_0(\rho F_{{\rm ext},1})\right|&\leq C(1+y)^{-s}\|\rho\|_{L^2(\Omega)}\|(1+y)^sF_{{\rm ext},1}\|_{L^2}\nonumber\\
	&\leq C(1+y)^{-s}\|\rho\|_{L^2(\Omega)}\|(F_{{\rm ext},1},F_{{\rm ext},2})\|_{w}.
	 \label{Po2}
\end{align}
Plugging \eqref{Po1} and \eqref{Po2} into \eqref{Po0}, we obtain
\begin{align}
	|\partial_y^{-1}N_{u,0}(\rho,u,v)|\leq& \left|\mathcal{P}_0[(1+\rho)uv]\right|+C\|\mathcal{P}_0(\rho v)\|_{L^\infty}(1+Y)^{-s+1}\nonumber\\
	&+C(1+y)^{-s}\|\rho\|_{L^2(\Omega)}\|(F_{{\rm ext},1},F_{{\rm ext},2})\|_{w}.\label{Po3}
\end{align}
Similarly, we have
\begin{align}
	|\partial_y^{-1}N_{v,0}(\rho,u,v)|\leq& C\left|\mathcal{P}_0[(1+\rho)v^2]\right|+\left|\mathcal{P}_0[P(1+\rho)-P'(1)\rho]\right|\nonumber\\
	&+C(1+y)^{-s}\|\rho\|_{L^2(\Omega)}\|\|(F_{{\rm ext},1},F_{{\rm ext},2})\|_{w}.\label{Po4}
\end{align}
Then from \eqref{Po3} and \eqref{Po4}, and using \eqref{A9.1}-\eqref{A9.3}, we obtain
\begin{align}
	\|\partial_y^{-1}N_{u,0}(\rho,u,v)\|_{L^1_y}&\leq C\|\mathcal{P}_0[(1+\rho)uv]\|_{L^1_y}+C\nu^{\frac12}\|\mathcal{P}_0(\rho v)\|_{L^\infty_y}+C\|\rho\|_{L^2(\Omega)}\|(F_{{\rm ext},1},F_{{\rm ext},2})\|_{w}\nonumber\\
	&\leq C\|(\rho,u,v)\|_{1}^2\left(1+\nu^{-\frac14}\|(\rho,u,v)\|_{1}\right)+C\|(\rho,u,v)\|_{1}\|(F_{{\rm ext},1},F_{{\rm ext},2})\|_{w}\label{Po5},\\
	\|\partial_y^{-1}N_{u,0}(\rho,u,v)\|_{L^2_y}&\leq C\|\mathcal{P}_0[(1+\rho)uv]\|_{L^2_y}+C\nu^{\frac14}\|\mathcal{P}_0(\rho v)\|_{L^\infty_y}+C\|\rho\|_{L^2(\Omega)}\|(F_{{\rm ext},1},F_{{\rm ext},2})\|_{w}\nonumber\\
	&\leq C\nu^{-\frac{1}{4}}\|(\rho,u,v)\|_{1}^2\left(1+\nu^{-\frac14}\|(\rho,u,v)\|_{1}\right)+C\|(\rho,u,v)\|_{1}\|(F_{{\rm ext},1},F_{{\rm ext},2})\|_{w}\label{Po6},\\
	\|\partial_y^{-1}N_{v,0}(\rho,u,v)\|_{L^2_y}&\leq C\|\mathcal{P}_0[(1+\rho)v^2]\|_{L^2_y}+\|\rho\|_{L^\infty(\Omega)}\|\rho\|_{L^2(\Omega)}+C\|\rho\|_{L^2(\Omega)}\|(F_{{\rm ext},1},F_{{\rm ext},2})\|_{w}\nonumber\\
	&\leq C\nu^{-\frac{1}{4}}\|(\rho,u,v)\|_{1}^2\left(1+\nu^{-\frac14}\|(\rho,u,v)\|_{1}\right)+C\|(\rho,u,v)\|_{1}\|(F_{{\rm ext},1},F_{{\rm ext},2})\|_{w}\label{Po7},\\	\|\partial_y^{-1}N_{v,0}(\rho,u,v)\|_{L^\infty_y}&\leq \|\mathcal{P}_0[(1+\rho)v^2]\|_{L^\infty_y}+\|\rho\|_{L^\infty(\Omega)}^2+C\|\rho\|_{L^2(\Omega)}\|(F_{{\rm ext},1},F_{{\rm ext},2})\|_{w}\nonumber\\
	&\leq C\nu^{-\frac{1}{2}}\|(\rho,u,v)\|_{1}^2\left(1+\nu^{-\frac{1}{2}}\|(\rho,u,v)\|_{1}\right)+C\|(\rho,u,v)\|_{1}\|(F_{{\rm ext},1},F_{{\rm ext},2})\|_{w}\label{Po8}.
\end{align}
 Moreover, we obtain the following pointwise estimate:
\begin{equation}
\begin{aligned}
	\left|(N_u,N_v)\right|\leq& |(u,v)|^2|\nabla_{x,y}\rho|+|(\rho,u,v)||\nabla_{x,y}(\rho,u,v)|+\nu^{-\frac12}|\rho v|+|(F_{{\rm ext},1},F_{{\rm ext},2})|,\\
	\left|\nabla_{x,y}(N_u,N_v)\right|\leq& |(u,v)|^2|\nabla_{x,y}^2\rho|+|(\rho,u,v)|(|\nabla^2_{x,y}(\rho,u,v)|+|\nabla_{x,y}(\rho,u,v)|^2)\\
	&+|\nabla_{x,y}(\rho,u,v)|^2+\nu^{-\frac12}|(\rho,u,v)||\nabla_{x,y}(\rho,u,v)|\\
	&+\nu^{-1}|\rho v|+
	|\nabla_{x,y}(F_{{\rm ext},1},F_{{\rm ext},2})|+|\nabla_{x,y}\rho||(F_{{\rm ext},1},F_{{\rm ext},2})|,\\
	\left|\nabla_{x,y}^2(N_u,N_v)\right|\leq& |(u,v)|^2|\nabla_{x,y}^3\rho|+|(\rho,u,v)|(|\nabla^3_{x,y}(\rho,u,v)|+|\nabla_{x,y}(\rho,u,v)||\nabla^2_{x,y}(\rho,u,v)|)\\
	&+|\nabla_{x,y}(\rho,u,v)|^3+|\nabla_{x,y}(\rho,u,v)||\nabla^2_{x,y}(\rho,u,v)|+\nu^{-1}|(\rho,u,v)||\nabla_{x,y}(\rho,u,v)|\\
	&+\nu^{-\frac12}|(\nabla^2_{x,y}(\rho,u,v)||\rho,u,v|+|\nabla_{x,y}(\rho,u,v)|^2)+\nu^{-\frac{3}{2}}|\rho v|\\
	&+
	|\nabla_{x,y}^2(F_{{\rm ext},1},F_{{\rm ext},2})|+|\nabla_{x,y}\rho||\nabla_{x,y}(F_{{\rm ext},1},F_{{\rm ext},2})|+|\nabla^2_{x,y}\rho||(F_{{\rm ext},1},F_{{\rm ext},2})|.\label{Po8-1}
\end{aligned}
\end{equation}
Then from \eqref{Po8-1} and $W^{2,\infty}$-estimates \eqref{Lw1}-\eqref{Lw3}, we get
\begin{align}
	\|(N_u,N_v)\|_{L^2(\Omega)}\leq& C\|(\rho,u,v)\|_{L^\infty(\Omega)}\bigg( \|\nabla_{x,y}(\rho,u,v)\|_{L^2(\Omega)}+\|(\rho,u,v)\|_{L^\infty(\Omega)}\|\nabla_{x,y}\rho\|_{L^2(\Omega)}\bigg)\nonumber\\
	&+C\nu^{-\frac12}\|\rho\|_{L^\infty(\Omega)}\|v\|_{L^2(\Omega)}+C\|(F_{{\rm ext},1},F_{{\rm ext},2})\|_{L^2(\Omega)},\nonumber\\
	\leq& C\nu^{-1^-}\|(\rho,u,v)\|_{2}^2\left(1+\nu^{-\frac{1}{2}}\|(\rho,u,v)\|_{2}\right)+C\|(F_{{\rm ext},1},F_{{\rm ext},2})\|_{L^2(\Omega)}.\label{P09}
	\end{align}
Similarly, we can deduce
\begin{align}
		\|\nabla_{x,y}(N_u,N_v)\|_{L^2(\Omega)}\leq&
	  \nu^{-(\frac{17}{8})^-}\|(\rho,u,v)\|_{2}^2\left(1+\nu^{-\frac{1}{2}}\|(\rho,u,v)\|_{2}\right)\nonumber\\
	  &+C\|\nabla_{x,y}\rho\|_{L^2(\Omega)}\|(F_{{\rm ext},1},F_{{\rm ext},2})\|_{L^\infty(\Omega)}+C\|\nabla_{x,y}(F_{{\rm ext},1},F_{{\rm ext},2})\|_{L^2(\Omega)}\nonumber\\
	  \leq & \nu^{-(\frac{17}{8})^-}\|(\rho,u,v)\|_{2}^2\left(1+\nu^{-\frac{1}{2}}\|(\rho,u,v)\|_{2}\right)\nonumber\\
	  &+C\nu^{-\frac{13}{8}}\left(1+\nu^{-\frac12}\|(\rho,u,v)\|_{2}\right)\|(F_{\rm ext,1},F_{\rm ext,2})\|_{w},
	  \label{P10}
	 \end{align}
 and
 \begin{align}
	 	\|\nabla_{x,y}^2(N_u,N_v)\|_{L^2(\Omega)}\leq&
	 \nu^{-(\frac{13}{4})^-}\|(\rho,u,v)\|_{\mathfrak{X}}^2\left(1+\nu^{-\frac{1}{2}}\|(\rho,u,v)\|_{\mathfrak{X}}\right)+C\|\nabla_{x,y}^2\rho\|_{L^2(\Omega)}\|(F_{{\rm ext},1},F_{{\rm ext},2})\|_{L^\infty(\Omega)}\nonumber\\
	 &+C\|\nabla_{x,y}\rho\|_{L^\infty(\Omega)}\|\nabla_{x,y}(F_{{\rm ext},1},F_{{\rm ext},2})\|_{L^2(\Omega)}+C\|\nabla_{x,y}^2(F_{{\rm ext},1},F_{{\rm ext},2})\|_{L^2(\Omega)}\nonumber\\
	 \leq&\nu^{-(\frac{13}{4})^-}\|(\rho,u,v)\|_{\mathfrak{X}}^2\left(1+\nu^{-\frac{1}{2}}\|(\rho,u,v)\|_{\mathfrak{X}}\right)\nonumber\\
	 &+\nu^{-\frac{13}{8}}\left(1+\nu^{-\frac58}\|(\rho,u,v)\|_{2}\right)\|(F_{\rm ext,1},F_{\rm ext,2})\|_{w}.
	 \label{P10-1}
\end{align}
Here we have used $L^\infty$ bound \eqref{Lw4} on $(F_{\rm ext,1},F_{\rm ext,2})$ in the last inequalities of \eqref{P10} and \eqref{P10-1}.

Plugging estimates \eqref{Po5}-\eqref{P10-1} into \eqref{P1}, we deduce that
\begin{align}
	\|(\rho^{i+1},u^{i+1},v^{i+1})\|_{\mathfrak{X}}\leq& C\nu^{-(\frac{9}{8})^-}	\|(\rho^{i},u^{i},v^{i})\|_{\mathfrak{X}}^2\left(1+\nu^{-\frac{1}{2}}	\|(\rho^{i},u^{i},v^{i})\|_{\mathfrak{X}}\right)\nonumber\\
	&+C\|(F_{{\rm ext},1},F_{{\rm ext},2})\|_{w}(1+\nu^{-1}\|(\rho^i,u^i,v^i)\|_{\mathfrak{X}}).\label{P11}
\end{align}
Therefore, if $\|(F_{{\rm ext},1} F_{{\rm ext},2})\|_{w}\ll \nu^{(\frac{9}{8})^+}$, we can inductively derive from \eqref{P11} the following uniform bounds on $(\rho^i,u^i,v^i):$
\begin{align}
	\|(\rho^{i},u^{i},v^{i})\|_{\mathfrak{X}}\leq C\|(F_{{\rm ext},1} F_{{\rm ext},2})\|_{w}\leq \nu^{(\frac{9}{8})^+}.\label{P12}
\end{align}

\underline{\it Step 2. Contraction in $\|\cdot\|_{2}$ norm.} We introduce the difference $(\rho^i_d,u^i_d,v^i_d)\eqdef(\rho^{i}-\rho^{i-1},u^i-u^{i-1},v^{i}-v^{i-1})$. Then the equations of $(\rho^i_d,u^i_d,v^i_d)$ reads
\begin{align*}
	\begin{cases}
		U_s\partial_x\rho^{i+1}_d+\div_{x,y}(u^{i+1}_d,v^{i+1}_d)+\partial_x(u^i\rho^{i+1}_d)+\partial_y(v^i\rho^{i+1}_d)\nonumber\\
		\qquad\qquad\qquad=-\partial_x(\rho^iu^i_d)-\partial_y(\rho^iv^i_d)\eqdef g_{\rho,d}^i,\\
		U_s\partial_xu^{i+1}_d+v^{i+1}_d\partial_yU_s+m^{-2}\partial_x\rho_d^{i+1}-\nu\Delta_{x,y} u_d^{i+1}-\lambda\nu\partial_x\div_{x,y}(u_d^{i+1},v_d^{i+1})+\nu\rho_d^{i+1}\partial_y^2U_s\\
		\qquad\qquad\qquad=\mathcal{N}_u(\rho^{i},u^i,v^i)-\mathcal{N}_u(\rho^{i-1},u^{i-1},v^{i-1})\eqdef g_{u,d}^i,\\
		U_s\partial_xv_d^{i+1}+m^{-2}\partial_y\rho_d^{i+1}-{\nu}\Delta_{x,y} v_d^{i+1}-\lambda\nu\partial_y\div_{x, y}(u_d^{i+1}, v_d^{i+1})\\
		\qquad\qquad\qquad=\mathcal{N}_v(\rho^i,u^i,v^i)-\mathcal{N}_v(\rho^{i-1},u^{i-1},v^{i-1})\eqdef g^i_{v,d},\\
		u^{i+1}_d|_{y=0}=v^{i+1}_d|_{y=0}=0,\\
		(\rho^0_d,u^0_d,v^0_d)=(0,0,0).
	\end{cases}
\end{align*}
By an inductive argument, we have 
\begin{align}
	\iint_{\Omega}g_{\rho,d}^idxdy=-\iint_\Omega \partial_x(\rho^iu^i_d)+\partial_y(\rho^iv^i_d)d x d y=0,\nonumber
\end{align}
$\div_{x,y}(u^i_d,v^i_d)|_{y=0}=0$, and
\begin{align}
[\partial_x(\rho^iu^i_d)+\partial_y(\rho^iv^i_d)]\big|_{y=0}=(u^{i}_d,v^{i}_d)\cdot \nabla_{x,y}\rho^i|_{y=0}+\rho^i\div_{x,y}(u^{i}_d,v^{i}_d)|_{y=0}=0.\nonumber
\end{align}
Then by \eqref{LS}, we have
\begin{align}
	\|(\rho^{i+1}_d,u_{d}^{i+1},v_{d}^{i+1})\|_{2}\leq C[[(g_{\rho,d}^i,g^{i}_{u,d},g^{i}_{v,d})]]_2\nonumber.
\end{align}
For $g_{d,0}^i$, we can obtain the following pointwise estimate
\begin{equation}
\begin{aligned}
	|g_{\rho,d}^i|&\leq |\rho^i||\nabla_{x,y}(u_d^i,v_d^i)|+|\nabla_{x,y}\rho^i||(u_d^i,v_d^i)|,\\
	|\nabla_{x,y}g_{\rho,d}^i|&\leq |\rho^i||\nabla_{x,y}^2(u_d^i,v_d^i)|+|\nabla_{x,y}\rho^i||\nabla_{x,y}(u_d^i,v_d^i)|+|\nabla_{x,y}^2\rho^i||(u_d^i,v_d^i)|,\\
	|\nabla_{x,y}^2g_{\rho,d}^i|&\leq |\rho^i||\nabla_{x,y}^3(u_d^i,v_d^i)|+|\nabla_{x,y}\rho^i||\nabla_{x,y}^2(u_d^i,v_d^i)|+|\nabla_{x,y}^2\rho^i||\nabla_{x,y}(u_d^i,v_d^i)|+|\nabla_{x,y}^3\rho^i||(u_d^i,v_d^i)|.\nonumber
\end{aligned}
\end{equation}
Therefore, we have
	\begin{align}
		\|g_{\rho,d}^i\|_{L^2(\Omega)}&\leq \|\rho^i\|_{L^\infty}\|\nabla_{x,y}(u_d^i,v_d^i)\|_{L^2(\Omega)}+\|\nabla_{x,y}\rho^i\|_{L^2(\Omega)}\|(u_d^i,v_d^i)\|_{L^\infty}\nonumber\\
		&\leq C\nu^{-1^-}\|(\rho^i,u^i,v^i)\|_{\mathfrak{X}}\|(\rho^i_d,u^i_d,v^i_d)\|_{2},\label{d1}\\
		\|\nabla_{x,y}g_{\rho,d}^i\|_{L^2(\Omega)}&\leq \|\rho^i\|_{L^\infty}\|\nabla_{x,y}^2(u_d^i,v_d^i)\|_{L^2(\Omega)}+\|\nabla_{x,y}\rho^i\|_{L^\infty}\|\nabla_{x,y}(u_d^i,v_d^i)\|_{L^2(\Omega)}+\|\nabla_{x,y}^2\rho^i\|_{L^2(\Omega)}\|(u_d^i,v_d^i)\|_{L^\infty}\nonumber\\
		&\leq C\nu^{-(\frac{17}{8})^-}\|(\rho^i,u^i,v^i)\|_{\mathfrak{X}}\|(\rho^i_d,u^i_d,v^i_d)\|_{2},\label{d2}\\
		\|\nabla_{x,y}^2g_{\rho,d}^i\|_{L^2(\Omega)}&\leq \|\rho^i\|_{L^\infty}\|\nabla_{x,y}^3(u_d^i,v_d^i)\|_{L^2(\Omega)}+\|\nabla_{x,y}\rho^i\|_{L^\infty}\|\nabla_{x,y}^2(u_d^i,v_d^i)\|_{L^2(\Omega)}\nonumber\\
		&\qquad+\|\nabla_{x,y}^2\rho^i\|_{L^2(\Omega)}\|\nabla_{x,y}(u_d^i,v_d^i)\|_{L^\infty}+\|\nabla_{x,y}^3\rho^i\|_{L^2(\Omega)}\|(u_d^i,v_d^i)\|_{L^\infty}\nonumber\\
		&\leq \nu^{-(\frac{13}{4})^-}\|(\rho^i,u^i,v^i)\|_{\mathfrak{X}}\|(\rho^i_d,u^i_d,v^i_d)\|_{2}.\label{d3}
	\end{align}
Other terms in $[[(g_{\rho,d}^i,g^{i}_{u,d},g^{i}_{v,d})]]_2$ can be bounded  similarly as  \eqref{L1}-\eqref{L6}, \eqref{Po5}-\eqref{Po8}, \eqref{P09}, and \eqref{P10}. We end up with the following bound
\begin{align}
\|(\rho^{i+1}_d,u_{d}^{i+1},v_{d}^{i+1})\|_{2}\leq& C\nu^{-(\frac{9}{8})^-}\|(\rho^{i}_d,u_{d}^{i},v_{d}^{i})\|_{2}\times\left(\|(\rho^i,u^i,v^i)\|_{\mathfrak{X}}+\|(\rho^{i-1},u^{i-1},v^{i-1})\|_{\mathfrak{X}}\right)\nonumber\\
&\qquad\times\left(1+\nu^{-\frac{9}{8}}\|(\rho^i,u^i,v^i)\|_{\mathfrak{X}}+\nu^{-\frac{9}{8}}\|(\rho^{i-1},u^{i-1},v^{i-1})\|_{\mathfrak{X}}\right).\label{P14}
\end{align}
Thanks to \eqref{P12}, we obtain from \eqref{P14} that   $\|(\rho^{i+1}_d,u_{d}^{i+1},v_{d}^{i+1})\|_{2}\ll\|(\rho^{i}_d,u_{d}^{i},v_{d}^{i})\|_{2}.$ Therefore, $(\rho^{i}_d,u_{d}^{i},v_{d}^{i})$ is a Cauchy sequence in $\|\cdot\|_2$ norm so that it has a limit $(\rho,u,v)$. It is straightforward to see that $(\rho,u,v)$ is the solution to the compressible Navier-Stokes system \eqref{1.5}. Moreover, following the argument in  Step 1, one can show that $(\rho,u,v)$ satisfies the estimate \eqref{1.4}. The proof of Theorem \ref{T1.1} is complete. \qed
\subsection{Proof of Theorem \ref{T1.2}}\label{S8.4} In this subsection, we prove the low Mach number limit of solutions $(\rho^m,u^m,v^m)$ obtained in Theorem \ref{T1.1}. By the uniform-in-$m$ estimate \eqref{1.4}, we have
$$
\begin{aligned}
	\rho^m &\rightarrow 0 \text{ in }H^3(\Omega),\\
	m^{-2}\rho^m&\stackrel{\ast}{\rightharpoonup} P^{\rm in} \text{ in } W^{1,\infty}(\Omega)\cap H^3(\Omega),\\
	u^m&\stackrel{\ast}{\rightharpoonup} u^{\rm in} \text{ in } W^{2,\infty}(\Omega),\\
	\nabla_{x,y}u^m&\rightharpoonup \nabla_{x,y}u^{\rm in} \text{ in }H^3(\Omega),\\
	v^m&\stackrel{\ast}{\rightharpoonup} v^{\rm in} \text{ in } W^{2,\infty}(\Omega)\cap H^4(\Omega),
\end{aligned}
$$
as the Mach number $m\rightarrow 0^+$. It is straightforward to check that $(P^{\rm in},u^{\rm in},v^{\rm in})$ is a solution to steady incompressible Navier-Stokes equations \eqref{1.5-1}. Moreover, it satisfies the following bounds
\begin{align}
	\|(P^{\rm in},u^{\rm in},v^{\rm in})\|_{\mathfrak{X}^{\rm in}}\leq C\|(F_{\rm ext,1},F_{\rm ext,2})\|_{w},\label{e1}
\end{align}
where 
\begin{align}
	\|(P^{\rm in},u^{\rm in},v^{\rm in})\|_{\mathfrak{X}^{\rm in}}\eqdef & ~\|v_0^{\rm in}\|_{L^1_y}+\|(P^{\rm in}_0,u_0^{\rm in},v_0^{\rm in})\|_{L^\infty_y}+\|(P^{\rm in}_0,v^{\rm in}_0)\|_{L^2_y}+\|(P^{\rm in}_{\neq},u^{\rm in}_{\neq},v^{\rm in}_{\neq})\|_{L^2(\Omega)}\nonumber\\
	& +\nu^{\frac12}\|\nabla_{x,y}(P^{\rm in},u^{\rm in},v^{\rm in})\|_{L^2(\Omega)}+\nu^{\frac{13}{8}}\|\nabla_{x,y}^2(P^{\rm in},u^{\rm in},v^{\rm in})\|_{L^2(\Omega)}\nonumber\\
	&+\nu^{\frac{21}{8}} \|\nabla_{x,y}^3(P^{\rm in },u^{\rm in},v^{\rm,in})\|_{L^2(\Omega)}+\nu^{\frac{29}{8}} \|\nabla_{x,y}^4(u^{\rm in},v^{\rm in})\|_{L^2(\Omega)}.\nonumber
\end{align}
Here $(P^{\rm in}_0,u^{\rm in}_0,v^{\rm,in}_0)$ and $(P^{\rm in}_{\neq},u^{\rm in}_{\neq }, v^{\rm in}_{\neq})$ are zero and non-zero modes of $(P^{\rm in},u^{\rm in},v^{\rm in})$ respectively.

To obtain convergence rate, we introduce the difference
\begin{align}
	\rho_{d}^m=\rho^m-m^2P^{\rm in},~u^m_d=u^m-u^{\rm in},~v^{m}_d=v^m-v^{\rm in},\nonumber
\end{align}
and derive equations of $(\rho^m_d,u^m_d,v^m_d)$ as follows:
\begin{align}\nonumber
	\begin{cases}
		U_s\partial_x\rho^m_d+\div_{x,y}(u_d^m,v_d^m)+\partial_x(u^m\rho^m_d)+\partial_y(v^m\rho_d^m)=G_{\rho},\\
		U_s\partial_xu^m_d+v^m_d\partial_yU_s+m^{-2}\partial_x\rho^m_d-\nu\Delta u^m_d-\lambda\nu\partial_x\div_{x,y}(u_d^m,v^m_d)+\nu\rho^m_d\partial_y^2U_s=G_{u,d}+G_u,\\
		U_s\partial_xv^m_d+m^{-2}\partial_y\rho^m_d-\nu\Delta v^m_d-\lambda\nu\partial_y\div_{x,y}(u^m_d,v^m_d)=G_{v,d}+G_v,\\
		u^m_d|_{y=0}=v^m_d|_{y=0}=0,
	\end{cases}
\end{align}
where  
\begin{align}
	G_\rho&=-m^2U_s\partial_xP^{\rm in}-m^2\partial_x(P^{\rm in}u^m)-m^2\partial_y(P^{\rm in}v^m),\nonumber\\
	G_{u}&=-\rho^m(U_s+u^m)\partial_xu^m-\rho^mv^m\partial_yu^m-\rho^mv^m\partial_yU_s\nonumber\\
	&\qquad-[P'(1+\rho^m)-P'(1)]\partial_x\rho^m+m^2\nu P^{\rm in}\partial_y^2U_s+\rho^mF_{\rm ext,1},\nonumber\\
	G_v&=-\rho^m(U_s+u^m)\partial_xv^m-\rho^mv^m\partial_yv^m-[P'(1+\rho^m)-P'(1)]\partial_y\rho^m+\rho^mF_{\rm ext,2},\nonumber\\
	G_{u,d}&=-u_{d}^m\partial_xu^m-u^{\rm in}\partial_xu^m_d-v^m_d\partial_yu^m-v^{\rm in}\partial_yu^m_d,\nonumber\\
	G_{v,d}&=-u_{d}^m\partial_xv^m-u^{\rm in}\partial_xv^m_d-v^m_d\partial_yv^m-v^{\rm in}\partial_yv^m_d.\nonumber
\end{align}
Applying linear bounds \eqref{LS} to $(\rho^m_d,u^m_d,v^m_d)$, we have
\begin{align}
	\|(\rho^m_d,u^m_d,v^m_d)\|_{2}\leq C[[(0,G_{u,d},G_{v,d})]]_2+C[[(G_\rho,G_u,G_v)]]_2.\label{e4}
\end{align}

For the first term, we have
$$
\begin{aligned}
\|(G_{u,d},G_{v,d})\|_{L^2(\Omega)}&\leq \|(u^m_d,v^m_d)\|_{L^\infty(\Omega)}\|\nabla_{x,y}(u^m,v^m)\|_{L^2(\Omega)}+\|(u^{\rm in},v^{\rm in})\|_{L^\infty(\Omega)}\|\nabla_{x,y}(u^m_d,v^m_d)\|_{L^2(\Omega)}\\
&\leq \nu^{-1-}\|(\rho^m_d,u^m_d,v^m_d)\|_{2}(\|(\rho^m,u^m,v^m)\|_{\mathfrak{X}}+\|(P^{\rm in},u^{\rm in},v^{\rm in})\|_{\mathfrak{X}^{\rm in}})\nonumber\\
&\leq \nu^{-1-}\|(\rho^m_d,u^m_d,v^m_d)\|_{2}\|(F_{\rm ext,1},F_{\rm ext,2})\|_{w},
\end{aligned}
$$
where we have used $L^\infty$ bounds \eqref{Lw1} and uniform estimates \eqref{1.4} and \eqref{e1}. Similarly, for derivative of $(G_{u,d},G_{v,d})$ we have
$$
\begin{aligned}
|\nabla_{x,y}(G_{u,d},G_{v,d})|\leq& |(u_d^m,v_d^m)||\nabla^2_{x,y}(u^m,v^m)|+|\nabla_{x,y}(u^m_d,v^m_d)||\nabla(u^m,v^m)|\nonumber\\
&+|(u^{\rm in},v^{\rm in})||\nabla^2(u^m_d,v^m_d)|+|\nabla_{x,y}(u^{\rm in},v^{\rm in})||\nabla_{x,y}(u_d^m,v_d^m)|,
\end{aligned}
$$
which implies
$$
\begin{aligned}
	\|\nabla_{x,y}(G_{u,d},G_{v,d})\|_{L^2(\Omega)}&\leq \|(u^m_d,v^m_d)\|_{L^\infty(\Omega)}\|\nabla_{x,y}^2(u^m,v^m)\|_{L^2(\Omega)}+\|(u^{\rm in},v^{\rm in})\|_{L^\infty(\Omega)}\|\nabla_{x,y}^2(u^m_d,v^m_d)\|_{L^2(\Omega)}\\
	&\qquad+(\|\nabla_{x,y}(u^{\rm in},v^{\rm in})\|_{L^\infty(\Omega)}+\|\nabla_{x,y}(u^m,v^m)\|_{L^\infty(\Omega)})\|\nabla_{x,y}(u_d^m,v^m_d)\|_{L^2(\Omega)}\\
	&\leq \nu^{-\frac{17}{8}-}\|(\rho^m_d,u^m_d,v^m_d)\|_{2}(\|(\rho^m,u^m,v^m)\|_{\mathfrak{X}}+\|(P^{\rm in},u^{\rm in},v^{\rm in})\|_{\mathfrak{X}^{\rm in}})\nonumber\\
	&\leq \nu^{-\frac{17}{8}-}\|(\rho^m_d,u^m_d,v^m_d)\|_{2}\|(F_{\rm ext,1},F_{\rm ext,2})\|_{w}.
\end{aligned}
$$
Here we have used $W^{1,\infty}$ bounds \eqref{Lw1}-\eqref{Lw3}. Other terms in $[[(0,G_{u,d},G_{v,d})]]_2$ can be treated in the same manner as \eqref{Po5}-\eqref{Po8}. Thus, we end up with the following bound:
\begin{align}
	[[(0,G_{u,d},G_{v,d})]]_2\leq C\nu^{-\frac{9}{8}-}\|(\rho^m_d,u^m_d,v^m_d)\|_{2}\|(F_{\rm ext,1},F_{\rm ext,2})\|_{w}.\label{e5}
\end{align}

Next we estimate the second term $[[(G_{\rho},G_u,G_v)]]_2$ on the right hand side of \eqref{e4}. Similar to \eqref{d1}-\eqref{d3}, we have
\begin{equation}
\begin{aligned}
	\|G_\rho\|_{L^2(\Omega)} &\leq Cm^2\left[\|\nabla_{x,y}P^{\rm in}\|_{L^2(\Omega)}(1+\|(u^m,v^m)\|_{L^\infty(\Omega)})+\|P^{\rm in}\|_{L^\infty(\Omega)}\|\nabla_{x,y}(u^m,v^m)\|_{L^2(\Omega)}\right]\\
	&\leq Cm^2\nu^{-\frac12-}\|(P^{\rm in},u^{\rm,in},v^{\rm,in})\|_{\mathfrak{X}^{\rm in}}(1+\nu^{-\frac12}\|(\rho^m,u^m,v^m)\|_{\mathfrak{X}})\nonumber\\
	&\leq Cm^2\nu^{-\frac12-}\|(F_{\rm ext,1},F_{\rm ext,2})\|_w,\\
	\|\nabla_{x,y}G_\rho\|_{L^2(\Omega)}&\leq Cm^2\bigg[\|\nabla^2_{x,y}P^{\rm in}\|_{L^2(\Omega)}+\nu^{-\frac12}\|\nabla_{x,y}P^{\rm in}\|_{L^2(\Omega)}+\|P^{\rm in}\|_{L^\infty(\Omega)}\|\nabla_{x,y}^2(u^m,v^m)\|_{L^2(\Omega)}\\
	&\qquad +\|\nabla_{x,y}P^{\rm in}\|_{L^2(\Omega)}\|\nabla_{x,y}(u^m,v^m)\|_{L^\infty(\Omega)}+\|\nabla_{x,y}^2P^{\rm in}\|_{L^2(\Omega)}\|(u^m,v^m)\|_{L^\infty(\Omega)}\bigg]\\
	&\leq Cm^2\nu^{-\frac{13}{8}-}\|(P^{\rm in},u^{\rm,in},v^{\rm,in})\|_{\mathfrak{X}^{\rm in}}(1+\nu^{-\frac12}\|(\rho^m,u^m,v^m)\|_{\mathfrak{X}})\nonumber\\
	&\leq Cm^2\nu^{-\frac{13}{8}-}\|(F_{\rm ext,1},F_{\rm ext,2})\|_w,\\
	\|\nabla_{x,y}^2G_\rho\|_{L^2(\Omega)}&\leq Cm^2\bigg[\sum_{j=1}^3\nu^{-\frac{3-j}{2}}\|\nabla^j_{x,y}P^{\rm in}\|_{L^2(\Omega)}+\|P^{\rm in}\|_{L^\infty(\Omega)}\|\nabla^3_{x,y}(u^m,v^m)\|_{L^2(\Omega)}\\
	&\qquad+\sum_{j=1}^3\|\nabla_{x,y}^jP^{\rm in}\|_{L^2(\Omega)}\|\nabla^{3-j}_{x,y}(u^m,v^m)\|_{L^\infty(\Omega)}\bigg]\\
	&\leq Cm^2\nu^{-\frac{21}{8}-}\|(P^{\rm in},u^{\rm,in},v^{\rm,in})\|_{\mathfrak{X}^{\rm in}}(1+\nu^{-\frac58}\|(\rho^m,u^m,v^m)\|_{\mathfrak{X}})\nonumber\\
	&\leq Cm^2\nu^{-\frac{21}{8}-}\|(F_{\rm ext,1},F_{\rm ext,2})\|_w.
\end{aligned}
\end{equation}
Here we have used $\|(\rho^m,u^m,v^m)\|_{\mathfrak{X}}+\|(\rho^{\rm in},u^{\rm in},v^{\rm in})\|_{\mathfrak{X}^{\rm in}}\leq C\|(F_{\rm ext,1},F_{\rm ext,2})\|_w\leq \nu^{\frac98+}$. Moreover, for $(G_{u},G_v)$ we have the following pointwise estimate
\begin{align}
	|(G_u,G_v)|\leq& |\rho^m|(1+|(u^m,v^m)|)|\nabla_{x,y}(u^m,v^m)|+|\rho^m||\nabla_{x,y}\rho^m|\nonumber\\
	&+\nu^{-\frac12}|\rho^mv^m|+|\rho^m||(F_{\rm ext,1},F_{\rm ext, 2})|+m^2|P^{\rm in}|.\nonumber
\end{align}
Note that each term involving $\rho^m$ in the solution norm $\|\cdot\|_{\mathfrak{X}}$ is accompanied by a factor of $m^{-2}$. Hence we have
$$
\begin{aligned}
\|(G_u,G_v)\|_{L^2(\Omega)}\leq& C\|\rho^m\|_{L^\infty(\Omega)}\bigg[\nu^{-\frac12}\|v^m\|_{L^2(\Omega)}+\|\nabla_{x,y}(\rho^m,u^m,v^m)\|_{L^2(\Omega)}+\|(F_{\rm ext,1},F_{\rm ext, 2})\|_{L^2(\Omega)}\nonumber\\
&+\|(u^m,v^m)\|_{L^\infty(\Omega)}\|\nabla_{x,y}(u^m,v^m)\|_{L^2(\Omega)}\bigg]+Cm^2\|P^{\rm in}\|_{L^2(\Omega)}\\
\leq& Cm^2\nu^{-\frac12-}\|(\rho^m,u^m,v^m)\|_{\mathfrak{X}}\bigg[\nu^{-\frac12}\|(\rho^m,u^m,v^m)\|_{\mathfrak{X}}+\nu^{-1}\|(\rho^m,u^m,v^m)\|_{\mathfrak{X}}^2\nonumber\\
&\qquad+\|(F_{\rm ext,1},F_{\rm ext, 2})\|_{w}
\bigg]+Cm^2\|(F_{\rm ext,1},F_{\rm ext, 2})\|_{w} \\
\leq &Cm^2\|(F_{\rm ext,1},F_{\rm ext, 2})\|_{w}.
\end{aligned}
$$
Similarly, we deduce
$$
\begin{aligned}
	\|\nabla_{x,y}(G_u,G_v)\|_{L^2(\Omega)}
	\leq& Cm^2\nu^{-\frac{13}{8}-}\|(\rho^m,u^m,v^m)\|_{\mathfrak{X}}\bigg[\nu^{-\frac12}\|(\rho^m,u^m,v^m)\|_{\mathfrak{X}}+\nu^{-1}\|(\rho^m,u^m,v^m)\|_{\mathfrak{X}}^2\nonumber\\
	&\qquad+\|(F_{\rm ext,1},F_{\rm ext, 2})\|_{w}
	\bigg]+Cm^2\nu^{-\frac12}\|(F_{\rm ext,1},F_{\rm ext, 2})\|_{w} \\
	\leq &Cm^2\nu^{-1-}\|(F_{\rm ext,1},F_{\rm ext, 2})\|_{w}.
\end{aligned}
$$
Other terms in $[[(G_\rho,G_u,G_v)]]_2$ can be bounded similarly as \eqref{L1}-\eqref{L6}, and \eqref{Po5}-\eqref{Po8}. Then we have
$$[[(G_\rho,G_u,G_v)]]_2\leq Cm^2\nu^{-\frac{5}{8}-}\|(F_{\rm ext,1},F_{\rm ext, 2})\|_{w}.
$$
Combining this with \eqref{e5} and \eqref{e4} we get
$$\|(\rho^m_d,u^m_d,v^m_d)\|_{2}\leq C\nu^{-\frac{9}{8}-}\|(\rho^m_d,u^m_d,v^m_d)\|_{2}\|(F_{\rm ext,1},F_{\rm ext, 2})\|_{w}+Cm^2\nu^{-\frac{5}{8}-}\|(F_{\rm ext,1},F_{\rm ext, 2})\|_{w}.
$$
Therefore, \eqref{T1.2-1} follows from $\|(F_{\rm ext,1},F_{\rm ext, 2})\|_{w}\leq \nu^{\frac98+}$ and \eqref{T1.2-2} follows from $L^\infty$ bound \eqref{Lw1}. The proof of Theorem \ref{T1.2} is complete. \qed

\section{Appendix}

\subsection{Higher-order estimates on Stokes equations}
In this section, we derive the higher-order estimates on solutions $(\rho,u,v)$ to the following Stokes equations in the unbounded strip $\Omega=\{(x,y)\mid x\in \mathbb{T}_L, y\in \mathbb{R}_+\}$, which has been used in the proof of Lemma \ref{lmHF3}.
\begin{equation}\label{10.1}
	\left\{
	\begin{aligned}
		&\partial_xu+\partial_yv=h_\rho,\\
&	-{\nu}\Delta_{x,y}u+m^{-2}\partial_x\rho=h_u,\\
&-{\nu}\Delta_{x,y}v+m^{-2}\partial_y\rho=h_v,\\
&u|_{y=0}=v|_{y=0}=0.
	\end{aligned}
\right.
\end{equation}
\begin{lemma}\label{lem9.1}
Let $k=2$ or $3$. Suppose that $(h_\rho,h_u,h_v)\in H^k(\Omega)\times H^{k-1}(\Omega)^2$ and satisfies
	\begin{align}
		\iint_{\Omega}h_\rho(x,y)d xd y=0.\nonumber
	\end{align}
 There exists a unique solution  $(\rho,u,v)
\in H^k(\Omega)\times \dot{H}^{k+1}(\Omega)^2$ to \eqref{10.1} 
satisfying the following estimate
\begin{align}
	&\nu\|\nabla^{k+1}_{x,y}(u,v)\|_{L^2(\Omega)}+m^{-2}\|\nabla^k_{x,y}\rho\|_{L^2(\Omega)}\nonumber\\ &\qquad\leq C\nu\|h_\rho\|_{H^k(\Omega)}+C\|(h_u,h_v)\|_{H^{k-1}(\Omega)}+Cm^{-2}\|\rho\|_{H^{k-1}(\Omega)}+C\nu\|\nabla_{x,y}(u,v)\|_{H^{k-1}(\Omega)}.\label{10.2}
\end{align}
\end{lemma}
\begin{proof}
%	To get uniform-in-$\nu$ bounds on the solutions, we rescale the Stokes equations by  $(\tilde{\rho},\tilde{u},\tilde{v})(x,y)=(\rho,u,v)(\nu^{-\frac12}x,\nu^{-\frac12}y)$. The equations for $(\tilde{\rho},\tilde{u},\tilde{v})(x,y)$ reads
%	\begin{equation}\label{10.3}
%		\left\{
%		\begin{aligned}
%			&\partial_x\tilde{u}+\partial_y\tilde{v}=\nu^{-\frac12}g_0,\\
%			&	-\nu\Delta_{x,y}\tilde{u}+m^{-2}\partial_x\tilde{\rho}=\nu^{-\frac12}g_1,\\
%			&-{\nu}\Delta_{x,y}v+m^{-2}\partial_y\tilde{\rho}=\nu^{-\frac12}g_2,\\
%			&u|_{y=0}=v|_{y=0}=0,
%		\end{aligned}
%		\right.
%	\end{equation}
%where $(x,y)\in \Omega= \mathbb{T}_L\times\mathbb{R}_+.$  
We introduce a smooth cut-off function $\chi(y)$, which satisfies $\chi\equiv 1$ for $y\in [0,1]$, and $\chi\equiv0$ for $y\geq 2.$ Then we define
$(\tilde{\rho}_1,\tilde{u}_1,\tilde{v}_1)=(\chi{\rho},\chi{u},\chi{v})$, and 
$(\tilde{\rho}_2,\tilde{u}_2,\tilde{v}_2)=\left((1-\chi){\rho},(1-\chi){u},(1-\chi){v}\right)$.
 The corresponding  equations read 
	\begin{equation}\label{10.4}
	\left\{
	\begin{aligned}
		&\partial_x\tilde{u}_1+\partial_y\tilde{v}_1=\chi'v+\chi h_\rho\eqdef \tilde{h}_{\rho,1},\\
		&	-\nu\Delta_{x,y}\tilde{u}_1+m^{-2}\partial_x\tilde{\rho}_1=-2\nu\chi'\partial_y{u}-\nu\chi''{u}+\chi h_u\eqdef \tilde{h}_{u,1},\\
		&-{\nu}\Delta_{x,y}\tilde{v}_{1}+m^{-2}\partial_y\tilde{\rho}_1=-2\nu\chi'\partial_y{v}-\nu\chi''{v}+m^{-2}\chi'{\rho}+\chi h_v\eqdef\tilde{h}_{v,1},\\
		&\tilde{u}_1|_{y=0,2}=\tilde{v}_1|_{y=0,2}=0,
	\end{aligned}
	\right.
\end{equation}
and
	\begin{equation}\label{10.5}
	\left\{
	\begin{aligned}
		&\partial_x\tilde{u}_2+\partial_y\tilde{v}_2=-\chi'v+(1-\chi) h_\rho\eqdef \tilde{h}_{\rho,2},\\
		&	-\nu\Delta_{x,y}\tilde{u}_2+m^{-2}\partial_x\tilde{\rho}_2=2\nu\chi'\partial_y{u}+\nu\chi''{u}+(1-\chi)h_u\eqdef \tilde{h}_{u,2},\\
		&-{\nu}\Delta_{x,y}\tilde{v}_{2}+m^{-2}\partial_y\tilde{\rho}_2=2\nu\chi'\partial_y{v}+\nu\chi''{v}-m^{-2}\chi'{\rho}+(1-\chi)h_v\eqdef\tilde{h}_{v,2},\\
		&\tilde{u}_2|_{y=0}=\tilde{v}_2|_{y=0}=0.
	\end{aligned}
	\right.
\end{equation}
Since 
$$\int_{0}^2\int_{0}^{2\pi L}\tilde{h}_{\rho,1} dxdy=\int_{0}^2\int_{0}^{2\pi L}\chi(y)\partial_xud xdy=0,
$$
then applying higher-order estimates (cf. monographs \cite{GA,SH}) for Stokes equations to  $(\tilde{u}_1,\tilde{v}_1)(x,y)$ and using Hardy's inequality
$\|\partial_y^j\chi(y)(u,v)\|_{L^2(\Omega)}\leq C\|\partial_y(u,v)\|_{L^2(\Omega)}~, j=0,1,2$, we have
\begin{align}
	&\nu\|\nabla^{k+1}_{x,y}(\tilde{u}_1,\tilde{v}_1)\|_{L^2(\Omega)}+m^{-2}\|\nabla_{x,y}^k\tilde{\rho}_{1}\|_{L^2(\Omega)}\nonumber\\
	&\qquad\quad\leq C\nu\|\tilde{h}_{\rho,1}\|_{H^k(\Omega)}+C\|(\tilde{h}_{u,1},\tilde{h}_{v,1})\|_{H^{k-1}(\Omega)}+C\nu\|(\tilde{u}_1,\tilde{v}_1)\|_{L^2(\Omega)}\nonumber\\
	&\qquad\quad\leq C\nu\|h_\rho\|_{H^k(\Omega)}+C\|(h_u,h_v)\|_{H^{k-1}(\Omega)}+C\nu\|\nabla_{x,y}({u},{v})\|_{H^{k-1}(\Omega)}+Cm^{-2}\|{\rho}\|_{H^{k-1}(\Omega)}.\label{10.6}
\end{align}
To estimate $(\tilde{\rho}_2,\tilde{u}_2,\tilde{v}_2)$, we consider vorticity $\tilde{\omega}_2=\partial_y\tilde{u}_2-\partial_x\tilde{v}_2$, which satisfies
$$\nu\Delta_{x,y}\tilde{\omega}_2=\partial_y\tilde{h}_{u,2}-\partial_x\tilde{h}_{v,2}.
$$
Since $\tilde{\omega}_{2}$ is supported away from the boundary,  we have, from the regularity of the Laplace equation, that
\begin{align}
	\nu\|\nabla^k_{x,y}\tilde{\omega}_2\|_{L^2(\Omega)}&\leq C\nu\|\Delta_{x,y}\tilde{\omega}_2\|_{H^{k-2}(\Omega)}\leq  C\|(\tilde{h}_{u,2},\tilde{h}_{v,2})\|_{H^{k-1}(\Omega)}\nonumber\\
	&\leq C\|(h_u,h_v)\|_{H^{k-1}(\Omega)}+\nu\|\nabla_{x,y}({u},{v})\|_{H^{k-1}(\Omega)}+m^{-2}\|{\rho}\|_{H^{k-1}(\Omega)}.\nonumber
\end{align}
Combining this with the estimates on divergence part
$$
\begin{aligned}
\nu\|\nabla^{k}_{x,y}(\partial_x\tilde{u}_2+\partial_y\tilde{v}_2)\|_{L^2(\Omega)}&\leq C\nu\|\tilde{h}_{\rho,2}\|_{H^k(\Omega)}\\
&\leq C\nu\|h_\rho\|_{H^k(\Omega)}+C\nu\|\nabla_{x,y}{v}\|_{H^{k-1}(\Omega)},
\end{aligned}
$$
we obtain
\begin{align}
\nu\|\nabla^{k+1}_{x,y}(\tilde{u}_2,\tilde{v}_2)\|_{L^2(\Omega)}
\leq& C\|(h_u,h_v)\|_{H^{k-1}(\Omega)}+C\nu\|h_\rho\|_{H^k(\Omega)}\nonumber\\
&+C\nu\|\nabla_{x,y}({u},{v})\|_{H^{k-1}(\Omega)}+Cm^{-2}\|{\rho}\|_{H^{k-1}(\Omega)}.\label{10.7}
\end{align}
From the equations $\eqref{10.5}_1$ and $\eqref{10.5}_2$ we can bound the density $\tilde{\rho}_2$ as follows. 
\begin{align}\label{10.8}
	m^{-2}\|\nabla^k_{x,y}\tilde{\rho}_2\|_{L^2(\Omega)}
	\leq& C\nu\|\nabla^{k+1}_{x,y}(\tilde{u}_2,\tilde{v}_2)\|_{L^2(\Omega)}+C\|(\tilde{h}_{u,2},\tilde{h}_{v,2})\|_{H^{k-1}(\Omega)}\nonumber\\
	\leq& C\|(h_u,h_v)\|_{H^{k-1}(\Omega)}+C\nu\|h_\rho\|_{H^k(\Omega)}\nonumber\\
	&+C\nu\|\nabla_{x,y}({u},{v})\|_{H^{k-1}(\Omega)}+Cm^{-2}\|{\rho}\|_{H^{k-1}(\Omega)}.
\end{align}
Conbining \eqref{10.6}, \eqref{10.7} and \eqref{10.8}, we get
\begin{align}
	&\nu\|\nabla^{k+1}_{x,y}({u},{v})\|_{L^2(\Omega)}+m^{-2}\|\nabla_{x,y}^k{\rho}\|_{L^2(\Omega)} \nonumber\\ &\qquad\leq C\nu\|h_\rho\|_{H^k(\Omega)}+C \|(h_u,h_v)\|_{H^{k-1}(\Omega)}+C\nu\|\nabla_{x,y}({u},{v})\|_{H^{k-1}(\Omega)}+Cm^{-2}\|{\rho}\|_{H^{k-1}(\Omega)},\nonumber
\end{align}
which is \eqref{10.2}. The proof of Lemma \ref{lem9.1} is complete.
\end{proof}
\subsection{Nonlinear estimates}
For any function $f$ defined on $\Omega$, we denote its zero and non-zero modes by $f_0$ and $f_{\neq}$ respectively.
\begin{lemma}[$L^\infty$-inequality]
	Let $\beta>0$. 
	There exists a constant $C_\beta$, such that
	\begin{align}
		\|f\|_{L^\infty(\Omega)}\leq \|f_0\|_{L^\infty_y}+C_\beta\|\partial_y f_{\neq}\|_{L^2(\Omega)}^{\frac12}\|\partial_x^{1+\beta}f_{\neq}\|_{L^2(\Omega)}^{\frac12}.\label{lw}
	\end{align}
\end{lemma}
\begin{proof}
	Recall the identity: $f=f_0+\sum_{n\neq 0}f_ne^{i\hat{n}x}$. Then by Sobolev inequality $\|f_n\|_{L^\infty_y}\leq C\|f_n\|_{L^2_y}^{\frac12}\|\partial_yf_n\|_{L^2_y}^{\frac12},$ we obtain
	\begin{align}
		\|f\|_{L^\infty(\Omega)}&\leq \|f_0\|_{L^\infty_y}+\sum_{n\neq 0}\|f_n\|_{L^\infty_y}\leq \|f_0\|_{L^\infty_y}+C\sum_{n\neq 0 }\|f_n\|_{L^2_y}^{\frac12}\|\partial_yf_n\|_{L^2_y}^{\frac12}\nonumber\\
		&\leq \|f_0\|_{L^\infty_y}+C\left(\sum_{n\neq 0} n^{2(1+\beta)}\|f_n\|_{L^2_y}^2\right)^{\frac14}\left(\sum_{n\neq 0} \|\partial_yf_n\|_{L^2_y}^2\right)^{\frac14}\left(\sum_{n\neq 0} n^{-(1+\beta)}\right)^{\frac12}\nonumber\\
		&\leq \|f_0\|_{L^\infty_y}+C_\beta\|\partial_x^{1+\beta}f_{\neq}\|_{L^2(\Omega)}^{\frac12}\|\partial_yf_{\neq}\|_{L^2(\Omega)}^{\frac12},\nonumber
	\end{align}
	which is \eqref{lw}. Thus, the proof is complete.
\end{proof}

\begin{lemma}\label{ne}
	Let $f$ and $g$ be two functions defined on $\Omega$. It holds that
	\begin{align}
		\|\mathcal{P}_0(fg)\|_{L^1_y}&\leq C\min\{\|f_0\|_{L^\infty_y}\|g_0\|_{L^1_y},\|f_0\|_{L^2_y}\|g_0\|_{L^2_y}\}+C\|f_{\neq}\|_{L^2(\Omega)}\|g_{\neq}\|_{L^2(\Omega)},\label{A9.1}\\
		\|\mathcal{P}_0(fg)\|_{L^2_y}&\leq C\|f_0\|_{L^\infty_y}\|g_0\|_{L^2_y}+ C\|f_{\neq}\|_{L^2(\Omega)}^{\frac12}\|\partial_yf_{\neq}\|_{L^2(\Omega)}^{\frac12}\|g_{\neq}\|_{L^2(\Omega)},\label{A9.2}\\
		\|\mathcal{P}_0(fg)\|_{L^\infty_y}&\leq C\|f_0\|_{L^\infty_y}\|g_0\|_{L^\infty_y}+C\|f_{\neq}\|_{L^2(\Omega)}^{\frac12}\|g_{\neq}\|_{L^2(\Omega)}^{\frac12}\|\partial_yf_{\neq}\|_{L^2(\Omega)}^{\frac12}\|\partial_yg_{\neq}\|_{L^2(\Omega)}^{\frac12}\label{A9.3},\\
		\|\mathcal{P}_{\neq}(fg)\|_{L^2(\Omega)}&\leq C\|f\|_{L^\infty(\Omega)}\|g_{\neq}\|_{L^2(\Omega)}+C\|g_0\|_{L^\infty_y}\|f_{\neq}\|_{L^2(\Omega)}.\label{A9.4}
	\end{align}
If $g=\partial_xh$ or $\partial_yh$, it holds that
\begin{align}
	\|\mathcal{P}_{\neq}(fg)\|_{L^2(\Omega)}\leq \|f_0\|_{L^\infty_y}\|g_{\neq}\|_{L^2(\Omega)}+\|\partial_yh_0\|_{L^2_y}\|f_{\neq}\|_{L^\infty(\Omega)}+\|f_{\neq}\|_{L^\infty(\Omega)}\|g_{\neq}\|_{L^2(\Omega)}\label{A9.5}.
\end{align}
\end{lemma}
\begin{proof}
	Note that $\mathcal{P}_0(fg)=f_0g_0+\sum_{n+m=0}f_ng_m$. A straightforward computation yields
	\begin{align}
		\|\mathcal{P}_0(fg)\|_{L^1_y}&\leq \|f_0g_0\|_{L^1_y}+\sum_{n+m=0}\|f_ng_m\|_{L^1_y}\leq \|f_0g_0\|_{L^1_y}+\sum_{n+m=0}\|f_n\|_{L^2_y}\|g_m\|_{L^2_y}\nonumber\\
		&\leq C\min\{\|f_0\|_{L^\infty_y}\|g_0\|_{L^1_y},\|f_0\|_{L^2_y}\|g_0\|_{L^2_y}\}+(\sum_{n\neq 0}\|f_n\|_{L^2_y}^2)^{\frac12}(\sum_{n\neq 0}\|g_n\|_{L^2_y}^2)^{\frac12}\nonumber\\
		&\leq C\min\{\|f_0\|_{L^\infty_y}\|g_0\|_{L^1_y},\|f_0\|_{L^2_y}\|g_0\|_{L^2_y}\}+\|f_{\neq}\|_{L^2(\Omega)}\|g_{\neq}\|_{L^2(\Omega)},\nonumber
	\end{align}
which is \eqref{A9.1}. Similarly, by Sobolev inequality $\|f\|_{L^\infty_y}\leq C\|f\|_{L^2_y}^{\frac12}\|\partial_yf\|_{L^2_y}^{\frac12}$, we can obtain
\begin{align}
	\|\mathcal{P}_0(fg)\|_{L^\infty_y}&\leq \|f_0\|_{L^\infty_y}\|g_0\|_{L^\infty_y}+\sum_{n\neq 0}\|f_n\|_{L^\infty_y}\|g_{-n}\|_{L^\infty_y}\nonumber\\
	&\leq \|f_0\|_{L^\infty_y}\|g_0\|_{L^\infty_y}+\sum_{n\neq 0}\|f_n\|_{L^2_y}^{\frac12}\|\partial_yf_n\|_{L^2_y}^{\frac12}\|g_{-n}\|_{L^2_y}^{\frac12}\|\partial_yg_{-n}\|_{L^2_y}^{\frac12}\nonumber\\
	&\leq \|f_0\|_{L^\infty_y}\|g_0\|_{L^\infty_y}+(\sum_{n\neq 0}\|f_n\|_{L^2_y}^2)^{\frac12}(\sum_{n\neq 0}\|g_n\|_{L^2_y}^2)^{\frac12}(\sum_{n\neq 0}\|\partial_yf_n\|_{L^2_y}^2)^{\frac12}(\sum_{n\neq 0}\|\partial_yg_n\|_{L^2_y}^2)^{\frac12}\nonumber\\
	&\leq \|f_0\|_{L^\infty_y}\|g_0\|_{L^\infty_y}+\|f_{\neq}\|_{L^2(\Omega)}^{\frac12}\|g_{\neq}\|_{L^2(\Omega)}^{\frac12}\|\partial_yf_{\neq}\|_{L^2(\Omega)}^{\frac12}\|\partial_yg_{\neq}\|_{L^2(\Omega)}^{\frac12},\nonumber
\end{align}
which yields \eqref{A9.3}. The estimate \eqref{A9.2} can be proved similarly. 

For non-zero modes, the estimate \eqref{A9.4} follows from the identity
$$\mathcal{P}_{\neq}(fg)=f_0g_{\neq}+g_0f_{\neq}+f_{\neq}g_{\neq}.
$$
For $g=\partial_xh$ or $\partial_yh$, note that $\mathcal{P}_0(\partial_xh)=0$. Thus,  we have
$$\|\mathcal{P}_{\neq}(fg)\|_{L^2(\Omega)}=\|f_0\|_{L^\infty_y}\|g_{\neq}\|_{L^2(\Omega)}+\|\partial_yh_0\|_{L^2_y}\|f_{\neq}\|_{L^\infty(\Omega)}+\|f_{\neq}\|_{L^\infty(\Omega)}\|g_{\neq}\|_{L^2(\Omega)},
$$
which is \eqref{A9.5}. The proof of Lemma \ref{ne} is complete.
\end{proof}
\bigbreak

\noindent{\bf Acknowledgements}. 
The research of T. Yang is supported by the General Research Fund of Hong Kong (Project No. 11302020). He would also like to thank the Kuok Group foundation for its generous support.  The research of Z. Zhang is supported by the General Research Fund of Hong Kong (Project No. 15300024), a start-up grant (Project No. P0043862), and Research Center for Nonlinear Analysis (Project No. P0046121), PolyU.\\

\noindent{\bf Data availability.} This manuscript has no associated data.

\end{document}